\documentclass[11pt,a4paper,twoside]{report}

\usepackage{graphicx}
\usepackage{amscd}
\usepackage{amsmath,latexsym,amssymb,amsthm}
\usepackage{setspace,cite}
\usepackage{lscape,fancyhdr,fancybox}
\usepackage{a4}
\usepackage{mystyle}
\usepackage[hmarginratio=1:1, vmarginratio =5:5,
textheight=22cm,bindingoffset=1.5cm, textwidth=14.6cm]{geometry}
\renewcommand{\baselinestretch}{1.2}

\numberwithin{equation}{section}

\begin{document}
\pagestyle{empty}
\vspace*{1cm}
\begin{center}
\begin{huge}
\textbf{QUANTUM ISOMETRY GROUPS\\}
\end{huge}
\end{center}

\vspace{2cm}

\begin{center}
\begin{large}
\textbf{JYOTISHMAN BHOWMICK\\}
\end{large}
\end{center}

\vspace{3cm}
\begin{figure}[h]
\begin{center}
\includegraphics[width=4cm]{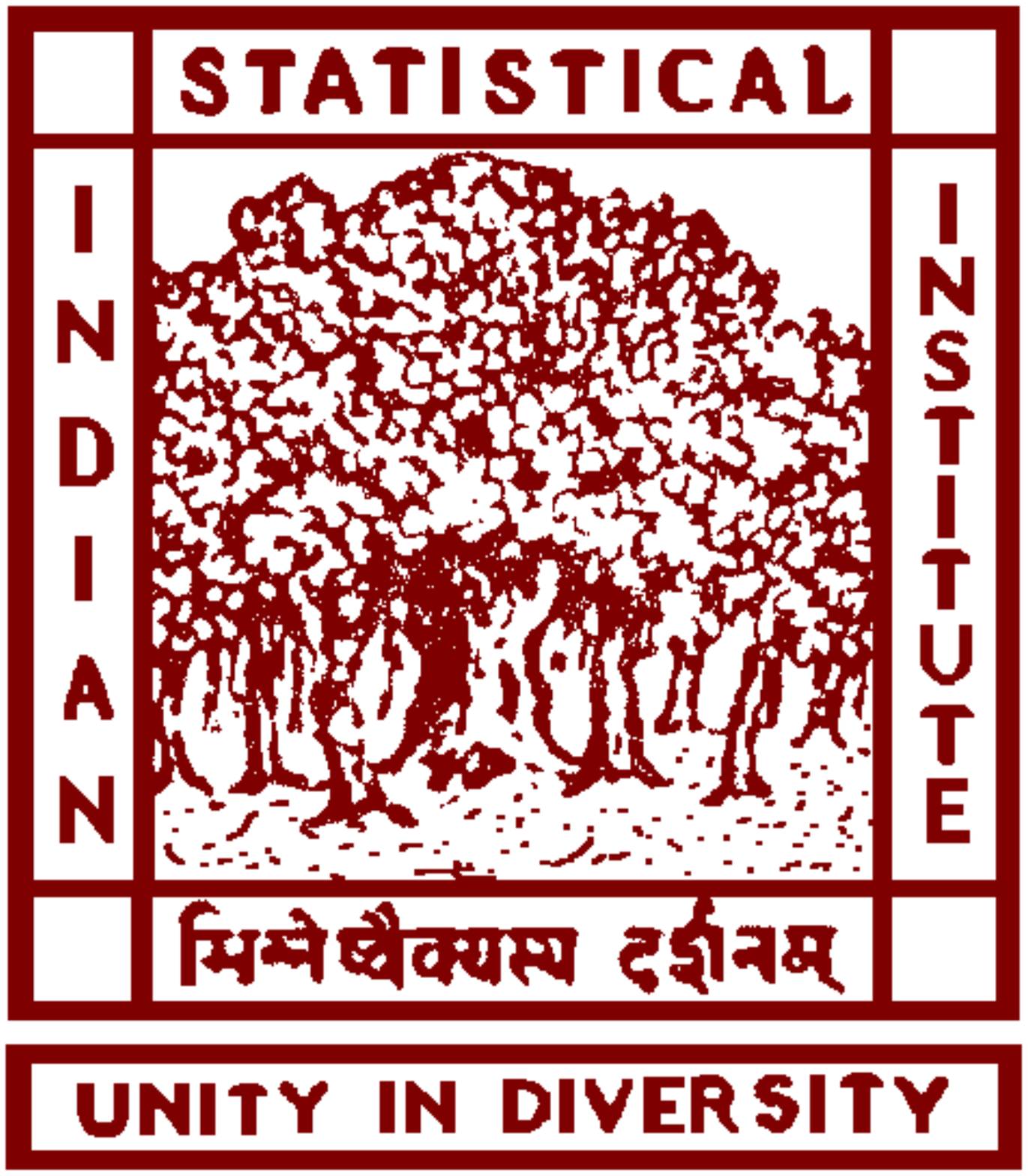}
\end{center}
\end{figure}
\vspace{2cm}

\begin{center}
\begin{large}
\textbf{Indian Statistical Institute, Kolkata\\
June,
2009}
\end{large}
\end{center}

\cleardoublepage

\pagestyle{empty}
\vspace*{1cm}
\begin{center}
\begin{huge}
\textbf{QUANTUM ISOMETRY GROUPS\\}
\end{huge}
\end{center}

\vspace{3cm}

\begin{center}
\begin{large}
\textbf{JYOTISHMAN BHOWMICK\\}
\end{large}
\end{center}

\vspace{4cm}
\begin{center}
\textbf{Thesis submitted to the Indian Statistical Institute\\
in partial fulfillment of the requirements\\
for the award of
the degree of\\
Doctor of Philosophy.\\
June,2009}\\
\textbf{Thesis Advisor: Debashish Goswami}
\end{center}

\vspace{3cm}
\begin{center}
\begin{large}
\textbf{Indian Statistical Institute\\
203, B.T. Road, Kolkata, India.}
\end{large}
\end{center}

\cleardoublepage
\pagenumbering{roman}
\pagestyle{plain}
\pagestyle{empty}
\vspace*{8cm}
\begin{center}
\em{In the memory of my grandmother\\}
\end{center}
\cleardoublepage
\vspace*{1cm}
\begin{center}
\begin{Large}
\textbf{Acknowledgements}\\
\end{Large}
\end{center}
\vspace{2cm}
\def\baselinestretch{1.2}

I would like to start by expressing my deepest gratitudes to my supervisor Debashish Goswami who introduced me to the theories of noncommutative geometry and compact quantum groups and whose mere presence acted as a psychological support to an extent unknown even to him. I thank each and every faculty member of Stat Math unit, ISI Kolkata from each of whom I learnt something or the other. I am grateful to Prof. Shuzhou Wang, whose valuable comments and suggestions helped me to have a better understanding about the contents of this thesis and also led to an improvement of the work. I thank the National Board for Higher Mathematics, India for providing me with partial financial support. I should also mention the names of Max Planck Institute fur Mathematics of Bonn and Chern Institute of mathematics of Nankai University for allowing me to attend two workshop and conference on Non commutative Geometry and Quantum Groups hosted by them, from where I had valuable exposures about the subjects. I thank my parents, sister, aunt and uncle for their continuous support during the time of working on this thesis. I would have been unable to fix some critical latex problem unless my friend Rajat Subhra Hazra  had spent his valuable time on it. I am grateful to Biswarup for the discussions I had with him on operator algebras. Many of my other friends including Koushik Saha, Abhijit da, Pusti da, Ashis da, Subhra, Subhajit and obviously Rajat were a continuous source of encouragement. I thank Subhajit for allowing me to use his room and providing me with a steady supply of 13 Tzameti et al. Lastly, I would like to mention the names of Subhra and Rajat again for bearing with my bhoyonkor and haayre during a critical part of my thesis work.

\cleardoublepage
\tableofcontents
\cleardoublepage

\addtolength{\headheight}{15pt}
\pagenumbering{arabic}
\setcounter{page}{1}
\pagestyle{myheadings}
\setcounter{chapter}{-1}
\pagestyle{fancy}
\renewcommand{\chaptermark}[1]{\markboth{\chaptername
\ \thechapter:\,\ #1}{}}
\renewcommand{\sectionmark}[1]{\markright{\ #1}}
\fancyhf{}
\fancyhead[LE]{\sl\leftmark}
\fancyhead[LO,RE]{\rm\thepage}
\fancyhead[RO]{\sl\rightmark}
\fancyfoot[C,L,E]{}

\addcontentsline{toc}{chapter} 
{\protect\numberline{Notations\hspace{-96pt}}} 
\cleardoublepage
\begin{Huge}
\textbf{Notations}
\end{Huge}

\begin{tabbing}
\hspace{2cm}\qquad\= \hspace{12cm} \\[1.2ex]

$ \IN $ \> The set of natural numbers \\

$ \IN_0 $\> $\IN \cup \{0 \} $\\

$ \IR $\>  The set of real numbers\\

$ \IC $ \>  The set of complex numbers\\

$ \clm_n ( \IC ) $\> The set of all $ n \times n $ complex matrices\\

$ S^1 $ \> The circle group\\

$ \IT^n $ \> The $ n $-torus\\

$ {\rm ev} $\>  Evaluation map\\

$ {\rm id}  $ \>  The identity map\\

$ l^2 ( \IN ) $ \>  The Hilbert space of square summable sequences\\

$ C^{\infty} ( M ) $ \>  The space of smooth functions on a smooth manifold $ M $\\

$ C^{\infty}_c ( M ) $ \>  The space of compactly supported smooth functions on $ M $\\

$ V_1 \otimes_{{\rm alg}} V_2 $ \>  Algebraic tensor product of two vector spaces $ V_1 $ and $ V_2 $\\

$ \cla \otimes \clb $ \>  Minimal tensor product of two $ C^* $ algebras $ \cla $ and $ \clb $\\

$ \clm ( \cla ) $ \>  The multiplier algebra of a $ C^* $ algebra $ \cla $\\

$ \cll ( E, F ) $ \> The space of adjointable maps from Hilbert modules $ E $ to $ F $\\

$ \cll ( E ) $ \> The space of adjointable maps from a Hilbert module $ E $ to itself\\

$ \clk ( E, F ) $ \> The space of compact operators from Hilbert modules $ E $ to $ F $\\

$ \clk ( E ) $ \> The space of compact operators from  a Hilbert module $ E $ to itself\\

$ \clb ( \clh ) $ \> The set of all bounded linear operators on a Hilbert space $ \clh $\\

$ \cla \ast \clb $ \> Free product of two $ C^* $ algebras $ \cla $ and $ \clb $\\

$ G \ast H $ \> Free product of two groups $ G $ and $ H $\\ 
 
$ G >\!\!\! \lhd H $ \>  Semi direct product of two groups $ G $ and $ H $\\   

\end{tabbing}

\cleardoublepage
\chapter{Introduction}

The theme of this thesis lies on the interface of two areas of the so called `` noncommutative mathematics '', namely noncommutative geometry (NCG) a la Connes, cf \cite{con} and the theory of ( $C^*$-algebraic ) compact quantum groups (CQG) a la Woronowicz, cf \cite{woro} which are generalizations of classical Riemannian spin geometry and that of compact topological groups respectively.

The root of NCG can be traced back to the Gelfand Naimark theorem which says that there is an anti-equivalence between the category of (locally) compact Hausdorff spaces and (proper, vanishing at infinity) continuous maps and the category of (not necessarily) unital $C^*$ algebras and $\ast$-homomorphisms. This means that the entire topological information of a locally compact Hausdorff space is encoded in the commutative $C^*$ algebra of continuous functions vanishing at infinity. This motivates one to view a possibly noncommutative $C^*$ algebra as the algebra of ``functions on some noncommutative space''. 

In classical Riemannian geometry on spin manifolds, the Dirac operator on the Hilbert space $L^2(S)$ of square integrable sections of the spinor bundle contains a lot of geometric information. For example, the metric, the volume form and the dimension of the manifold can be captured from the Dirac operator. This motivated Alain Connes to define his noncommutative geometry with the central object as the spectral triple which is a triplet $(\cla, \clh, D)$ where $\clh$ is a separable Hilbert space, $\cla$ is a (not necessarily closed) $\ast$-algebra of $\clb(\clh), ~ D $ is a self adjoint (typically unbounded) operator ( sometimes called the Dirac operator of the spectral triple ) such that $[D, a]$ admits a bounded extension. This generalizes the classical spectral triple $( C^{\infty}(M), L^{2}(S), D )$ on any Riemannian spin manifold $ M $, where $ D $ denotes the usual Dirac operator.

 On the other hand, quantum groups have their origin in different problems in mathematical physics as well as the theory of classical locally compact groups. It was S.L. Woronowicz, who in \cite{woro1} and \cite{woro} was able to pinpoint a set of axioms for defining compact quantum groups (CQG for short) as the correct generalization of compact topological groups.

The idea of a group acting on a space was extended to the idea of a CQG co-acting on a noncommutative space (that is, a possibly noncommutative $C^*$ algebra ). Following suggestions of Alain Connes, Shuzhou Wang in \cite{wang} defined and proved the existence of quantum automorphism groups on finite dimensional $C^*$ algebras. Since then, many interesting examples of such quantum groups, particularly the quantum permutation groups of finite sets and finite graphs, were extensively studied by a number of mathematicians (see for example \cite{ban1}, \cite{ban2}, \cite{ban3}, \cite{ban4}, \cite{ban5}, \cite{univ1} and references therein ). The underlying basic principle of defining a quantum automorphism group corresponding to some given mathematical structure (for example a finite set, a graph, a $C^*$ or von Neumann algebra) consists of two steps: first, to identify (if possible) the group of automorphism of the structure as a universal object in a suitable category, and then, try to look for the universal object in a similar but bigger category, replacing groups by quantum groups of appropriate type. However, most of the work done by them concerned some kind of quantum automorphism group of a `finite' structure,  for example, of finite sets or finite dimensional matrix algebras. It was thus quite natural to try to extend these ideas to the `infinite' or `continuous' mathematical structures, for example classical and noncommutative manifolds. With this motivation, Goswami ( \cite{goswami} )  formulated and studied the quantum analogues of the group of Riemannian isometries called the quantum isometry group. Classically, an isometry is characterized by the fact that its action commutes with the Laplacian. Therefore, to define the quantum isometry group, it is reasonable to consider a category of compact quantum groups which act on the manifold (or more generally on a noncommutative manifold given by a spectral triple) with the action commuting with the Laplacian, say $\cll,$ coming from the spectral triple. It is proven in \cite{goswami} that a universal object in the category ( denoted by $ {\bf Q^{\prime}_{\cll}} $ ) of such quantum groups does exist (denoted by $ QISO^{\cll}) $ if one makes some mild assumptions on the spectral triple all of which are valid for a compact Riemannian spin manifold. The work of this thesis starts with the computation of the quantum isometry group in several commutative and noncommutative examples (\cite{jyotish}, \cite{jyotish_torus}).

However, the formulation of quantum isometry groups in \cite{goswami} had a major drawback from the view point of noncommutative geometry, since it needed a `good' Laplacian to exist. In noncommutative geometry it is not always easy to verify such an assumption about the Laplacian, and thus it would be more appropriate to have a formulation in terms of the Dirac operator directly. This is what is done in \cite{goswami2} where the notion of a quantum group analogue of the group of orientation preserving isometries was given and its existence as the universal object in a suitable category was proved. Then, a number of computations for this quantum group were done in \cite{goswami2}, \cite{AF} and \cite{sphere}.

Now we try to give an idea of the contents of each of the chapters. In chapter 1, we discuss the concepts and results needed in the later chapters of the thesis. For the sake of completeness, we begin with a glimpse of operator algebras and Hilbert modules, free product and tensor products of $C^*$ algebras and some examples. The next section is on quantum groups which we start with the basics of Hopf algebras and then define compact quantum groups (CQG) and give relevant definitions and properties including a brief review of Peter Weyl theory. After that, we introduce the quantum groups $U_{\mu} (2), SU_{\mu} (2), A_{u} (Q) $ and $\clu_{\mu} (su(2)). $ We end this section by introducing the notion of a $C^*$ action of a compact quantum group on a $C^*$ algebra and giving an account of Shuzhou Wang's work in \cite{wang}. The next section is on Rieffel deformation where we recall a part of the work done in \cite{rieffel} and \cite{wang_def} which are relevant to us. We describe some important examples which are going to appear in chapter \ref{deformation}. The 4th section is on classical Riemannian geometry which includes, among other things, the definition and properties of Dirac operator which will serve as a motivation for the definition of ` spectral triple' in the 5th section. This section also contains a subsection on isometry groups of classical Riemannian manifold in which the characterizing property of an isometry, in the form given in \cite{goswami}, is given. In the 5th section, we define spectral triples, give examples of them, introduce the Hilbert space of forms (as in \cite{fro}), noncommutative volume form and the notion of Laplacian in noncommutative geometry. 

Chapter \ref{qisol} is on the Laplacian-based approach to quantum isometry groups as proposed in \cite{goswami}. Here we recall the formulation of quantum isometry groups from \cite{goswami} and then compute them for the space of continuous functions on the classical 2-spheres, the circle and the n-tori. In each of these cases, the quantum isometry group turns out to be the same as the classical ones, that is, $ C ( O ( 3 ) ), ~ C ( S^{1}  >\!\!\! \lhd Z_{2} ) $ and $ C ( \IT^n >\!\!\! \lhd ( \IZ^{n}_2 >\!\!\! \lhd S_n ) ) $ respectively ( $ S_n $ being the permutation group on $ n $ symbols ).

Chapter \ref{qorient} deals with the quantum group of orientation preserving isometries. The classical situation is stated clearly, which will serve as a motivation for the quantum formulation. Then, the quantum group of orientation-preserving isometries of an $ R $-twisted spectral triple is defined ( see \cite{goswami_rmp} for the definition of an $R$-twisted spectral triple ) and its existence is proven. Given an $R$-twisted spectral triple  $ ( \cla^{\infty}, \clh, D ) $ of compact type, we consider a category $ \bf Q^{\prime}$ of pairs $ ( \clq, U ) $ where $ \clq $ is a  compact quantum group which has  a unitary (co)-representation $ U $ on $ \clh $ commuting with $ D, $ and such that for all state $ \phi $ on $ \clq, ~ ( {\rm id} \otimes \phi ) ad_{U} $ maps $ \cla^{\infty} $ inside $ \cla^{\prime \prime}_{\infty} .$ Moreover, let $ {\bf Q^{\prime}_{R}} $ be a subcategory of $ {\bf Q^{\prime}} $ consisting of those $  ( \clq, U ) $  for which $ ad_U $ preserves the $ R $-twisted volume form. In section \ref{qorient_section_defn_existence_qtm_group_oreintation_pres_isometries}, we have proved that  $ {\bf Q^{\prime}}_{R}$ has a universal object to be denoted by $ \widetilde{QISO^+_R} ( D ).$ The Woronowicz $ C^{\ast} $ subalgebra of $ \widetilde{QISO^+_R} ( D ) $ generated by elements of the form $ \left\langle ad_{U} ( a ) ( \eta \otimes 1 ) , ~ \eta^{\prime} \otimes 1 \right\rangle_{\widetilde{QISO^+_R} ( D )} $ where $ \eta, \eta^{\prime} $ are in $ \clh,~ a $ is in $ \cla^{\infty} $ and $ \left\langle . ~,~ . \right\rangle_{\widetilde{QISO^+_R} ( D )} $ denotes the $ \widetilde{QISO^+_R} ( D ) $ valued inner product of $ \clh \otimes \widetilde{QISO^+_R} ( D ) $ is defined to be the quantum group of orientation and volume preserving isometries of the spectral triple $ ( \cla^{\infty}, \clh, D )$ and is denoted by $ QISO^+_R ( D ).$  The next section explores the conditions under which the action of this compact quantum group keeps the $ C^* $ algebra invariant and is a $ C^* $ action. Moreover, we have given some sufficient conditions under which the universal object in the bigger category $ {\bf Q^{\prime}} $ exists which is denoted by $ \widetilde{ QISO^{+} ( D )}$ and the corresponding Woronowicz $C^*$ subalgebra as above is denoted by $QISO^{+} ( D ).$   After this, we compare this approach with the Laplacian-based approach in \cite{goswami}. We obtain the following results:

( 1 ) Under some reasonable conditions, $ QISO^{+}_{I} ( D )  $  is a sub-object of $ QISO^{\cll} $ in the category  $ {\bf Q^{\prime}_{\cll}}.$

( 2 ) $ QISO^{\cll} $ is isomorphic to a quantum subgroup of $ QISO^{+}_{I} $ corresponding to the Hodge Dirac operator coming from $ D. $

( 3 ) Moreover, under some conditions which are valid for compact spin manifolds, $ QISO^{\cll} $ and the $ QISO^+_{I} $ of the Hodge Dirac operator are isomorphic.

The next section is on examples and computations. To begin with, the compact quantum group $ U_{\mu} ( 2 ) $ is identified as the  $   QISO^+  $ of $ SU_{\mu} ( 2 )  $ corresponding to the spectral triple constructed by Chakraborty and Pal in \cite{partha}. Then we derive that $ QISO^{+} $ for the classical spectral triple on $ C ( \IT^2 ) $ is $ C ( \IT^2 ) $ itself. We end the chapter by showing that $ QISO^+ $ of spectral triples associated with some approximately finite dimensional $ C^* $ algebras arise as the inductive limit of $ QISO^+ $ of the constituent finite dimensional algebras. The results of this chapter are taken from \cite{goswami2} and \cite{AF}.

Chapter \ref{deformation} is about the $ QISO^{\cll} $ and $ QISO^{+}_R $ of a Rieffel deformed noncommutative manifold. We first discuss the isospectral deformation of a spectral triple, followed by the proof of some preparatory technical results which will be needed later. Then in the final section we prove that $ QISO^{+}_{R} $ and $ QISO^{\cll} $ of a Rieffel deformed ( noncommutative ) manifold is a Rieffel deformation of the $ QISO^{+}_R $ and $ QISO^{\cll} $ ( respectively ) of the original ( undeformed ) manifold.

In chapter \ref{quantumsphere}, we compute the quantum group of orientation preserving isometries for two different families of spectral triples on the Podles spheres, one constructed by Dabrowski et al in \cite{Dabrowski_et_al} and the other by Chakraborty and Pal in \cite{chak_pal}. We start by giving the different descriptions of the Podles spheres ( as in \cite{podles}, \cite{Dabrowski_et_al}, \cite{klimyk}, and \cite{Schmudgen_wagner_crossproduct} ) and the formula for the Haar functional on it. Then we introduce the spectral triples on the Podles spheres as in \cite{Dabrowski_et_al} and show that it is indeed $ SU_{\mu} ( 2 ) $ equivariant and R-twisted ( for a suitable $ R $ ). After this, the compact quantum group $ SO_{\mu} ( 3 ) $ is defined and its action on the Podles sphere is discussed. In the 3rd section, the computation for identifying $ SO_{\mu} ( 3 ) $ as $ QISO^+_R $ for this spectral triple is given. In the 4th section, the spectral triple defined in \cite{chak_pal} is introduced and then the corresponding $ QISO^{+} $ is computed. In particular, it follows that $ QISO^+ $ in general may not be a matrix quantum group and that it may not have a $ C^* $ action.

\cleardoublepage
\chapter{Preliminaries}

\label{preliminaries} 
 
 \section{Operator algebras and Hilbert modules}

 \label{preliminaries_section_OA_&_HM}
 
 We presume the reader's familiarity with the theory of operator algebras and Hilbert modules. However, for the sake of completeness, we give a sketchy review of some basic definitions and facts and refer to \cite{takesaki}, \cite{lance} for the details. Throughout this thesis, all algebras will be over $\IC$ unless otherwise mentioned.
 
\subsection{$C^*$ algebras}

\label{preliminaries_subsection_C*algebras}

A $ C^* $ algebra $ \cla $ is a Banach $ \ast $-algebra satisfying the $ C^* $ property : $ \left\| x^* x \right\| = \left\| x \right\|^2 $ for all $ x $ in $ \cla. $ The algebra $\cla$ is said to be unital or non-unital depending on whether it has an identity or not. Every commutative $C^*$ algebra $\cla$ is isometrically isomorphic to the $C^*$ algebra $C_0 ( X )$ consisting of complex valued functions on a locally compact Hausdorff space $X$ vanishing at infinity( Gelfand's theorem ). An arbitrary ( possibly noncommutative   ) $ C^* $ algebra is isometrically isomorphic to a $ C^* $-subalgebra of $ \clb ( \clh ), $ the set of all bounded operators on a  Hilbert space $ \clh. $

For $x$ in $\cla,$ the {\bf spectrum} of $x,$ denoted by $\sigma ( x )  $ is defined as the complement of the set $ \{ z \in \IC: {( z1 - x  )}^{- 1} \in \cla \}.$ An element $ x $ in $ \cla $ is called {\bf self adjoint} if $ x = x^*, $ {\bf normal} if $x^* x = x x^*,$ {\bf unitary} if $ x^* = x^{- 1} ,$ {\bf projection} if $ x = x^* = x^2 $ and {\bf positive} if $ x = y^* y $ for some $ y $ in $ \cla .$  When $x$ is normal, there  is a continuous functional calculus sending $f$ in $C ( \sigma (x) )$ to $ f (x) $ in $\cla.$ where $f \mapsto f (x) $ is a $\ast$ isometric isomorphism from $C ( \sigma (x) ) $ onto $C^* (x).$ 

A linear map between two $ C^* $ algebras is said to be {\bf positive} if it maps positive elements to positive elements.  A positive linear functional $\phi$ such that $\phi (1) = 1$ is called a {\bf state} on $\cla.$  A state $\phi$ is called a {\bf trace} if $\phi ( ab) = \phi (b a) $ for all $a,b  $ in $\cla$ and {\bf faithful} if $\phi (x^* x) = 0 $ implies $x = 0.$ Given a state $\phi$ on a $C^*$ algebra $\cla,$ there exists a triple (called the GNS triple) $( \clh_{\phi}, \pi_{\phi}, \xi_{\phi} )$ consisting of a Hilbert space $\clh_{\phi},$ a $\ast$ representation $\pi_{\phi} $ of $\cla$ into $\clb (\clh_{\phi})$ and a vector $\xi_{\phi}$ in $\clh_{\phi}$ which is cyclic in the sense that $\{ \pi_{\phi} (x) \xi_{\phi} : x \in \cla \}$ is total in $\clh_{\phi}$ satisfying 
$$ \phi (x) = \left\langle  \xi_{\phi}, ~ \pi_{\phi} (x) \xi_{\phi} \right\rangle. $$

For a two-sided norm closed ideal $ \cli $ of a $ C^* $ algebra $ \cla, $ the canonical quotient norm on the Banach space $ \cla/\cli $ is in fact the unique $ C^* $ norm making $ \cla /\cli $ into a $ C^* $ algebra. Here we prove two results which we are going to need later on.

\blmma

\label{preliminary_lemma_goswami_paper}

\cite{goswami} ~ Let $ \clc $ be a  $ C^* $ algebra and $ \clf $ be a nonempty collection of $ C^* $-ideals ( closed two-sided ideals ) of $ \clc. $ Then for any $ x $ in $ \clc, $ we have 
  $$ {\rm sup}_{I \in \clf} \left\| x + I \right\| = \left\| x + I_0 \right\| ,$$  
  where  $ I_0 $ denotes the intersection of all $ I $ in $ \clf $ and $ \left\| x + I \right\| = {\inf} \{ \left\| x - y \right\| : y \in I \} $ denotes the norm in  $ \clc/ I.$    

\elmma

{\it Proof :} It is clear that $ {\rm sup}_{I \in \clf} \left\| x + I \right\| $ defines a norm on $ \clc/ I_0,$ which is in fact a $ C^* $ norm since each of the quotient norms $ \left\| . + I \right\| $ is so. Thus the lemma follows from the uniqueness of $ C^* $ norm on the $ C^* $ algebra $ \clc / I_0. $ \qed 

\blmma 
  
  \label{lim} 
  
  Let $\clc$ be a unital $C^*$ algebra and $\clf$ be a nonempty collection of $C^*$-ideals (closed two-sided ideals) of
$\clc$. Let  $\cli_0$ denote the intersection of all $\cli$ in $\clf$, and  let $\rho_\cli$ denote the map $\clc /\cli_0 \ni x+\cli_0 \mapsto x+I  \in \clc/\cli$ for $\cli $ in $ \clf$.  Denote by $\Omega$  the set  $\{ \omega \circ \rho_\cli, \cli \in \clf,~\omega ~{\rm state}~{\rm on}~ \clc/\cli \}$, and let $K$ be the weak-$\ast$ closure of the convex hull of $\Omega \bigcup (- \Omega)$. Then $K$ coincides with the set of bounded linear functionals $\omega$ on $\clc/\cli_0$ satisfying $\| \omega \|=1$ and $\omega(x^*+ \cli_0)=\overline{\omega(x+ \cli_0)}$.  

\elmma

 {\it Proof :}
We will use Lemma \ref{preliminary_lemma_goswami_paper}. Clearly, $K$ is a weak-$\ast$ compact, convex subset of the unit ball $(\clc/\cli_0)^*_1$ of the dual of $\clc/\cli_0$, satisfying $-K=K$.  If $K$ is strictly smaller than the self-adjoint part of unit ball of the dual of $\clc/\cli_0$, we can find a state $\omega$ on $\clc/\cli_0$ which is not in $K$. Considering the real Banach space $X=(\clc/\cli_0)^*_{\rm s.a.}$ and using standard separation theorems for real Banach spaces (for example, Theorem 3.4 of \cite{rudin}, page 58), we can find a self-adjoint element $x$ of $\clc$ such that $\| x+\cli_0 \|=1,$ and   $$ \sup_{\omega^\prime \in K} \omega^\prime(x+\cli_0) < \omega(x+\cli_0).$$ Let $\gamma $ belonging to $ \IR$ be such that $\sup_{\omega^\prime \in K} \omega^\prime(x+\cli_0)< \gamma  < \omega(x+\cli_0).$ Fix  $0<\epsilon<\omega(x+\cli_0)-\gamma$, and let $\cli $ be an element of $ \clf$ be such that $\| x+\cli_0\| -\frac{\epsilon}{2} \leq \| x+\cli \| \leq \| x+\cli_0\|.$ Let $\phi$ be a state on $\clc/\cli$ such that $\| x+ \cli \|=| \phi(x+\cli)|.$ Since $x$ is self-adjoint, either $\phi(x+\cli)$ or $-\phi(x+\cli)$ equals $\| x+\cli \|$, and $\phi^\prime:=\pm \phi \circ \rho_\cli $, where the sign is chosen so that $\phi^\prime(x+\cli_0)=\| x+\cli\|$. Thus, $\phi^\prime $ is in $ K$, so $\| x+\cli_0\|=\phi^\prime(x+ \cli) \leq \gamma < \omega(x+\cli_0)-\epsilon.$ But this implies $\|x+\cli_0\| \leq \| x+\cli \| +\frac{\epsilon}{2}<\omega(x+\cli_0)-\frac{\epsilon}{2}\leq \| x+\cli_0\| -\epsilon$ (since $\omega$ is a state), which is a contradiction completing the proof.
 \qed 
 
 \vspace{4mm}

For a $C^*$ algebra $\cla$ ( possibly non unital ), its multiplier algebra, denoted by $\clm ( \cla ),$ is defined as the maximal $C^*$ algebra which contains $\cla$ as an essential two sided ideal, that is, $\cla$ is an ideal in $\clm ( \cla )$ and for $ y $ in $\clm ( \cla), ~ y a = 0 $ for all $ a $ in $\cla$ implies $y = 0.$  The norm of $\clm ( \cla )$ is given by $ \left\| x \right\| = {\rm sup}_{a \in \cla, \left\| a \right\| \leq 1} \{ \left\| xa \right\|,~ \left\| ax \right\| \}.$ There is a locally convex topology called the strict topology on $\clm (\cla),$ which is given by the family of seminorms $ \{ \left\| . \right\|_a, a \in \cla  \},$ where $ \left\| x \right\|_a = {\rm Max} ( \left\| x a \right\|, \left\| a x \right\| )  ,$ for $ x $ in $ \clm ( \cla ).$ $ \clm ( \cla ) $ is the completion of $\cla$ in the strict topology.

\vspace{4mm}

We now come to the {\bf inductive limit of $ C^* $ algebras
.} Let $ I $ be a directed set and $ \{ \cla_i \}_{i \in I} $ be a family of $ C^* $ algebras equipped with a family of $ C^* $  homomorphisms $ \Phi_{ij} : \cla_j \rightarrow \cla_i $ ( when $ j < i $ ) such that $ \Phi_{ij} = \Phi_{ik} \Phi_{kj} $ when $ j <  k < i .$ Then there exists a unique $ C^* $ algebra denoted by $ {\rm lim}_{i} \cla_i $ and $ C^* $ homomorphisms $ \phi_i : \cla_i \rightarrow {\rm lim}_{i} \cla_i $ with the universal property that given any other $ C^* $ algebra $ \cla^{\prime}$ and $ C^* $ homomorphisms $ \psi_i : \cla_i \rightarrow \cla^{\prime}$ satisfying $ \psi_j = \psi_i \Phi_{ij}$ for $ j < i, $ then there exists unique $ C^* $ homomorphism  $ \chi : {\rm lim}_{i} \cla_i \rightarrow \cla^{\prime}  $ satisfying $ \chi \phi_i = \psi_i.~ {\rm lim}_{i} \cla_i $ is called the inductive limit $ C^* $ algebra corresponding to the inductive system $ ( \cla_i, \Phi_{ij} ).$  Inductive limit of a sequence of finite dimensional $C^*$ algebras are called approximately finite dimensional $C^*$ algebras or AF algebras.

\vspace{4mm}

A large class of $ C^* $ algebras are obtained by the following construction. Let $ \cla_0 $ be an associative $ \ast $-algebra without any a-priori norm such that the set $ \clf = \{ \pi : \cla_0 \rightarrow \clb ( \clh_{\pi} ) ~ \ast{\rm - homomorphism,}~ \clh_{\pi} {\rm ~ a ~ Hilbert ~ space} \} $ is non empty and $ \| . \|_{u} $  defined by $  \| a \|_{u} = {\rm sup} \{ \| \pi (a ) \| : \pi \in \clf \} $ is finite for all $a.$  Then the completion of $ \cla_0 $ in  $ \| . \|_{u} $ is a $ C^* $ algebra known as the universal $C^*$ algebra corresponding to $\cla_0.$

\vspace{2mm}

{\bf Example 1:}  {\bf Noncommutative two-torus}

Let $\theta$ belongs to $[0,1].$ Consider the $\ast$ algebra $\cla_0$ generated by two unitary symbols $U$ and $V$ satisfying the relation $UV = e^{2 \pi i \theta} VU.$  It has a representation $ \pi $ on  the Hilbert space $L^2 ( S^{1} ) $ defined by $\pi ( U ) (f) (z) = f ( e^{2 \pi i \theta} z ),~ \pi ( V ) (f) (z) = z f(z) $ where $f$ is in $L^2 (S^{1}), ~ z  $ is in $S^1.$ Then $\| a \|_u $ is finite for all $a$ in $\cla_0.$ The resulting $C^*$ algebra is called noncommutative two-torus and denoted by $\cla_{\theta}.$

\vspace{2mm}

{\bf Example 2:} {\bf Group $C^*$ algebra }

Let $ G $ be a locally compact group with left Haar measure $ \mu .$ One can make $ L^{1} ( G ) $ into a Banach $ \ast $-algebra by defining 
$$ ( f \ast g ) ( t ) = \int_{G} f ( s )  g ( s^{- 1} t ) d \mu ( s ), $$

$$ f^* ( t ) = {\Delta ( t )}^{- 1} \overline{f ( t^{- 1} ) } .$$

Here $ f, ~ g $ are in  $ L^{1} ( G ), ~ \Delta $ is the modular homomorphism of $ G. $

$ L^{1} ( G ) $ has a distinguished representation $ \pi_{reg} $ on $ L^{2} ( G ) $ defined by $ \pi_{reg} ( f ) = \int f ( t ) \pi ( t ) d \mu ( t ) $ where $ \pi ( t ) $ is a unitary operator on $ L^{2} ( G ) $  defined by $ ( \pi ( t ) f ) ( g ) = f ( t^{- 1} g ) ~ ( f \in L^{2} ( G ), ~ t,g \in G  ) .$  The reduced group $ C^* $ algebra of $ G $ is defined to be $ C^*_{r} ( G ): = {\overline{\pi_{reg} ( L^{1} ( G ) )}}^{\clb ( L^{2} ( G ) )}.$

\brmrk

\label{preliminaries_reduced_dual}

 For $ G $ abelian, we have $ C^*_{r} ( G ) \cong C_0 ( \widehat{G} ) $ where $ \widehat{G}  $ is the group of characters on $ G. $
 
 \ermrk

One can also consider the universal $ C^* $ algebra described before corresponding to the Banach $ \ast $-algebra $ L^{1} ( G ).$ This is called the free or full group $ C^* $ algebra and denoted by $ C^* ( G ). $ 
 
 \brmrk
 
 \label{preliminaries_amenable_reduced=full}
 
 For the so-called amenable groups ( which include compact and abelian groups ) we have $ C^* ( G ) \cong C^*_r ( G ). $ 
 
 \ermrk

\subsection{von Neumann algebras}  

\label{preliminaries_subsection_vN_algebras}  

We recall that for a Hilbert space $ \clh, $ the {\bf strong operator topology ( SOT ) }, the {\bf weak operator topology ( WOT ) } and the {\bf ultra weak topology} are the locally convex topologies on $ \clb ( \clh ) $ given by  families of seminorms  $ \clf_1, ~ \clf_2, ~ \clf_3 $ respectively  where  $ \clf_1 =  \{ p_{\xi} : \xi \in \clh \}, ~ \clf_2 =  \{ p_{\xi, \eta} : \xi, \eta \in \clh \}, ~ \clf_3 = \{ p_{\rho} :  \rho ~ {\rm is ~ a ~ trace ~ class ~ operator ~ on} ~ \clh \}  $ and  $ p_{\xi} ( x ) = \left\| x \xi \right\| ,~ p_{\xi, \eta} ( x ) = \left| \left\langle x \xi , \eta \right\rangle \right|, ~ p_{\rho} ( x ) = \left| {\rm Tr} ( x \rho ) \right| $ ( where $ {\rm Tr} $ denotes the usual trace on $ \clb ( \clh ) ) .$   



Now we state a well known fact.

\blmma

\label{preliminaries_SOT_convergence_Trace}

If $ T_n $ is a sequence of bounded operators converging to zero in SOT, then for any trace class operator $ W, ~ Tr ( T_n W ) \rightarrow 0 $ as $ n \rightarrow \infty.$

\elmma

       \vspace{4mm}

For any subset $\clb$ of $ \clb ( \clh ), $ we denote by $\clb^{\prime}$ the commutant of $\clb,$ that is, $\clb^{\prime} = \{ x \in \clb ( \clh ) : xb = bx {\rm ~ for ~ all} ~ b \in \clb \} .$ A unital $C^*$ subalgebra $\cla \subseteq \clb ( \clh )  $ is called a {\bf von Neumann algebra} if $\cla = \cla^{\prime \prime}$ which is equivalent to being closed in any of the three topologies mentioned above.

A state $ \phi $ on a von Neumann algebra $ \cla $ is called normal if $ \phi ( x_{\alpha} ) $ increases to $ \phi ( x ) $ whenever $ x_{\alpha} $ increases to $ x. $ A state $ \phi $ on $ \cla $ is normal if and only if there is a trace class operator $ \rho $ on $ \clh $ such that $ \phi ( x ) = {\rm Tr } ( \rho x ) $ for all $ x $ in  $ \cla. $  More generally, we call a linear map $ \Phi : \cla \rightarrow \clb $ ( where $ \clb $ is a von Neumann algebra ) normal if whenever $ x_{\alpha}  $ increases to $ x $ for a net $ x_{\alpha}  $ of positive elements from $ \cla, $ one has that $ \Phi ( x_{\alpha} ) $ increases to $ \Phi ( x ) $ in $ \clb. $ It is known that a positive linear map is normal if and only if it is continuous with respect to the ultra-weak-topology. In view of this fact, we shall say that a bounded linear map between two von Neumann algebras is normal if it is continuous with respect to the respective ultra-weak topologies.

\subsection{Free product and tensor product}

\label{preliminaries_subsection_free_product_tensor_product}

If $(\cla_i)_{i \in I}$ is a family of unital $C^*$ algebras, then their unital $C^*$ algebra free product $\ast_{i \in I}A_i$ is the unique $C^*$ algebra $\cla$ together with unital $\ast$-homomorphism $\psi_i:\cla_i \rightarrow \cla $ such that given any unital $C^*$ algebra $\clb$ and unital $\ast$-homomorphisms $\phi_i: \cla_i \rightarrow \clb$ there exists a unique unital $\ast$-homomorphism $\Phi : \cla \rightarrow \clb $ such that $\phi_i = \Phi \circ \psi_i.$ 

\brmrk

\label{preliminaries_free_product_of_homomorphisms}

It is a direct consequence of the above definition that given a family of $ C^* $ homomorphisms $ \phi_i $ from $ \cla_i $ to $ \clb,$ there exists a $ C^* $ homomorphism $ \ast_i \phi_i $ such that $ ( \ast_i \phi_i ) \circ \psi_i = \phi_i $ for all $ i.$

\ermrk

\brmrk

\label{preliminaries_free_product_C*_distributive}

We recall that for discrete groups $ \{ G_i \}_{i \in I }, ~ C^* ( \ast_{i \in I} G_i ) \cong \ast_{i \in I} C^* ( G_i ).$

\ermrk









 For $\cla$ and $\clb$ two algebras, we will denote the algebraic tensor product of $\cla$ and $\clb$ by the symbol $\cla \otimes_{alg} \clb.$  When $\cla$ and $\clb$ are $C^*$ algebras, there is more than one norm on $\cla \otimes_{alg} \clb$ so that the completion with respect to that norm is a $C^*$ algebra. Throughout this thesis, we will work with the so called injective tensor product, that is, the completion of $\cla \otimes_{alg} \clb $ with respect to the norm given on  $\cla \otimes_{alg} \clb $  by $ \| \sum^n_{i = 1} a_i \otimes b_i \| = {\rm sup} \| \sum^n_{i = 1} \pi_1 (a_i) \otimes \pi_2 (b_i) \|_{\clb ( \clh_1  \otimes \clh_2 )}  $ where $a_i$ is in $\cla, ~  b_i$ is in  $\clb$ and the supremum runs over all possible choices of $ ( \pi_1, \clh_1 ), ~ ( \pi_2, \clh_2 ) $ where $ \clh_1, \clh_2 $ are Hilbert spaces and  $ \pi_1 : \cla \rightarrow \clb ( \clh_1 ) $ and $ \pi_2 : \cla_2 \rightarrow \clb ( \clh_2 ) $ are $ \ast $-homomorphisms. When $ \cla \subseteq \clb ( \clh_1 ), ~ \clb \subseteq \clb ( \clh_2 ) $ are von Neumann algebras, then by the notation $ \cla \otimes \clb, $ we mean the von Neumann algebra tensor product, that is, the WOT closure of $ \cla \otimes_{\rm{alg}} \clb $ in $ \clb ( \clh_1 \otimes \clh_2 ).$  We refer to \cite{takesaki} for more details.
   
We now prove a useful general fact.
 
 \blmma

\label{classical_case_family_ofmaps2}

Let $\cla$ be a $C^*$ algebra and $\omega, \omega_j $ ($j=1,2,...$) be states on $\cla$ such that $\omega_j \raro \omega$ in the weak-$\ast$ topology of $\cla^*$. Then for any separable Hilbert space $\clh$ and 
for all $  Y $ in $ \clm ( \clk ( \clh ) \otimes \cla ),$ we have  $ ({\rm  id}  \otimes \omega_j ) (Y) \rightarrow ({\rm  id} \otimes \omega ) (Y) $ in the strong operator topology.

\elmma

{\it Proof:} Clearly,  $({\rm id} \ot \omega_j)(Y) \raro ({\rm id} \ot \omega)(Y)$ (in the strong operator topology)  for all $Y $ in $ {\rm Fin}(\clh) \ot_{\rm alg} \cla$, where ${\rm Fin}(\clh)$ denotes the set of finite rank operators on $\clh$.  Using the strict density of ${\rm Fin}(\clh) \ot_{\rm alg} \cla $ in $\clm(\clk(\clh) \ot \cla)$, we choose, for a given $Y $ in $  \clm(\clk(\clh) \ot \cla)$, $\xi $ in $ \clh$ with $\| \xi \| =1$,  and $\delta>0$, an element $Y_0 $ in $ {\rm Fin}(\clh) \ot_{\rm alg} \cla$ such that $\|(Y-Y_0) (|\xi><\xi| \ot 1) \| < \delta$. Thus, 
\bean \lefteqn{\|({\rm id} \ot \omega_j)(Y) \xi- ({\rm id} \ot \omega)(Y)\xi \|}\\
&=& \|({\rm id} \ot \omega_j)(Y (|\xi><\xi| \ot 1))\xi- ({\rm id} \ot \omega)(Y ( |\xi><\xi| \ot 1))\xi \|\\
&\leq & \|({\rm id} \ot \omega_j)(Y_0 (|\xi><\xi| \ot 1))\xi- ({\rm id} \ot \omega)(Y_0 ( |\xi><\xi| \ot 1))\xi \|\\
&+& 2  \| (Y-Y_0)(|\xi><\xi| \ot 1)\|\\
&\leq & \|({\rm id} \ot \omega_j)(Y_0 (|\xi><\xi| \ot 1))\xi- ({\rm id} \ot \omega)(Y_0 ( |\xi><\xi| \ot 1))\xi \|
+2  \delta, \eean  from which it follows that $({\rm id} \ot \omega_j)(Y) \raro ({\rm id} \ot \omega)(Y)$ in the strong operator topology. \qed

\vspace{4mm}

Let $ \cla $ and $ \clb $ be two unital $ \ast $-algebras. Then a linear map $ T $ from $ \cla $ to $ \clb $ is called $ n$-positive if $ T \otimes {\rm Id_n}: \cla \otimes \clm_{n} ( \IC ) \rightarrow \clb \otimes \clm_{n} ( \IC ) $ is positive for all $ k \leq n $ but not necessarily for $ k = n + 1. $ $ T $ is said to be completely positive ( CP for short ) if it is n-positive for all $ n. $ It is a well known result that for a CP map $ T : \cla \rightarrow \clb ( \clh ) ,$ one has the following operator inequality for all $a$ in $ \cla $: 
\be \label{preliminaries_CP_maps_result1} {T ( a )}^* T ( a ) \leq \left\| T ( 1 ) \right\| T ( a^* a )  .\ee 
Tensor product and composition of two CP maps are CP. The following is an useful result about CP maps.
\bppsn

\label{preliminaries_CP_maps_result2}

If $A$ and $B$ are $C^*$ algebras with $A$ commutative, $\phi $ is a positive map from $A$ to $B,$ then $\phi$ is CP. The same holds if $\phi$ is from $B$ to $A.$

\eppsn

\subsection{Hilbert modules}

\label{preliminaries_subsection_Hilbert_C*_modules}

Given a $ \ast $-subalgebra $ \cla \subseteq \clb ( \clh ) $ ( where $ \clh $ is a Hilbert space ), a semi-Hilbert $ \cla $ module $ E $ is a right $ \cla $-module equipped with a sesquilinear map $ \left\langle . ~, ~ . \right\rangle : E \times E \rightarrow \cla $ satisfying $ {\left\langle x, y \right\rangle}^* = \left\langle y,x \right\rangle, ~ \left\langle x, y a \right\rangle = \left\langle x,y \right\rangle a $ and $ \left\langle  x, x \right\rangle \geq 0 $ for $ x,y $ in $ E $ and $ a  $ in $ \cla. $ A semi Hilbert module is called a pre-Hilbert module if $ \left\langle x, x \right\rangle = 0  $ if and only if $ x = 0. $ It is called a Hilbert module if furthermore, $ \cla $ is a $ C^* $ algebra and $ E $ is complete in the norm $ x \mapsto \left\| \left\langle x,x \right\rangle \right\|^{\frac{1}{2}} $ where  $ \left\| . \right\| $ is the $ C^* $ norm of $ \cla. $ 

The simplest examples of Hilbert $ \cla $ modules  are the so called trivial $ \cla $ modules of the form $ \clh \otimes \cla $ where $ \clh $ is a Hilbert space with an $\cla$ valued sesquilinear form defined on $\clh \otimes_{{\rm alg}} \cla$ by :  $ \left\langle  \xi \otimes a, ~ \xi^{\prime} \otimes a^{\prime} \right\rangle = \left\langle \xi, \xi^{\prime} \right\rangle a^* a^{\prime}.$ The completion of  $\clh \otimes_{{\rm alg}} \cla$ with respect to this pre Hilbert module structure is a Hilbert $ \cla $ module and is denoted by $\clh \otimes \cla.$

We recall that for a pre Hilbert $ \cla $ module $ E $ ( $ \cla $ is a $ C^* $ algebra ), the Cauchy Schwarz inequality holds in the following form: $ 0 \leq \left\langle x, y \right\rangle \left\langle y, x \right\rangle \leq \left\langle x,x \right\rangle \left\| \left\langle y,y \right\rangle \right\|. $  

Let $ E  $ and $  F $ be two Hilbert $ \cla $ modules. We say that a $ \IC $ linear map $ L $ from $ E $ to $ F $ is adjointable if there exists a $ \IC $ linear map $ L^* $ from $ F $ to $ E $ such that $ \left\langle L ( x ), y \right\rangle = \left\langle x, L^* ( y ) \right\rangle $ for all $ x $ in $ E, ~ y $ in $ F. $ We call $ L^* $ the adjoint of $ L. $ The set of all adjointable maps from $ E $ to $ F $ is denoted by $ \cll ( E, F ). $ In case, $ E = F, $ we write $ \cll ( E )  $ for $ \cll ( E, E ). $ For an adjointable map  $ L, $ both $ L $ and $ L^* $ are automatically $ \cla $-linear and norm bounded maps between Banach spaces.  We say that an element $ L $ in $ \cll ( E, F ) $ is an isometry if $ \left\langle L x, L y \right\rangle = \left\langle x, y \right\rangle $ for all $ x, y $ in $ E. ~ L $ is said to be a unitary if $ L $ is isometry and its range is the whole of $ F. $ One defines a norm on  $ \cll ( E, F ) $ by $ \left\| L \right\| = {\rm sup}_{x \in R, ~ \left\| x \right\| \leq 1 } \left\| L ( x ) \right\|. ~ \cll ( E ) $ is a $ C^* $ algebra with this norm.

There is a topology on $ \cll ( E, F ) $ given by a family of seminorms $ \{ \left\| . \right\|_{x}, \left\| . \right\|_{y} : x \in E, ~ y \in F \} $ ( where $ \left\| t \right\|_{x} = \left\| {\left\langle  tx, tx \right\rangle}^{\frac{1}{2}} \right\| $ and $ \left\| t \right\|_{y} = \left\| {\left\langle  t^* y, t^* y \right\rangle}^{\frac{1}{2}} \right\|  $ )  known as the strict topology. For $ x $ in $ E, ~ y $ in $ F, $ we denote by $ \theta_{x,y} $ the element of $ \cll ( E, F ) $ defined by $ \theta_{x,y} ( z ) = y \left\langle x, z \right\rangle $  ( where  $ z $ is in $ E $ ). The norm closure ( in $ \cll ( E, F  ) $ ) of the $ \cla $ linear span of $ \{ \theta_{x,y} : x \in E, ~ y \in F \} $ is called the set of compact operators and denoted by $ \clk ( E, F ) $  and we denote $ \clk ( E, E )  $ by $ \clk ( E ).$ These are not necessarily compact in the sense of compact operators between two Banach spaces. One has the following important result:

\bppsn 

The multiplier algebra $ \clm ( \clk ( E ) )  $ of $ \clk ( E ) $ is isomorphic with $ \cll ( E ) $ for any Hilbert module $ E. $

\eppsn                                    

Using this, for a possibly non-unital $ C^* $ algebra $ \clb, $  we often identify an element $ V $ of $ \clm ( \clk ( \clh ) \otimes \clb ) $ with the map from $ \clh $ to $ \clh \otimes \clb $ which sends a vector $ \xi $ of $ \clh $ to $ V ( \xi \otimes 1_{\clb} ) $ in $ \clh \otimes \clb. $

 Given a Hilbert space $ \clh $ and a $ C^* $ algebra $ \cla, $ and a unitary element $U$ of $ \clm ( \clk ( \clh ) \otimes \cla ), $ we shall denote by $\alpha_U ( \equiv ad_{U} ) $ the $\ast$-homomorphism $\alpha_U(X)= \widetilde{U}(X \ot 1){\widetilde{U}}^*$ for $X $ belonging to $ \clb(\clh)$. For a  not necessarily bounded, densely defined (in the weak operator topology)  linear functional $\tau$ on $\clb(\clh)$,  we say that $\alpha_U$ preserves $\tau$ if $\alpha_U$ maps a suitable (weakly) dense $\ast$-subalgebra   (say $\cld$) in the domain of $\tau$ into $\cld \ot_{\rm alg} \cla$ and $( \tau \ot {\rm id}) (\alpha_U(x))=\tau(x).1_\cla$  for all $ x $ in $ \cld$. When $\tau$ is bounded and normal, this is equivalent to $(\tau \ot {\rm id}) (\alpha_U(x))=\tau(x) 1_\cla$ for all $x$ belonging to $ \clb(\clh)$.

We say that a (possibly unbounded) operator $T$ on $\clh$ commutes with $U$ if $T \ot I$ (with the natural domain) commutes with $\widetilde{U}$. Sometimes such an operator will be called $U$-equivariant.

 \section{Quantum Groups}

\label{preliminaries_section_quantum_groups}
 
 In this section, we will recall the basics of Hopf algebras and then define compact quantum groups ( as in  \cite{woro}, \cite{woro1} ). After that, we will discuss a few examples of quantum groups and the concept of an action of a  compact quantum group on a $ C^* $ algebra. For more detailed discussion, we refer to \cite{klimyk}, \cite{Kassel}, \cite{chari} \cite{masuda_nakagami_wor} and references therein. In this thesis, we will be concerned about compact quantum groups only. For other types of quantum groups, we refer to \cite{klimyk}, \cite{kustermans}, \cite{vandaele_discrete} etc. 
 
 \subsection{Hopf algebras}

\label{preliminaries_subsection_Hopf_algebras}
 
 We recall that an associative algebra with an unit is a vector space $A$ over $\IC$ together with two linear maps $m: A \otimes A \rightarrow A $ called the multiplication or the product and $\eta: \IC \rightarrow A$ called the unit such that $ m ( m \otimes {\rm id} ) = m ({\rm id} \otimes m ) $ and $ m (\eta \otimes {\rm id}) = {\rm id} = m ({\rm id} \otimes \eta).$ Dualizing this, we get the following definition.

\bdfn

A {\bf coalgebra} $A$ is a vector space over $\IC$ equipped with two linear maps $ \Delta: A \rightarrow A \otimes A $ called the comultiplication or coproduct and $\epsilon: A \rightarrow \IC$ such that 

$$ (\Delta \otimes {\rm id}) \Delta = ({\rm id} \otimes \Delta) \Delta,  $$
$$ (\epsilon \otimes {\rm id}) \Delta = {\rm id} = ({\rm id} \otimes \epsilon) \Delta. $$

\edfn

\bdfn  

Let $ (A, \Delta_{A}, \epsilon_A) $ and $ (B, \Delta_{B}, \epsilon_B) $ be co algebras. A $\IC$ linear mapping $\phi: A \rightarrow B$ is said to be a {\bf cohomomorphism} if

$$ \Delta_B \circ \phi = (\phi \otimes \phi) \Delta_A $$
$$\epsilon_A = \epsilon_B \circ \phi  $$

\edfn

\vspace{4mm}

Let $\sigma$ denote the flip map : $A \otimes A \rightarrow A \otimes A$ given by $\sigma (a \otimes b) = b \otimes a.$

\bdfn

A coalgebra is said to be {\bf cocommutative} if $\sigma \circ \Delta = \Delta.$

\edfn

\bdfn

A linear subspace $B$ of $A$ is a {\bf subcoalgebra} of $A$ if $\Delta ( B ) \subseteq B \otimes B.$ 

\edfn
 
 \bdfn

A $\IC$ linear subspace $\cli$ of $A$ is called a {\bf coideal} if $ \Delta (\cli)  \subseteq A \otimes \cli + \cli \otimes A $ and $\epsilon (\cli) = \{ 0\}.$ 

\edfn

If $\cli$ is a coideal of $A,$ the quotient vector space $A/\cli$ becomes a coalgebra with comultiplication and counit induced from $A.$

\vspace{4mm}

{\bf Sweedler notation}

We introduce the so called Sweedler notation for comultiplication. If  $ a $ is an element of a coalgebra $ \cla, $
 the element $ \Delta ( a ) $ in $ \cla \otimes \cla $ is a finite sum  $ \Delta ( a ) = \sum_i a_{1i} \otimes a_{2i} $ where $ a_{1i}, ~ a_{2i} $ belongs to $ \cla.$  Moreover, the representation of $ \Delta ( a ) $ is not unique. For notational simplicity we shall suppress the index $ i $ and write the above sum symbolically as $ \Delta ( a ) =  a_{(1)} \otimes a_{(2)} .$ Here the subscripts $ ( 1 ) $ and $ ( 2 ) $ refer to the corresponding tensor factors. 
 
 \bdfn
  
  A vector space $ \cla $ is called a bialgebra if it is an algebra and a coalgebra along with the condition that $ \Delta : \cla \rightarrow \cla \otimes \cla $ and $ \epsilon: \cla \rightarrow \IC $ are algebra homomorphisms
 ( equivalently, $ m : \cla \otimes \cla \rightarrow \cla $ and $ \eta : \IC \rightarrow \cla $ are coalgebra co-homomorphisms ).   
 \edfn  
 
 \bdfn
  
  A bialgebra $ \cla $ is called a Hopf algebra if there exists a linear map $ \kappa : \cla \rightarrow \cla $ called the antipode or the coinverse of $ \cla, $ such that $ m \circ ( \kappa \otimes {\rm id}  ) \Delta = \eta \circ \epsilon = m \circ ( {\rm id} \otimes \kappa ) \circ \Delta.$
   
 \edfn
 
 {\bf Dual Hopf algebra}
 
 Let us consider a finite dimensional Hopf algebra $ \cla. $ Then the dual vector space $ \cla^{\prime} $ is an algebra with respect to the multiplication $ fg ( a ) = ( f \otimes g ) \Delta ( a ).$ Moreover, for $ f $ in $ \cla^{\prime} ,$ one defines the functional $ \Delta ( f ) \in ( \cla \otimes \cla )^{\prime} $ by $ \Delta ( f ) ( a \otimes b ) = f ( a b ), ~ a,~ b $ in $ \cla.$ Since $ \cla $ is finite dimensional, $ ( \cla \otimes \cla )^{\prime} = \cla^{\prime} \otimes \cla^{\prime} $ and so $ \Delta ( f ) $ belongs to $ \cla^{\prime} \otimes \cla^{\prime}.$ Then the algebra $ \cla^{\prime} $ equipped with the comultiplication $ \Delta $, antipode $ \kappa $ defined by  $ ( \kappa  f ) ( a ) = f ( \kappa ( a ) ) ,$ counit $ \epsilon_{\cla^{\prime}} $ defined by $ \epsilon_{\cla^{\prime}} ( f )  = f ( 1 ) $ and $ 1_{\cla^{\prime} } ( a ) = \epsilon ( a ) $ gives a Hopf algebra. This is called the dual Hopf algebra of $ \cla.$ 
 
 \vspace{4mm}

 
 
 
 
 
 
 \bdfn
 
 A Hopf $ \ast $-algebra is a Hopf algebra ( $ \cla, \Delta, \kappa, \epsilon  $ )  endowed with an involution $ \ast $ which maps $ a $ to an element denoted by $ a^*$ satisfying the following properties:
 
 1. For all $ a, b $ in $ \cla ,~ \alpha, \beta $ in $ \IC, ~ ( \alpha a + \beta b  )^* = \overline{\alpha} a^* + \overline{\beta} b^* ,~ (  a^* )^* = a,~ {( a . b )}^* = b^* a^*. $
 
 2. $ \Delta : \cla \rightarrow \cla \otimes \cla $ is a $ \ast $-homomorphism which means that $ \Delta ( a^* ) = {\Delta ( a )}^{*} $ where the involution on $ \cla \otimes \cla  $ is defined by $ ( a \otimes b )^* = a^* \otimes b^*. $

 \edfn
 
 \bppsn
 
  In any Hopf $\ast$-algebra  ( $ \cla, \Delta, \kappa, \epsilon  $ ), we have
  
  1. $ \epsilon ( a^* ) = \overline{\epsilon ( a ) } $ for all $ a $ in $ \cla.$ 
  
  2. $ \kappa ( {\kappa ( a^* )}^* ) = a $ for all $ a $ in $ \cla.$  
 
 \eppsn
 
 \vspace{4mm}
 
 We recall that the dual algebra $ \cla^{\prime} $ of a Hopf $ \ast $-algebra $ \cla $ is a $ \ast $-algebra with involution defined by 
 $$ f^* ( a ) = \overline{f (\kappa{(a)}^* )}, ~ {\rm for} ~  f ~ {\rm in} ~ \cla^{\prime}.$$      

{\bf Dual Pairing}





A left action of a Hopf $ \ast $-algebra  ( $ \clu, \Delta_{\clu}, \kappa_{\clu}, \epsilon_{\clu} $ )  on another Hopf $\ast$-algebra  ( $ \cla, \Delta_{\cla}, \kappa_{\cla}, \epsilon_{\cla} $ ) is a   bilinear form $ \triangleright: \clu \times \cla \rightarrow \IC $ if the following conditions hold:

$$ {\rm ( 1 )}  f \triangleright ( a_1 a_2 ) = \Delta_{\clu} (f) \triangleright ( a_1 \otimes a_2 ) , ~ ( f_1 f_2 ) \triangleright a  = ( f_1 \otimes f_2 ) \triangleright  \Delta_{\cla} ( a ) ; $$

$$ {\rm ( 2 )}  f \triangleright 1_{\cla} = \epsilon_{\clu} ( f ), ~  1_{\clu} \triangleright a = \epsilon_{\cla} ( a ); $$

$$ {\rm ( 3 )}  f^* \triangleright a = \overline{ f \triangleright {\kappa_{\cla} ( a )}^*} ~ ( {\rm equivalently} ~ f \triangleright a^* = \overline{ {\kappa_{\clu}( f )}^* \triangleright a }  ) $$
for all $ f, f_1, f_2 $ in $ \clu $ and $ a, a_1, a_2 $ in $ \cla.$

Similarly, a right action of a Hopf $ \ast $-algebra  ( $ \clu, \Delta_{\clu}, \kappa_{\clu}, \epsilon_{\clu} $ )  on another Hopf $ \ast $-algebra  ( $ \cla, \Delta_{\cla}, \kappa_{\cla}, \epsilon_{\cla} $ ) is a   bilinear form $ \triangleleft: \cla \times  \clu \rightarrow \IC $ if the following conditions hold:  $  a_1 a_2  \triangleleft f = ( a_1 \triangleleft f_{(1)} ) ( a_2 \triangleleft f_{(2)} ), ~ a \triangleleft ( f_1 f_2 ) = \Delta_{\cla} ( a ) \triangleleft ( f_1 \otimes f_2 ), 1_{\cla} \triangleleft f = \epsilon_{\clu} ( f ), ~ a \triangleleft 1_{\clu} = \epsilon_{\cla} ( a ),  a \triangleleft f^* = \overline{ {\kappa_{\cla} ( a )}^* \triangleleft f} $ ( equivalently $ a^* \triangleleft f = \overline{ a \triangleleft  {\kappa_{\clu}( f )}^*}  )$
for all $ f, f_1, f_2 $ in $ \clu $ and $ a, a_1, a_2 $ in $ \cla.$

 $ \clu = \cla^{\prime} $ gives a particular case of this duality pairing.

 \subsection{Compact Quantum Groups: basic definitions and examples}

 \label{preliminaries_subsection_CQG}

\bdfn  

A {\bf compact quantum group} (to be abbreviated as CQG from now on)  is given by a pair $(\cls, \Delta)$, where $\cls$ is a unital separable $C^*$ algebra equipped  with a unital $C^*$-homomorphism $\Delta : \cls \raro \cls \otimes \cls$ (where $\otimes$ denotes the injective tensor product)   satisfying \\
  (ai) $(\Delta \ot {\rm id}) \circ \Delta=({\rm id} \ot \Delta) \circ \Delta$ (co-associativity), and \\
  (aii) each of the linear spans of $ \Delta(\cls)(\cls \ot 1)$ and of  $\Delta(\cls)(1 \ot \cls)$ are norm-dense in $\cls \ot \cls$. \\
  
  \edfn
  
  It is well-known (see \cite{woro}, \cite{woro1}) that there is a canonical dense $\ast$-subalgebra $\cls_0$ of $\cls$, consisting of the matrix elements of
   the finite dimensional unitary (co)-representations (to be defined shortly) of $\cls$, and maps $\epsilon : \cls_0 \raro \IC$ (co-unit) and
   $\kappa : \cls_0 \raro \cls_0$ (antipode)  defined
    on $\cls_0$ which make $\cls_0$ a Hopf $\ast$-algebra.
    
    \vspace{4mm}
    
    The following theorem is the analogue of Gelfand Naimark duality for commutative CQG s.
    
    \bppsn
    
    Let $ G $ be a compact group. Let $ C ( G ) $ be the algebra of continuous functions on $ G .$ If we define $ \Delta $ by $ \Delta ( f ) ( g, h ) = f ( g.h ) $ for $ f $ in $ C ( G ), ~ g,h  $ in $ G, $ then this defines a CQG structure on $ C ( G ).$
    
    Conversely, let $ ( \cls, \Delta ) $ be a commutative CQG. Let $ H ( \cls ) $ denote the Gelfand spectrum of $ \cls $ and endow it with the product structure given by  $  \chi \chi^{\prime} = ( \chi \otimes \chi^{\prime} ) \Delta $ where $\chi, ~ \chi^{\prime} $ are in $ H ( \cls ). $ Then $ H (  \cls ) $ is a compact group.   
     
    \eppsn
    
    \brmrk
    
    \label{T_n_Van_Daele}
      
 In \cite{VanDaele2}, A Van Daele removed  Woronowicz's separability assumption (in \cite{woro} ) for the $ C^* $ algebra of the underlying compact quantum group. We remark that although we assume that CQG s are separable, most of the results in this thesis go through in the non separable case also. 
 
 \ermrk

 \bdfn

\label{preliminaries_CQG_comodule}

 Let $ (\cls, \Delta_\cls) $ be a compact quantum group.  
   A vector space $ M $  is said to be an algebraic  $ \cls $ co-module (or $\cls$ co-module) if  there exists a linear map $ \alpha : M \rightarrow M \otimes_{\rm alg} \cls_0  $ such that\\
 1. $ \alpha \otimes {\rm id} ) \alpha = ({\rm  id}  \otimes \Delta_{\cls} ) \alpha, $\\
  2. $ ( {\rm id} \otimes \epsilon ) \alpha(m) = m $ for all $ m $ in $ M.$
  
  \edfn

  In the notations as above, let us define $ \widetilde{\alpha} : M \otimes \cls \rightarrow M \otimes \cls $ by $ \widetilde{\alpha} = ( {\rm id} \otimes m ) ( \alpha \otimes {\rm id} ).$ Then we claim that $ \widetilde{\alpha} $ is invertible with the inverse given by $ T ( m \otimes q ) = ( {\rm id} \otimes \kappa ) \alpha ( m ) ( 1 \otimes q ),$ where $ m $ is in $ M, ~ q $ is in $ \cls.$ 
  As $ T $ is defined to be $ \cls_0 $ linear, it is enough to check that $ \widetilde{\alpha} T ( m \otimes 1 ) = m \otimes 1.$
  \bean \lefteqn{\widetilde{\alpha} T ( m \otimes 1 ) }\\
        &=& \widetilde{\alpha} ( m_{(1)} \otimes \kappa ( m_{(2)} ) 1 )\\
        &=& m_{(1)(1)} \otimes m_{(1)(2)} \kappa ( m_{(2)} )\\
        &=& ( {\rm id} \otimes m ( {\rm id} \otimes \kappa ) \Delta ) \alpha ( m )\\
        &=& ( {\rm id} \otimes \epsilon ( ).1 ) \alpha ( m )\\
        &=& m \otimes 1 .\eean
        Similarly, $ T \widetilde{\alpha} = {\rm id}. $
        Thus, \be \label{preliminaries_comodule_inverse} T = {\widetilde{\alpha}}^{- 1}.\ee

 \bdfn
 
 A morphism from a CQG  $ ( \cls_1, \Delta_1 ) $ to another CQG $ ( \cls_2, \Delta_2 ) $ is a unital $ C^* $ homomorphism $ \pi : \cls_1 \rightarrow \cls_2 $ such that 
 
               $$ ( \pi \otimes \pi ) \Delta_1 = \Delta_2 \pi.$$  
 
 \edfn     
  
  It follows that in such a case, $ \pi $ preserves the Hopf $ \ast $-algebra structures, that is, we have  
  $$ \pi ( {( \cls_1 )}_0 ) \subseteq  {( \cls_2 )}_0, ~ \pi \kappa_1 = \kappa_2 \pi, ~ \epsilon_2 \pi = \epsilon_1 , $$
  where $ \kappa_1, \epsilon_1 $ denotes the antipode and counit of $ \cls_1 $ respectively while $ \kappa_2, \epsilon_2 $ denotes those of $ \cls_2.$

 \bdfn
 
A {\bf Woronowicz $C^*$-subalgebra} of a CQG $ ( \cls_1, \Delta ) $ is a $ C^* $ subalgebra $ \cls_2 $  of $ \cls_1 $ such that $ ( \cls_2, \Delta|_{\cls_2} ) $ is a CQG such that the inclusion map from $ \cls_2 \rightarrow \cls_1 $ is a morphism  of CQG s.
 
 \edfn
 
 \bdfn
 
A {\bf  Woronowicz $C^*$-ideal} of a CQG $ ( \cls, \Delta  ) $ is a $ C^* $ ideal $ J $ of $ \cls $ such that $ \Delta ( J ) \subseteq {\rm Ker} ( \pi \otimes \pi ), $ where $ \pi $ is the quotient map from $ \cls $ to $ \cls / J. $ 
 
 \edfn

It can be easily seen that a kernel of a CQG morphism is a Woronowicz $C^*$-ideal. 
 
 We recall the following isomorphism theorem.
 
 \bppsn
 
 The quotient of a CQG $ ( \cls, \Delta ) $ by a Woronowicz $ C^* $-ideal $ \cli $ has  a unique CQG structure such that the quotient map $ \pi $  is a morphism of CQG s. More precisely, the coproduct $ \widetilde{\Delta} $   on $ \cls / \cli $ is given by $ \widetilde{\Delta} ( s + \cli ) = ( \pi \otimes \pi ) \Delta ( s ).$
 
 \eppsn
 
 \bdfn
 
 A CQG $ ( \cls^{\prime}, \Delta^{\prime} ) $ is called a quantum subgroup of another CQG $ ( \cls, \Delta )  $ if there is a Woronowicz $ C^* $-ideal $ J $ of $ \cls $ such that $ ( \cls^{\prime}, \Delta^{\prime} ) \cong ( \cls, \Delta )/J. $
 
 \edfn

Let us mention a convention which we are going to follow. We shall use most of the terminologies of
\cite{free}, for example Woronowicz $C^*$ -subalgebra, Woronowicz $C^*$-ideal etc, however with the exception that we shall call the Woronowicz $C^*$ algebras just compact quantum groups, and not use the term compact quantum groups for the dual objects as done in \cite{free}.

Let $(\cls, \Delta)$ be a compact quantum group. Then there exists a state $h$ on $\cls$, to be called a {\bf Haar state} on $\cls$ such that $( h \otimes {\rm id} ) \Delta (s) = ( {\rm id} \otimes h ) \Delta ( s ) = h (s).1.$  We recall that unlike the group case, $h$ may not be faithful. But on the dense Hopf $\ast$-algebra $\cls_0$ mentioned above, it is faithful. We have the following result.

\bppsn

Let $ i : \cls_1 \rightarrow \cls_2 $ be an injective morphism of CQG s. Then the Haar state on $ \cls_1 $ is the restriction of that of $ \cls_2 $ on $ \cls_1.$ 

\eppsn

\brmrk

\label{preliminaries_haar_state_not_tracial}

In general, the Haar state might not be tracial. In fact, there exists a multiplicative linear functional denoted by $ f_1 $ in \cite{woro1} such that $ h ( a b ) = h ( b ( f_{1} \triangleleft a \triangleright f_1 ) ).$ Moreover, from Theorem 1.5 of \cite{woro}, it follows that the Haar state of a CQG is tracial if and only if $ \kappa^2 = {\rm id}.$ 

\ermrk  
    
    \vspace{4mm}
    
{\bf Co-Representations of a compact quantum group}
    
\bdfn
  
   A co-representation of a compact quantum group $ ( S, \Delta ) $ on a Hilbert space $ \clh $ is a map $ U $ from $ \clh $ to the Hilbert $\cls$ module $ \clh \otimes \cls $  such that the  element $ \widetilde{U} $ belonging to $ \clm ( \clk ( \clh ) \otimes \cls ) $ given by $\widetilde{U}( \xi \ot b)=U(\xi)(1 \ot b)$ (  $\xi $ in $ \clh, b $ in $ \cls $) satisfies  $$ ({\rm  id} \otimes \Delta ) \widetilde{U} = {\widetilde{U}}_{(12)} {\widetilde{U}}_{(13)},$$ where for an operator $X $ in $ \clb(\clh_1 \ot \clh_2)$ we have denoted by $X_{(12)}$ and $X_{(13)}$ the operators $X \ot I_{\clh_2} $ in $ \clb(\clh_1 \ot \clh_2 \ot \clh_2)$, and $\Sigma_{23} X_{(12)} \Sigma_{23}$ respectively and $\Sigma_{23}$ is the unitary on $\clh_1 \ot \clh_2 \ot \clh_2$ which flips the two copies of $\clh_2$.
   
 If $\widetilde{U}$ is an unitary element of  $ \clm ( \clk ( \clh ) \otimes \cls ) ,$ then $U$ is called a unitary co-representation.  

\edfn

From now on, we will drop the term co in the word co-representation unless there is any confusion.

 
 \brmrk
 
 Let $ \pi $ be a CQG  morphism from a CQG $ ( \cls_1, \Delta_1 ) $ to another CQG $ ( \cls_2, \Delta_2 ).$ Then for every unitary representation $ U $ of $ \cls_1, ~ ( {\rm id} \otimes  \pi ) U $ is a unitary representation of $ \cls_2.$
 
 \ermrk   
  
 Following the definitions given in the last part of subsection \ref{preliminaries_subsection_Hilbert_C*_modules} and a unitary representation $U$ of a CQG  on a Hilbert space $ \clh, $ and a  not necessarily bounded, densely defined (in the weak operator topology)  linear functional $\tau$ on $\clb(\clh)$,  we will use the notation $ \alpha_U $ and the terms ``$ \alpha_U $ preserves $ \tau $" and `` $ U $ equivariant " throughout this thesis. 

A CQG  $( \cls, \Delta )$ has a distinguished representation which corresponds to the right regular representation in the group case. Let $ \clh $ be the GNS space of $ \cls $ associated with the Haar state $ h, ~ \xi_0 $ be the associated cyclic vector and $ \clk $ be a Hilbert space on which $ \cls $ acts faithfully and non-degenerately. There is a unitary operator $ u $ on $ \clh  \otimes \clk $ defined by $ u ( a \xi_0 \otimes \eta ) = \Delta ( a ) ( \xi_0 \otimes \eta ) $ when $ a $ is in $ \cls, ~ \eta $ is in $ \clk. $ Then $ u $ can be shown to be an element of multiplier of $ \clk ( \clh ) \otimes \cls$ and called the right regular representation of  $ \cls. $    

Let $v$ be a representation of a CQG $(\cls, \Delta)$ on a Hilbert space $\clh.$ A closed subspace $\clh_1$ of $\clh$ is said to be invariant if $ ( e \otimes 1) v (e \otimes 1) = v (e \otimes 1),$ where $e$ is the orthogonal projection onto this subspace. The representation $v$ is called irreducible if the only invariant subspaces are $\{0\}$ and $\clh.$ It is clear that one can make sense of direct sum of (co)-representations in this case also. Moreover, for two  representations $v$ and $w$ of  a CQG $(\cls, \Delta)$ on Hilbert spaces $\clh_1$ and $\clh_2,$ the tensor product of $v$ and $w$ is given by the element $ v_{(13)} w_{(23)}.$ The intertwiner between $v$ and $w$ is an element $x$ in $\clb (\clh_1, \clh_2)$ such that $ ( x \otimes 1 ) v = w (x \otimes 1).$ The set of intertwiners between $v$ and $w$ is denoted by ${\rm Mor} ( v,w).$ Two representations are said to be equivalent if there is an invertible intertwiner. They are unitarily equivalent if the intertwiner can be chosen to be unitary. 

Just like the case of compact groups, CQG s have an analogous Peter Weyl theory which corresponds to the usual Peter Weyl theory in the group case. We will give a sketch of it by mentioning the main results and refer to \cite{vandaelenotes}, \cite{woro1} and \cite{woro} for the details.

Let $v$ be a unitary representation of $(\cls, \Delta)$ on $\clh.$ If $ \clh_1 $ is an invariant subspace, then the orthogonal complement of $\clh_1$ is also invariant. Any non degenerate finite dimensional representation is equivalent with a unitary representation. 




Every irreducible unitary representation of a CQG is contained in the regular representation. Let $v$ be a representation on a finite dimensional Hilbert space $\clh.$ If we denote the matrix units in $\clb (\clh) $ by $(e_{pq}),$ we can write $v = \sum e_{pq} \otimes v_{pq}. $ $v_{pq}$ are called the matrix elements of the finite dimensional representation $v.$  Define $\overline{v} = \sum e_{pq} \otimes v^*_{pq}. $ Then $\overline{v}$ is a representation and is called the adjoint of $v.$ It can be shown that if $v$ is a finite dimensional irreducible representation, then $\overline{v}$ is also irreducible. Moreover, for an irreducible unitary representation, its adjoint is equivalent with a unitary representation.

 The subspace spanned by the matrix elements of finite dimensional unitary representations is denoted by $\cls_0 .$      Firstly, $ \cls_0 $ is a subalgebra as the product of two matrix elements of finite dimensional unitary representations is a matrix element of the tensor product of these representations. Moreover, as the adjoint of a finite dimensional unitary representation is equivalent with a unitary representation, $ \cls_0 $ is $ \ast $ invariant. We note that $ 1 $ is in $ \cls_0 $ as $ 1 $ is a representation. Now, we will  recall some basic facts about the subalgebra $ \cls_0. $ We will denote the Haar state of $ \cls $ by $ h. $ 
 
 \bppsn
 
 \label{preliminaries_CQG_Peter_Weyl}
 
 (1) $ \cls_0 $ is a dense $ \ast $-subalgebra of $ \cls. $ 
 
 \vspace{2mm}  
 
 (2)Let $ \{ u^{\alpha} : \alpha \in I \} $ be a complete set of mutually inequivalent, irreducible unitary representations. We will denote the representation space and  dimension of $ u^{\alpha} $ by $ \clh_{\alpha}  $ and $ n ( \alpha ) $ respectively.  Then the Schur's orthogonality relation takes the following form:
  
 For any $ \alpha $ in $ I, $ there is a positive invertible operator $ F^{\alpha} $ acting on $ \clh_{\alpha} $ such that for any $ \alpha, \beta $ in $ I $ and $ 1 \leq j,  q \leq n ( \alpha ), ~ 1 \leq i, p \leq n ( \beta ) $
    $$ h ( {( u^{\beta}_{ip} )}^* u^{\alpha}_{jq}  ) = \delta_{\alpha \beta} \delta_{pq} F^{\alpha}_{ij}.   $$

\vspace{2mm}

(3) $ \{ u^{\alpha}_{pq} : \alpha \in I, 1 \leq p,q \leq n ( \alpha ) \} $ form a basis for  $ \cls_0.$ 

\vspace{2mm}

(4) Moreover, $ \Delta $ maps $ \cls_0 $ into $ \cls_0 \otimes \cls_0. $ In fact, $ \Delta $ is given by $ \Delta ( u^{\alpha}_{pq} ) = \sum^{n_{\alpha}}_{k = 1} u^{\alpha}_{pk} \otimes u^{\alpha}_{kq}.$ A counit and an antipode are defined on $ \cls_0 $ respectively by the formulae,
$$ \epsilon ( u^{\alpha}_{pq} ) = \delta_{pq}, ~ \kappa ( u^{\alpha}_{pq} ) = {( u^{\alpha}_{qp} )}^* .$$ 
It follows that $ \cls_0 $ becomes a Hopf $ \ast$-algebra.  
 
 \eppsn 
 
 \vspace{4mm}

    
    
    
    
    
    

    A {\bf compact matrix quantum group}  is a CQG such that there exists a distinguished unitary irreducible representation called the fundamental representation such that the $\ast$-algebra spanned by its matrix elements is a dense Hopf $\ast$-subalgebra of the CQG.

 We now discuss the {\bf free product of CQG} s which were developed in \cite{free}. 
 Let $ (  \cls_1, \Delta_1 ) $ and $ ( \cls_2, \Delta_2 ) $ be two CQG s. Let $ i_1 $ and $ i_2 $ denote the canonical injections of $ \cls_1 $ and $ \cls_2 $ into the $ C^* $ algebra $ \cls_1 \ast \cls_2. $ Put $ \rho_1 = ( i_1 \otimes i_1 ) \Delta_1 $ and $ \rho_2 = ( i_2 \otimes i_2 ) \Delta_2.$ By the universal property of $ \cls_1 \ast \cls_2 ,$ there exists a map $ \Delta : \cls_1 \ast \cls_2 \rightarrow ( \cls_1 \ast \cls_2  ) \otimes ( \cls_1 \ast \cls_2 )$ such that $ \Delta i_1 = \rho_1  $ and $ \Delta i_2 = \rho_2.$ It can be shown that $ \Delta $ indeed has the required properties so that $ ( \cls, \Delta  ) $ is a CQG.
 
 Let $ \{ \cls_n \}_{n \in \IN} $ be an inductive sequence of CQG s, where the connecting morphisms $ \pi_{mn} $ from $ \cls_n $ to $ \cls_m ~ ( n < m ) $ are injective morphisms of CQG s. Then from Proposition 3.1 of \cite{free}, we have that the inductive limit $ \cls_0  $ of $ \cls_n $ s has a unique CQG structure with the following property: for any CQG $ \cls^{\prime}$ and any family of CQG morphisms $ \phi_n : \cls_n \rightarrow \cls^{\prime} $ such that  $ \phi_m \pi_{mn} = \phi_n, $ the uniquely defined morphism $ {\rm lim}_n \phi_n $ in the category of unital $ C^* $ algebras is a morphism in the category of CQG s. 
 
 Combining the above two results, it follows that the free product $ C^* $ algebra  of an arbitrary sequence of CQG s  has a natural CQG structure.  
 
 Moreover, the following result was derived in \cite{free}. 
 
 \bppsn

 Let $ \Gamma_1, \Gamma_2 $ be a discrete abelian groups. Then the natural isomorphisms $ C^* ( \Gamma_1 ) \cong C ( \widehat{\Gamma_1} ) $ and $ C^* ( \Gamma_1 ) \ast C^* ( \Gamma_2 ) \cong C^* ( \Gamma_1 \ast \Gamma_2 ) $  are isomorphism of CQG s. 
  
 \eppsn

Let $ i_1 $ and $ i_2 $ be the inclusion of CQG s $ \cls_1 $ and $ \cls_2 $ into $ \cls_1 \ast \cls_2.$  If $U_1$ and $U_2$ are unitary representations of CQG s $\cls_1$ and $\cls_2$ on Hilbert spaces $\clh_1$ and $\clh_2$ respectively, then the {\bf free product representation} of $U_1$ and $U_2$ is a representation of the CQG $\cls_1 \ast \cls_2$ on the Hilbert space $\clh_1 \oplus \clh_2$ given by the $\cls_1 \ast \cls_2$ valued matrix $  \left ( \begin {array} {cccc}
 ( {\rm id} \otimes i_1 )  U_1 & 0  \\ 0 & ( {\rm id} \otimes i_2 ) U_2 \end {array} \right ) .$

Similarly, the free product representation of an arbitrary sequence of CQG representations are defined.

\vspace{4mm} 
 
The inductive limit of an arbitrary sequence of CQG s has the structure of a CQG. The following lemma is probably known, but we include the proof ( taken from \cite{AF} ) for the sake of completeness.

\blmma
\label{basic}
Suppose that $(\cls_n)_{n \in \IN}$ is a sequence of CQG s and for each $n,m $ in $ \IN,$ $n \leq m$ there is a CQG morphism $\pi_{n,m}:\cls_n \to \cls_m$ with the compatibility property
\[ \pi_{m,k} \circ \pi_{n,m} = \pi_{n,k}, \;\;\; n\leq m \leq k.\]
Then the inductive limit of $C^*$-algebras $(\cls_n)_{n \in \IN}$ has a canonical structure of a CQG. It will be
denoted $\cls_\infty$ or $\lim_{n \in \IN} \cls_n$. It has the following universality property:\\
for any CQG $(\cls, \Delta)$ such that there are CQG morphisms $\pi_n : \cls_n \to \cls$ satisfying for all $m,n \in \IN$, $m \geq n$ the equality $\pi_m \circ \pi_{n,m}=\pi_n$, there exists a unique CQG morphism $\pi_\infty : \cls_\infty \to \cls$ such that $\pi_n=\pi_\infty \circ \pi_{n, \infty}$ for all $n\in \IN$, where we have denoted by $\pi_{n,\infty}$ the canonical unital $C^*$-homomorphism from $\cls_n$ into $\cls_\infty$.
\elmma

{\it Proof:}\\ 
Let us denote the coproduct on $\cls_n$ by $\Delta_n$. We consider the unital $C^*$-homomorphism $\rho_n : \cls_n \to \cls_\infty
\ot \cls_\infty$ given by $ \rho_n=(\pi_{n,\infty} \ot \pi_{n,\infty}) \circ \Delta_n $, and observe that these maps  do satisfy
the compatibility property: $$ \rho_m \circ \pi_{n,m}=\rho_n~~~\forall n \leq m.$$ Thus, by the general properties of the
$C^*$-algebraic inductive limit, we have a unique unital $C^*$-homomorphism $\Delta_\infty : \cls_\infty \to \cls_\infty \ot
\cls_\infty$ satisfying $\Delta_\infty \circ \pi_{n,\infty}=\rho_n$ for all $n$. We claim that $(\cls_\infty, \Delta_\infty)$ is a
CQG.

We first  check that $ \Delta_{\infty} $ is coassociative.
It is enough to verify the coassociativity on the dense set $\cup_n \pi_{n,\infty}(\cls_n)$. Indeed, for $s=\pi_{n,\infty}(a)$ $(a \in \cls_n$),  by using $\Delta_\infty \circ \pi_{n,\infty}=(\pi_{n,\infty} \otimes \pi_{n,\infty}) \circ \Delta_n,$ we have the following:
\bean \lefteqn{(\Delta_\infty \otimes {\rm id})\Delta_\infty(\pi_{n,\infty}(a))}\\
&=& (\Delta_\infty \otimes {\rm id})(\pi_{n,\infty }\otimes \pi_{n,\infty})(\Delta_n(a))\\
&=&(\pi_{n,\infty} \otimes \pi_{n,\infty} \otimes \pi_{n,\infty})(\Delta_n \otimes {\rm id})(\Delta_n(a))\\
&=& (\pi_{n,\infty} \otimes \pi_{n,\infty} \otimes \pi_{n,\infty})( {\rm id} \otimes \Delta_n)(\Delta_n(a))\\
&=&(\pi_{n,\infty} \otimes (\pi_{n,\infty} \otimes \pi_{n,\infty})\circ \Delta_n)(\Delta_n(a))\\
&=&(\pi_{n,\infty} \otimes \Delta_\infty \circ \pi_{n,\infty})(\Delta_n(a))\\
&=&({\rm id} \otimes \Delta_\infty)((\pi_{n,\infty} \otimes \pi_{n,\infty})(\Delta_n(a)))\\
&=&({\rm id} \otimes \Delta_\infty)(\Delta_\infty(\pi_{n,\infty}(a)))\eean
 which proves the coassociativity.

Finally, we need to verify the quantum cancellation properties. Note that to show that  $ \Delta_{\infty} ( \cls_{\infty} ) ( 1 \otimes  \cls_{\infty} )  $ is dense in $ \cls_{\infty} \otimes  \cls_{\infty}$
it is enough to show that the above assertion is true with $\cls_{\infty}$ replaced by a dense subalgebra $ \bigcup_{n} \pi_{n, \infty} ( \cls_{n} ) $.

Using the density of $ \Delta_{n} ( \cls_{n} ) ( 1 \otimes \cls_{n} ) $ in $\cls_{n} \otimes \cls_{n} $ and the contractivity of the map $ \pi_{n, \infty} $ we note that $ ( \pi_{n, \infty} \otimes  \pi_{n, \infty} ) ( \Delta_{n} ( S_{n} ) ( 1 \otimes \cls_{n} ))$ is dense in $ ( \pi_{n, \infty} \otimes  \pi_{n, \infty} ) ( \cls_{n} \otimes \cls_{n} ).$
 This implies that  $ ( \pi_{n, \infty} \otimes  \pi_{n, \infty} ) ( \Delta_{n} ( \cls_{n} ) ) ( 1 \otimes   \pi_{n, \infty} ( \cls_{n} ) ) $ is dense in $  \pi_{n, \infty} ( \cls_{n} ) \otimes \pi_{n, \infty} ( \cls_{n} ) $ and hence $ \Delta_{\infty} ( \pi_{n, \infty} ( \cls_{n} ) ) ( 1 \otimes \pi_{n, \infty} ( \cls_{n} ) ) $ is dense in  $ \pi_{n, \infty} ( \cls_{n} ) \otimes  \pi_{n, \infty} ( \cls_{n} ).$ The proof of the claim now follows by  noting that $ \pi_{n, \infty} ( \cls_{n} ) =  \pi_{m, \infty} \pi_{n, m} ( \cls_{n} ) \subseteq \pi_{m, \infty} ( \cls_{m} ) $ for any $m \geq n$, along with the above observations. The right quantum cancellation property can be shown in the same way.

The proof of the universality property is routine and hence omitted.

\qed

\vspace{4mm}

We note that the proof remains valid for any other indexing set for the net, not necessarily $\IN.$

\subsection{The CQG $ U_{\mu} ( 2 ) $}

\label{preliminaries_subsection_U_mu_2}

We now introduce the compact quantum group $U_{\mu}( 2 ).$ We refer to \cite{klimyk} for more details.

As a unital $ C^{*} $ algebra, $ U_{\mu}( 2 ) $ is generated by $4$ elements $ u_{11}, u_{12}, u_{21}, u_{22} $ satisfying:

\be \label{qorient_Umu2_def_1} u_{11}u_{12} = \mu u_{12}u_{11} \ee 
\be \label{qorient_Umu2_def_2} u_{11}u_{21} = \mu u_{21}u_{11} \ee 
\be \label{qorient_Umu2_def_3} u_{12}u_{22} = \mu u_{22}u_{12} \ee 
\be \label{qorient_Umu2_def_4} u_{21}u_{22} = \mu u_{22}u_{21} \ee 
\be \label{qorient_Umu2_def_5} u_{12}u_{21} = u_{21}u_{12} \ee 
\be \label{qorient_Umu2_def_6} u_{11}u_{22} - u_{22}u_{11} = ( \mu - {\mu}^{-1} )u_{12}u_{21} \ee
and the condition that the matrix  $ u = \left(  \begin {array} {cccc}
     u_{11}   & u_{12}  \\ u_{21} & u_{22} \end {array} \right ) $ is a unitary. Thus, the above matrix $ u $ is the fundamental unitary for $ U_{\mu} ( 2 ).$ 
 
 The  CQG structure is given by  
 \be \label{qorient_U_mu_2_CQG_structure} \Delta ( u_{ij} ) = \sum_{k = 1,2} u_{ik} \otimes u_{kj} , ~ \kappa ( u_{ij} ) = {u_{ji}}^{*}, \epsilon ( u_{ij} ) = \delta_{ij} .\ee 
The quantum determinant $ D_{\mu} $ is defined by
\be \label{qorient_U_mu_2_quantum_determinant} D_{\mu} = u_{11}u_{22} - \mu u_{12}u_{21} = u_{22}u_{11} - {\mu}^{-1} u_{12}u_{21} .\ee
 Then, $ {D_{\mu}}^{*}D_{\mu} = D_{\mu}{D_{\mu}}^{*} = 1 .$ Moreover, $ D_{\mu} $ belongs to the centre of $ U_{\mu} ( 2 ) .$

We mention the following result for future use. 

\bppsn

\label{qorient_U_mu_2_k(u_ij)}

$$ \kappa ( u_{11} ) = u_{22} D^{- 1}_{\mu}, ~ \kappa ( u_{12} ) = - \mu^{- 1} u_{12} D^{- 1}_{\mu}, ~ \kappa ( u_{21} ) = - \mu u_{21} D^{- 1}_{\mu}, ~ \kappa ( u_{22} ) = u_{11} D^{- 1}_{\mu}. $$

\eppsn

{\it Proof :}  By \cite{klimyk} ( Proposition 10, Page 314 ), we have that  $ \kappa ( u_{ij} ) = \widetilde{u_{ij}} D^{- 1}_{\mu} $ where $ \widetilde{u_{ij}} $ is the $ (i,~ j )  $th entry of a matrix $ \widetilde{u} $ satisfying $ \widetilde{u} u = u \widetilde{u} = D_{\mu} I_{2}. $ One can easily check that the matrix $ \left(  \begin {array} {cccc}
     u_{22}   & - \mu^{- 1} u_{12}  \\ - \mu u_{21} & u_{11} \end {array} \right )  $ satisfies this equation and this proves the Proposition. \qed

\subsection{The CQG $ SU_{\mu} ( 2 ) $}

\label{preliminaries_subsection_SU_mu_2}

Let $ \mu $ belongs to $ [ -1, 1 ].$
 The $ C^* $ algebra $ SU_{\mu}( 2 ) $ is defined as the universal unital $ C^{*} $ algebra generated by $\alpha, \gamma $ satisfying: 
  \be \label{su2def1}  {\alpha}^{*} \alpha + {\gamma}^{*}\gamma = 1, \ee  
  \be \label{su2def2}  \alpha {\alpha}^{*} + {\mu}^{2}{\gamma}{\gamma}^{*} = 1, \ee  
  \be \label{su2def3} {\gamma}{\gamma}^{*} = {\gamma}^{*}\gamma, \ee  
  \be \label{su2def4} \mu \gamma \alpha = \alpha \gamma, \ee  
  \be \label{su2def5} \mu {\gamma}^{*} \alpha = \alpha {\gamma}^{*}. \ee

  The fundamental representation of $ SU_{\mu} ( 2 ) $ is given by : $  \left ( \begin {array} {cccc}
   \alpha & - \mu {\gamma}^{*}  \\ \gamma & {\alpha}^{*} \end {array} \right ) .$
  
  There is a coproduct $ \Delta $ of $ SU_{\mu} ( 2 ) $ given by :  
  \bean \Delta ( \alpha ) = \alpha \otimes \alpha - \mu \gamma^* \otimes \gamma, ~~ \Delta ( \gamma ) = \gamma \otimes \alpha + \alpha^* \otimes \gamma \eean  
  which makes it into a CQG. Let $ h $ denote the Haar state and $ \clh = L^{2}( SU_{\mu}( 2 ) ) $ be the corresponding G.N.S space.
  
  {\bf Haar state on $SU_{\mu} (2)$}

We restate the content of Theorem 14, Chapter 4 ( page 113 ) of \cite{klimyk} in a convenient form below.
For all $ m \geq 1, n,l,k \geq 0, k^{\prime} \neq k^{\prime \prime},$ 
\be \label{sphere_haar_state_su_mu_2_original}h ( (\gamma^* \gamma)^{k}) = \frac{1 - \mu^2}{1 - \mu^{2k + 2}},~ h ( \alpha^{m} \gamma^{*n}\gamma^{l}) = 0, ~ h ( \alpha^{*m} \gamma^{*n}\gamma^{l}) = 0, ~ h ( \gamma^{*k^{\prime}} \gamma^{*k^{\prime \prime}} ) = 0  .\ee

\vspace{4mm}

{\bf ( Co )- representations of $ SU_{\mu} ( 2 ) $}
   
   For each $ n $ in $ \{ 0, 1/2,1,..... \},$ there is a unique irreducible representation $ T^{n} $ of dimension $ 2n + 1.$ Denote by $ t^{n}_{ij} $ the $ ij $ th entry of $ T^{n}.$ They form an orthogonal basis of  $ \clh.$ Denote by $e^{n}_{ij} $ the normalized $ t^{n}_{ij} $ s so that $ \{ e^{n}_{ij} : n = 0, 1/2,1,......,i,j = -n, -n +1,.....n \}$ is an orthonormal basis.
   
   We recall from \cite{klimyk} that 

\be t^{1/2}_{-1/2, -1/2} = \alpha,~ t^{1/2}_{-1/2, 1/2} = - \mu {\gamma}^{*} ,~ t^{1/2}_{1/2, -1/2} = \gamma,~ t^{1/2}_{1/2, 1/2} =  {\alpha}^{*} . \ee

Moreover, if we define

 \be \label{qorient_SU_mu_2_recursive_lemma_equation_0} f_{n,i} = a ( n,i ) \alpha^{n - i} \gamma^{* n + i} \ee 
 
  where $ a( n,i ) $ s are some constants as in \cite{klimyk}, then  $ \{ f_{n,i} : n = 0, \frac{1}{2}, 1, \frac{3}{2},....., ~ -n \leq i \leq n \} $ is an orthonormal basis of $ SU_{\mu}( 2 ) $ and $ \Delta ( f_{n,i} ) = \sum^{n}_{k = - n} f_{n,k} \otimes t^{n}_{k,i}.$  

The following recursive relations will be useful to us.

\bppsn

\label{qorient_SU_mu_2_recursive_relations}

\bea 
\label{rec_1} \lefteqn{ t^{l + 1/2}_{i, l + 1/2} }{\nonumber}\\
 &=&  c_{11}(i,l) \gamma^* t^{l}_{i + 1/2, l} + c_{12}(i,l) \alpha^* t^{l}_{i- 1/2, l} ~~~~~~~~ -l + 1/2 \leq i \leq l - 1/2 ,{\nonumber}\\
&=&  c_{21}(i,l) {\gamma}^{*} t^{l}_{i+ 1/2, l}   ~~~~~~~~~~         i = -l - 1/2, {\nonumber}\\ 
&=&  c_{31}(i,l) {\alpha}^{*} t^{l}_{i - 1/2, l} ~~~~~~~~          i= l + 1/2 , {\nonumber}\\
\eea

and for $ j \le l,$
\bea   \lefteqn{ t^{l + 1/2}_{i, j}  }{\nonumber}\\
 & =&  c(l,i,j) \alpha t^{l}_{i+ 1/2, j + 1/2}  +c^\prime(l,i,j)  \gamma t^{l}_{i - 1/2, j + 1/2} ~~~~~~~    - l + 1/2 \leq i \leq l - 1/2, {\nonumber}\\
& = & d(l,j) \alpha  t^{l}_{- l, j + 1/2} + d^{\prime}(l,j) \gamma^* t^{l}_{- l, j - \frac{1}{2}} ~~~~~~~   i = -l - 1/2, ~ - l + \frac{1}{2} \leq j \leq l - \frac{1}{2}, {\nonumber}\\
& = & d^{\prime \prime}(l,j) \alpha  t^{l}_{i+ 1/2, j + 1/2} ~~~~~~~   i = -l - 1/2, ~ j = - l - \frac{1}{2}, {\nonumber}\\
&=&  e(l,j)  \gamma t^{l}_{i - 1/2, j+ 1/2} +  e^{\prime}( l,j ) \alpha^* t^{l}_{i - \frac{1}{2}, j - \frac{1}{2}} ~~~~~~ i= l + 1/2,{\nonumber} \label{rec_2} \\
\eea
where $C_{pq}(il), c(l,i,j), d(l,j), d^\prime_{l,j}, d^{\prime \prime}(l,j), e(l,j), e^{\prime}( l,j )  $ are all complex numbers.

\eppsn

{\it Proof :} It can be easily seen that 
\be \label{qorient_SU_mu_2_recursive_lemma_equation_1} f_{l + \frac{1}{2},i} = c ( l,i ) \alpha f_{l, i + \frac{1}{2}} \ee
 for some constants $ c ( l,i ) .$ 
 
 Moreover, from ( \ref{qorient_SU_mu_2_recursive_lemma_equation_0} ) we have 
$  \gamma^* f_{l,i} = a^{\prime} ( l,i ) \alpha^{l - i} {\gamma^{*}}^{l + i + 1} $  for some constant $ a^{\prime} ( l,i ). $ This means that 
\be \label{qorient_SU_mu_2_recursive_lemma_equation_2} \gamma^* f_{l,i} = \frac{a^{\prime} ( l,i )}{a ( l+ \frac{1}{2}, i + \frac{1}{2} )} f_{l + \frac{1}{2}, i + \frac{1}{2}}.\ee

We have $ f_{l + \frac{1}{2}, l + \frac{1}{2}} = a ( l + \frac{1}{2}, l + \frac{1}{2} ) {\gamma^*}^{2l + 1}$ and $ f_{l,l} = a ( l,l ) {\gamma^*}^{2l} $ which implies that 

\be \label{qorient_SU_mu_2_recursive_lemma_equation_3} f_{l + \frac{1}{2}, l + \frac{1}{2}} = \frac{a ( l + \frac{1}{2}, l + \frac{1}{2} )}{a ( l,l )} \gamma^* f_{l,l}.\ee

Now, we proceed to prove (  \ref{rec_1}  ). Applying coproduct on (  \ref{qorient_SU_mu_2_recursive_lemma_equation_3}  ) and using ( \ref{qorient_SU_mu_2_recursive_lemma_equation_0}  ) and (  \ref{qorient_SU_mu_2_recursive_lemma_equation_2}    ), we have
\bean \lefteqn{\sum^{l + \frac{1}{2}}_{k = - ( l + \frac{1}{2} ) } f_{l + \frac{1}{2},k} \otimes t^{l + \frac{1}{2}}_{k, l + \frac{1}{2}} }\\
&=& \frac{a ( l + \frac{1}{2},l + \frac{1}{2} )}{a ( l,l )} ( \gamma^* \otimes \alpha^* + \alpha \otimes \gamma^* ) ( \sum^{l}_{k = - l} f_{l,k} \otimes t^{l}_{k,l}  )\\
&=& \frac{a ( l + \frac{1}{2},l + \frac{1}{2} )}{a ( l,l )} \left[ \sum^{ l}_{k = - l } \gamma^* f_{l,k} \otimes \alpha^* t^{l}_{k,l} + \sum^{l}_{k = - l} \alpha f_{l,k} \otimes \gamma^* t^{l}_{k,l} \right]\\
&=& \frac{a ( l + \frac{1}{2},l + \frac{1}{2} )}{a ( l,l )} \left[ \sum^{l}_{k = - l} \frac{a^{\prime}( l,k )f_{l + \frac{1}{2}, k + \frac{1}{2}} \otimes \alpha^* t^{l}_{k,l}}{a ( l + \frac{1}{2},k + \frac{1}{2} )}  + \sum^{l}_{k = - l} \frac{f_{l + \frac{1}{2}, k - \frac{1}{2}} \otimes \gamma^* t^{l}_{k,l}}{{c ( l, k - \frac{1}{2} )}} \right]\\
&=& \frac{a ( l + \frac{1}{2},l + \frac{1}{2} )}{a ( l,l )} \left[ \sum^{l + \frac{1}{2}}_{k = - l + \frac{1}{2}} \frac{a^{\prime} ( l, k - \frac{1}{2} )f_{l + \frac{1}{2},k} \otimes \alpha^* t^{l}_{k - \frac{1}{2},l}}{a ( l + \frac{1}{2},k )}  + \sum^{l - \frac{1}{2}}_{k = - l - \frac{1}{2}} \frac{f_{l + \frac{1}{2},k} \otimes \gamma^* t^{l}_{k + \frac{1}{2},l}}{c ( l,k )}\right] . 
\eean

Let $ - l + \frac{1}{2} \leq k \leq l - \frac{1}{2}.$ Then comparing coefficient of $ f_{l + \frac{1}{2},k} $ we have $ t^{l + \frac{1}{2}}_{k, l + \frac{1}{2}} = \frac{a ( l + \frac{1}{2},l + \frac{1}{2}  )}{a ( l,l )} [ \frac{a^{\prime} ( l,k - \frac{1}{2} )}{a ( l + \frac{1}{2},k )} \alpha^* t^{l}_{k - \frac{1}{2},l} + \frac{1}{c( l,k )} \gamma^* t^{l}_{k + \frac{1}{2},l}   ]$ which proves the first equation of (  \ref{rec_1}  ).

Applying the same procedure for $ k = - l - \frac{1}{2}, $ we have $ t^{l + \frac{1}{2}}_{k, l + \frac{1}{2}} = \frac{a_{l + \frac{1}{2},l + \frac{1}{2}}}{a_{l,l}} [ \frac{1}{c( l,k )} \gamma^* t^{l}_{k + \frac{1}{2},l}  ]$ which proves the second equation of (  \ref{rec_1}  ).

Similarly, for $ k = l + \frac{1}{2}, $ we have $ t^{l + \frac{1}{2}}_{k, l + \frac{1}{2}} = \frac{a ( l + \frac{1}{2},l + \frac{1}{2} )}{a( l,l )} [ \frac{a^{\prime} ( l,k - \frac{1}{2} )}{a ( l + \frac{1}{2},k )} \alpha^* t^{l}_{k - \frac{1}{2},l}   ] $ which proves the third equation of (  \ref{rec_1}  ). This completes the proof of  (  \ref{rec_1}  ).    

Next, to prove (  \ref{rec_2}  ), we apply coproduct on ( \ref{qorient_SU_mu_2_recursive_lemma_equation_1} ) and use ( \ref{qorient_SU_mu_2_recursive_lemma_equation_0}  ) and (  \ref{qorient_SU_mu_2_recursive_lemma_equation_2}    ) to have
 
 \bean \lefteqn{\sum^{l + \frac{1}{2}}_{k = -( l + \frac{1}{2}) } f_{l + \frac{1}{2},k} \otimes t^{l + \frac{1}{2}}_{k,i} }\\
 &=& c(l,i) ( \alpha \otimes \alpha - \mu \gamma^* \otimes \gamma ) ( \sum^{l}_{k = - l} f_{l,k} \otimes t^{l}_{k, i + \frac{1}{2}} )\\
 &=& c ( l,i ) ( \sum^{l}_{k = - l} \alpha f_{l,k} \otimes \alpha t^{l}_{k, i + \frac{1}{2}} - \sum^{n}_{k = - n} \gamma^* f_{l,k} \otimes \mu \gamma t^{l}_{k, i + \frac{1}{2}}    )\\
 &=& c( l,i ) \sum^{l}_{k = - l} \frac{1}{c ( l, k - \frac{1}{2} )} f_{l + \frac{1}{2}, k - \frac{1}{2}} \otimes \alpha t^{l}_{k, i + \frac{1}{2}} - c ( l,i ) \sum^{l}_{k = - l}  \frac{a^{\prime}( l,k )}{a ( l + \frac{1}{2}, k + \frac{1}{2})} f_{l + \frac{1}{2}, k + \frac{1}{2}} \otimes\\
 &&  \mu \gamma t^{l}_{k, i + \frac{1}{2}}\eean
 $$ = c( l,i ) \sum^{l - \frac{1}{2}}_{k = - l - \frac{1}{2}} \frac{1}{c ( l,k )} f_{l + \frac{1}{2},k} \otimes \alpha t^{l}_{k + \frac{1}{2}, i + \frac{1}{2}} - c( l,i ) \sum^{l + \frac{1}{2}}_{k = - l + \frac{1}{2}} \frac{a ( l, k - \frac{1}{2} )}{a ( l + \frac{1}{2}, k )} f_{l + \frac{1}{2},k} \otimes $$ 
 $ \mu \gamma t^{l}_{k - \frac{1}{2}, i + \frac{1}{2}} .$
 
For $ - l + \frac{1}{2} \leq k \leq l - \frac{1}{2} ,$ by comparing coefficient of $ f_{l + \frac{1}{2},k} $ we have $ t^{l + \frac{1}{2}}_{k,i} = \frac{c ( l,i  )}{c ( l, k )} \alpha t^{l}_{k + \frac{1}{2}, i + \frac{1}{2}} - \frac{c ( l,i ) a ( l, k - \frac{1}{2} ) }{a ( l + \frac{1}{2}, k  )} \mu \gamma t^{l}_{k - \frac{1}{2}, i + \frac{1}{2}}   $ which proves the first equation of (  \ref{rec_2}  ).
 
 Comparing coefficient of $ f_{l + \frac{1}{2}, - l - \frac{1}{2}}, $ we have $ t^{l + \frac{1}{2}}_{- l - \frac{1}{2},i} = \frac{c ( l,i )}{c ( l, - l - \frac{1}{2})} \alpha t^{l}_{- l, i + \frac{1}{2}} $ from which we get the second and the third equation of (  \ref{rec_2}  ).
 
 Comparing coefficient of $ f_{l + \frac{1}{2}, l + \frac{1}{2}}, $ we have $ t^{l + \frac{1}{2}}_{l + \frac{1}{2},i} = - \frac{c ( l,i ) a ( l,l) }{a ( l + \frac{1}{2},l + \frac{1}{2} )} \mu \gamma t^{l}_{l, i + \frac{1}{2}} $ from which we get the last equation of (  \ref{rec_2}  ).  \qed
 
 \vspace{4mm}

We recall the following multiplication rule from Page 74, \cite{klimyk} which we are going to need :

 \be \label{qorient_SUmu2_multiplication_rule_klimyk} t^{l}_{i,j} t^{1/2}_{{i}^{'},{j}^{'}} = \sum_{k = \left| l - 1/2 \right|,..... l + 1/2 } c_{k}( l,i,j,{{i}^{'},{j}^{'}} ) t^{k}_{ i + {i}^{'},j + {j}^{'}} \ee  
 ($c_k(l,i,j,i^\prime, j^\prime)$ are scalars).
 
 \subsection{The Hopf $\ast$-algebras $ \clo ( SU_{\mu} ( 2 ) ) $ and  $\clu_{\mu} ( su ( 2 ) ) $}

\label{preliminaries_subsection_U_mu_su2}

 We define the Hopf $ \ast $-algebra $ \clo ( SU_{\mu} ( 2 ) ) $  following the notations of \cite{klimyk}.
 
 $ {\clo} ( SL_{\mu} ( 2 ) ) $ is the  complex associative algebra with generators $ a,b,c,d $ such that 
 \be \label{SL_mu_2_definition} ab = \mu ba,~ ac = \mu ca,~ bd = \mu db,~ cd = \mu dc,~ bc = cb,~ ad - \mu bc = da - \mu^{- 1} bc = 1 . \ee
 
 The coproduct is given by 
 $$ \Delta ( a ) = a \otimes a + b \otimes c,~ \Delta ( b ) = a \otimes b + b \otimes d , $$ 
 $$ \Delta ( c ) = c \otimes a + d \otimes c, ~ \Delta ( d ) = c \otimes b + d \otimes d . $$ 
 The antipode is $$ \kappa ( a ) = d,~ \kappa ( b ) = - b ,~ \kappa ( c ) = - c,~ \kappa ( d ) = a $$ 
 Finally, the counit is $$ \epsilon ( a ) = \epsilon ( d ) = 1,~ \epsilon ( b ) = \epsilon ( c ) = 0 . $$ 
 For all $ \mu $ in $ \IR ,$ there is an involution of the algebra $ {\cal O} ( SL_{\mu} ( 2 ) ) $ determined by 
 \be \label{SL_mu_2_real form} a^* = d, ~~ b^* = - \mu c . \ee 
 The corresponding Hopf $ \ast $-algebra is denoted by $ {\cal O} ( SU_{\mu} ( 2 ) ). $
 
 \bppsn
 
 $ {\cal O} ( SU_{\mu} ( 2 ) ) $ can be identified with  $ {( SU_{\mu} ( 2 ) )}_0, $ i.e the Hopf $ \ast $-algebra generated by the matrix elements of irreducible unitary representations of $ SU_{\mu} ( 2 ) ,$ via the isomorphism given on the generators by 
 \be \label{SU_mu_2_woronowicz_klimyk_correspondence} \alpha \mapsto a, ~ \gamma \mapsto c, ~ \alpha^* \mapsto d,~ \gamma^*  \mapsto - \mu^{ - 1} b. \ee
 
 \eppsn

{\it Proof :} $ {( SU_{\mu} ( 2 ) )}_0 $ is generated by the matrix elements  of the fundamental unitary of  $ SU_{\mu} ( 2 ), $ that is, the $ \ast $-algebra generated by $ \alpha  $ and $ \gamma.$ On the other hand, inserting ( \ref{SL_mu_2_real form}  ) in (  \ref{SL_mu_2_definition} ), we have that $ {\cal O} ( SU_{\mu} ( 2 ) ) $ is generated by 4 elements $a,b,c,d$ such that  $ ac = \mu ca,~ ac^* = \mu c^* a,~ c c^* = c^* c,~ a^* a + c^* c = 1,~ a a^* + \mu^2 c^* c = 1 .$ 
Comparing with the defining equations of $ SU_{\mu} ( 2 ),$ that is, ( \ref{su2def1} ) - ( \ref{su2def5} ), it is clear that the above correspondence gives the required isomorphism. \qed 

\vspace{4mm}

Next, we recall from \cite{wagner} the  Hopf $ \ast $ algebra $ \clu_{\mu} ( su( 2 ) ) $ which is the dual Hopf $ \ast $-algebra of $ {\cal O} ( SU_{\mu} ( 2 ) ) .$ It is generated by elements $ F, E, K, K^{-1} $ with defining relations:
$$ K K^{-1} = K^{-1} K = 1,~ K E = \mu E K,~ F K = \mu K F, ~ E F - F E = ( \mu - \mu^{-1} )^{-1} ( K^2 - K^{-2} ) $$
with involution $ E^* = F, ~ K^* = K $ and comultiplication :
$$ \Delta ( E ) = E \otimes K + K^{-1} \otimes E,~ \Delta ( F ) = F \otimes K + K^{-1} \otimes F, ~ \Delta ( K ) = K \otimes K .$$
The counit is given by $ \epsilon ( E ) = \epsilon ( F ) = \epsilon ( K - 1 ) = 0 $ and antipode $ \kappa ( K ) = K^{- 1}, ~ \kappa ( E ) = - \mu E, ~ \kappa ( F ) = - \mu^{- 1} F. $

There is a dual pairing $ \left\langle . , . \right\rangle $ of $ \clu_{\mu} ( su( 2 ) ) $ and $ {\cal O} ( SU_{\mu} ( 2 ) ) $ given on the generators by :

 $ \left\langle K^{\pm 1} , \alpha^* \right\rangle = \left\langle K^{\mp 1} , \alpha \right\rangle = \mu^{\pm \frac{1}{2}}, ~  \left\langle E , \gamma \right\rangle = \left\langle F , - \mu \gamma^* \right\rangle = 1 $
 
 and zero otherwise.

The left action $ \triangleright $ and right action $ \triangleleft $ of $ \clu_{\mu} ( su( 2 ) ) $ on $ SU_{\mu} ( 2 ) $ are given by:

$ f \triangleright x =  \left\langle f , x_{(2)} \right\rangle x_{(1)},~ x \triangleleft f =  \left\langle f , x_{(1)} \right\rangle x_{(2)}, ~ x ~ \in ~ {\clo} ( SU_{\mu} ( 2 ) ), ~ f ~ \in ~ \clu_{\mu} ( su( 2 ) ) $ where we use the Sweedler notation $ \Delta ( x ) = x_{( 1 )} \otimes x_{( 2 )}.$

The actions satisfy :

$ {( f \triangleright x )}^{*} = {\kappa ( f )}^{*} \triangleright x^{*}, ~ {( x \triangleleft f )}^{*} = x^{*} \triangleleft {\kappa ( f )}^{*}, ~ f \triangleright xy = ( f_{(1)} \triangleright x ) ( f_{(2)} \triangleright y ), ~ xy \triangleleft f = ( x \triangleleft f_{(1)} ) ( y \triangleleft f_{(2)} ). $

The action on generators is given by :

\be \label{sphere_Uq_su2_left_action_su2}
 \left\{
 \begin{array}{cccc}
 E\triangleright \alpha = - \mu \gamma^* &  E  \triangleright \gamma = \alpha^*, & E \triangleright \gamma^* = 0, &  E        \triangleright \alpha^* = 0, \\ F \triangleright ( - \mu \gamma^* ) = \alpha & F \triangleright \alpha^* = \gamma, &  F  \triangleright \alpha = 0,  & F \triangleright \gamma = 0,\\  K \triangleright \alpha = \mu^{- \frac{1}{2}} \alpha, & K \triangleright ( \gamma^* ) = \mu^{\frac{1}{2}} \gamma^*,&  K \triangleright \gamma = \mu^{- \frac{1}{2}} \gamma,&  K \triangleright \alpha^* = \mu^{\frac{1}{2}} \alpha^* .
 \end{array}
 \right\} 
 \ee
 
\be \label{sphere_Uq_su2_right_action_su2}
 \left\{
 \begin{array}{cccc}
 \gamma \triangleleft E = \alpha, & ~ \alpha^{*}  \triangleleft E = - \mu \gamma^{*},&  \alpha \triangleleft E = 0, & \gamma^{*} \triangleleft E = 0 \\  \alpha \triangleleft F = \gamma,&   - \mu \gamma^* \triangleleft F =  {\alpha}^{*},& \gamma \triangleleft F = 0,& \alpha^* \triangleleft F = 0,\\  \alpha \triangleleft K = \mu^{- \frac{1}{2}} \alpha,& \gamma^* \triangleleft K = \mu^{ - \frac{1}{2}} \gamma^{*},& \gamma \triangleleft K = \mu^{ \frac{1}{2}} \gamma,&  \alpha^{*} \triangleleft K = \mu^{\frac{1}{2}} \alpha^{*} .\end{array}
 \right\} 
 \ee

 \subsection{The Wang algebras} 

 \label{preliminaries_subsection_Wang_algebras}
 
  Let us now recall the universal quantum groups as in
\cite{univ1}, \cite{free}
and references therein. 
For an $n \times n$ positive invertible matrix $Q=(Q_{ij})$. let $A_{u,n}(Q)$ be the compact quantum group defined and  studied in \cite{wang}, \cite{univ1}, which is the universal $C^{*}$-algebra generated by $ \{ u^{Q}_{kj}, k,j=1,...,n \}$ such that $u:=(( u^{Q}_{kj} ))$ satisfies \be \label{wangalg} u u^*=I_n =u^{*}u, ~~u^{\prime} Q  \overline{u} Q^{-1}=I_n=Q{\overline{u}} Q^{-1} u^{\prime}.\ee Here 
$u^{\prime} =(( u_{ji} ))$ and $\overline{u}=(( u_{ij}^{*} ))$. The coproduct, say $\tilde{\Delta}$, is given by, $$ \tilde{\Delta}(u_{ij})=\sum^n_{k = 1} u^Q_{ik} \ot u^Q_{kj}.$$ It may be noted that $A_{u,n}(Q)$ is the universal object in the category of compact quantum groups which admit an action on the finite dimensional $C^*$ algebra $M_n(\IC)$ which preserves  the functional $M_n \ni x \mapsto {\rm Tr({Q}^{T} x)}$,( see \cite{wangergodic} ) where the notion of a CQG and that of preservation of a functional by an action are as in subsection \ref{preliminaries_subsection_action_CQG}. We
refer the reader to \cite{univ1}  for a detailed  discussion on the structure and classification of
such quantum groups.

\brmrk

\label{preliminaries_wang_algebra_I_k2=I}

It was proved in \cite{free} that in the case where $ Q = I, ~ \kappa ( u^{I}_{ij} ) = u^{I*}_{ji} $ and hence $ \kappa^2 = {\rm id} $ holds for $ A_{u,n} ( I ). $  
                            
\ermrk

  \subsection{Action of a compact quantum group on a $ C^* $ algebra}

   \label{preliminaries_subsection_action_CQG}
   
   We say that  the compact quantum group $(\cls,\Delta)$ (co)-acts on a unital $C^*$ algebra $\clb$,
    if there is a  unital $C^*$-homomorphism (called an action) $\alpha : \clb \raro \clb \ot \cls$ satisfying the following :\\
    (bi) $(\alpha \ot {\rm id}) \circ \alpha=({\rm id} \ot \Delta) \circ \alpha$, and \\
    (bii) the linear span of $\alpha(\clb)(1 \ot \cls)$ is norm-dense in $\clb \ot \cls$.\\
     
It is known ( see, for example,  \cite{wang} , \cite{podles_subgroup} ) that  (bii) is equivalent to the existence of a norm-dense, unital $\ast$-subalgebra $\clb_0$ of $\clb$ such that $\alpha(\clb_0) \subseteq \clb_0 \ot_{\rm alg} \cls_0$ and on $\clb_0$, $({\rm id }\ \otimes \epsilon) \circ \alpha={\rm id}$.

We shall sometimes say that $ \alpha $ is a  `topological' or $ C^* $ action to distinguish it from a normal action of von Neumann algebraic quantum group.

    \bdfn
   
   Let $ ( \cls, \alpha ) $ has a $ C^* $ action $ \alpha $ on the $ C^* $ algebra $ \clb.$ We say that the action $ \alpha $ is {\bf faithful} if there is no proper Woronowicz $ C^* $-subalgebra $ \cls_1 $ of $ \cls $ such that $ \alpha $ is a $ C^* $ action of $ \cls_1 $ on $ \clb.$   
      
   \edfn
   
   \bdfn
   
   Let $ ( \cls, \alpha ) $ has a $ C^* $ action $ \alpha $ on the $ C^* $ algebra $ \clb.$ A continuous linear functional $ \phi $ on $ \clb $ is said to be {\bf invariant under $ \alpha $} if 
    $$ ( \phi \otimes {\rm id} ) \alpha ( b ) = \phi ( b ).1_{\cls}. $$ 
   
   \edfn
   
   Now, we recall the work of Shuzhou Wang done in \cite{wang}. One can also see \cite{ban4}, \cite{ban5}.
   
   The {\bf quantum permutation group} $ \clq \clu_{n}  $ is defined to be the $ C^* $ algebra generated by $ a_{ij} $ ( $ i,j = 1,2,...n $ ) satisfying the following relations:
   $$ a^2_{ij} = a_{ij} = a^*_{ij}, ~  i,j = 1,2,...n, $$
   $$ \sum^n_{j = 1} a_{ij} = 1, ~ i = 1,2,...n,  $$
   $$ \sum^n_{i = 1} a_{ij} = 1, ~ i = 1,2, ...n. $$
   
   The name comes from the fact that the universal commutative $ C^* $ algebra generated by the above set of relations is isomorphic to $ C ( S_n ) $ where $ S_n $ denotes the permutation group on $ n $  symbols.
   
 Let us consider the category with objects as compact groups acting on on a $ n $-point set $ X_n = \{ x_1, x_2, ..., x_n \} .$ If two groups $ G_1 $ and $ G_2 $ have actions $ \alpha_1 $ and $ \alpha_2 $ respectively, then a morphism from $ G_1 $ to $ G_2 $ is a group homomorphism $ \phi $ such that $ \alpha_2 ( \phi \times {\rm id} ) = \alpha_1.$ Then $ C ( S_n  ) $ is the universal object in this category. It is proved in \cite{wang} that the quantum permutation group enjoys a similar property.
 
We have that $ C ( X_n ) = C^* \{ e_i : e^2_i = e_i = e^*_i, ~ \sum^n_{r = 1} e_r = 1, i = 1,2, ...,n \}.$ Then $  \clq \clu_{n}   $ has a $ C^* $ action on $ C ( X_n ) $ via the formula:
 $$ \alpha ( e_j ) = \sum^n_{i = 1} e_i \otimes a_{ij}, j = 1,2, ...n. $$

   \bppsn
   
  Consider the category with objects as CQG s having a $ C^* $ action on $ C ( X_n ) $ and morphisms as CQG morphisms intertwining the actions as above. Then $ \clq \clu_{n}  $ is the universal object in this category.
   
   \eppsn 
   
   Now we note down a simple fact for future use.
   
   \blmma

\label{fact}

 Let $ \alpha $ be an action of a CQG $ \cls $ on $C(X)$ where $ X $ is a finite set.  Then $\alpha$ automatically preserves the functional $\tau$   corresponding to the counting measure:
$$ ( \tau \otimes {\rm id} )( \alpha ( f )) = \tau ( f ).1_{\cls}.$$  

   \elmma   

  {\it Proof:}\\  Let $X=\{ 1,...,n\}$ for some $n \in \IN$ and denote by $\delta_i$ the characteristic function of the point $i$. Let $\alpha(\delta_i)=\sum_j \delta_j \ot q_{ij}$ where $ \{q_{ij}:i,j=1\ldots n\}$  are the images of the canonical  generators of the quantum permutation group as above. Then  $\tau$-preservation of $\alpha$ follows from the properties of the generators of the quantum permutation group, which in particular imply that $\sum_j q_{ij}=1=\sum_i q_{ij}$. \qed   
    
   \vspace{4mm}

   Wang also identified the universal object in the category of all CQG s having a $ C^* $ action $ \alpha_1  $ on $ M_n ( \IC ) $ ( with morphisms as before ) such that the functional  $ \frac{1}{n} {\rm Tr} $ is kept invariant under $ \alpha_1 .$ However, no such universal object exists if the invariance of the functional is not assumed. The precise statement is contained in the following theorem.
   
   Before that, we recall that $ M_n ( \IC ) = C^* \{ e_{ij} : e_{ij} e_{kl} = \delta_{jk} e_{il}, ~ e^*_{ij} = e_{ji}, ~ \sum^n_{r = 1} e_{rr} = 1, i,j,k,l = 1,2,...n\}.$
   
   \bppsn
   
   \label{preliminaries_action_the_work_of_wang}
   
  Let $ \clq \clu_{M_n (\IC), {\rm \frac{1}{n} Tr } } $ be the $ C^* $ algebra with generators $ a^{kl}_{ij} $ and the following defining relations:
   $$ \sum^n_{v = 1} a^{kv}_{ij} a^{vl}_{rs} = \delta_{jr} a^{kl}_{is}, ~ i,j,k,l,r,s = 1,2,...,n, $$
   $$ \sum^n_{v = 1} a^{sr}_{lv} a^{ji}_{vk} = \delta_{jr} a^{si}_{lk}, i,j,k,l,r,s = 1,2,...,n, $$
   $$ {a^{kl}_{ij}}^* = a^{lk}_{ji}, ~ i,j,k,l =1,2,...,n, $$
   $$ \sum^n_{r = 1} a^{kl}_{rr} = \delta_{kl}, ~ k,l = 1,2, ...,n, $$
   $$ \sum^n_{r = 1} a^{rr}_{kl} = \delta_{kl}, ~ k,l = 1,...,n. $$
   Then,
   
  ( 1 )  $  \clq \clu_{M_n (\IC), {\rm \frac{1}{n} Tr } } $ is a CQG with coproduct $ \Delta $ defined by $ \Delta ( a^{kl}_{ij} ) = \sum^{n}_{r,s = 1} a^{kl}_{rs} \otimes a^{rs}_{ij}, ~ i,j,k,l = 1,2,...,n.  $
  
  ( 2 ) $  \clq \clu_{M_n (\IC), {\rm \frac{1}{n} Tr } } $ has a $ C^* $ action $ \alpha_1 $ on $ M_n ( \IC ) $  given by $ \alpha_1 ( e_{ij} ) = \sum^n_{k,l = 1} e_{kl} \otimes a^{kl}_{ij}, ~ i,j = 1,2, ...,n. $ Moreover, $  \clq \clu_{M_n (\IC), {\rm \frac{1}{n} Tr } } $ is the universal object in the category of all CQG s having $ C^* $ action on $ M_n ( \IC ) $ such that the functional $ \frac{1}{n} {\rm Tr} $ is kept invariant under the action.
  
  ( 3 ) There does not exist any  universal object in the category of all CQG s having $ C^* $ action on $ M_n ( \IC ) .$     
   
   \eppsn

   \bppsn
   
Since, any faithful state on a finite dimensional $ C^* $ algebra $ \cla $ is of the form $ {\rm Tr} ( R x ) $ for some operator $ R, $  it follows from  Theorem 6.1, ( 2 ) of  \cite{wang} that the universal CQG acting on any finite dimensional $ C^* $ algebra  preserving a faithful state $ \phi $ exists and is going to be denoted by $ \clq \clu_{\cla, \phi}.$
   
   \eppsn  
   
{\bf Notations:}  

 We conclude this section on quantum groups by fixing some notations which will be used throughout this thesis. In particular, given a compact quantum group $(\cls,\Delta)$, the dense unital Hopf $\ast$-subalgebra  of $\cls$ generated by the matrix elements of the irreducible unitary representations will be denoted by $ \cls_0. $ Moreover, given an action $\gamma : \clb \raro \clb \ot \cls$ of the compact quantum group $(\cls, \Delta)$  on a unital $C^*$-algebra $\clb$, the dense, unital  $\ast$-subalgebra of $\clb$ on which the action becomes an action by the Hopf $\ast$-algebra $\cls_{0}$ will be denoted by $ \clb_0. $ We shall use the Sweedler convention of abbreviating $\gamma(b) \in \clb_0 \ot_{\rm alg} \cls_{0}$ by $b_{(1)} \ot b_{(2)}$, for $b $ in $ \clb_0.$ This applies in particular to the canonical action of the quantum group $\cls$ on itself, by taking $\gamma=\Delta$. 
 
 Moreover, for a linear functional $f$ on $\cls$ and an element $c $ in $ \cls_0$ we recall the `convolution' maps $f \triangleleft c :=(f \otimes {\rm id} ) \Delta ( c )$ and $c \triangleright f := ({\rm id} \otimes f) \Delta ( c)$. We also define convolution of two functionals $f$ and $g$ by $(f \diamond g)(c)=(f \ot g)(\Delta(c)).$

   \section{Rieffel deformation}

   \label{preliminaries_section_Rieffel_deformation}
   
   In this section, we recall the notions of Rieffel's formulation of deformation quantization ( \cite{rieffel} ) as well as Rieffel type deformation of CQG s due to Rieffel and Wang ( as in \cite{toral} and \cite{wang_def} ).
   
  We begin with Rieffel deformation ( as in \cite{rieffel} ) from action of $ \IR^{n} $ on a $ C^* $ algebra. In the following discussion and henceforth, the symbol $ e ( x ) $ will stand for $ e^{2 \pi i x}.$
     Let V be a real vector space of dimension n and $ \alpha $ be its strongly continuous isometric action on a complex Frechet space $ \cla .$ Let $ \{ \left\|  \right\|_{j} \} $ denote the family of seminorms which determine the topology of $\cla.$ It is assumed that $ \alpha $ is isometric for each of the given seminorms on $\cla.$
     
     Let $ \cla^{\infty} $ denote the space of smooth vectors for the action $ \alpha ,$ that is, $ \cla^{\infty} = \{ a \in \cla : v \rightarrow \alpha_v ( a ) ~ \rm{is} ~ C^{\infty} \}.$
     
     Let $ \{ X_1, X_2, ........,X_n \} $ be a basis of V and $ \delta_k $ denotes the operator of partial differentiation on $ \cla^{\infty} $ in the direction of $ X_k.$
     
     For any multi index $ \mu = ( \mu_1, \mu_2,.....\mu_n  ) ,$ we will let $ \delta^{\mu} = \delta^{\mu_{1}} \delta^{\mu_{2}}...\delta^{\mu_{n}}, ~ \mu ! = \mu_1 ! ......\mu_n !, ~ \left| \mu \right| = \sum^{n}_{i = 1} \mu_i.$  We will equip $ \cla^{\infty} $ with the semi norms: $ \left\| a \right\|_{jk} = sup_{i \leq j} \sum_{ \left| \mu \right| \leq k } \frac{\left\| \delta^{\mu} a \right\|_{i} }{\mu !} .$ 
     
     Let $ C_b ( V, \cla ) $ denote the Frechet space of continuous bounded functions from $ V $ to $ \cla ,$ equipped with the semi norms  $ \left\| f \right\|_k = sup_{v \in V} \left\| f ( v ) \right\|_k.$ There is also a natural action of $ V $ on $ C_b ( V, \cla ) $ by translation and let $ C_u ( V, \cla ) $ denote the largest subspace of $ C_b ( V, \cla ) $  on which this action is strongly continuous and let $ \clb^{\cla} ( V ) $ denote the space of smooth vectors with respect to this action.

     Let $ W = V \times V.$ Then $ \clb^{\cla} ( W ) $ makes sense and for $ F $ in  $ \clb^{\cla} ( W ) ,$ one can define the oscillatory integral $ \int \int F ( u, v ) e ( u . v ) du dv $ ( where $ u.v $ denote the usual inner product  ) in the following way: 
     
     We choose a basis of $ W $ and let $ L $ denote the lattice of points of $ W $ which have integer co-ordinates w.r.t this basis. Moreover, choose a positive $ \phi_0 $ in $ C^{\infty}_c ( W ) $ such that $ \Phi = \sum_{p \in L} \phi_p $ vanishes nowhere on $ W$ where $\phi_p$ denotes the translate of $\phi$ by $p $ belonging to  $ L.$ Let $ \phi = \frac{\phi_0}{\Phi} .$

It can be shown ( \cite{rieffel} ) that for $ F $ in $ \clb^{\cla} ( W ),~  \sum_{p \in L} \int ( F \phi_p ) ( u, v ) e ( u. v ) du dv  $ converges absolutely in $ \cla $ and $ \int \int F ( u, v ) e ( u. v ) du dv $ is defined to be this sum. Moreover, this sum is independent of the choice of lattice and of $ \phi.$
Thus, \be \label{preliminaries_oscillatory_definition} \int \int F ( u, v ) e ( u. v ) du dv = \sum_{p \in L} \int ( F \phi_p ) ( u, v ) e ( u. v ) du dv . \ee
     
     For more details of oscillatory integral, we refer to \cite{rieffel} and references therein.
     
     We will need the following results from \cite{rieffel}.
     
     \bppsn
     
     \label{preliminiaries_rieffel_1.12}
     
  ( Corollary 1.12, \cite{rieffel} ) Let $ F $ be a function in $ \clb^{\cla} ( V \times V ) $ which depends only on the first variable, so that it is essentially an element of $ \clb^{\cla} ( V ) .$ Then $ \int \int F ( u ) e ( u. v ) du dv = F ( 0 ).$ The same is true if instead $ F $ depends only on the second variable.

     \eppsn

     \bppsn
     
     \label{preliminiaries_rieffel_1.14}
     
    ( Proposition 1.14, \cite{rieffel} ) Let $ S $ be a continuous linear transformation from $ \cla $ into a Frechet space $ C .$ Let $ F $ belongs to $ \clb^{\cla} ( W ).$ Then $ S \circ F $ belongs to $ \clb^{C} ( W ) $ and $ S ( \int \int F ( u, v ) e ( u.v ) du dv  ) = \int \int S ( F ( u, v ) ) e ( u.v ) du dv .$

     \eppsn

     Now, let $ \cla $ be a Frechet algebra. Fix a skew symmetric matrix $ J $ on $ V .$ Then  for all $ a, b $  in $ \cla^{\infty},~ \alpha_{Ju} ( a ) \alpha_v ( b ) $ belongs to $ \clb^{A} ( W ) $ and a new product $ \times_J $ is defined on $ \cla^{\infty} $ by declaring $ a \times_J b = \int_{V} \int_{V} \alpha_{Ju} ( a ) \alpha_v ( b ) e ( u.v ) du dv .$
     
     If the Frechet algebra has an involution $ \ast $ which is continuous and if $ \alpha $ acts by $ \ast $ automorphisms, then $ \ast $ is also an involution for the deformed product $ \times_J.$
     
     \bppsn
     
     \label{preliminiaries_rieffel_2.20}
     
  ( Lemma 2.20, \cite{rieffel} ) Let $ f, g $ belongs to $ \clb^{\cla} $ and let $ g $ have the lattice $ L $ as a period lattice so that $ g $ can be viewed as a smooth function on the compact group $ H = V / L .$ Then $ \int \int f ( u ) g ( v ) e ( u.v ) du dv = \sum_{L} f ( p ) (  \int_{H} g ( v ) e ( p.v ) dv  ) .$ A similar statement holds if instead it is $ f $ which is periodic.
     
     \eppsn

   \bcrlre
   
   \label{preliminiaries_rieffel_2.20_corollary}
   
   $$ \int \int e( \theta z_1 ) e(z_2) e(z_1. z_2) dz_1 dz_2 = e ( - \theta ). $$

   \ecrlre
   
   {\it Proof :} We have, \bean \lefteqn{ \int \int e( \theta z_1 ) e(z_2) e(z_1 . z_2) dz_1 dz_2}\\
                &=& \sum_{p \in \IZ} e( \theta p) ( \int_{S^{1}} e(z_2) e( p.z_2 ) dz_2  )\\
                &=&  \sum_{p \in \IZ} e ( \theta p ) \delta_{p,- 1}\\
                &=& e ( - \theta ).\eean \qed

   \vspace{4mm}
   
   Now, we will define the $ C^* $ algebra constructed by Rieffel corresponding to the data $ ( \cla,~ V, ~ \alpha, ~ J )$ where $\cla$ is also assumed to be a $C^*$ algebra and $\alpha$ a $C^*$ automorphism.
   
   Let $ \cls^{\cla} $ be the space of $ \cla $ valued smooth functions on $ V $ such that the product of their derivatives with any complex valued polynomials on $ V $ are bounded under the supremum norm of $ \cls^{\cla} .$  Then $ \cls^{\cla} $ is a pre Hilbert right  $ \cla $ module with $ \cla $ valued inner product defined by   
   $$  \left\langle f, ~ g  \right\rangle_{A} = \int_{V} {f ( v )}^* g ( v ) dv, $$   
   for $ f, g $ belonging to $ \cls^{\cla}.$
   
   Then, for $ a $ in $ \cla,$ one defines the operator $ L_{\widetilde{a}} $ on $ \cls^{\cla}  $ by   
   $$ L_{\widetilde{a}} ( f ) ( x ) = \int_{V} \int_{V} \alpha_{x + Ju} ( a ) f ( x + v ) e ( u.v ) du dv,$$   
   where $  f $ belongs to $ \cls^{\cla}.$ Then 
   
   \bppsn
   
   \label{deformation_content_of_theorem4.6_rieffel}
   
( Theorem 4.6 of \cite{rieffel} )  $ L_{\widetilde{a}}  $ is a bounded operator having an adjoint on the pre Hilbert module $ \cls^{\cla} $ and $ a \mapsto L_{\widetilde{a}} $ is a  $ \ast  $ representation of the algebra $ ( \cla^{\infty}, ~ \times_{J}  )  $  into the $ C^* $ algebra of bounded operators on $ \cls^{\cla}.$

    \eppsn
   
   Now, by defining   
        $$ \left\| a \right\|_{J} = \left\|  L_{\widetilde{a}} \right\|, $$        
   we have a pre-$ C^* $ norm $ \left\| \right\|_{J} $ on $ \cla^{\infty} $ endowed with the new product $ \times_{J}.$
   
   The completion of this pre $ C^* $ algebra is the deformed $ C^* $ algebra and is  denoted by $ \cla_J.$ 
  
  One has a natural Frechet topology on $ \cla^{\infty}_{J} ,$ given by a family of seminorms $ \{ \left\|  \right\|_{n,J} \} $ where $  \left\| a \right\|_{n,J} = \sum_{\left| \mu \right| \leq n} \frac{\left\| \alpha_{X^{\mu}} ( a ) \right\|_{J} }{\mu !} $
  
  We recall the following Proposition from \cite{rieffel}.
  
  \bppsn
  
  \label{preliminiaries_rieffel_4.10}
  
 ( Proposition 4.10, \cite{rieffel} ) Let $ J $ be fixed. Then for large enough $ k $ there is a constant $ c_k $ such that for all $ a $ in $ \cla^{\infty},$ we have $ \left\| a \right\|_{J} \leq c_k \left\| a \right\|_{2k}. $
  
  \eppsn
  
  \bppsn
  
  \label{preliminiaries_rieffel_7.1}
  
  ( Proposition 7.1, \cite{rieffel} ) Let $ \alpha $ be an action of $ V $ on the $ C^* $ algebra $ \cla, $ with $ \cla^{\infty} $ its subalgebra of smooth vectors. Let $ J $ be a skew-symmetric operator on $ V, $ and let $ \alpha $ also denote the corresponding action of $ V $ on $ \cla_J.$ Then the subalgebra of smooth vectors in $ \cla_{J} $ for $ \alpha $ is exactly $ \cla^{\infty}.$ Moreover, $ {( \cla_J )}_{- J} \cong \cla.$
  
  \eppsn

 \bcrlre
 
 \label{preliminiaries_rieffel_4.10_corollary}
 
 $ \cla^{\infty} $ and $ \cla^{\infty}_{J}$ coincide as  topological ( Frechet ) spaces.
 
 \ecrlre

{\it Proof :} The proof is essentially contained in the proof of Proposition 7.1 in \cite{rieffel} ( Proposition \ref{preliminiaries_rieffel_7.1} above ). By Proposition \ref{preliminiaries_rieffel_4.10}, we know that there is a constant $ c_k $ such that for any $ a $ in $ \cla^{\infty} $ and for any $ \mu ,$
$$ \left\| X^{\mu} a \right\|_{J} \leq c_k \left\| X^{\mu} a \right\|_{2k} \leq c^{\prime}_k \left\| a \right\|_{j}  $$
for $ j = \left| \mu \right| + 2k $ and a new constant $ c^{\prime}_k. $ Thus, the inclusion of $ \cla^{\infty} $ into $ \cla^{\infty}_{J} $ is continuous for their Frechet topologies. Similarly, using $ {( \cla_{J} )}_{- J} = \cla, $ we deduce that the inclusion of $ \cla^{\infty}_{J} $ into $ \cla^{\infty} $ is continuous. This proves the result. \qed   

\vspace{4mm} 

{\bf Examples}

\vspace{4mm}

 {\bf  The Noncommutative Torus}
 
 Let $ \cla = C ( \IT^n ).$ For $ v = ( v_1, v_2, ..., v_n ) $ in $ \IR^n, ~ x = ( x_1, x_2,..., x_n ) $ in $ \IT^n, ~ f $ in  $ C ( \IT^n ), $ the action $ \alpha $ of $ \IR^n $ on $ \cla $ is given by $ \alpha_v f ( x ) = f  ( x_1 e( v_1 ), x_2 e( v_2 ), ... x_n e ( v_n ) ).$ 
  Let $ \theta $ be a $ n \times n $ skew symmetric matrix and $ J = \frac{\theta}{2}.$ Then $ \cla_{J} $ can be seen to be equal to the noncommutative $ n $ tori $ \IT^n_{\theta}, $ that is the universal $ C^* $ algebra generated by unitaries $ U_i,~ i = 1,2,...,n $ satisfying $ U_i U_j = e ( \theta_{ij} ) U_j U_i $ where $ \theta_{ij} $ denotes the $ ( i,j ) $ th entry of the matrix $ \theta.$ We will denote $ \IT^2_{\theta} $ by the notation $ \cla_{\theta}.$
 
 \vspace{4mm}
 
 {\bf The Rieffel deformed spheres}
 
For a skew symmetric matrix $ \theta ,$ we recall from \cite{connes_etal}, the definition of $ S^{n}_{\theta} .$

Let $ \mu,~ \nu = 1,2,.... n.$ Let $ \lambda^{\mu \nu} = e ^{i \theta_{\mu \nu}} $ where $ \theta_{\mu \nu} $ is the $ ( \mu, \nu )  $ th entry of the matrix $ \theta.$ 


 
 
 

  
  
  
  
   
   
   

$ S^{2n - 1}_{\theta} $ is the universal $C^*$ algebra generated by $ 2n $ elements $ z^{\mu}, ~ {\overline{z}}^{\mu}$ with relations:

 \be \label{preliminaries_Rieffel_deform_eg_S_2n-1_theta_1} z^{\mu} z^{\nu} = \lambda^{\mu \nu} z^{\nu} z^{\mu}, ~~~ {\overline{z}}^{\mu} {\overline{z}}^{\nu} = \lambda^{\mu \nu} {\overline{z}}^{\nu} {\overline{z}}^{\mu} ,\ee 
 \be \label{preliminaries_Rieffel_deform_eg_S_2n-1_theta_2} {\overline{z}}^{\mu} z^{\nu} = \lambda^{\nu \mu} z^{\nu} {\overline{z}}^{\mu},\ee 
 \be  \label{preliminaries_Rieffel_deform_eg_S_2n-1_theta_3} {( z^{\mu} )}^{*} = {\overline{z}}^{\mu}, \ee
 \be \label{preliminaries_Rieffel_deform_eg_S_2n-1_theta_4} \sum^{n}_{\mu = 1} z^{\mu} {\overline{z}}^{\mu} = 1.\ee
It can be easily seen that $ S^{2n - 1}_{\theta}  $ is obtained by the Rieffel deformation of  $ C ( S^{2n - 1} ) $ using the $ 2n \times 2n $ skew symmetric matrix $ J $ whose $ ( \mu, \nu ) $ th entry is $  \frac{\lambda^{\mu,\nu}}{2}$ and the $ \IR^{2n} $ action on $ C ( S^{2n - 1} ) $ given by $ \alpha_v f ( x_1,..., x_{2n} ) = f ( x_1 e( v_1 ),..., x_{2n} e ( v_{2n} ) ) $ ( $ v = ( v_1,..., v_{2n} ) $ is in $ \IR^{2n}, ~ f $ is in  $ C^{\infty} ( S^{2n - 1} ) $ ).
  
 $  S^{2n}_{\theta}  $ is the universal $C^*$ algebra generated by $ 2n + 1 $ elements $ \{ z^{\mu}, ~ {\overline{z}}^{\mu}, x : \mu = 1,2,..,n \} $ where $z^{\mu}, ~ {\overline{z}}^{\mu} $ satisfy ( \ref{preliminaries_Rieffel_deform_eg_S_2n-1_theta_1} ) - ( \ref{preliminaries_Rieffel_deform_eg_S_2n-1_theta_3} ),    $x$ is a self adjoint element satisfying the relations $ x z^{\mu} = z^{\mu} x ~  $ for all $ \mu = 1,2,...n$ and $ \sum^{n}_{\mu = 1} z^{\mu} {\overline{z}}^{\mu} + x^2 = 1.$
 
$S^{2n}_{\theta}$ is the Rieffel deformation of $ C ( S^{2n} ) $ by the action of $\IR^{2n + 1}$ on $ C ( S^{2n} ) $ similar to above and a $ ( 2n + 1 ) \times ( 2n + 1 ) $ matrix $J^{\prime}$ such that $ {( J^{\prime} )}_{\mu, \nu} = \frac{\lambda^{\mu, \nu}}{2} $ if $\mu \leq 2n, ~ \nu \leq 2n$ and $0$ otherwise.


     
   \subsection{Rieffel Deformation  of compact quantum group}

   \label{preliminaries_subsection_Rieffel_deformation_CQG}
   
   Here we describe the Rieffel deformation of a CQG as in \cite{wang_def}.
   
   Let $ ( \cla, \Delta ) $ be a CQG with $ C ( \IT^{n} ) $ as a quantum subgroup. Let $ \pi $ be the corresponding CQG morphism from $ \cla $ to $ C ( \IT^{n} ).$
   
   Let $ \eta $ be the canonical homomorphism from $ \IR^n $ to $ \IT^n $ given by $ \eta ( x_1, x_2,......, x_n ) = ( e ( x_1 ), e ( x_2 ),..... e ( x_n ) ) $ and $ ev_x $ be the state on $  C ( \IT^{n} ) $ obtained by evaluation of a function at the point $ x $ in $ \IT^{n}.$  
   
   Now, put   
   \be \label{preliminiaries_deformation_wang_lambda} \lambda_{\eta ( s )} = ( ev_{\eta( - s )} \pi \otimes {\rm id} ) \Delta,   \ee 
   \be \label{preliminiaries_deformation_wang_rho} \rho_{\eta( u )} = ( {\rm id} \otimes ev_{\eta( u )} \pi ) \Delta .  \ee 
We will use the notation $ \Omega ( u ) $ for $ ev_{\eta( u )} \pi.$   
      
   Then there is a $ \IR^{2n} $ action on $ \cla $ defined by    
   \be \label{preliminiaries_deformation_wang_formula}  \chi_{(s,u)} = \lambda_{\eta ( s )} \rho_{\eta ( u )}.   \ee   
   Fix a skew symmetric matrix $ J $ on $ \IR^{n}$ and put     
    \bean \widetilde{J} = J \oplus ( - J ). \eean    
    Then, by the prescription of Rieffel as described above, we have a $ C^* $ algebra $ \cla_{\widetilde{J}}.$  Shuzhou Wang showed in \cite{wang_def} that $ \cla_{\widetilde{J}}$ can be made into a CQG.
    
     The $ \ast $-algebra generated by the matrix elements of unitary irreducible representations of $\cla $ ( denoted by $ \cla_0 $ ) is dense in the space $ \cla^{\infty} $ of smooth vectors of the action $ \chi $ under the Frechet topology and hence is dense in the $ C^* $ algebra $ \cla_{\widetilde{J}} $ under the $ C^* $ norm of $ \cla_{\widetilde{J}} .$ On  $ \cla_0,$ the Hopf  $ \ast $-algebra structure remains unchanged and this extends to a CQG structure on $ \cla_{\widetilde{J}}.$  
     
     We quote the following result ( Remark 3.10 ( 2 ), \cite{wang_def} ) which will be used later.
     
     \bppsn
     
     \label{preliminiaries_deformation_wang_haar_state_same}
     
     The Haar measure $ h_{\widetilde{J}} $ of $ \cla_{\widetilde{J}} $ is still the same as the Haar measure on the common subspace $ \cla_0 .$

     \eppsn
     
     \blmma           
         \label{deformationLemma7}
         
     The Haar state (say $h$) of $\cla$ coincides with the Haar state on $\cla_{\widetilde{J}}$ ( say $ h_{J} $ ) on the common subspace $\cla^{\infty}$, and moreover, $h(a \times_{\widetilde{J}} b)=h(ab)$ for $a,b $ in $ \cla^{\infty}.$
         \elmma
         
      {\it Proof :} We recall ( Proposition \ref{preliminiaries_deformation_wang_haar_state_same} ) that $h$ = $h_{J}$ on $\cla_{0}$. 
                       By using  $ h ( \Omega ( - s ) \otimes {\rm id} ) = \Omega ( - s )( {\rm id} \otimes h ) $ and $ h( {\rm id} \otimes \Omega( u ) ) = \Omega ( u ) ( h \otimes {\rm id} ) $, we  have for $a $ in $ \clq_0,$  
 \bean \lefteqn { h  ( \chi_{s,u} ( a ) ) }\\
            & = & \Omega ( - s ) ( {\rm id} \otimes h ) \Delta ( {\rm id} \otimes \Omega ( u ) ) \Delta ( a ) \\
            & = & \Omega ( -s ) ( h ( ( {\rm id} \otimes \Omega ( u ) ) \Delta ( a ) )1 ) \\
            & = & h ( ( {\rm id} \otimes \Omega ( u ) ) \Delta ( a ) ) \\
            & = & \Omega ( u ) ( h( a ).1 ) \\
            & = & h( a ) .\eean             
    Therefore,  \be \label{preliminaries_deformation_h_chi=h} h \chi_{s,u} ( b ) = h( b ) ~ {\rm ~ for ~ all ~  b  ~  in} ~  \clq_0 .\ee  
      Now,
      \bean
      \lefteqn{ h( a {\times}_{\widetilde{J}} b )}\\
      & =& \int\int h( \chi_{\widetilde{J}u} ( a ) \chi_{v} ( b ) ) e( u . v ) du dv \\
      & =& \int\int h( \chi_{v} ( \chi_{\widetilde{J}u - v} ( a ) b ) ) e( u . v ) du dv \\
      &=& \int\int h( \chi_{t} ( a ) b ) e( s.t ) ds dt ,\eean where $ s = - u, t = \widetilde{J}u - v ,$
   which by Proposition \ref{preliminiaries_rieffel_1.12} equals $ h( \chi_{0} ( a ) b ) = h( a b ). $  That is, we have proved 
   \be \label{etaseta} \lgl a, b \rgl_J=\lgl a, b \rgl ~~\forall a,b \in \clq_0, \ee
      where $\lgl \cdot, \cdot \rgl_J$ and $\lgl \cdot, \cdot \rgl$ respectively denote the inner products of $L^2(h_J)$ and $L^2(h)$. 
    We now complete the proof of the lemma by extending ( \ref{etaseta} ) from   $\clq_0$ to $\clq^\infty$, by using the fact that $\clq^\infty$ is a common  subspace of the Hilbert spaces $L^2(h)$ and $L^2(h_J)$ and moreover, $\clq_0$ is  dense in both these Hilbert spaces. In particular, taking $a=1 $ in $ \clq_0,$ we have $h=h_J$ on $\clq^\infty$.   
      \qed

    \brmrk
    
    \label{haarrem}
    
    Lemma \ref{deformationLemma7} implies in particular that for every fixed $a_1,a_2 $ in $ \clq_0,$ the functional $\clq_0 \ni b \mapsto h(a_{1} \times_{\widetilde{J}} b \times_{\widetilde{J}} a_2)= h ( b \times_{\widetilde{J}} a_2 \times_{\widetilde{J}} ( f_1 \triangleleft a_1 \triangleright f_1 ) ) $ ( where $ f_1 $ is as in Remark \ref{preliminaries_haar_state_not_tracial} ) $ =  h( b ( a_2\times_{\widetilde{J}}( f_1 \triangleleft a_1 \triangleright f_1 )) )$ extends to a bounded linear functional on $\clq.$

\ermrk 

Let $e$ be the identity of $ \IT^{2n} $ and $U_{n}$ be a sequence of neighbourhoods of $e$ shrinking to $e$, $f_{n}$  smooth, positive functions with support contained inside $U_{n}$ such that $ \int_{\IT^{2n}} f_{n}(z) dz = 1 $  for all $n$.  
    
Let us denote the action of $\IT^{2n}$ action on $\cla^{\infty}$ induced by $ \chi $ by $\widetilde{\chi}.$  Define $ \lambda_{f_{n}} ( a ) = \int_{\IT^{2n}} {\widetilde{\chi}}_{z} ( a ) f_{n}( z ) dz .$ Then, we have the following result:

\blmma

\label{deformation_laplacian_haar state faithful_0}

$ \lambda_{f_n} ( a ) $ belongs to $ \clq^{\infty} $ and
\bean \int_{\IT^{2n}} {\widetilde{\chi}}_{z} ( a ) f_n ( z ) dz \rightarrow a  ~{\rm as} ~  n \rightarrow \infty.\eean

\elmma   
    
{\it Proof :} We note that, by using the translation invariance of Haar measure, for all $ g $ in $ \IT^{2n}, ~ {\widetilde{\chi}}_g ( \lambda_{f_n} ( a ) ) = \int_{\IT^{2n}} f_n ( g^{- 1} h ) {\widetilde{\chi}}_h ( a ) dh. $ Therefore, $ {\widetilde{\chi}}_g ( \lambda_{f_n} ( a ) ) - \lambda_{f_n} ( a ) = \int_{\IT^{2n}} ( f_n ( g^{- 1} h ) - f_n ( h )  ) \widetilde{\chi}_h ( a ) dh $ which proves the first part.

Now we prove the second part. As $ \int_{\IT^{2n}} f_{n} ( z ) = 1 $  for all $ n $ and supp $ ( f_n ) \subseteq U_n ,$ we have                
           \bean \lefteqn{ \left\| \int_{\IT^{2n}} {\widetilde{\chi}}_{z} ( a ) f_n ( z ) - a  \right\|}\\
                  &=& \left\| \int_{\IT^{2n}} {\widetilde{\chi}}_{z} ( a ) f_n ( z ) - a \int_{\IT^{2n}} f_n ( z ) dz  \right\| \\
                  &=& \left\|  \int_{\IT^{2n}} ( {\widetilde{\chi}}_{z} ( a ) - a ) f_n( z ) dz \right\| \\
                  &=& \left\|  \int_{U_n} ( {\widetilde{\chi}}_{z} ( a ) - a ) f_n( z ) dz  \right\| .\eean 

 Now, using the fact that the map  $ z \mapsto {\widetilde{\chi}}_{z} ( a ) $ is continuous for all $ a ,$ we deduce that  for all $ \epsilon >0, $ there exists $ n $ such that  for all $ z $ in $ U_n,~ \left\| {\widetilde{\chi}}_z ( a ) - \widetilde{\chi}_0 ( a ) \right\| < \epsilon,$ that is, $ \left\| {\widetilde{\chi}}_z ( a ) - a  \right\| < \epsilon.$
 
 Hence,$ \left\| \int_{\IT^{2n}} \lambda_{z} ( a ) f_n ( z ) dz  - a \right\| \leq \epsilon  \int_{\IT^{2n}} f_{n} ( z ) = \epsilon $ which proves the lemma. \qed

 \blmma
     \label{deformation_laplacian_haar state faithful_1}
     
   If $h$ is faithful on $ \clq $, then $ h_{J} $ is faithful on $\clq_{\widetilde{J}}$.
   
   \elmma
   
     {\it Proof :} Let $a \geq 0, \in \clq_{\widetilde{J}} $ be such that $ h_J( a ) = 0$. Let $ \lambda_{f_n} $ be as defined above.
     
 Then,  \bean \lefteqn {h_{J}( \lambda_{f_{n}} ( a ) )} \\
                       & = & \int_{\IT^{2n}} h_{J} ( \widetilde{\chi}_{z} ( a ) ) f_{n} ( z ) dz \\
                       & = & \int_{\IT^{2n}} h_{J} ( a ) f_{n} ( z ) dz \\
                       && ( ~ by ~ ( \ref{preliminaries_deformation_h_chi=h} ) )\\ 
                       & = & 0, \eean
           so we have
                      $ h( \lambda_{f_{n}} ( a ) ) = 0 $, since $h$ and $h_J$ coincide on $\clq^\infty$ by Lemma \ref{deformationLemma7} and $ \lambda_{f_{n}} ( a )  $ belongs to $ \clq^{\infty}.$
         
 Now we fix some notation which we are going to use in the rest of the proof. Let $ L^{2}( h ) $ and $ L^{2}( h_{J} ) $ denote the G.N.S spaces of $ \clq $ and $ \clq_{\widetilde{J}} $ respectively with respect to the Haar states. Let $ i $ and $ i_{J} $ be the canonical maps from $ \clq $ and $ \clq_{\widetilde{J}} $ to  $ L^{2}( h ) $ and $ L^{2}( h_{J} ) $ respectively. Also, let $ \Pi_{J} $ denote the G.N.S representation of $ \clq_{J} $.
  Using the facts  $ h ( b^{\ast} \times_{\widetilde{J}} b ) = h ( b^{\ast} b ) $  for all $ b $ in $ \clq^{\infty} $ and $ h = h_{J} $ on $ \clq^{\infty} $ ( by Lemma \ref{deformationLemma7} ) ,  we get $ \left\|i_{J}( b ) \right\|^{2}_{L^{2}( h_{J} )} = \left\|i( b ) \right\|^{2}_{L^{2}( h ) } $ for all $ b $ in $ \clq^{\infty} .$ So the map sending $ i ( b ) $ to $ i_{J} ( b ) $ is an isometry from a dense subspace of $ L^{2}( h )$ onto a dense subspace of $ L^{2}( h_{J} ) $, hence it extends to a unitary, say $ \Gamma : L^{2} ( h ) \rightarrow L^{2} ( h_{J} ) $. We also note that the maps $  i $ and $ i_{J} $ agree on $ \clq^{\infty} $.
  
  Now, $ a \geq 0 $ means that $ \lambda_{f_n} ( a ) $ is positive in $ \clq_{\widetilde{J}}$ and therefore, $ \lambda_{f_{n}} ( a ) = b^{\ast} \times_{\widetilde{J}} b $ for some $ b $ in $ \clq_{\widetilde{J}} $. So $ h( \lambda_{f_{n}} ( a ) ) = 0 $ implies $ h_J ( \lambda_{f_{n}} ( a ) ) = 0 $ and therefore $ \left\|i_{J} ( b ) \right\|^{2}_{L^{2} ( h_{J} ) } = 0 $. Therefore, one has  $ \Pi_{J}( b^{\ast} ) i_{J} ( b ) = 0$, and hence $ i_{J}( b^{\ast} b ) = i_{J} ( \lambda_{f_{n}} ( a ) ) = 0$. It thus follows that $ \Gamma ( i ( \lambda_{f_{n}}( a ) ) ) = 0$, which implies  $ i ( \lambda_{f_{n}}( a ) ) = 0$. But the faithfulness of $h$  means that $ i $ is one one, hence $ \lambda_{f_{n}}( a ) = 0$ for all $n$.
 Thus, recalling Lemma \ref{deformation_laplacian_haar state faithful_0}, we have    $a=lim_{n \rightarrow \infty} \lambda_{f_{n}} ( a ) =0 $, which proves the faithfulness of $h_J$. \qed

 At this point, we note a useful implication of the Lemma \ref{deformationLemma7}. Let us make use of the identification of $\clq_0$ as a common vector-subspace of all $\clq_{\widetilde{J}}$. To be precise, we shall sometimes denote this identification map from $\clq_0 $ to $\clq_{\widetilde{J}}$ by $\rho_J$. 
\bcrlre
\label{repdef}
Let $W$ be a finite-dimensional (say, $n$-dimensional) unitary representation of $\clq$, with $\widetilde{W} $ belonging to $ M_ n(\IC) \ot \clq_0$ be the corresponding unitary. Then, for any $J$, we have that  $\widetilde{W}_J:=({\rm id} \ot \rho_J)(\widetilde{W}) $ is unitary in $\clq_{\widetilde{J}}$, giving a unitary $n$-dimensional representation of $\clq_{\widetilde{J}}$.   In other words, any finite dimensional unitary representation of $\clq$ is also a unitary representation of $\clq_{\widetilde{J}}$.
\ecrlre
     {\it Proof:}\\
 Since the coalgebra structures of $\clq$ and $\clq_{\widetilde{J}}$ are identical, and $\widetilde{W}_J$ is identical with $\widetilde{W}$ as a linear map, it is obvious that $\widetilde{W}_J$ gives a nondegenerate representation of $\clq_{\widetilde{J}}$. Let $y=({\rm id} \ot h)(\widetilde{W}_J^* \widetilde{W}_J)$. 
It follows from the proof of Proposition 6.4 of \cite{vandaelenotes} that $y$ is invertible positive element of $M_n$ and $(y^{\frac{1}{2}} \ot 1) \widetilde{W}_J (y^{-\frac{1}{2}} \ot 1)$ gives a unitary representation of $\clq_{\widetilde{J}}$. We claim that $y=1$, which will complete the proof of the corollary.
For convenience, let us write $W$ in the Sweedler notation: $W=w_{(1)} \ot w_{(2)}$. We note that by Lemma \ref{deformationLemma7}, we have \bean \lefteqn{({\rm id} \ot h)(\widetilde{W}_J^*\widetilde{W}_J)}\\
&=& w_{(1)}^*w_{(1)} h(w_{(2)}^* \times_{\widetilde{J}} w_{(2)})\\
&=& w_{(1)}^*w_{(1)} h(w_{(2)}^* w_{(2)})\\
&=& ({\rm id} \ot h)(\widetilde{W}^* \widetilde{W})=({\rm id} \ot h)(1 \ot 1)=1.\eean
\qed

   {\bf Example ~ The Rieffel deformed orthogonal groups}

Let $ \theta $ be a $ n \times n $ skew symmetric matrix. $ C ( \IT^n ) $ sits inside $ C ( O ( n ) ) $ as a quantum subgroup. It can be easily seen that $ O_{\theta} ( n ) $  is obtained by Rieffel deformation from $ C ( O ( n ) ) $ by using the induced $ \IR^{2n} $ action as given in the equation ( \ref{preliminiaries_deformation_wang_formula} ) and considering the matrix $ \widetilde{J} =  - J \oplus  J$ when $ n $ is even and $- J^{\prime} \oplus J^{\prime} $ when $n$ is odd where $J$ and $ J^{\prime} $ are the matrices introduced while giving the definition of the $\theta$ deformed spheres.

\section{Classical Riemannian geometry}

\label{preliminaries_section_Riemannian_geometry}

In this section we recall some classical facts regarding manifolds which will be useful to us later on.

\subsection{Classical Hilbert space of forms}

\label{preliminaries_subsection_classical_Hilbert_space_forms}

Let $ M $ be an $ n $ dimensional Riemannian manifold and $ \Omega^{k} ( M ) $ ( $ k = 0,1,2,...n $ ) be the space of smooth $k$-forms. Set $ \Omega^{k} ( M ) = \{0\} $ for $ k > n.$ The de-Rham differential $ d $ maps $ \Omega^k ( M ) $ to $ \Omega^{k + 1} ( M ).$ Let $ \Omega \equiv \Omega ( M )  = \oplus_k \Omega^k ( M ).$ We will denote the Riemannian volume element by dvol. We recall that the Hilbert space $ L^2 ( M ) $ is obtained by completing the space $ \{ f \in C^{\infty}_c ( M ) \} $ with respect to the pre-inner product given by $ \left\langle f_1, f_2 \right\rangle = \int_{M} \overline{f_1} f_2 dvol.$ 

In an analogous way, one can construct a canonical Hilbert space of forms. The Riemannian metric $ \left\langle . ~ , ~ . \right\rangle_{m} $ ( for $ m $ in $M$ ) on $ T_m M $ induces an inner product on the vector space $ T^*_m M $ and hence also $ \Lambda^k T^*_m M ,$ which will be again denoted by $ \left\langle . ~, ~ . \right\rangle_{m}.$ This gives a natural pre-inner product on the space of compactly supported $k$-forms by integrating the compactly supported smooth function $ m \mapsto \left\langle \omega( m ), \eta ( m ) \right\rangle_{m} $ over $ M. $ We will denote the completion of this space by $ \clh^k ( M ).$ Let $ \clh = \oplus_k \clh^k ( M ).$ 

Then, one can view $ d : \Omega \rightarrow \Omega $ as an unbounded, densely defined operator ( again denoted by $ d $ ) on the Hilbert space $ \clh $ with the domain $ \Omega.$ It can be verified that it is closable.

\subsection{Isometry groups of classical manifolds}

\label{preliminaries_subsection_Isometry_group_manifolds}
 
 Let $ M $ be a Riemannian manifold of dimension $n.$ Then the collection of all isometries of $ M $ has a natural group structure and is denoted by $ISO ( M ).$ 
 Let $ C $ and $ U $ be respectively a compact and open subset of $ M $ and let $ W ( C, U ) = \{ h \in ISO ( M ) : h.C \subseteq U \}.$ The compact open topology on $ ISO ( M ) $ is the smallest topology on $ ISO ( M ) $ for which the sets $  W ( C, U ) $ are open. It follows ( see \cite{helgason_symmetric_space} ) that under this topology, $ ISO ( M ) $ is a closed locally compact topological group. Moreover, if $ M $ is compact, $ ISO ( M ) $ is also compact.
 
We recall that the Laplacian $ \cll $ on $ M $ is an unbounded densely defined self adjoint operator $ - d^* d $ on the space of zero forms $ \clh^0 ( D ) = L^2 ( M, {\rm dvol} )$ which has the local expression 
$$ \cll ( f ) = - \frac{1}{\sqrt{det( g )}} \sum^n_{i,j = 1} \frac{\partial}{\partial x_j} (  g^{ij} \sqrt{{\rm det} ( g )} \frac{\partial}{\partial x_i} f ) $$ 
for $ f $ in $ C^{\infty} ( M ) $ and where $ g = (( g_{ij} )) $ is the Riemannian metric and $ g^{- 1} = (( g^{ij} )). $ It is well known that on a compact manifold, the Laplacian has compact resolvents. Thus, the set of eigenvalues of $ \cll $ is countable, each having finite multiplicities, and accumulating only at infinity. Moreover, there exists an orthonormal basis of $ L^2 ( M ) $ consisting of eigenvectors of $ \cll$ which  belong to $ C^{\infty} ( M ) .$ It can be shown ( Lemma 2.3 of \cite{goswami} ) that for a compact manifold, the complex linear span of the eigenvectors of  $ \cll $ is dense in $ C^{\infty} ( M )$ in the sup norm.

\vspace{2mm}
 
 The following result is in the form in which it has been stated and proved in \cite{goswami} ( Proposition 2.1 ).
     
     \bppsn
     
     \label{preliminiaries_isometry_laplacian_characterization}
     
     Let $ M $ be a compact Riemannian manifold. Let $ \cll $ be the Laplacian of $ M.$ A smooth map $ \gamma: M \rightarrow M $ is a Riemannian isometry if and only if $ \gamma $ commutes with  $ \cll $ in the sense that $ \cll ( f \circ \gamma ) = ( \cll ( f ) ) \circ \gamma  $ for all  $ f $ in $ C^{\infty} ( M ).$    
      
     \eppsn
     
Using this fact, we give an operator theoretic proof of the fact that for a compact manifold, $ ISO( M ) $ is compact. As the action of $ ISO ( M ) $  commutes with the Laplacian, it has a unitary representation on $ L^2 ( M ).$ As the action preserves the finite dimensional eigenspaces of the Laplacian, $ ISO ( M ) $ is a subgroup of $ U ( d_1 ) \times U ( d_2 ) \times ..........  $( where $ \{ d_i: i \geq 0 \} $ denote the dimension of the eigenspaces of the Laplacian and $ U ( d ) $ denotes the group of unitary operators on a Hilbert space of dimension $ d $  )   which is a compact group. As $ ISO ( M ) $ is closed, it is a closed subgroup of a compact group, hence compact. We will see that this technique can be generalized in the noncommutative set-up in the chapters \ref{qisol} and \ref{qorient}.
     
  Proposition \ref{preliminiaries_isometry_laplacian_characterization} has the generalization in a more general context.
     
    Let us fix some notations. Let $ Y $ be a compact metrizable space and $ \theta: M \times Y \rightarrow M .$ Let $ \xi_y : M \rightarrow M $ defined by $ \xi_y ( m ) = \theta ( m,y ).$ Let $\alpha : C ( M ) \rightarrow C ( M ) \otimes C ( Y ) \cong C ( M \times Y )  $ be defined by $ \alpha ( f ) ( m, y ) = f ( \theta ( m,y )  ) $ for all $ y $ in $ Y, ~ m $ in $ M.$ For a state $ \phi $ on $ C ( Y ), $ denote by  $ \alpha_{\phi},$ the map: $ ( {\rm id} \otimes \phi ) \alpha : C ( M ) \rightarrow C ( M ).$ Lastly, let $ \cla^{\infty}_0$ be the span of eigenvectors of the Laplacian $ \cll $ of $ M .$
    
    Then, we have the following( Lemma 2.5 of \cite{goswami} ):
    
    \bppsn
    
    \label{preliminiaries_family_of_isometry_characterization}
    
  The following are equivalent:
  
  a. For every $ y $ in $ Y,$ $ \xi_y $ is smooth isometric.
  
  b. For every state $ \phi $ on $ C ( Y ),$ we have $ \alpha_{\phi} ( \cla^{\infty}_0 ) \subseteq \cla^{\infty}_0 $ and $ \alpha_{\phi} \cll = \cll \alpha_{\phi} $ on $ \cla^{\infty}_0.$

    \eppsn
    
\bxmpl 

1. The isometry group of the $n$-sphere  $ S^n  $ is $ O ( n + 1 ) $ where the action is given by the usual action of $ O ( n + 1 ) $ on $ \IR^{n + 1}.$ The subgroup of $ O ( n + 1 ) $ consisting of all orientation preserving isometries on $ S^{n} $ is $ S O ( n + 1 ).$

\vspace{1mm}
 
2.  The isometry group of the circle $S^1$ is $ S^{1}  >\!\!\! \lhd Z_{2}  .$ Here the $ \IZ_2 ~ ( ~ = \{ 0, 1\} ~ ) $ action on $ S^{1} $ is given by $ 1. z = \overline{z} $ where $ z  $ is in $ S^{1}$ while the action of $ S^{1} $ is its action on itself.

 \vspace{1mm}

3. $ ISO ( \IT^n ) \cong  \IT^n >\!\!\! \lhd ( \IZ^{n}_2 >\!\!\! \lhd S_n )  $  where  $ S_n $ is the permutation group on $ n $ symbols. Here an element of $ S_n $ acts on an element $ ( z_1, z_2, ..., z_n ) \in \IT^n $ by permutation. If the generator of $ i $-th copy of $ \IZ^n_2 $ is denoted by $ 1_i, $ then the action of $ 1_i $ is given by $ 1_i ( z_1, z_2,..., z_n ) = ( z_1,..., z_{i - 1}, \overline{z_i}, z_{i + 1},..., z_n ) $ where $ ( z_1, z_2,..., z_n ) \in \IT^n .$ Lastly, the action of $ \IT^n $ on itself is its usual action.
 
 \exmpl

\subsection{Spin Groups and Spin manifolds}

\label{preliminaries_subsection_Spin_groups_Spin_manifolds}

We begin with Clifford algebras. Let $ Q $ be a quadratic form on an $ n $ dimensional vector space $ V. $ Then $ Cl ( V, Q ) $ will denote the universal associative algebra $ \clc $ equipped with a linear map $ i : V \rightarrow \clc, $ such that $  i ( V ) $ generates $ \clc $ as a unital algebra satisfying $ {i ( V )}^2  = Q ( V ).1$ 

Let $ \beta : V \rightarrow Cl ( V, Q ) $ be defined by $ \beta ( x ) = - i ( x ).$ Then, $ Cl ( V, Q ) = Cl^0 ( V, Q ) ~ \oplus ~ Cl^{1} ( V, Q ) $ where $ Cl^0 ( V, Q ) = \{ x \in Cl ( V, Q ) : \beta ( x ) = x \}, ~ Cl^{1} ( V, Q ) = \{ x \in Cl ( V, Q ) : \beta ( x ) = - x \}.$

We will denote by $ \clc_n $ and $ \clc^{\IC}_n $ the Clifford algebras $ Cl ( \IR^n, - x^2_1 - ... - x^2_n ) $ and $ Cl ( \IC^n, z^2_1 + ... + z^2_n ) $ respectively.

We will denote the vector space $ {{\IC}^2}^{[ \frac{n}{2}]} $ by the symbol $ \Delta_n. $ It follows that $ \clc^{\IC}_n = {\rm End} ( \Delta_n ) $ if $ n   $ is even and equals $ {\rm End} ( \Delta_n ) \oplus {\rm End} ( \Delta_n ) $ is $ n $ is odd. There is a representation $  \clc^{\IC}_n \rightarrow {\rm End} ( \Delta_n ) $ which is the isomorphism with $ {\rm End} ( \Delta_n ) $ when $ n $ is even and in the odd case, it is the isomorphism with $ {\rm End} ( \Delta_n ) \oplus {\rm End} ( \Delta_n ) $ followed by the projection onto the first component. This representation restricts to $ \clc_n,$ to be denoted by $ \kappa_n$ and called the spin representation. This representation is  irreducible when $ n $ is odd and for $ n $ even, it decomposes into two irreducible representations which decomposes $ \Delta_n $ into a direct sum of two vector spaces $ \Delta^+_n $ and $ \Delta^{-}_n.$ 

 $ \rm{Pin} ( n ) $ is defined to be the subgroup of $ \clc_n $ generated by elements of the form $\{ x: \left\| x \right\| = 1, x \in \IR^n \}.$ Spin ( n ) is the group given by $\rm{ Pin} ( n )  \cap ~ \clc^0_n.$  There exists a continuous group homomorphism from $ {\rm Pin ( n )} $ to $ O ( n ) $ which restricts to a two covering map $ \lambda : {\rm Spin ( n )} \rightarrow SO ( n ).$
 
 \vspace{2mm}
 
 Let $ M $ be an $n$-dimensional orientable Riemannian manifold.  Then we have the oriented orthonormal bundle of frames over $ M $ ( which is a principal $ S O ( n ) $ bundle ) which we will denote by $ F.$
 
  Such a manifold $ M $ is said to be a {\bf spin manifold} if there exists a pair $ ( P, \Lambda )  $ ( called a spin structure ) where
   
   ( 1 ) P is a Spin( n ) principal bundle over $ M. $
   
   ( 2 ) $ \Lambda $ is a map from $ P $ to $ F $ such that it is a 2-covering as well as a bundle map over $ M. $
   
   ( 3 ) $ \Lambda ( p. \widehat{g} ) = \Lambda ( p ).g  $ where $ \lambda ( \widehat{g} ) = g.$
   
   Given such a spin structure, we consider the associated bundle $ S = P \times_{{\rm Spin}( n )} \Delta_n$ called the {\bf ` bundle of spinors '}.

\subsection{Dirac operators}

\label{preliminaries_subsection_Dirac_operators}

We follow the notations of the previous subsection. On the space of smooth sections of the bundle of spinors $ S, $ one can define an inner product by 
$$ \left\langle s_1, s_2 \right\rangle_{S} = \int_{M} \left\langle s_{1} ( x ), s_{2} ( x ) \right\rangle dvol ( x ) $$
 
 The Hilbert space obtained by completing the space of smooth sections with respect to this inner product is denoted by $ L^{2} ( S ) $ and its members are called square integrable spinors. The Levi Civita connection on $ M $ induces a canonical connection on $ S $ which we will denote by $ \nabla^{S}. $
 
 \bdfn
 
 The {\bf Dirac operator} on $ M $ is the self-adjoint extension of the following operator $ D $ defined on the space of smooth sections of $ S : $
 $$ ( D s ) ( m ) = \sum^{n}_{i = 1} \kappa_{n} ( X_{i} ( m ) ) (  \nabla^{S}_{X_i} s ) ( m ) , $$
 where $ ( X_1, ...X_n )  $ are local orthonormal ( with respect to the Riemannian metric ) vector fields defined in a neighborhood of $ m.$ In this definition, we have viewed $ X_i ( m ) $ belonging to $ T_m ( M ) $ as an element of the Clifford algebra $ Cl_{\IC} ( T_m M ), $ hence $ \kappa_n ( X_i ( m ) ) $ is a map on the fibre of $ S $ at $ m, $ which is isomorphic with $ \Delta_n. $ The self-adjoint extension of $ D $ is again denoted by the same symbol. 
 
 \edfn
 
  We recall three important facts about the Dirac operator:

\bppsn

\label{preliminaries_classical_Dirac_properties}

( 1 ) $ C^{\infty} ( M ) $ acts on $ S $ by multiplication and this action extends to a representation, say $ \pi,$ of the $ C^* $ algebra $C ( M )$ on the Hilbert space $L^{2} ( S ).$ 

( 2 ) For $ f  $ in $ C^{\infty} ( M ), ~ [ D, \pi ( f ) ] $ has a bounded extension. 

( 3 ) Furthermore, the Dirac operator on a compact manifold has compact resolvents. 

 \eppsn

As the action of an element $ f $ in $ C^{\infty} ( M ) $ on $ L^2 ( S ) $ is  by multiplication operator, we will use the symbol $ M_f $ in place of $ \pi ( f ).$

The Dirac operator carries a lot of geometric and topological information. We give two examples. 


( a ) The Riemannian metric of the manifold is recovered by 
 \be \label{preliminaries_Dirac_operator_metric_expression}  d ( p, q ) = {\rm sup}_{\phi \in C^{\infty} ( M ), ~ \left\| [ D, M_{\phi}  ] \right\| \leq 1 } \left| \phi ( P ) - \phi ( q ) \right|. \ee

( b ) For a compact manifold, the operator $ e^{- t D^2} $ is trace class for all $ t > 0 .$ Then the volume form of the manifold can be recovered by the formula  $$ \int_{M} f dvol = c ( n ) {\rm lim_{t \rightarrow 0}} \frac{ {\rm Tr} ( M_f e^{- t D^2} )}{ {\rm Tr} ( e^{- t D^2} )}    $$ 
where $ {\rm dim}  M = n, ~  c ( n ) $ is a constant depending on the dimension.

 \section{Noncommutative Geometry}

 \label{preliminaries_section_NCG}
 
 In this section, we recall those basic concepts of noncommutative geometry which we are going to need. We refer to \cite{con}, \cite{landi_book}, \cite{connes_marcolli} for more details.

 \subsection{Spectral triples}

  \label{preliminaries_subsection_spectral_triples}
 
 Motivated by the facts in Proposition \ref{preliminaries_classical_Dirac_properties}, Alain Connes defined his formulation of noncommutative manifold based on the idea of a spectral triple:
 
 \bdfn
 
 A {\bf spectral triple} or {\bf spectral data} is a triple $ ( \cla^{\infty}, \clh, D )  $ where $ \clh $ is a separable Hilbert space, $ \cla^{\infty} $ is a $ \ast $ subalgebra of $ \clb ( \clh ), $ ( not necessarily norm closed ) and $ D $ is a self adjoint ( typically unbounded ) operator such that for all $ a $ in $ \cla^{\infty},$ the operator $ [ D, a ] $ has a bounded extension. Such a spectral triple is also called an odd spectral triple. If in addition, we have $ \gamma  $ in $ \clb ( \clh ) $ satisfying $ \gamma = \gamma^* = \gamma^{- 1}, ~ D \gamma = - \gamma D $ and $ [ a, \gamma ] = 0 $ for all $a$ in $ \cla^{\infty}, $  then we say that the quadruplet $ ( \cla^{\infty}, \clh, D, \gamma )  $ is an even spectral triple. The operator $ D $ is called the Dirac operator corresponding to the spectral triple.   
  
 \edfn
 
 Furthermore, given an abstract $ \ast $-algebra $ \clb ,$ an odd ( even ) spectral triple on $ \clb $ is an odd ( even ) spectral triple ( $ \pi( \clb ), \clh, D   $ ) ( respectively, ( $ \pi ( \clb ), \clh, D, \gamma $  )  ) where $ \pi: \clb \rightarrow \clb ( \clh ) $ is a  $ \ast $-homomorphism. 

Since in the classical case, the Dirac operator has compact resolvent if the manifold is compact, we say that the spectral triple is of {\bf compact type} if $ \cla^{\infty} $ is unital and  $ D $ has compact resolvent. 

\bdfn

We say that two spectral triples ( $ \pi_1 ( \cla ), \clh_1, D_1 $ ) and ( $ \pi_2 ( \cla ), \clh_2, D_2   $ ) are said to be unitarily equivalent if there is a unitary operator $ U : \clh_1 \rightarrow \clh_2 $ such that $D_2 = U D_1 U^*  $ and $ \pi_2 ( . ) = U \pi_1 ( . ) U^* $ where $ \pi_j, j = 1,2 $ are the representations of $ \cla $ in $ \clh_j, $ respectively. 

\edfn

Next, we will give two examples of spectral triples in classical geometry and a non-classical example. We will  give more examples in chapters \ref{qorient} and \ref{quantumsphere}. 
 
 \bxmpl
 
 Let $ M $ be a smooth spin manifold. Then from Proposition \ref{preliminaries_classical_Dirac_properties},  we see that ( $ C^{\infty} ( M ), \clh, D $ ) is a spectral triple over $ C^{\infty} ( M )  $ and is of compact type if $ M $ is compact. 
 
 We recall that when the dimension of the manifold is even, $ \Delta_n = \Delta^+_n \oplus \Delta^{-}_n.$  An $ L^2 $ section $ s $ has a decomposition $ s = s_1 + s_2 $ where $ s_1 ( m ), s_2 ( m ) $ belongs to $ \Delta^+_n ( m ) $ and $ \Delta^{-}_{n} ( m )$ ( for all $ m $ ) respectively where $ \Delta^{\pm}_n ( m ) $ denotes the subspace of the fibre over $ m. $ This decomposition of $ L^2 ( S ) $ induces a grading operator $ \gamma $ on $ L^2 ( S ). $  It can be seen that $ D $ anticommutes with $ \gamma.$  
 
 \exmpl  
 
\bxmpl

 This example comes from the classical Hilbert space of forms discussed in subsection \ref{preliminaries_subsection_spaceofforms_NCG}. One considers the self adjoint extension of the operator $ d + d^* $ on $ \clh = \oplus_k \clh^k ( M ) $ which is again denoted by $ d + d^*.$ $ C^{\infty} ( M ) $ has a representation on each $ \clh^k ( M ) $ which gives a representation, say $ \pi $ on $ \clh.$ Then it can be seen that ( $ C^{\infty} ( M ), \clh, d + d^*  $ ) is a spectral triple and $ d + d^* $ is called the Hodge Dirac operator. When $ M $ is compact, this spectral triple is of compact type.
 
 \exmpl
 
 \brmrk
 
 \label{preliminaries_spectral_triples_spinors_caution}
 
 Let us make it clear that by a `classical spectral triple' we always mean the spectral triple obtained by the Dirac operator on the spinors (so, in particular, manifolds are assumed to be Riemannian spin manifolds), and not just any spectral triple on the commutative algebra $C^\infty(M)$.
 
 \ermrk
 
 \bxmpl
 
 The Noncommutative torus
 
 We recall from subsection \ref{preliminaries_subsection_C*algebras} that the noncommutative 2-torus $ \cla_{\theta} $ is the universal $ C^* $ algebra generated by two unitaries $ U $ and $ V $ satisfying $ U V = e^{2 \pi i \theta} V U $ where $ \theta $ is a number in $ [ 0, 1 ].$ 
 
 There are two derivations $ d_1 $ and $ d_2 $ on $ \cla_{\theta} $ obtained by extending linearly the rule: 
      $$ d_1 ( U ) = U, ~ d_1 ( V ) = 0,  $$
      $$ d_2 ( U ) = 0, ~ d_2 ( V ) = V. $$  
      
 Then $ d_1 $ and $ d_2 $ are well defined on the following dense $ \ast $-subalgebra of $ \cla_{\theta} $ : 
 $$ \cla^{\infty}_{\theta} = \{ \sum_{m,n \in \IZ} a_{mn} U^m V^n : {\rm sup}_{m,n} \left| m^k n^l a_{mn} \right| < \infty ~ {\rm ~ for ~ all} ~ k,l ~ {\rm in} ~ \IN \}. $$ There is a unique faithful trace on $ \cla_{\theta} $ defined as follows:
 $$ \tau ( \sum a_{mn} U^m V^n ) = a_{00}.$$
 
 Let $ \clh = L^2 ( \tau ) \oplus L^2 ( \tau ) $ where $ L^2 ( \tau ) $ denotes the GNS Hilbert space of $ \cla_{\theta} $ with respect to the state $ \tau.$ We note that $ \cla^{\infty}_{\theta} $ is embedded as a subalgebra of $ \clb ( \clh ) $ by $ a \mapsto \left(  \begin {array} {cccc}
     a   &  0  \\ 0 & a \end {array} \right ) .$
     
     Now, we define $ D = \left(  \begin {array} {cccc}
     0   & d_1 + i d_2  \\ d_1 - i d_2 & 0 \end {array} \right )  .$
     
     Then, $( \cla^{\infty}_{\theta}, \clh, D )  $ is a spectral triple of compact type.  In particular, for $ \theta = 0, $ this coincides with the classical spectral triple on $  C ( \IT^2 ). $

 \exmpl 
 
 \subsection{The space of forms in noncommutative geometry}

\label{preliminaries_subsection_spaceofforms_NCG}
 
 We start this subsection by recalling the universal space of one forms corresponding to an algebra.
 
 \bppsn
 
 \label{preliminaries_Hilbert_space_forms_universal_one_form}
 
 Given an algebra $ \clb ,$ there is a ( unique upto isomorphism ) $ \clb-\clb $ bimodule $ \Omega^{1} ( \clb ) $ and a derivation $ \delta: \clb \rightarrow \Omega^{1} ( \clb ) $ ( that is, $ \delta ( a b ) = \delta ( a ) b + a \delta ( b )  $ for all $ a,b $ in $ \clb $  ), satisfying the following properties:
 
 (i) $ \Omega^{1} ( \clb ) $ is spanned as a vector space by elements of the form $ a \delta ( b ) $ with $ a,b $ belonging to $ \clb; $ and
 
 (ii) for any $ \clb-\clb $ bimodule $ E $ and a derivation $ d: \clb \rightarrow E, $ there is an unique $ \clb-\clb $ linear map $ \eta: \Omega^{1} ( \clb ) \rightarrow E $ such that $ d = \eta \circ \delta. $

 \eppsn
 
 \vspace{2mm}
 
 The bimodule $ \Omega^{1} ( \clb ) $ is called the space of universal 1-forms an $ \clb $ and $ \delta $ is called the universal derivation.
 
 We can also introduce universal space of higher forms on  $ \clb,~ \Omega^{k} ( \clb ), $ say, for $ k = 2,3,..., $ by defining them recursively as follows: $ \Omega^{k + 1} ( \clb ) = \Omega^{k} ( \clb ) \otimes_{\clb} \Omega^{1} ( \clb ) $ and also set $ \Omega^{0} ( \clb ) = \clb. $ 
 
 \vspace{2mm}
 
 Now we briefly discuss  the notion of the noncommutative Hilbert space of forms which will need noncommutative volume form for a spectral triple of compact type. We refer to \cite{fro} ( page 124 -127 ) and the references therein for more details.


\bdfn

A spectral triple $ ( \cla^{\infty}, \clh, D ) $ of compact type is said to be $ \Theta $-summable if $ e^{- t D^2} $ is of trace class for all $ t > 0.$ A $ \Theta $-summable spectral triple is called finitely summable  when there is some $p>0$ such that $ t^{\frac{p}{2}} ~ {\rm Tr} ( e^{- t D^2} ) $ is bounded on $ ( 0, \delta ] $ for some $ \delta > 0.$ The infimum of all such $ p, $ say $p^{\prime}$ is called the dimension of the spectral triple and the spectral triple is called $p^{\prime}$-summable.
  
  \edfn 

  \brmrk
  
  We remark that the definition of $ \Theta $-summability to be used in this thesis is stronger than the one in \cite{con} ( page 390, definition 1. ) in which a spectral triple is called $ \Theta $-summable if $ {\rm Tr} ( e^{-  D^2} ) < \infty.$  
  
  \ermrk 







For a $ \Theta $-summable spectral triple, let $ \sigma_{\lambda} ( T ) = \frac{{\rm Tr} ( T e^{- \frac{1}{\lambda} D^2} )}{ {\rm Tr} ( e^{- \frac{1}{\lambda} D^2} ) } $ for $ \lambda > 0.$ We note that $ \lambda \mapsto \sigma_{\lambda} ( T ) $ is bounded.  



 


Let $$ \tau_{\lambda} (T) = \frac{1}{\log \lambda} \int^{\lambda}_{a} \sigma_u(T) \frac{du}{u} {\rm for~ } \lambda \geq a \geq e. $$ 


Now consider the quotient $C^*$ algebra $\clb_{\infty} = C_b ([a, \infty))/ C_0 ([a, \infty)).$ Let  for $T$ in $ \clb ( \clh ) , \tau (T)  $ in $\clb_{\infty} $  be the class of $\lambda \rightarrow \tau_{\lambda} (T) .$ 





For any state $\omega$ on the $C^*$ algebra $B_{\infty}, ~ Tr_{\omega} (T) = \omega (\tau(T)) $ for all $T$ in $\clb ( \clh )$ defines a functional on $ \clb ( \clh ) .$ As we are not going to need the choice of $ \omega $ in this thesis, we will suppress the suffix  $ \omega $ and simply write $ {\rm Lim}_{t \rightarrow 0^+} \frac{{\rm Tr} ( T e^{- t D^2} )}{{\rm Tr} ( e^{ - t D^2 } ) }$ for $ Tr_{\omega} (T) .$ This is a kind of Banach limit because if $ {\rm lim}_{t \rightarrow 0^+} \frac{{\rm Tr} ( T e^{- t D^2} )}{{\rm Tr} ( e^{ - t D^2 } ) }$ exists, then it agrees with the functional $ {\rm Lim}_{t \rightarrow 0^+} .$ Moreover, $ {\rm Tr}_{\omega} ( T ) $ coincides ( upto a constant ) with the Dixmier trace ( see chapter IV, \cite{con} ) of the operator $  T {\left|  D \right|}^{- p}  $ when the spectral triple has  a finite dimension $ p > 0,$ where $ {\left| D \right|}^{- p} $ is to be interpreted as the inverse of the restriction of $|D|^p$ on the closure of its range. In particular, this functional gives back the volume form for the classical spectral triple on a compact Riemannian manifold.

 
 Let $ \Omega^k ( \cla^{\infty} ) $ be the space of universal k-forms on the algebra $ \cla^{\infty} $ which is spanned by  $  a_0 \delta ( a_1 ) \cdots  \delta ( a_k ) $, $a_i $ belonging to $ \cla^\infty$,  where $ \delta $ is as in Proposition \ref{preliminaries_Hilbert_space_forms_universal_one_form}. There is a natural graded algebra structure on $\Omega \equiv \bigoplus_{k \geq 0} \Omega^k(\cla^\infty),$ 
 which also has a natural involution given by $(\delta(a))^*=-\delta(a^*),$ and using the spectral triple, we get a  $\ast$-representation 
$\Pi : \Omega \raro \clb(\clh)$ which sends $a_0 \delta(a_1) \cdots \delta(a_k)$ to $a_0d_D(a_1) \cdots d_D(a_k)$, where $d_D(a)=[D,a]$. 
 Consider the state $\tau$ on $\clb(\clh)$ given by, $ \tau ( X ) = {\rm Lim}_{t \rightarrow 0^{+}} \frac{{\rm Tr}( X  e^{- t D^2})}{{\rm Tr} ( e^{- t D^2} )}, $ where $ {\rm Lim} $ is as above. Using $\tau $,  we define a positive semi definite sesquilinear form on $ \Omega^k ( \cla^{\infty} ) $ by setting $ \left\langle  w, \eta \right\rangle = \tau ( {\Pi ( w )}^* \Pi ( \eta ) ).$ Let $ K^k = \{ w \in \Omega^k ( \cla^{\infty} ) : \left\langle w, ~ w \right\rangle = 0\},$ for $k \geq 0$, and $K^{-1} :=(0)$. Let $\overline{\Omega^k_D}$ be the Hilbert space obtained by completing the quotient $\Omega^k(\cla^\infty)/K^k$ with respect to the inner product mentioned above, and we  define $ \clh^k_{D}:=
 P^\perp_k \overline{\Omega^k_D},$ where $P_k$ denotes the projection onto the closed subspace generated by $\delta(K^{k-1})$. The map $D^\prime:=d+d^* \equiv d_D+d_D^*$ on $\clh_{d + d^*}:=\bigoplus_{k \geq 0} \clh^k_{D}$ has a self-adjoint extension (which is again denoted by $d+d^*$).  
Clearly, $\clh^k_{D}$ has a total set consisting of elements of the form $[a_0\delta(a_1) \cdots \delta(a_k)]$, with $a_i $ in $ \cla^\infty$ and where $[\omega]$ denotes the equivalence class $P_k^\perp (w+K^k)$ for $\omega $ belonging to $ \Omega^k(\cla^\infty)$. There is a $\ast$-representation $\pi_{d + d^*} : \cla \raro \clb(\clh_{d + d^*})$, given by $\pi_{d + d^*}(a)([a_0  \delta(a_1) \cdots \delta(a_k)])=[aa_0 \delta(a_1) \cdots \delta(a_k)]$.
  Then it is easy to see that
  
  \bppsn
  
   $( \cla^{\infty}, \clh_{d + d^*}, d + d^* ) $ is a spectral triple.
   
   \eppsn

 \subsection{Laplacian in Noncommutative geometry}

  \label{preliminaries_subsection_Laplacian_NCG}

   
  
  Now, we discuss the notion of Laplacian in noncommutative geometry as introduced in \cite{goswami}. Let $ ( \cla^{\infty}, \clh, D ) $ be a spectral triple of compact type. To define the Laplacian in the noncommutative case ( as in \cite{goswami} ), we need the following assumptions on the spectral triple.
    
    {\bf Assumptions}
    
  {\bf 1.} $( \cla^{\infty}, \clh, D ) $ is a compact type spectral triple.
  
  {\bf 2.} It is $QC^{\infty},$ that is, ${\cla^\infty}$ and $\{ [D,a], ~a \in {\cla^\infty} \}$ are contained in the domains of all powers of the derivation $[|D|, \cdot].$ 
  
  {\bf 3.} The  unbounded densely defined map $d_D$ from  $\clh^0_D$ to $\clh^1_D$  given by $d_D(a)=[D,a]$ for $a $ in $ \cla^\infty$, is  closable. 
  
  {\bf 4.} $\cll:=-d^*_D d_D$ has $\cla^\infty$ in its domain,  and it is left invariant by $\cll.$
  
  Under assumption {\bf 2.}, $\tau$ defined by $ \tau ( X ) = {\rm Lim}_{t \rightarrow 0} \frac{{\rm Tr} ( X e^{- t D^2} )}{{\rm Tr} ( e^{- t D^2} )} $ is a positive trace on the  $C^*$-subalgebra generated by $\cla^\infty$ and $\{ [D,a]: ~a \in {\cla^\infty} \}.$ 
  
  {\bf 5.} We assume that it is also faithful on this subalgebra. 
  
  \vspace{2mm}

 Then,  $ \cll = - d^*_D d_D $ is defined to be the Laplacian for the spectral triple $ ( \cla^{\infty}, \clh,  D ).$ It coincides with  the Hodge Laplacian $-d^\ast d$ (restricted on space of smooth functions) in the classical case, where $d$ denotes the de-Rham differential.

  The linear span of eigenvectors of $\cll$, which is a subspace of $\cla^\infty$, is denoted by $\cla^\infty_0$, and the $\ast$-subalgebra  of $\cla^\infty$ generated by $\cla^\infty_0$ is denoted by $\cla_0.$   
     

 
 

\cleardoublepage

\chapter{Quantum isometry groups: approach based on Laplacian}

\label{qisol}
 
The idea of quantum isometry group of a noncommutative manifold (given by a  spectral triple), which has been defined by Goswami, is motivated by the definition and study of quantum permutation groups of finite sets and finite graphs by a number of mathematicians (see, e.g.\cite{ban1}, \cite{ban2}, \cite{wang}, \cite{univ1} and references therein).
  
   In this chapter, we first recall the definition of quantum isometry groups as proposed in \cite{goswami} and then compute it for some examples.

     \section{Formulation of the quantum isometry group}
     
     \subsection{Characterization of isometry group for a compact Riemannian manifold}
     
     Let $ M $ be a compact Riemannian manifold. Consider the category with objects being the pairs $ ( G, \alpha ) $ where $ G $ is a compact metrizable group acting on $M$ by the smooth and isometric action $ \alpha .$ If $ ( G_1, \alpha  )$ and $ ( G_2, \beta ) $ are two objects in this category, $ \rm{Mor} ( ( G_1, \alpha  ), ~ ( G_2, \beta )  ) $ consists of group homomorphisms $ \pi $ from $ G_1 $ to $ G_2 $ such that $ \beta \circ \pi = \alpha.$ Then the isometry group of $ M $ is the universal object in this category.
 
More generally, the isometry group of a classical compact Riemannian  manifold, viewed as a compact metrizable space ( forgetting the group structure ), can be seen to be the universal object of a category whose object class consists of subsets ( not generally subgroups ) of the set of smooth isometries of the manifold. Then it can be proved that this universal compact set has a canonical group structure. Thus, motivated by the ideas of Woronowicz and Soltan ( \cite{soltan}, \cite{woro_pseudo} ), Goswami considered in \cite{goswami} a bigger category with objects as the pair $ ( S, f ) $ where $ S $ is a compact metrizable space and $ f : S \times M \rightarrow M $ such that the map from $ M $ to itself defined by $ m \mapsto f ( s, m )  $ is a smooth isometry  for all $ s $ in $ S .$ The morphism set is defined as above ( replacing group homomorphisms by continuous set maps ).  
  
  Therefore, to define the quantum isometry group, it is reasonable to  consider a category of compact quantum groups which act on the manifold (or more generally, on a noncommutative manifold given by spectral triple) in a `nice' way, preserving the Riemannian structure in some suitable sense, which is precisely formulated in \cite{goswami}, where it is also proven that a universal object in the category of such quantum groups does exist if one makes some natural regularity assumptions on the spectral triple.

    \subsection{The definition and existence of the quantum isometry group}

     \label{qisol_subsection_definition_existence_of_quantum_isometry_group}

 Let $ ( \cla^{\infty}, \clh, D ) $ be a $\Theta$-summable spectral triple of compact type. We recall from section \ref{preliminaries_section_NCG} the Hilbert spaces of $ k $-forms $ \clh^k_{D}, k = 0,1,2,... $ and also the Laplacian $ \cll = - d^*_D d_D.$

     To define the quantum isometry group, we need the following assumptions:
     
     {\bf Assumptions}
     
     {\bf 1.} $ d_{D} $ is closable and $ \cla^{\infty} \subseteq \rm{Dom} ( \cll ) $ where  $ \cla^{\infty} $ is viewed as a dense subspace of $ \clh^{0}_{D}.$
     
     {\bf 2.} $ \cll $ has compact resolvents.
     
     {\bf 3.} $ \cll ( \cla^{\infty} ) \subseteq \cla^{\infty}. $
     
     {\bf 4.} Each eigenvector of $ \cll $ ( which has a discrete spectrum, hence a complete set of eigenvectors ) belongs to  $ \cla^{\infty}.$
       
     {\bf 5.} ( connectedness assumption ) The kernel of $ \cll $ is one dimensional, spanned by the identity $ 1 $ of $ \cla^{\infty} ,$ viewed as a unit vector in $ \clh^{0}_{D}.$
     
     {\bf 6.} The complex linear span of the eigenvectors of $ \cll, $ denoted by  $ \cla^{\infty}_{0} $ is norm dense in  $ \cla^{\infty}.$
     
     \bdfn
     
     We say that a spectral triple satisfying the assumptions {\bf 1.} - {\bf 6.} admissible.
     
     \edfn 
     
     The following result is contained in Remark 2.16 of \cite{goswami}.
  
  \bppsn
  
  \label{QISO_Laplacian_formulation_remark_2.16}
  
  If an admissible spectral triple $ ( \cla^{\infty}, \clh, D ) $ satisfies the condition $\bigcap  \rm{Dom} ( \cll^{n} ) = \cla^{\infty},$ and if $ \alpha : \bar{\cla} \rightarrow \bar{\cla} \otimes S  $ is a smooth isometric action on $ \cla^{\infty} $ by a CQG $ S,$ then for all state $ \phi $ on $ S, ~ \alpha_{\phi} ( = ( {\rm id} \otimes \phi ) \alpha  ) $ keeps $ \cla^{\infty} $ invariant. 
  
  \eppsn
  
  
  

     In view of the characterization of smooth isometric action on a classical compact manifold ( Proposition \ref{preliminiaries_isometry_laplacian_characterization} and Proposition \ref{preliminiaries_family_of_isometry_characterization} in Chapter \ref{preliminaries} ), Goswami gave the following definition in \cite{goswami}.
     
     \bdfn
     
     A quantum family of smooth isometries of the noncommutative manifold $ \cla^{\infty} $ ( or more precisely on the corresponding spectral triple ) is a pair $( \cls, \alpha ) $ where $ \cls $ is a separable unital $ C^* $ algebra, $ \alpha: \overline{\cla} \rightarrow \overline{\cla} \otimes \cls $ ( where $ \overline{\cla} $ denotes the $ C^* $ algebra obtained by completing $ \cla^{\infty} $ in the norm of $ \clb ( \clh^{0}_{D} ) ) $ is a unital $ C^* $ homomorphism, satisfying the following:
     
     a. $ \overline{\rm{Sp}} ( \alpha ( \overline{\cla} ) ( 1 \otimes \cls ) = \overline{\cla} \otimes \cls $
     
     b. $ \alpha_{\phi} = ( {\rm id} \otimes \phi ) \alpha $ maps $ \cla^{\infty}_0 $ into itself and commutes with  $ \cll $ on $ \cla^{\infty}_0,$ for every state $ \phi $ on $ \cls .$
     
     In case, the $ C^* $ algebra has a coproduct $ \Delta $ such that $( \cls, \Delta  )$ is a compact quantum group and $ \alpha $ is an action of $( \cls, \Delta )$ on $ \overline{\cla},$ we say that $ ( \cls, \Delta )  $ acts smoothly and isometrically on the noncommutative manifold. 
     
     \edfn
     
   {\bf Notations}
   
  {\bf 1.} We will denote by $ {\bf Q^{\cll}} $ the category with the object class consisting of all quantum families of isometries $( \cls, \alpha )$ of the given noncommutative manifold, and the set of morphisms $ {\rm Mor} ( ( \cls, \alpha ), ( \cls^{\prime}, \alpha^{\prime}) ) $ being the set of unital $ C^* $ homomorphisms $ \phi: \cls \rightarrow \cls^{\prime} $ satisfying $ ( {\rm id} \otimes \phi ) \alpha = \alpha^{\prime}.$   
      
   {\bf 2.} We will denote by $ {\bf Q^{\prime}_{\cll}} $ the category whose objects are triplets  $ ( \cls, \Delta, \alpha )$ where $ ( \cls, \Delta ) $ is a CQG acting smoothly and isometrically on the given noncommutative manifold, with $ \alpha $ being the corresponding action. The morphisms are the homomorphisms of compact quantum groups which are also morphisms of the underlying quantum families.

     \vspace{4mm}
  
 Let $ \{ \lambda_1, \lambda_2,...\} $ be the set of eigenvalues of $ \cll, $ with $ V_i $ being the corresponding ( finite dimensional ) eigenspace. We will denote by $ \clu_i $ the Wang algebra $ A_{u,d_i} ( I ) $ ( as introduced in the  chapter \ref{preliminaries} ) where $ d_i $ is the dimension of the subspace $ V_i.$ We fix a representation $ \beta_i : V_i \rightarrow V_i \otimes \clu_i $ on the Hilbert space $ V_i, $ given by $ \beta_i ( e_{ij} ) = \sum_k e_{ik} \otimes u^{(i)}_{kj} $ for $ j = 1,2, ... d_i, $ where $ \{ e_{ij} \} $ is an orthonormal basis for $ V_i, $ and $ u^{(i)} \equiv u^{(i)}_{kj} $ are the generators of $ \clu_i .$  Thus, both $ u^{(i)} $ and $ \overline{u^{(i)}} $ are unitaries. The representations $ \beta_i $ canonically induce the free product representation $ \beta = \ast_i \beta_i $ of the free product CQG $ \clu = \ast_i \clu_i $ on the Hilbert space  $ \clh^0_D $ such that the restriction of $ \beta $ on $ V_i $ coincides with $ \beta_i $ for all $i.$ 
  
   The following Lemma ( Lemma 2.12 of \cite{goswami} ) will be needed later and hence we record it. 
  
  \blmma
  
  \label{qiso_L_formulation_Lemma2.12}
  
  Consider an admissible spectral triple $( \cla^{\infty}, ~ \clh, ~ D )$ and let  $( \cls, \alpha )$ be a quantum family of smooth isometries of the spectral triple. Moreover, assume that the action is faithful in the sense that there is no proper $ C^* $ subalgebra $\cls_1$ of $ \cls $ such that $ \alpha ( \cla^{\infty} ) \subseteq \cla^{\infty} \otimes \cls_1. $  Then $ \widetilde{\alpha}: \cla^{\infty} \otimes \cls \rightarrow \cla^{\infty} \otimes \cls $ defined by $\widetilde{\alpha} ( a \otimes b ) = \alpha ( a ) ( 1 \otimes b )$ extends to an $ \cls $ linear unitary on the Hilbert $ \cls $ module $ \clh^{0}_{D} \otimes \cls ,$ denoted again by $ \widetilde{\alpha}.$ Moreover, we can find a $ C^* $ isomorphism $ \phi : \clu / \cli \rightarrow \cls $ between $\cls $ and a quotient of $ \clu $ by a $ C^* $ ideal $ \cli $ of $ \clu ,$ such that $ \alpha = ( {\rm id} \otimes \phi ) \circ ( {\rm id} \otimes \Pi_{\cli} ) \circ \beta $ on $ \cla^{\infty} \subseteq \clh^{0}_{D}, $ where $\Pi_{\cli} $ denotes the quotient map from $ \clu $ to  $ \clu / \cli.$
  
  If furthermore, there is a CQG structure on $ \cls $ given by a coproduct $ \Delta $ such that $ \alpha $ is a $ C^* $ action of a CQG on $ \overline{\cla},$ then the map $ \alpha: \cla^{\infty} \rightarrow \cla^{\infty} \otimes \cls $ extends to a unitary representation ( denoted again by $ \alpha $  ) of the CQG $ ( \cls, \Delta ) $ on $ \clh^{0}_{D}.$ In this case, the ideal $ \cli $ is a Woronowicz $C^*$ ideal and the $ C^*$ isomorphism $ \phi: \clu/\cli \rightarrow \cls $ is a morphism of CQG s.

  \elmma
  
 Using this, the following result has been proved in \cite{goswami}, which defines and gives the existence of $ QISO^{\cll}.$
  
  \bthm
  
  \label{qisol_main_existence_theorem}
  
  For any admissible spectral triple $ ( \cla^{\infty}, \clh, D ) ,$ the category  $ {\bf Q^{\cll}} $ has a universal object denoted by $ ( QISO^{\cll}, \alpha_0 ).$ Moreover, $ QISO^{\cll} $ has a coproduct $ \Delta_0 $ such that $ ( QISO^{\cll}, \Delta_0 ) $ is a CQG and $ ( QISO^{\cll}, \Delta_0, \alpha_0 ) $ is a universal object in the category $ {\bf Q^{ \prime}_{ \cll}} .$ The action $ \alpha_0 $ is faithful.    
  
  \ethm
  
We very briefly outline the main ideas of the proof. The universal object $ QISO^{\cll} $ is constructed as a suitable quotient of $ \clu. $ Let $ \clf $ be the collection of all those $ C^* $-ideals $ \cli $ of $ \clu $ such that the composition $ \Gamma_{\cli} = ( {\rm id} \otimes \Pi_{\cli} ) \circ \beta : \cla^{\infty}_0 \rightarrow \cla^{\infty}_0 \otimes_{{\rm alg}} ( \clu / \cli ) $ extends to a $ C^* $-homomorphism from $ \overline{\cla} $ to $ \overline{\cla} \otimes ( \clu / \cli ).$ Then it can be shown that $ \cli_0 ~ ( = \cap_{\cli \in \clf} \cli ~ ) $ is again a member of $ \clf $ and $ ( \clu / \cli_0, ~ \Gamma_{\cli_0} ) $ is the required universal object.  Thus, 
 \brmrk
 
 \label{qiso_laplacian_qiao_tracial_haar_state}
 
 $ QISO^{\cll} $ is a quantum subgroup of the CQG $ \clu = \ast_i A_{u,d_i} ( I ) .$ As $ A_{u,d_i} ( I ) $ satisfies $ \kappa^2 = {\rm id}, $ ( by Remark \ref{preliminaries_wang_algebra_I_k2=I} ) the same is satisfied by  $ QISO^{\cll} $ so that by Remark \ref{preliminaries_haar_state_not_tracial}, $ QISO^{\cll} $ has tracial Haar state.
 
 \ermrk

\brmrk

It is proved in \cite{goswami}  that to ensure the existence of $ QISO^{\cll},$  the assumption ({\bf 5}) can be replaced by the condition that the action  $ \alpha $ is $ \tau $ preserving, that is, $ ( \tau \otimes {\rm id} ) \alpha ( a ) = \tau ( a ).1.$  
In \cite{goswami} it was also shown ( Lemma 2.5, $ b \Rightarrow a $ ) that for an isometric group action on a not necessarily connected classical manifold, the  volume functional is automatically preserved. It can be easily seen that the proof goes verbatim for a quantum group action, and consequently we get the existence of  $ {QISO}^{\cll} $ for a ( not necessarily connected ) compact Riemannian manifold.

\ermrk

\vspace{4mm}

  {\bf Unitary representation  of $ QISO^{\cll} $ on a spectral triple}

We shall also need the following result proved in section 2.4 of \cite{goswami}.

\bppsn

\label{qisol_unitary_rep_of_qisol}

 $ QISO^{\cll} $ has a unitary representation $ U \equiv U_\cll $ on
 $ \clh_{D} $ such that $ U $ commutes with $ d + d^* .$ Let $\delta$ be as in subsection \ref{preliminaries_subsection_spaceofforms_NCG}.  On the Hilbert space of $ k$-forms, that is. $\clh^k_D$, $U$ is defined by:
 $$  U([a_{0}\delta(a_{1})\cdots \delta(a_{k})] \ot q) = [ a^{(1)}_0 \delta(a^{(1)}_1 )\cdots \delta(a^{(1)}_k )] \otimes (a^{(2)}_0 a^{(2)}_1  \cdots a^{(2)}_k )q, $$
 where $q $ belongs to $ QISO^{\cll},$ $ a_i $ belongs to $ \cla^\infty_0,$ and  for  $ x  $ in $ \cla_0, $ ( the $\ast$-subalgebra generated by the eigenvectors of $\cll$ ) we write in Sweedler notation  $ \alpha ( x ) =  x^{(1)} \otimes x^{(2)} \in \cla_0 \otimes (QISO^{\cll})_{0}$ ($\alpha$ denotes the action of $QISO^{\cll}$). 
 
 \eppsn


       


  \section{Computation of $ QISO^{\cll} $}
  
  Here we compute $ QISO^{\cll} $ for three commutative examples, viz: the sphere, the circle and the n tori. In Chapter \ref{deformation}, we will be able to compute it for two noncommutative examples, namely $ \cla_{\theta} $ and  $ S^{n}_{\theta}$ by using Theorem \ref{abcd}.

  \subsection{The commutative spheres}
  
  Let $QISO^{\cll}$ be the quantum isometry group of $S^2$ and let $\alpha$ be the action of $QISO^{\cll}$ on $C(S^2)$.
      Let $\cll$ be the Laplacian on $S^2$ defined as $$\cll=\frac{\partial^2}{\partial \theta^2}+
{\rm cot}(\theta) \frac{\partial}{\partial \theta}+\frac{1}{{\rm sin}^2(\theta)}\frac{\partial^2}{\partial \psi^2},$$ where the cartesian coordinates $x_1$, $x_2$, $x_3$ for $S^2$ are given by $x_1=r\cos{\psi}  \sin{\theta}$, $x_2=r\sin{\psi} \sin{\theta}$, $x_3=r\cos{\theta}$. In the cartesian coordinates, $\cll=\sum_{i=1}^3 \frac{\partial^2}{\partial x_i^2}.$
      
      The eigenspaces of $\cll$ on $S^2$ are of the form $$E_k={\rm Sp}\{(c_1X_1+c_2X_2+c_3X_3)^k~:~c_i\in \IC,i=1,2,3,~ \sum c_i^2=0\},$$ where $k \geq 1$. $E_k$ consists of harmonic homogeneous polynomials of degree $k$ on $\IR^{3}$ restricted to $S^2$( See \cite{Helgason}, page 29-30 ).
      
      We begin with the following lemma, which says that any smooth isometric action by a quantum group must be `linear'. 
      
      \blmma
      
      The action  $\alpha$ satisfies $\alpha(x_i)=\sum_{j=1}^{3} x_j\otimes Q_{ij}$ where $Q_{ij} $ belongs to $ QISO^{\cll}, i=1,2,3.$
      \elmma
    {\it Proof :}  Since $\alpha$ is a smooth isometric action of $QISO^{\cll}$ on $C(S^2)$, $\alpha$ has to preserve the eigenspaces of the Laplacian $\cll$. In particular, it has to preserve $E_1={\rm Sp}\{ c_1x_1+c_2x_2+c_3x_3~:~c_i \in \IC,i=1,2,3,~      \sum_{i=1}^{3}c^2_i=0\}.$
   
   Now note that $x_1+ix_2,~  x_1-ix_2 $ are in $ E_1$, hence $x_1,x_2 $ are in $ E_1.$
          Similarly $x_3 $ belongs to $ E_1$ too. 
     Therefore $E_1={\rm Sp}\{ x_1,x_2,x_3 \}$, which completes the proof of the lemma.
     \qed

\vspace{4mm}
     
     Now, we state and prove the main result of this section, which identifies $QISO^{\cll}$ with the commutative $C^*$ algebra of continuous functions on the isometry group of $S^2$, that is $O(3)$.
     
     \bthm
     
     \label{QISO_Laplacian_computations_sphere}
     
     The quantum isometry group $QISO^{\cll}$ is commutative as a $C^*$ algebra, and hence $ QISO^{\cll} \cong C ( O ( 3 ) ) .$
     
     \ethm
     {\it Proof :}          
     We begin with the expression $$ \alpha(x_i)=\sum_{j=1}^3 x_j \otimes Q_{ij},~i=1,2,3,$$ and also note that $x_1,x_2,x_3$ form a basis of $E_1$ and $\{ x_1^2,x_2^2,x_3^2,x_1x_2,x_1x_3, x_2x_3 \}$ is a basis of $E_2.$ 
          Since $x_i^*=x_i$ for each $i$ and $\alpha$ is a $\ast$-homomorphism, we must have $Q_{ij}^*=Q_{ij} $ for all $ i,j=1,2,3.$ 
           Moreover, the condition $x^2_1+x^2_2+x^2_3=1$ and the fact that $\alpha$ is a homomorphism gives:
           
                     $$ Q^2_{1j}+Q^2_{2j}+Q^2_{3j}=1,~\forall j=1,2,3.$$                 
                  Again, the condition that $x_i$,$x_j$ commute for all $ i,j$ gives \be \label{QISO_Laplacian_sphere_2c} Q_{ij}Q_{kj}=Q_{kj}Q_{ij} ~ \forall i,j,k, \ee
          \be \label{QISO_Laplacian_sphere_3c} Q_{ik}Q_{jl}+Q_{il}Q_{jk}=Q_{jk}Q_{il}+Q_{jl}Q_{ik}.\ee

               Now, it follows from the Lemma \ref{qiso_L_formulation_Lemma2.12} 
                that $\tilde{\alpha}: C(S^2) \otimes QISO^{\cll} \raro C(S^2) \otimes QISO^{\cll}$ defined by $\tilde{\alpha}(X \otimes Y)=\alpha(X)(1 \otimes Y)$ extends to  a unitary of the Hilbert $QISO^{\cll}$-module $L^2 ( S^2 ) \otimes QISO^{\cll}$ (or in other words, $\alpha$ extends to a unitary representation of $QISO^{\cll}$ on $L^2(S^2)$). 
               But  $\alpha$ keeps $V={\rm Sp}\{x_1,x_2,x_3\}$ invariant.
               So $\alpha$ is a unitary representation  of $QISO^{\cll}$ on $V$, that is $Q = (( Q_{ij} )) $  belonging to $ M_3 ( QISO^{\cll} )$ is a unitary,      hence $Q^{-1}=Q^*=Q^T$, since in this case entries of $Q$ are self-adjoint elements.

               Clearly, the matrix $Q$ is a $3$-dimensional unitary representation of $QISO^{\cll}$. From ( 4 ) of Proposition \ref{preliminaries_CQG_Peter_Weyl}, 
            the antipode $\kappa$ on the matrix elements of a finite-dimensional   unitary representation $U^\alpha \equiv ( u_{pq}^\alpha)$ is given by $\kappa (u_{pq}^\alpha ) =( u_{qp}^\alpha )^* .$
               
               So we obtain \be \label{QISO_Laplacian_sphere_5c} \kappa( Q_{ij} )= Q^{-1}_{ij}=Q^T_{ij}=Q_{ji}.\ee         Now from ( \ref{QISO_Laplacian_sphere_2c} )  , we have $Q_{ij}Q_{kj} = Q_{kj}Q_{ij}.$
               Applying $ \kappa$ on this equation and using the fact that $\kappa$ is an antihomomorphism along with ( \ref{QISO_Laplacian_sphere_5c} ) , 
            we have $Q_{jk}Q_{ji} = Q_{ji}Q_{jk}$
               Similarly , applying $\kappa$ on ( \ref{QISO_Laplacian_sphere_3c} ), we get               
               $$ Q_{lj}Q_{ki} + Q_{kj}Q_{li} = Q_{li}Q_{kj} +Q_{ki}Q_{lj}~ \forall i,j,k,l.$$
         Interchanging between $k$ and $i$ and also between $l,j$ gives 
               \be \label{QISO_Laplacian_sphere_6c} Q_{jl}Q_{ik} +Q_{il}Q_{jk} =Q_{jk}Q_{il} +Q_{ik}Q_{jl}~ \forall i,j,k,l.\ee               
               Now, by (\ref{QISO_Laplacian_sphere_3c} ) - ( \ref{QISO_Laplacian_sphere_6c} ), we have               
               $$ [ Q_{ik},Q_{jl} ] =[ Q_{jl},Q_{ik} ],$$               
               hence $$ [ Q_{ik},Q_{jl} ] = 0.$$
               
                            Therefore the entries of  the matrix $Q$ commute among themselves. 
            However, by faithfulness of the action of $QISO^{\cll}$, it is clear that the $C^*$-subalgebra generated by entries of $Q$ (which forms a quantum subgroup of $QISO^{\cll}$ acting on $C(S^2)$ isometrically) must be the same as  $QISO^{\cll}$,  so $QISO^{\cll}$ is commutative.
            
  So $QISO^{\cll}=C( G )$ for some compact group $G$ acting by isometry on $C (S^2 )$ and $ G $ is clearly universal in the category of compact metrizable groups acting on $ S^2 $ isometrically, so that $ G \cong  O( 3 ).$ \qed

     \brmrk
      
        Similarly, it can be shown that $ QISO ( S^n ) $ is commutative for all $ n \geq 2. $
        
        \ermrk

      \subsection {The commutative one-torus}
       
      Let $\clc=C(S^1)$ be the $C^*$-algebra of continuous functions on the one-torus $S^1$. Let us denote by $z$ and $\overline{z}$ the identity function (which is the generator of $C(S^1)$) and its conjugate respectively. The Laplacian coming from the standard Riemannian metric is given by $\cll(z^n)=-n^2 z^n$, for $n $ in $ \IZ,$ hence the eigenspace corresponding to the eigenvalue $-1$ is spanned by $z$ and $\overline{z}.$ Thus, the action of a compact quantum group acting smoothly and isometrically (and faithfully) on $C(S^1)$ must be {\it linear} in the sense that its action must map $z$ into an element of the form $z \ot A +\overline{z} \ot B$. However, we show below that this forces the quantum group to be commutative as a $C^*$ algebra, that is it must be the function algebra of some compact group .
       
       \bthm
       
       \label{QISO_Laplacian_computations_circle}
       
       Let $\alpha $ be a faithful, smooth and  linear action of a compact quantum group $(\clq,\Delta)$ on $C ( S^1 )$ defined by $ \alpha ( z ) = z \otimes A + \overline{z} \otimes B$. 
      Then $\clq$ is a commutative $C^*$ algebra.
      
      \ethm 
      {\it Proof :}
    By the assumption of faithfulness, it is clear that $\clq$ is generated (as a unital $C^*$ algebra) by $A$ and $B$. Moreover, recall that smoothness in particular means that $A$ and $B$ must belong to the algebra $\clq_0$ spanned by matrix elements of irreducible representations of $\clq$ . Since $z \overline{z} =\overline{z}z =1$ and $\alpha$ is a $\ast$-homomorphism, we have 
        $ \alpha ( z ) \alpha (  \overline{z} ) = \alpha ( \overline{z} ) \alpha ( z ) = 1 \otimes 1 $.

    Comparing coefficients of $z^2,{ \overline {z}}^2$ and $1$ in both hand sides of the relation  $ \alpha ( z ) \alpha (  \overline{z} ) =1 \otimes 1$, we get
    
    \be \label{QISO_Laplacian_circle_1d} 
    AB^* = BA^* = 0,~~~
         AA^* + BB^* =1.\ee
    
    Similarly, $ \alpha (\overline{z} ) \alpha ( z ) = 1 \otimes 1$ gives 
    
    \be \label{QISO_Laplacian_circle_2d} B^*A = A^*B =0,~~~
         A^*A +B^*B =1. \ee
    
    Let $U =A+B$ ,
           $P=A^*A$ ,         
        $Q=AA^*$.      
        Then it follows from  (\ref{QISO_Laplacian_circle_1d}) and (\ref{QISO_Laplacian_circle_2d}) that $U$ is a unitary  and $P$ is a projection since  $P$ is self adjoint and  \bean \lefteqn{ P^2} ~~
        & =& A^*AA^*A ~~       
           = A^*A( 1-B^*B )~~
           = A^*A - A^*AB^*B ~~
            = A^*A  ~~                           
           =P.\eean 
       Moreover ,  \bean \lefteqn{ UP}\\
          & =& ( A + B ) A^*A~~
                             =AA^*A + BA^*A ~~                   
                   =AA^*A \\&& (~{\rm since}~ BA^* =0~{\rm from}~(\ref{QISO_Laplacian_circle_1d}) )\\ 
                   &=& A ( 1-B^*B )
                    ~~=A-AB^*B ~~=A .\eean             
       Thus, $A=UP$ , $B=U-UP=U(1-P)\equiv UP^\perp$ , so $\clq=C^*(A,B)=C^*(U,P)$.

  We can rewrite the action $\alpha$ as follows:\\    
    $$ \alpha ( z ) =z \otimes UP + \overline {z} \otimes UP^\bot.$$    
   The coproduct $\Delta$ can easily be calculated from the requirement 
     $ ({\rm id} \otimes \Delta )\alpha =( \alpha \otimes {\rm id} ) \alpha$ , and it is given by :
    
    
        
    
   \be \label{QISO_Laplacian_circle_6d} \Delta ( UP ) = UP \otimes UP +P^ \bot U^{-1 } \otimes UP^ \bot,  \ee   
    \be \label{QISO_Laplacian_circle_7d} \Delta ( UP^ \bot ) = UP^\bot \otimes UP + PU^{-1} \otimes UP^ \bot. \ee   
  From this, we get
   \be \label{QISO_Laplacian_circle_8d}  \Delta ( U ) =U \otimes UP +U^{-1} \otimes UP^ \bot, \ee   
 \be \label{QISO_Laplacian_circle_9d} \Delta ( P ) = \Delta (U^{ -1 }) \Delta ( UP )
                          =P \otimes P +UP^\bot U^{-1} \otimes P^\bot.\ee 
                       
   It can be checked that $ \Delta $ given by  the above expression is coassociative.
   
   Let $h$ denote the right-invariant Haar state on $\clq$. By the general theory of compact quantum groups, $h$ must be faithful on $\clq_0$.  
      We have (by right-invariance of $h$): $$ ({\rm id} \otimes h) ( P \otimes P + UP^\bot U^{-1} \otimes P^\bot ) =h( P )1.$$
   
  That is, we have \be \label{QISO_Laplacian_circle_10d} h( P^\bot )UP^ \bot U^{-1} = h ( P )P^ \bot.\ee
   
   Since $P$ is a positive element in $\clq_0$ and $h$ is faithful on $\clq_0$, $h(P)=0$ if and only if $P=0$. Similarly , $h(P^\bot)=0$, that is $h(P)=1$, if and only if $P=1$. However, if $P$ is either $0$ or $1$, clearly  $\clq =C^*( U,P )=C^*( U )$, which is commutative. On the other hand, if we assume that $P$ is not a trivial projection, then $h(P)$ is strictly between $0$ and $1$, and we have 
 from  ( \ref{QISO_Laplacian_circle_10d} ) $$  UP^\bot U^{-1} = \frac{h( P )}{ 1-h( P )} P^\bot .$$
   
   Since both $UP^ \bot U^{-1} $ and $P^ \bot $ are nontrivial projections, they can be scalar multiples of each other if and only if they are equal, so we conclude that 
   $ UP^\bot U^{-1}=P^\bot$, that is $U$ commutes with $P^\bot$, hence with $P$, and $\clq$ is commutative.
  \qed

   \subsection{The commutative  n-tori} 
   
   Consider  $ C ( \IT^n ) $ as the universal commutative $ C^{*} $ algebra generated by n commuting unitaries $ U_{1}, U_{2},.....U_{n}  $. It is clear that  the set $ \{ U_{i}^{m}U_{j}^{n} : m,n \in \IZ \}  $ is  an orthonormal basis for  $ L^{2} ( C ( \IT^n ) , \tau_0 ),     $ where $ \tau_0$ denotes the unique faithful normalized trace on $ C ( \IT^n ) $ given by, $\tau_0 ( \sum a_{m n} U_{i}^{m} U_{j}^{n} ) = a_{0 0} $ which is just the integration against the Haar measure.     
      We shall denote by  $ \left\langle A , B \right\rangle = \tau_0 ( A^{*} B )$    the inner product on $\clh_0:=L^2(C ( \IT^n ),\tau_0)$. Let $ {C ( \IT^n )}^{\rm fin}$ be the unital  $\ast$-subalgebra generated by finite complex linear combinations of $U^mV^n$, $m,n \in \IZ$. The Laplacian $\cll$ is given by $\cll(U_{1}^{m_{1}} ......U_{n}^{m_{n}})=-(m_{1}^2 + ...m_{n}^2) U_{1}^{m_{1}}......U_{n}^{m_{n}},$ and it is also easy to see that the algebraic span of eigenvectors of $\cll$ is nothing but the space ${C ( \IT^n )}^{\rm fin}$, and moreover, all the assumptions {\bf 1.} - {\bf 6.} in subsection \ref{qisol_subsection_definition_existence_of_quantum_isometry_group} required for defining the quantum isometry group are satisfied. 
         
         Let $QISO^{\cll}$ be the quantum isometry group coming from the above Laplacian, with the 
  smooth isometric action of $ QISO^{\cll} $  on $ C ( \IT^n ) $ given by $\alpha : C ( \IT^n ) \raro C ( \IT^n ) \ot QISO^{\cll}$. 
     By definition,  $ \alpha $  must keep invariant the eigenspace of $ \cll $ corresponding to the eigenvalue $- 1, $ spanned by $ U_{1},.....U_{n},U_{1}^{-1},.......,U_{n}^{-1} $. Thus, the action $ \alpha $ is given by:   
    $$  \alpha ( U_{i} ) = \sum_{j = 1}^{n} U_{j} \otimes A_{ij} + \sum_{j=1}^{n} U_{j}^{-1} \otimes B_{ij} ,$$  
   where $ A_{ij},B_{ij} $ are in $ QISO^{\cll} ,~ i,j =1,2....n $. By faithfulness of the action of quantum isometry group, the norm-closure of the unital $\ast$-algebra generated by $ \{ A_{ij},B_{ij} : ~ i,j=1,2,....n \} $ must be the whole of $QISO^{\cll}$.
   
   Next we derive a number of conditions on $ A_{ij},B_{ij}, i,j = 1,2,...n $ using the fact that $ \alpha $ is a $ \ast $   homomorphism.
  \blmma    
    \label{Lemma 1}

   The condition $ U^{*} U = 1 = U U^{*} $ gives:   
  \be \label{lem1.1a}  \sum_{j} ( A^{*}_{ij} A_{ij} +  B^{*}_{ij} B _{ij} ) = 1,    \ee 
  \be \label{lem1.2} B^{*}_{ij} A_{ik} +  B^{*}_{ik} A_{ij} = 0 ~ \forall j \neq k, \ee  
  \be \label{lem1.3} A^{*}_{ij} B_{ik} +  A^{*}_{ik} B_{ij} = 0 ~ \forall j \neq k, \ee
   \be \label{lem1.4} A^{*}_{ij} B_{ij} = B^{*}_{ij} A_{ij} = 0, \ee 
  \be \label{lem1.5} \sum_{j} ( A_{ij} A^{*}_{ij} + B_{ij} B^{*}_{ij} ) = 1, \ee
  \be \label{lem1.6}A_{ik} B^{*}_{ij} + A_{ij}B^{*}_{ik} = 0 ~ \forall j \neq k,  \ee
 \be \label{lem1.7}B_{ik} A^{*}_{ij} + B_{ij}A^{*}_{ik} =  0 ~ \forall j \neq k, \ee
  \be \label{lem1.8} A_{ij} B^{*}_{ij} = B_{ij}A^{*}_{ij} =  0. \ee

   \elmma
 {\it Proof :}
  We get ( \ref{lem1.1a} ) - ( \ref{lem1.4} ) by using the condition $ U_{i}^{*} U_{i} = 1 $ along with the fact that $ \alpha $ is a homomorphism and then comparing the coefficients of $ 1, U_{j}U_{k}, U_{j}^{-1}U_{k}^{-1},$ ( for $ j \neq k ),$ $ U_{j}^{-2}, U_{k}^{-2} . $
   
   Similarly the condition $ U_{i} U_{i}^{*} = 1 $ gives ( \ref{lem1.5} ) - ( \ref{lem1.8} ).\qed \vspace{4mm}

  Now,  for all $ i \neq j, ~ U_{i}^{*}U_{j}, U_{i}U_{j}^{*} $ and $ U_{i}U_{j} $ belong to the eigenspace of the Laplacian with eigenvalue $-2$, while $ U_{k}^{2}, U_{k}^{-2} $ belong to the eigenspace corresponding to the eigenvalue $-4$. As $ \alpha $ preserves the eigenspaces of the Laplacian, the coefficients of $ U_{k}^{2}, U_{k}^{-2} $ are zero  for all $ k $ in $ \alpha( U_{i}^{*}U_{j} ), ~ \alpha( U_{i}U_{j}^{*} ), ~ \alpha( U_{i}U_{j} ) $ for all $ i \neq j. $

   We use this observation in the next lemma.
   
   \blmma
   \label{Lemma 3}
   
    For all $ k $ and  for all $ i \neq j ,$
   
   \be \label{lem3.1}  B^{*}_{ik} A_{jk} = A^{*}_{ik} B_{jk} = 0,  \ee
   \be \label{lem3.2}  A_{ik}B_{jk} = B_{ik}A^{*}_{jk} = 0,  \ee
   \be \label{lem3.3}  A_{ik} A_{jk} = B_{ik} B_{jk} = 0.  \ee

   \elmma
   
    {\it Proof :}   
   The equation ( \ref{lem3.1} ) is  obtained from the coefficients of $  U_{k}^{2} $ and $ U_{k}^{-2} $  in $ \alpha ( U_{i}^{*} U_{j} ) $ while ( \ref{lem3.2} ) and ( \ref{lem3.3} ) are obtained from the same coefficients in $ \alpha ( U_{i} U_{j}^{*} ) $ and $ \alpha ( U_{i} U_{j} )  $ respectively. \qed
   
  \vspace{4mm}

      Now, by  Lemma \ref{qiso_L_formulation_Lemma2.12}  it follows that
            $ \tilde{\alpha}: C( \IT^{n} ) \otimes QISO^{\cll} \raro C( \IT^{n} ) \otimes QISO^{\cll} $ defined by                $\tilde{\alpha}(X \otimes Y)=\alpha(X)(1 \otimes Y)$ extends to  a unitary of the Hilbert $ QISO^{\cll}$-module $L^2 ( C( \IT^{n} ) ,\tau ) \otimes  QISO^{\cll}$ ( or in other words, $\alpha$ extends to a unitary representation of $ QISO^{\cll}$ on  $L^2( C( \IT^{n} ) ,\tau)$). 
               But $\alpha$ keeps $W = {\rm Sp}\{U_{i},U_{i}^{*} : 1\leq i \leq n \}$ invariant.
               So $\alpha$ is a unitary representation of $ QISO^{\cll}$ on $W$. Hence, the matrix ( say M ) corresponding to the $ 2n $ dimensional representation of $  QISO^{\cll} $ on $W$ is a unitary in $ M_{2n} (  QISO^{\cll} )$.
               
      From the definition of the action it follows that $ M =\left(  \begin {array} {cccc}
       A_{ij} &  B_{ij}^{*}  \\ B_{ij} & A_{ij}^{*} \end {array} \right ).   $
      
      Since $ M $ is the matrix corresponding to a finite dimensional unitary representation, using ( 4 ) of Proposition \ref{preliminaries_CQG_Peter_Weyl}, we have  $  ( k ( M_{k l} ) ) = \left ( \begin {array} {cccc}
   A_{ji}^{*} & B_{ji}^{*}  \\ B_{ji} & A_{ji} \end {array} \right ) $ ( $ \kappa $ is the antipode of $ QISO^{\cll} $ ).
   
   \vspace{4mm}
      
  Now we apply the antipode $ \kappa $ to get some more relations.
  
  \blmma
    \label{Lemma 2}:
    
     For all $ k $ and $ i \neq j ,$
        \be A_{kj}^{*}A_{ki}^{*} = B_{kj} B_{ki} = A_{kj}^{*} B_{ki}^{*} = B_{kj} A_{ki} = B_{kj} A_{ki}^{*} = A_{kj} B_{ki} = 0 . \ee
        \elmma
        
   {\it Proof :} The result follows by applying $ \kappa $ on the equations $ A_{ik} A_{jk} = B_{ik} B_{jk} =  B^{*}_{ik} A_{jk} = A^{*}_{ik} B_{jk} = A_{ik}B_{jk} = B_{ik}A^{*}_{jk} = 0 $ obtained from Lemma  \ref{Lemma 3}.\qed
      
       \blmma
   \label{Lemma 5a}:

      $  A_{li} $ is a normal partial isometry  for all $ l,i $  and hence has same domain and range.
      \elmma
      
   {\it Proof :} From the relation ( \ref{lem1.1a} )  in Lemma \ref{Lemma 1}, we have by applying $ \kappa $, $ \sum (  A^{*}_{ji} A_{ji} +  B_{ji}B _{ji}^{*} ) = 1 .$       
      Applying $ A_{li} $ on the right of this equation, we have 
      
    $  A_{li}^{*} A_{li} A_{li} + \sum_{j \neq l} ( A_{ji}^{*} A_{ji} A_{li} + B_{li} B_{li}^{*}A_{li} ) + \sum_{j \neq l} B_{ji} B_{ji}^{*}A_{li} = A_{li} .$     
   
   From Lemma \ref{Lemma 3}, we have $ A_{ji} A_{li} = 0 $ and $ B_{ji}^{*} A_{li} = 0 $ for all $ j \neq l.$
   Moreover, from Lemma \ref{Lemma 1}, we have $ B_{li}^{*} A_{li} = 0 .$
   Applying these to the above equation, we have    
  \be \label {i}    A^{*}_{li} A_{li} A_{li} = A_{li}.  \ee     
     Again, from the relation $  \sum_{j}( A_{ij} A^{*}_{ij} +  B_{ij}B^{*}_{ij} ) = 1 $  for all $ i $ in Lemma  \ref{Lemma 1}, applying $ \kappa $ and multiplying by $ A^{*}_{li} $ on the right, we have 
    $  A_{li} A_{li}^{*} A_{li}^{*} + \sum_{j \neq l}  A_{ji} A^{*}_{ji} A^{*}_{li}  + B^{*}_{li} B_{li}A^{*}_{li}   + \sum_{j \neq l} B^{*}_{ji} B_{ji}A^{*}_{li} = A^{*}_{li} .$   
  From Lemma \ref{Lemma 3}, we have  $ A_{li} A_{ji} = 0 $  for all $ j \neq l $ ( hence $  A^{*}_{ji} A^{*}_{li} = 0 $ ) and $ B_{ji} A^{*}_{li} = 0 .$ Moreover, we have $ B_{li} A^{*}_{li} = 0 $ from Lemma \ref{Lemma 1}.
  Hence, we have
     \be \label {ii}    A_{li} A^{*}_{li} A^{*}_{li} = A^{*}_{li} . \ee      
     From (\ref {i}), we have  \be \label {iii}  ( A^{*}_{li} A_{li} ) ( A_{li} A^{*}_{li} ) = A_{li} A^{*}_{li}.  \ee      By taking $*$ on (\ref {ii}), we have \be \label {iv}  A_{li} A_{li} A^{*}_{li} = A_{li} . \ee     
    Using this in (\ref {iii}), we have \be A_{li} A^{*}_{li} A_{li} = A_{li} A^{*}_{li}, \ee
      and hence $ A_{li} $ is normal.
          
     So, $ A_{li} = A^{*}_{li} A_{li} A_{li} $ ( from ( \ref {i} ) )     
                $ = A_{li} A^{*}_{li} A_{li}. $
                
      Therefore, $ A_{li} $ is a partial isometry which is normal and hence has same domain and range.\qed

      \blmma
   \label{Lemma 5b}:

   $ B_{li} $ is a normal partial isometry and hence has same domain and range.
   \elmma
   
   {\it Proof :} First we note that $ A_{ji} $ is a normal partial isometry and $ A_{ji} B_{li} = 0 $  for all  $ j \neq l $ ( obtained from Lemma \ref{Lemma 3} ) implies that  $ {\rm Ran} ( A^{*}_{ji} ) \subseteq {\rm Ker} ( B^{*}_{li} ) $ and hence $ {\rm Ran} ( A_{ji} ) \subseteq {\rm Ker} ( B^{*}_{li} ) $ which means  $ B^{*}_{li} A_{ji} = 0 $ for all $ j \neq l .$
   
   To obtain $ B^{*}_{li} B_{li} B_{li} = B_{li} $, we apply $ \kappa $ and multiply by $ B_{li} $ on the right of  ( \ref{lem1.5} ) and then use $ A^{*}_{li} B_{li} = 0 $ from Lemma \ref{Lemma 1}, $ A_{ji} B_{li} = 0 $  for all $ j \neq l $ ( from Lemma \ref{Lemma 3}  which implies  $ B^{*}_{li} A_{ji} = 0 $  for all $ j \neq l $ as above ) and $ B_{ji} B_{li} = 0 $  for all $ j \neq l $ from Lemma \ref{Lemma 3} .
 
 Similarly, we have $ B_{li} B^{*}_{li} B^{*}_{li} = B^{*}_{li} $ by applying $ \kappa $ and multiplying by $ B^{*}_{li} $ on the right of ( \ref{lem1.1a} ) obtained from Lemma \ref{Lemma 1} and then use  $ A_{li} B^{*}_{li} = 0 $  ( Lemma \ref{Lemma 1} ), $ B_{li} A^{*}_{ji} = 0 $  for all $ j \neq l $ and  $ B_{li} B_{ji} = 0 $  for all $ j \neq l $ ( Lemma \ref{Lemma 3} ).

     Using   $ B^{*}_{li} B_{li} B_{li} = B_{li} $ and $ B_{li} B^{*}_{li} B^{*}_{li} = B^{*}_{li} $  as in Lemma \ref{Lemma 5a}, we have $ B_{li} $ is a normal partial isometry. \qed \vspace{4mm}

  Now, we use the condition  $ \alpha( U_{i} ) \alpha ( U_{j} ) = \alpha ( U_{j} ) \alpha ( U_{i} ) $  for all $  i,j .$
  
  \blmma
   \label{Lemma 6}:

       For all $ k \neq l,$
  \be \label{lem6.1} A_{ik}A_{jl} + A_{il}A_{jk} = A_{jl}A_{ik} + A_{jk}A_{il}, \ee   
   \be A_{ik}B_{jl} + B_{il}A_{jk} = B_{jl}A_{ik} + A_{jk}B_{il}, \ee   
   \be B_{ik}A_{jl} + A_{il}B_{jk} = A_{jl}B_{ik} + B_{jk}A_{il}, \ee   
   \be B_{ik}B_{jl} + B_{il}B_{jk} = B_{jl}B_{ik} + B_{jk}B_{il}.  \ee
   
  \elmma
  
 {\it Proof :} The result follows by equating the coefficients of $ U_{k}U_{l}, U_{k}U^{-1}_{l}, U^{-1}_{k}U_{l} $ and $ U^{-1}_{k}U^{-1}_{l} $ ( where $  k \neq l $ ) in $ \alpha( U_{i} ) \alpha ( U_{j} ) = \alpha ( U_{j} ) \alpha ( U_{i} ) $ for all $ i,j .$

    \qed
   
   \blmma
   \label{Lemma 7}:

   $ A_{ik}B_{jl} = B_{jl}A_{ik} $  for all $ i \neq j, k \neq l. $
   \elmma
   
   {\it Proof :} From Lemma \ref{Lemma 6},  we have  for all $  k \neq l, A_{ik}B_{jl} - B_{jl}A_{ik} = A_{jk}B_{il} - B_{il}A_{jk} .$
  We consider the case where $ i \neq j.$ 
              
 We have,   $ {\rm Ran}( A_{ik}B_{jl} - B_{jl}A_{ik} ) \subseteq {\rm Ran}( A_{ik} ) + {\rm Ran}(  B_{jl} ) \subseteq  {\rm Ran}( B^{*}_{jl}B_{jl} + A^{*}_{ik}A_{ik} ) $ ( using the facts that $ A_{ik} $ and $ B_{jl} $ are normal partial isometries by Lemma \ref{Lemma 5a} and \ref{Lemma 5b} and also that  $ B^{*}_{jl}B_{jl} $ and $ A^{*}_{ik}A_{ik} $ are projections ).  
    
    Similarly, $ {\rm Ran} ( A_{jk}B_{il} - B_{il}A_{jk} ) \subseteq  {\rm Ran}( B^{*}_{il}B_{il} + A^{*}_{jk}A_{jk} ).$   
    
    Let \be T_{1} = A_{ik}B_{jl} - B_{jl}A_{ik}, \ee                   
         \be T_{2} = A_{jk}B_{il} - B_{il}A_{jk}, \ee          
        \be  T_{3} = B^{*}_{jl}B_{jl} + A^{*}_{ik}A_{ik}, \ee                   
         \be T_{4} =  B^{*}_{il}B_{il} + A^{*}_{jk}A_{jk}. \ee          
      Hence, $ T_{1} = T_{2}, {\rm Ran} T_{1} \subseteq {\rm Ran} T_{3}~, {\rm Ran} T_{2}\subseteq {\rm Ran} T_{4} .$             
                   
      We claim that $ T_{4}T_{3} = 0.$
      
      Then $ {\rm Ran} ( T_{3} )  \subseteq {\rm Ker} ( T_{4} ).$  
      
      But $ {\rm Ran} T_{1} \subseteq {\rm Ran} T_{3} $ will imply that $ {\rm Ran} T_{1} \subseteq {\rm Ker} T_{4}.$  Hence, $ {\rm Ran} ( T_{2} ) \subseteq {\rm Ker} ( T_{4} ) = { \overline{{\rm Ran} ( T^{*}_{4})}}^{\bot} = { \overline{{\rm Ran} ( T_{4})} }^{\bot} $
   But $ {\rm Ran} ( T_{2} ) \subseteq {\rm Ran} ( T_{4} ) $ which implies that $ {\rm Ran} ( T_{2} ) = 0$ and hence both $ T_{2} $ and $ T_{1} $ are zero.
   Thus, the proof of the lemma will be complete if we can prove the claim.
   \be   T_{4} T_{3} = ( B^{*}_{il}B_{il} + A^{*}_{jk}A_{jk} ) ( B^{*}_{jl}B_{jl} + A^{*}_{ik}A_{ik} ) \ee   
   \be = B^{*}_{il}B_{il}B^{*}_{jl}B_{jl} + B^{*}_{il}B_{il}A^{*}_{ik}A_{ik} + A^{*}_{jk}A_{jk}B^{*}_{jl}B_{jl} + A^{*}_{jk}A_{jk}A^{*}_{ik}A_{ik}. \ee
   
   From Lemma \ref{Lemma 3}, we have  for all $ i \neq j , B_{il}B_{jl} = 0 $ implying $ B_{il} B^{*}_{jl} = 0 $ as $ B_{jl} $ is a normal partial isometry. 
  
  Again, from Lemma \ref{Lemma 2}  for all $ k \neq l, B_{il} A_{ik} = 0.$ Then $ A_{ik} $ is a normal partial isometry implies that $ B_{il} A^{*}_{ik} = 0 $  for all $ k \neq l .$
  
  Similarly, by taking adjoint of the relation $ B_{jl} A^{*}_{jk} = 0 $  for all $ k \neq l $ obtained from Lemma \ref{Lemma 2}, we have $ A_{jk} B^{*}_{jl} = 0.$
  
  From Lemma \ref{Lemma 3}, we have $ A_{jk} A_{ik} = 0 $  for all $ i \neq j . ~ A_{ik} $ is a normal partial isometry implies that $ A_{jk} A^{*}_{ik} = 0 $  for all $ i \neq j .$
  
  Using these, we note that $ T_{4} T_{3} = 0 $ which proves the claim and hence the lemma.\qed

     \blmma
   \label{Lemma 8}:   
     
        \be A_{ik} B_{jk}= 0 =  B_{jk} A_{ik}, \ee         
        \be A_{ki} B_{kj} = 0 = B_{kj} A_{ki} \ee   for all $ i \neq j   $ and  for all $ k. $
   \elmma

    {\it Proof :}  By Lemma \ref{Lemma 3}, we have $ A_{ik} B_{jk} = 0 $ and $ B_{jk} A^{*}_{ik} = 0 $  for all $ i \neq j .$  
  The second relation along with the fact that $ A_{ik} $ is a normal partial isometry implies that $ B_{jk} A_{ik} = 0 $ for all $ i \neq j .$
  
  Thus, $ A_{ik} B_{jk}= 0 =  B_{jk} A_{ik} $  for all $ i \neq j .$
  
  Applying $ \kappa $ on the above equation and using $ B_{kj} $ and $ A_{ki} $ are normal partial isometries, we have $ A_{ki} B_{kj} = 0 = B_{kj} A_{ki} .$

   \qed \vspace{4mm}

        \blmma :
      \label{Lemma 9} 
     $ A_{ik} B_{ik} = B_{ik} A_{ik} $  for all $ i,k. $
      \elmma 
      
     {\it Proof :} We have $ A^{*}_{ij} B_{ij} = 0 = B^{*}_{ij} A_{ij} $ from Lemma \ref{Lemma 1}. Using the fact that $ B_{ij} $ and $ A_{ij} $ are normal partial isometry we have $ A^{*}_{ij} B^{*}_{ij} = 0 = B^{*}_{ij} A^{*}_{ij} $ and hence $ A_{ij} B_{ij} = B_{ij} A_{ij} .$ \qed

    \blmma :
      \label{Lemma 10}

 $ A_{ik} A_{jl} = A_{jl} A_{ik} $  for all $ i \neq j, k \neq l .$
  
   \elmma

  {\it Proof :} Using ( \ref{lem6.1} ) in Lemma \ref{Lemma 6},  we proceed as in Lemma \ref{Lemma 7} to get $ {\rm Ran}( A_{ik}A_{jl} - A_{jl}A_{ik} ) \subseteq   {\rm Ran}( A_{jl}A^{*}_{jl} + A_{ik}A^{*}_{ik} ) $ and  $ {\rm Ran} ( A_{jk}A_{il} - A_{il}A_{jk} ) \subseteq  {\rm Ran}( A_{il}A^{*}_{il} + A_{jk}A^{*}_{jk} ).$

  We claim that $ ( A_{ik}A^{*}_{ik} + A_{jl}A^{*}_{jl} )( A_{jk}A^{*}_{jk} + A_{il}A^{*}_{il} ) = 0.$
   
   Then by the same reasonings as given in Lemma \ref{Lemma 7} we will have : $  A_{jk} A_{il} = A_{il} A_{jk} .$

   To prove the claim,  we use $  A_{ik}A_{jk} = 0 $  for all $ i \neq j $ from Lemma \ref{Lemma 3} ( which implies $ A^{*}_{jk} A_{ik} = 0 $  for all $ i \neq j  $ as $ A_{ik} $ is a normal partial isometry ), $ A^{*}_{il} A^{*}_{ik} = 0 $  for all $ k \neq l $ from Lemma \ref{Lemma 2} ( which implies $ A^{*}_{il} A_{ik} = 0 $  for all $ k \neq l $ as $ A_{ik} $ is a normal partial isometry ) and  $ A_{il} A_{jl} = 0 $  for all $ i \neq j $ from Lemma \ref{Lemma 3} ( which implies  $ A^{*}_{jl} A_{il} = 0 $  for all $ i \neq j $ as $ A^{*}_{il} $ is a normal partial isometry ).

  \qed

    \blmma :
      \label{Lemma 11}

 \be A_{ik} A_{il} = A_{il} A_{ik} {\rm ~ for ~ all} ~  k \neq l, \ee
 \be A_{ik} A_{jk} = A_{jk} A_{ik} {\rm for ~ all} ~ i \neq j. \ee
   
   \elmma

  {\it Proof :} From Lemma \ref{Lemma 3}, we have $ A_{ki} A_{li} = 0 $  for all $ k \neq l .$
        
        Applying $ \kappa $ and taking adjoint, we have $ A_{ik} A_{il} = 0 $  for all $ k \neq l .$
        Interchanging $ k $ and $ l,$ we get $ A_{il}A_{ik} = 0 $  for all $ k \neq l .$
        Hence, $ A_{ik} A_{il} = A_{il} A_{ik} $ for all $  k \neq l .$ 
       
       From Lemma \ref{Lemma 3}, we have $ A_{ik} A_{jk} = 0 $  for all $ i \neq j .$ Interchanging $ i $ and  $ j,$  we have $ A_{jk} A_{ik} = 0 $  for all $ i \neq j .$ \qed
       
     \brmrk
      Proceeding in an exactly similar way, we have that $ B_{ij} $'s commute among themselves. 
                                                 \ermrk 
                                                 
     \bthm
     
     \label{QISO_Laplacian_computations_Torus_finaltheorem}
     
      The Quantum isometry group of $ C ( \IT^{n} ) $ is commutative as a $ C^{\ast} $ algebra and hence coincides with the classical isometry group $ C( {\IT}^{n} >\!\!\! \lhd ( \IZ^{n}_2 >\!\!\! \lhd S_n ) ).$
      
      \ethm
                                                  
    {\it Proof :}                                             
    Follows from the results in lemma \ref{Lemma 7} - \ref{Lemma 11} and the remark following them.  \qed

\cleardoublepage

 \chapter{Quantum group of orientation preserving Riemannian isometries}
 
  \label{qorient}
 
\section{Introduction}
The formulation of quantum isometry groups in \cite{goswami} had a major drawback from the viewpoint of noncommutative geometry since it needed a `good' Laplacian to exist.  In noncommutative geometry it is not always easy to verify such an assumption about the Laplacian, and thus it would be more appropriate to have a formulation in terms of the Dirac operator directly. This is what we aim to achieve in the present chapter.

   The group of Riemannian isometries of a compact Riemannian manifold $M$ can be viewed as the universal object in the
      category of all compact metrizable groups acting on $M$, with smooth and isometric action. Moreover, let us assume that the manifold has a spin structure (hence in particular orientable, so we can fix a choice of orientation) and $D$ denotes the conventional Dirac operator acting as an unbounded self-adjoint operator on the Hilbert space $\clh$ of square integrable spinors. Then, it can be proved that a group action on the manifold lifts as a unitary representation on the Hilbert space $\clh$ which commutes with $D$ if and only if   the action on the manifold is an orientation preserving isometric action. Therefore, to define the quantum analogue of the group of orientation-preserving Riemannian isometries of a possibly noncommutative manifold given by a spectral triple $(\cla^\infty, \clh, D)$, it is reasonable to  consider a category  ${\bf Q}^\prime$ of compact quantum groups having unitary (co-) representation, say $U$, on $\clh$,  which commutes with $D$, and  the action on $\clb(\clh)$ obtained by conjugation maps $\cla^\infty$ into its weak closure. 
A universal object in this category, if it exists, should define the `quantum group of orientation preserving Riemannian isometries' of the underlying spectral triple. It is easy to see that any object $(\cls, U)$  of the category ${\bf Q}^\prime$ gives an equivariant spectral triple $(\cla^\infty, \clh, D)$  with respect to the action of $\cls$ implemented by $U$. It  may be noted that recently there has been a lot of interest and work (see, for example, \cite{partha}, \cite{con2}, \cite{landiqgp}) towards construction of quantum group equivariant spectral triples. In all these works, given a $C^*$-subalgebra $\cla $ of $ \clb(\clh)$ and a CQG $\clq$ having a unitary representation $U$ on $\clh$ such that $ \alpha_U ( \equiv {\rm ad}_U ) $ gives an action of $\clq$ on $\cla$,  the authors investigate the possibility of constructing a (nontrivial) spectral triple $(\cla^\infty, \clh, D)$ on a suitable dense subalgebra $\cla^\infty$ of   $\cla$ such that $ U $ commutes with $ D \ot I ,$ that is, $D$ is equivariant. Our interest here is  in the (sort of) 
 converse direction: given a spectral triple, we want to consider all possible CQG representations  with respect to which the spectral triple is equivariant; and if there exists a universal object in the corresponding category, that is, ${\bf Q}^\prime$, we should call it the quantum group of orientation preserving isometries. 
 
  Unfortunately, even in the finite-dimensional (but with noncommutative $\cla$) situation  this category may often fail to have a universal object, as will be discussed later. It turns out, however, that if we fix a suitable densely defined ( in the WOT ) functional on $\clb(\clh)$ (to be interpreted as the choice of a `volume form') then there exists a universal object in the subcategory of ${\bf Q}^\prime$ obtained by restricting the object-class to the quantum group actions which also preserve the given functional.  The  subtle point to note here is that unlike the classical group actions on $\clb(\clh)$ which always preserve the usual trace, a quantum group action may not do so. In fact, it was proved by Goswami in \cite{goswami_rmp} that given an object $(\clq, U)$ of ${\bf Q}^\prime$ (where $\clq$ is a compact quantum group and $U$ denotes its unitary co-representation on $\clh$), we can find a  positive invertible operator $R$ in $\clh$ so that the given spectral triple is $R$-twisted in the sense of \cite{goswami_rmp} and the corresponding functional $\tau_R$ (which typically differs from the usual trace of $\clb(\clh)$ and can have a nontrivial modularity) is preserved by the action of $\clq$. This makes it quite natural to work in the setting  of twisted spectral data (as defined in \cite{goswami_rmp}).

       Motivated by the ideas of Woronowicz and Soltan ( \cite{woro_pseudo}, \cite{soltan} ), we actually consider a bigger category called the category of `quantum families of smooth orientation preserving Riemannian isometries'.  
       The underlying $C^*$-algebra of the quantum  group of orientation preserving isometries (whenever exists)  has been identified with the universal object in this bigger category and moreover, it is shown to be equipped with a canonical coproduct making it into a compact quantum group.  

In this chapter, we discuss a number of examples, covering the classical spectral triple on $ C^{\infty} ( \IT^2 ) $ as well as the equivariant spectral triples constructed recently on $SU_\mu(2).$ It may be relevant to point out here that it was not clear whether one could accommodate the spectral triples on $SU_\mu(2)$ and the Podles' spheres $S^2_{\mu ,  c}$ in the framework of \cite{goswami}, since it is very difficult to give a nice description of the space of `noncommutative' forms and the Laplacian for these examples. However, the present formulation in terms of the Dirac operator makes it easy to accommodate them, and we have been able to identify $U_\mu(2)$ and $SO_\mu(3)$ as the  universal quantum group of orientation preserving isometries for the spectral triples on $SU_\mu(2)$ and $S^2_{\mu , c}$ respectively (the computations for $S^2_{\mu c}$ have been presented in Chapter \ref{quantumsphere}). 


We conclude this section with an important remark about the use of the phrase `orientation -preserving' in our terminology. We recall from Remark \ref{preliminaries_spectral_triples_spinors_caution} that by a `classical spectral triple' we always mean the spectral triple obtained by the Dirac operator on the spinors. This is absolutely crucial in view of the fact that the Hodge Dirac operator $d+d^*$ on the $L^2$-space of differential forms also gives a spectral triple of compact type on any compact Riemannian (not necessarily with a spin structure) manifold $M$, but the action of the full isometry group $ISO(M)$ (and not just the subgroup of orientation preserving isometries $ISO^+(M)$, even when $M$ is orientable) lifts to a  canonical unitary representation on this space commuting with $d+d^*$. In fact, the category of groups acting on $M$ such that the action comes from a unitary representation commuting with $d+d^*$,  has $ISO(M)$, and not $ISO^+(M)$, as its universal object.  So, one must stick to the Dirac operator on spinors to obtain the group of orientation preserving isometries in the usual geometric sense. This also has a natural quantum generalization, as we shall see in section \ref{qorient_section_comparison_with_QISO_L}.  

\section{Definition and existence of the quantum  group of orientation-preserving isometries}

\label{qorient_section_defn_existence_qtm_group_oreintation_pres_isometries}

\subsection{The classical case}
We first discuss the classical situation clearly, which will serve as a motivation for our quantum formulation. 

We begin with  a few basic facts about topologizing the space $C^\infty(M,N)$ where $M,N$ are smooth manifolds. 
Let $ \Omega $ be an open set of $\IR^n$.  We endow  $ C^{\infty} ( \Omega ) $ with the usual Frechet topology coming from  uniform convergence (over compact subsets) of  partial derivatives of all orders. The space $C^\infty(\Omega)$ is complete with respect to this topology, so is a Polish space in particular. Moreover, by the Sobolev imbedding theorem ( Corollary 1.21, \cite{rosenberg} ), $\cap_{k \geq 0} H_k ( \Omega ) = C^{\infty} ( {\Omega} ) $ as a set, where $H_k(\Omega)$ denotes the $k$-th Sobolev space. Thus, $C^\infty(\Omega)$ has also the Hilbertian seminorms coming from the Sobolev spaces, hence the corresponding Frechet topology. We claim that these two topologies on $C^\infty(\Omega)$ coincide.   Indeed, the inclusion map from $ C^{\infty} ( {\Omega} ) $  into $  \cap_k H_k ( \Omega )  $,  is continuous and surjective, so by the open mapping theorem for Frechet space, the inverse is also continuous, proving our claim. 

Given  two  second countable smooth manifolds $M,N$, we shall equip  $C^\infty(M,N)$ with  the weakest locally convex topology making $C^\infty(M,N) \ni \phi \mapsto f \circ \phi \in C^\infty(M)$ Frechet continuous for every $f $ in $C^\infty(N)$. 



For topological or smooth fibre or principal  bundles $E,F$ over a second countable smooth manifold $M$, we shall denote by ${\rm Hom}(E,F)$ the set of bundle  morphisms from $E$ to $F$.   We remark that 
 the total space of a locally trivial topological bundle such that the base and the fibre spaces are locally compact Hausdorff second countable must itself be so, hence in particular Polish ( that is, a complete separable metric space ).

 In particular,  if $E, F$ are locally trivial principal $G$-bundles over a common base, such that the (common) base as well as structure group $G$ are locally compact Hausdorff and second countable,  then ${\rm Hom}(E,F)$ is a  Polish space.

We need a standard fact, stated below as  Lemma \ref{selection}, about the measurable lift of Polish space valued functions.

Before that, we introduce some notions.

A multifunction $ G : X \rightarrow Y $ is a map with domain $ X $ and whose values are nonempty subsets of $ Y.$ For $ A \subseteq Y, $ we put $ G^{- 1} ( A ) = \{ x \in X : G ( x ) \cap A \neq \phi \}. $

A selection of a multifunction $ G :  X \rightarrow Y  $ is a point map $ s: X \rightarrow Y $ such that $ s ( x )  $ belongs to $G ( x ) $ for all $ x $ in $ X. $ Now let $ Y $ be a Polish space and $ \sigma_X $ a $ \sigma $-algebra on $ X. $ A multifunction $ G : X \rightarrow Y $ is called $ \sigma_X $ measurable if $ G^{- 1} ( U ) $ belongs to $ \sigma_X $ for every open set $ U $ in $ Y.$

 The following well known selection theorem is  Theorem 5.2.1 of \cite{sms} and was proved by Kuratowski and Ryll-Nardzewski. 

\bppsn

Let $ \sigma_X $ be a $ \sigma $ algebra on $ X $ and $ Y $ a Polish space. Then, every $ \sigma_X $ measurable, closed valued multifunction $ F : X \rightarrow Y $ admits a $ \sigma_X $ measurable selection.     

\eppsn

A trivial consequence of this result is the following:       
 
\blmma
\label{selection}
Let $M$ be a compact metrizable space, $B, \tilde{B}$ Polish spaces  such that  there is an $n$-covering map $\Lambda : \tilde{B} \raro B$.  Then any continuous map $\xi : M \raro B$ admits a lifting $\tilde{\xi} : M \raro \tilde{B}$ which is Borel measurable  and $\Lambda \circ \tilde{\xi}=\xi$.  In particular, if $\tilde{B}$ and $B$ are topological bundles over $M$, with $\Lambda$ being a bundle map, any continuous section of $B$ admits a lifting which is a measurable section of $\tilde{B}$. 
\elmma

We shall now give an operator-theoretic characterization of the classical group of orientation-preserving Riemannian isometries, which will be the motivation of our definition of its quantum counterpart. Let $M$ be a compact Riemannian $ n $ dimensional spin manifold, with a fixed choice of orientation. We recall  the notations as in subsection \ref{preliminaries_subsection_Spin_groups_Spin_manifolds}. 
 In particular, the spinor bundle $S$ is the associated bundle of a principal $Spin(n)$-bundle, say $P$, on $M$  which has a canonical $2$-covering bundle-map $\Lambda$ from $ P$  to the frame-bundle $ F$ (which is an $SO(n)$-principal bundle), such that locally $\Lambda$ is of the form $({\rm id}_M \ot \lambda)$ where $ \lambda $ is the two covering map from $ Spin ( n ) $ to $ SO ( n ).$ Moreover, the spinor space will be denoted by $ \Delta_n.$
   Let $f$ be a smooth orientation preserving Riemannian isometry of $M,$ and consider the bundles $E={\rm Hom}(F, f^*(F))$ and $\tilde{E}={\rm Hom}(P, f^*(P))$ (where ${\rm Hom}$ denotes the set of bundle maps). We view $df$ as a section of the bundle $E$ in the natural way.  By the Lemma \ref{selection} we obtain a measurable lift $\tilde{df} : M \raro \tilde{E}$, which is a measurable section of $\tilde{E}$. Using this, we define a map on the space of measurable section of $S=P \times_{Spin(n)} \Delta_n$ as follows: given a (measurable) section $\xi$ of S, say of the form $\xi(m)=[p(m), v]$, with $p(m) $ in $ P_m, v $ in $ \Delta_n$,   we define $U \xi$ by $ ( U \xi )(m)=[\tilde{df}(f^{-1}(m))(p(f^{-1}(m))), v] $. Note that sections of the above form constitute a total subset in $L^2(S)$, and the map $\xi \mapsto U \xi$   is clearly a densely defined linear map on $L^2(S),$ whose  fibre-wise action is unitary since the $Spin(n)$ action is so on $\Delta_n$. Thus it extends to a unitary $U$ on $\clh=L^2(S)$.   Any such $U$, induced by the map $f$, will be denoted by $U_f$ ( it is not unique since the choice of the lifting used in its construction is not unique ). 
\bthm
\label{classical_case}
Let $M$ be a compact Riemannian spin manifold (hence orientable , and fix a choice of orientation)  with the usual Dirac operator $D$ acting as  an unbounded self-adjoint operator on the Hilbert space $\clh$ of the square integrable spinors, and let $S$ denote the spinor bundle, with $\Gamma(S)$ being the $C^\infty(M)$ module of smooth sections of  $S$. Let $f: M \raro M$ be a smooth one-to-one  map which is a Riemannian orientation preserving isometry. Then the unitary $U_f$ on $\clh$  commutes with $D$ and  $U_f M_\phi U_f^*=M_{\phi \circ f}$, for any $\phi $ in $ C(M)$, where $M_\phi$ denotes the operator of multiplication by $\phi$ on $L^2(S)$. Moreover, when the dimension of $ M $ is even, $ U_f $ commutes with the canonical grading $ \gamma $ on $ L^{2} ( S ).$ 

Conversely, suppose that $U$ is a unitary on $\clh$ such that $ UD = DU $ and the map $\alpha_U(X)=UXU^{-1}$ for $X $ in $ \clb(\clh)$ maps $\cla=C(M)$ into $L^\infty(M)=\cla^{\prime \prime}$, then there is a smooth one-to-one orientation-preserving  Riemannian isometry $f$ on $M$ such that $U=U_f$.  We have the same result in the even case, if we assume furthermore that $ U \gamma = \gamma U.$


\ethm
{\it Proof:} From the construction of $U_f$, it is clear that $U_f M_\phi U_f^{-1}=M_{\phi \circ f}$. Moreover, since the Dirac operator $D$ commutes  with the $Spin(n)$-action on $S$, we have $U_f D=D U_f$ on each fibre, hence on $L^2(S)$. In the even dimensional case, it is easy to see that the $ Spin ( n ) $ action commutes with $ \gamma $ ( the grading operator ), hence $ U_f $ does so.

For the converse, first note that $\alpha_U$ is a unital $\ast$-homomorphism on $L^\infty(M, dvol)$ and thus must be of the form $\psi \mapsto \psi \circ f$ for some measurable $f$. We claim that $f$ must be smooth. Fix any smooth $g$ on $M$ and consider $\phi=g \circ f$. We have to argue that $\phi$ is smooth. Let $\delta_D$ denote the generator of the strongly continuous one-parameter group of automorphism $\beta_t(X)=e^{itD} X e^{-itD}$ on $\clb(\clh)$ (with respect to the weak operator topology, say). From the assumption that $D$ and $U$ commute it is clear that  $\alpha_U$ maps $\cld:=\bigcap_{n \geq 1} {\rm Dom}(\delta_D^n)$ into itself and since $C^\infty(M) \subset \cld$, we conclude that $\alpha_U(M_\phi)=M_{  \phi \circ g }$ belongs to $\cld$. We claim that this implies the smoothness of $\phi$. Let $m $ belongs to $ M$ and choose a local chart $(V, \psi)$   at $m$,  with the coordinates $(x_1,...,x_n)$, such that $\Omega=\psi(V) \subseteq \IR^n$ has compact closure, $S|_{V}$ is trivial and $D$ has the local  expression $D= i \sum_{j=1}^n \mu(e_j) \nabla_j$, where $\nabla_j=\nabla_{\frac{\partial}{\partial x_j}}$ denotes the covariant derivative (with respect to the canonical Levi Civita connection) operator along the vector field $\frac{\partial}{\partial x_j}$ on $L^2(\Omega)$ and $\mu(v)$ denotes the Clifford multiplication by a vector $v$. Now, $\phi \circ \psi^{-1} \in L^\infty(\Omega) \subseteq L^2(\Omega)$ and it is easy to observe from the above  local structure of $D$  that 
$[D, M_\phi]$ has the local expression $\sum_j i M_{\frac{\partial}{\partial x_j}\phi} \ot \mu(e_j)$. Thus,  
 the fact $M_\phi \in \bigcap_{n \geq 1} {\rm Dom}(\delta_D^n)$ implies  $\phi \circ \psi^{-1} $ is in $ {\rm Dom}(d_{j_1}...d_{j_k})$ for every integer tuples $(j_1,...,j_k)$, $j_i \in \{ 1,..., n\}$, where $d_j:=\frac{\partial}{\partial x_j}$.  In other words, $\phi \circ \psi^{-1} $ is in $ H^k(\Omega)$ for all $k \geq 1$, where $H^k(\Omega)$ denotes the $k$-th Sobolev space on $\Omega$ (see \cite{rosenberg}).  By Sobolev's Theorem (see, for example. \cite{rosenberg},Corollary 1.21, page 24)  it follows that $\phi \circ \psi^{-1} $ is in $ C^\infty(\Omega)$. 

We note that $ f $ is one-to-one as $ \phi \rightarrow \phi \circ f $ is an automorphism of $ L^{\infty}.$
 Now, we shall show that $f$ is an isometry of the metric space $(M,d)$, where $d$ is the metric coming from the Riemannian structure, and we have the explicit formula (  \ref{preliminaries_Dirac_operator_metric_expression} )
$$ d(p,q)={\rm sup}_{\phi \in C^\infty(M) , \| [D, M_\phi] \| \leq 1} |\phi(p)-\phi(q)|.$$ Since $U$ commutes with $D$, we have $\| [D, M_{\phi \circ f}]\|=\|[D, UM_\phi U^*]\|=\| U[D, M_\phi]U^* \|=\|[D, M_\phi]\|$ for every $\phi$, from which it follows that $d(f(p), f(q))=d(p,q)$.  Finally, $f$ is orientation preserving if and only if the volume form (say $\omega$), which defines the choice of orientation, is preserved by the natural action of $df$, that is, $(df \wedge .... \wedge df)(\omega)=\omega$. This will follow from the explicit description of $\omega$ in terms of $D$, given by (see \cite{Varilly_book} Page 26, also see \cite{connes_characterization} )

$$ \omega(\phi_0d\phi_1...d \phi_n)=\tau(\epsilon M_{\phi_0} [D, M_{\phi_1}]...[D, M_{\phi_n}]),$$

 where $\phi_0,...,\phi_n $ belong to $C^{\infty}(M), ~ \epsilon = 1 $ in the odd case and $ \epsilon = \gamma $ ( the grading operator ) in the even case and $\tau$ denotes the volume integral. In fact, $\tau(X)={\rm Lim}_{t \raro 0+} \frac{{\rm Tr}(e^{-tD^2}X)}{{\rm Tr}(e^{-tD^2})}$ ( where ${\rm Lim}$ is as in subsection \ref{preliminaries_subsection_spaceofforms_NCG} ), which implies $\tau(UXU^*)=\tau(X)$ for all $X $ in $ \clb(\clh)$ (using the fact that $D$ and $U$ commute). Thus, 
 
$$ \omega(\phi_0 \circ f ~d(\phi_1\circ f) \ldots d( \phi_n\circ f)) $$
$$ = \tau( \epsilon U M_{\phi_{0}}U^{*} U [D, M_{\phi_{1}}]U^{*}...U[D, M_{\phi_{n}}]U^{*}) $$
$$ =  \tau(U \epsilon M_{\phi_{0}} [D, M_{\phi_1}]...[D, M_{\phi_n}]U^*) $$
$$ = \tau( \epsilon M_{\phi_{0}} [D, M_{\phi_1}]...[D, M_{\phi_n}]) $$
$$ = \omega(\phi_0 d\phi_1...d \phi_n).$$ \qed

\vspace{4mm}


Now we turn to the case of a family of maps. We are ready to state and prove the operator-theoretic characterization of `set of orientation-preserving isometries'. 

\bthm

\label{classical_characterization_set_orientation_preserving_isometries}

Let $ X $ be a compact metrizable space and $ \psi: X \times M \rightarrow M  $ is a map such that
 $ \psi_{x} $  defined by $ \psi_{x} ( m ) = \psi ( x, m )$  is a smooth orientation preserving Riemannian isometry and $x \mapsto \psi_x \in C^\infty(M,M)$ is continuous with respect to the locally convex topology of $C^\infty(M,M)$ mentioned before.

Then 
there exists a ( $C(X)$-linear ) unitary $U_\psi$  on the  Hilbert $C(X)$-module $\clh \ot C(X)$ (where $\clh = L^2 ( S )$ as in Theorem \ref{classical_case})  such that for all $ x $ belonging to $ X,~U_x:= ( id \otimes ev_x ) U_{\psi} $ is a unitary of the form $U_{\psi_x}$  on the Hilbert space $\clh$ commuting with $ D $ and $ U_x M_{\phi} U^{- 1}_{x} = M_{\phi \circ \psi^{- 1}_{x}} $. If in addition, the manifold is even dimensional, then $U_{\psi_x}$ commutes with the grading operator $ \gamma. $

Conversely, if there exists a $C(X)$-linear  unitary $ U $ on $\clh \ot C(X)$ such that $ U_x:=({\rm id} \ot {\rm ev}_x)(U) $ is a unitary commuting with $ D $ for all $ x ,$ ( and $ U_x $ commutes with the grading operator $ \gamma $ if the manifold is even dimensional )  and $ ({\rm id} \ot {\rm ev}_x) \alpha_{U} ( L^{\infty}(M) ) \subseteq L^{\infty}(M) $ for all $x $ in $ X$,  then there exists a map $ \psi : X \times M \rightarrow M $ satisfying the conditions mentioned above such that $U=U_\psi$.

\ethm

{\it Proof:} Consider the bundles $\hat{F}=X \times F$ and $\hat{P}=X \times P$ over $X \times M$, with fibres at $(x,m)$ isomorphic with $F_m$ and $P_m$ respectively, and where  $ F $ and $ P $ are respectively the bundles of orthonormal frames and the Spin( n ) bundle discussed before. Moreover, denote by $\Psi$ the map from $X \times M$ to itself given by $(x,m) \mapsto (x, \psi(x,m)).$ Let $\pi_X:{\rm Hom}(\hat{F}, \Psi^*(\hat{F})) \raro X$ be the obvious
 map obtained by composing the projection map of the $X \times M$ bundle with the projection from $X \times M$ to $X$, and let us denote by  $B$ 
 the  closed subset of the Polish space  $ C(X, {\rm Hom}(\hat{F}, \Psi^*(\hat{F})))$ consisting of those $f$ such that for all $x$, $\pi_X(f(x))=x$. Define   $\tilde{B}$ in a similar way replacing $\hat{F}$ by $\hat{P}$. The covering map from $P$ to $F$ induces a covering map from $\tilde{B}$ to $B$ as well.  Let $ d^{'}_\psi: M  \raro B$ be  the map given by  $ d^{\prime}_{\psi} ( m ) ( x ) \equiv d^{\prime}_\psi(x,m) = d \psi_{x}|_{m}.$ Then by Lemma \ref{selection} there exists a measurable lift of $ d^{'}_\psi,$ say $ \widetilde{d^{'}_\psi} $ from $M$ into $\tilde{B} $. Since $d^{\prime}_\psi(x,m) \in {\rm Hom}(F_m, F_{\psi(x,m)})$, it is clear that the lift $\widetilde{d^\prime_\psi}(x,m)$ will be an element of ${\rm Hom}(P_m, P_{\psi(x,m)})$.

We can identify  $  \clh \otimes C ( X )  $ with $ C ( X \rightarrow  \clh ) ,$ and since $\clh$ has a total set $\clf$ (say)  consisting of sections of the form $[p(\cdot), v]$, where $p: M \raro P$ is a measurable section of $P$ and $v $ belongs to $ \Delta_n$, we have a total set $\tilde{\clf}$ of $ \clh \ot C(X)$ consisting of $\clf$ valued continuous functions from $X$. Any such function can be written as   $[\Xi,v]$ with $\Xi : X \times M  \raro P$, $v \in \Delta_n$, and $\Xi(x, m) \in P_m$,  and  we define $U$ on $\tilde{\clf}$ by $U[\Xi, v]=[\Theta, v]$, where
$$   \Theta( x,m ) =  \widetilde{d^{'}_\psi} ( x,\psi^{- 1}_{x} ( m ) ) ( \Xi (x, \psi^{- 1}_{x} ( m ) ) ).$$

 It is clear from  the construction of the lift that $U$ is indeed a $C(X)$-linear isometry which maps the total set $\tilde{\clf}$ onto itself, so extends to a unitary on the whole of $ \clh \ot C(X)$ with the desired properties.


Conversely, given $U$ as in the statement of the converse part of the theorem, we observe that for each $x $ in $ X$, by Theorem \ref{classical_case},  $ ( {\rm  id} \otimes ev_x ) U = U_{\psi_{x}} $ for some $ \psi_x $ such that $ \psi_x $ is a smooth  orientation preserving Riemannian isometry. This defines the map $\psi$ by setting $\psi(x,m)=\psi_x(m)$. The proof will be complete if we can  show that $x \mapsto \psi_x \in C^\infty(M,M)$ is continuous, which is equivalent to showing that whenever $x_n \raro x$ in the topology of $X$, we must have $ \phi \circ \psi_{x_n} \rightarrow \phi \circ \psi_x $ in the Frechet topology of $ C^{\infty}(M)$,  for any $\phi \in C^\infty(M)$. However, by  Lemma \ref{classical_case_family_ofmaps2}, we have $ ( {\rm  id} \otimes ev_{x_n} ) \alpha_U ( [ D, ~M_\phi ] ) \rightarrow ( {\rm  id} \otimes ev_x ) \alpha_U ( [ D, ~ M_\phi ] ) $ in the strong operator topology where  $ \alpha_U ( X ) = U X U^{- 1} $. Since  $ U $ commutes with $ D ,$ this implies $$ ( {\rm  id} \otimes ev_{x_n} ) [ D \otimes {\rm  id}, ~ \alpha_{U} ( M_\phi ) ] \rightarrow ( {\rm  id} \otimes ev_x ) [ D \otimes {\rm  id}, ~ \alpha_{U} ( M_\phi ) ] ,$$ that is, 
for all $ \xi $ in $ L^2 ( S ),$
$$[ D, M_{\phi \circ \psi_{x_n}} ] \xi \stackrel{L^{2}}{\rightarrow} [ D, ~ M_{\phi \circ \psi_x} ] \xi .$$ 

By choosing $\phi$ with support in a local trivializing coordinate neighborhood for $S$, and then using the local expression of $D$ used  in the proof of Theorem \ref{classical_case},
we conclude that 
$ d_k ( \phi \circ \psi_{x_n} ) \stackrel{L^2}{\rightarrow}  d_k( \phi \circ \psi_x ) $ (where $d_k$ is as in the proof of Theorem \ref{classical_case}).
Similarly, by taking repeated commutators with $D$, we can show the $ L^{2} $ convergence with $d_k$ replaced by $d_{k_1}...d_{k_m}$ for any finite tuple $(k_1,...,k_m)$.  In other words, $\phi \circ \psi_{x_n} \raro \phi \circ \psi_x$ in the topology of $C^\infty(M)$ described before. \qed

\subsection{Quantum group of orientation-preserving isometries of an $R$-twisted spectral triple}

\label{qorient_subsection_proof_of_main_theorem}
 
In view of the  characterization of orientation-preserving isometric action on  a classical manifold ( Theorem \ref{classical_characterization_set_orientation_preserving_isometries} ), we give the following definitions.
         \bdfn 
         \label{def_q_fam}
A quantum family of orientation preserving  isometries for the ( odd, compact type ) spectral triple $({\cla^\infty}, \clh, D)$ is given by a pair $(\cls, U)$ where $\cls$ is a separable unital $C^*$-algebra and  $U$ is a linear map from $\clh$ to $\clh \ot \cls$ such that $\widetilde{U}$ given by $\widetilde{U}( \xi \ot b)=U(\xi) (1 \ot b)$ ( $ \xi $ in $ \clh, ~ b $ in $ \cls $ ) extends to a unitary element of  $ \clm(\clk(\clh) \ot \cls)$ satisfying the following:
         
(i) for every state $\phi$ on $\cls$ we have $U_\phi D=DU_\phi$, where $U_\phi:=({\rm id} \ot \phi)(\widetilde{U})$;

(ii) $({\rm id} \ot \phi) \circ \alpha_U(a) \in ({\cla^\infty})^{\prime \prime}$ for all $ a $ in $ \cla^\infty$ and state $\phi$ on $\cls$, where $\alpha_U(x):=\widetilde{U}( x \ot 1) {\widetilde{U}}^* $ for $x $ belonging to $ \clb(\clh)$.

In case the $C^*$-algebra $\cls$ has a coproduct $\Delta$ such that $(\cls,\Delta)$ is a compact quantum group and $U$ is a unitary representation  of $(\cls, \Delta)$ on $\clh$, we say that $(\cls, \Delta)$ acts  by orientation-preserving isometries  on the spectral triple.

In case the spectral triple is even with the grading operator $ \gamma ,$ a quantum family of orientation preserving isometries $ ( \cla^{\infty}, \clh, D, \gamma ) $ will be defined exactly as above, with the only extra condition being that $ U $ commutes with $ \gamma. $

\edfn

From now on, we will mostly consider odd spectral triples. However let us remark that in the even case,  all the definitions and results obtained by us will go through with some obvious modifications. We also remark that all our spectral triples are of compact type.

Consider the category ${\bf Q}\equiv {\bf Q}(\cla^\infty, \clh, D)\equiv {\bf Q}(D)$ with the object-class consisting of all quantum families of 
orientation preserving isometries $(\cls, U)$ of the given spectral triple, and the set of morphisms ${\rm Mor}((\cls,U),(\cls^\prime,U^\prime))$ 
being the set of unital $C^*$-homomorphisms $\Phi : \cls \raro \cls^\prime$ satisfying $({\rm id} \ot \Phi) (U)=U^\prime$. 
We also consider another category ${\bf Q}^\prime \equiv {\bf Q}^\prime(\cla^\infty, \clh, D) \equiv {\bf Q}^\prime(D)$ whose objects are triplets $(\cls, \Delta, U)$,  where $(\cls,\Delta)$ is a compact quantum group acting by orientation preserving isometries on the given spectral triple, with $U$ being the corresponding unitary representation. The morphisms  are the homomorphisms of compact quantum groups which are also morphisms of the underlying quantum families of orientation preserving isometries. The forgetful functor $F: {\bf Q}^\prime \raro {\bf Q}$ is clearly faithful, and we can view $F({\bf Q}^\prime)$ as a subcategory of ${\bf Q}$.

 Unfortunately, in general ${\bf Q}^\prime$ or ${\bf Q}$ will not have a universal object. It is easily seen by taking the standard example $\cla^\infty=M_n(\IC)$, $\clh=\IC^n$, $D=I.$ Any CQG having a unitary representation on $ \IC^n $ is an object of $ {\bf Q}^\prime ( M_n ( \IC ), \IC^n, I ).$ 
 But by Proposition \ref{preliminaries_action_the_work_of_wang}, there is no universal object in this category. However, the fact that comes to our rescue is that a universal object exists in each of the subcategories which correspond to the CQG actions preserving a given faithful functional on $M_n$.
 
  On the other hand, given any equivariant spectral triple, it has been shown in \cite{goswami_rmp} that there is a (not necessarily unique) canonical faithful functional which is preserved by the CQG action.  For readers' convenience, we state this result  (in a form suitable to us) briefly here. Before that, let us recall the definition of an $R$-twisted spectral data from \cite{goswami_rmp}.

 \bdfn
 
 \label{qorient_R_twisted_data}
 
  An $R$-twisted spectral data ( of compact type ) is given by a quadruplet  ($\cla^{\infty}, \clh, D, R$) where
  
  1. ( $\cla^{\infty}, \clh, D$ ) is a spectral triple of compact type.
  
  2. $R$ a positive (possibly unbounded) invertible operator such that $R$ commutes with $D$.
  
  3. For all $s \in \IR,$ the map $a \mapsto \sigma_s(a):=R^{-s} a R^s$ gives an automorphism of $\cla^\infty$ (not necessarily $\ast$-preserving) satisfying  $\sup_{s \in [-n,n]} \| \sigma_s(a) \| < \infty$ for all positive integer $n.$
  
  \edfn
  
  We shall also sometimes refer to $ ( \cla^{\infty}, \clh, D )  $ as an $ R $-twisted spectral triple.

\bppsn

\label{5678}

Given a spectral triple $(\cla^\infty, \clh, D)$ (of compact type)  which is $\clq$-equivariant with respect to a representation of a CQG $\clq$ on $\clh$, we can construct a positive (possibly unbounded) invertible operator $R$ on $\clh$ such that $(\cla^\infty, \clh, D, R)$ is  a twisted spectral data,
 and moreover, we have 

 $\alpha_U$ preserves the functional $\tau_R$ defined at least on a weakly dense $\ast$-subalgebra $\cle_D$ of $\clb(\clh)$ generated by the rank-one operators of the form $|\xi><\eta|$ where $\xi, \eta$ are eigenvectors of $D,$ given  by $$ \tau_R(x)= Tr ( R x ),~~x \in \cle_D.$$
   
\eppsn

\brmrk

\label{qorient_haarstate_tracial_implies_R=I}

When the Haar state of $ \clq $ is tracial, then it follows from the definition of $ R $ in Lemma 3.1 of \cite{goswami_rmp} and Theorem 1.5 part 1. of \cite{woro} that $ R $ can be chosen to be $ I. $

\ermrk

\brmrk
If $V_\lambda$ denotes the eigenspace  of $D$ corresponding to the eigenvalue, say $ \lambda ,$  it is clear that $ \tau_{R} (X) = e^{t \lambda^2} {\rm Tr}(Re^{-tD^2}X) $ for all $X =|\xi><\eta|$ with $\xi, \eta $ belonging to $ V_\lambda$ and for any $t>0$. Thus, the $\alpha_U$-invariance of the functional $\tau_R$ on $\cle_D$ is equivalent to the $\alpha_U$-invariance of the functional $X \mapsto {\rm Tr}(X Re^{-tD^2})$ on $\cle_D$ for each $t>0$.
 If, furthermore, the $R$-twisted spectral triple is $\Theta$-summable in the sense that  $Re^{-tD^2}$ is trace class for every $t>0$, the above is also equivalent to the $\alpha_U$-invariance of the bounded normal functional $X \mapsto {\rm Tr}(X Re^{-tD^2})$ on the whole of $\clb(\clh)$.  In particular, this implies that  $\alpha_U$ preserves the  functional $ \clb(\clh) \ni x \mapsto {\rm Lim}_{t \raro 0+} \frac{{\rm Tr}(xRe^{-tD^2})}{{\rm Tr}(Re^{-tD^2})},$ where ${\rm Lim}$ is as defined in subsection \ref{preliminaries_subsection_spaceofforms_NCG}. 
\ermrk

This motivates the following definition:
\bdfn
Given an  $R$-twisted spectral data $(\cla^{\infty}, \clh, D, R)$ of compact type, a quantum family  of orientation preserving isometries $(\cls, U)$ of $(\cla^\infty, \clh, D)$ is said to preserve the $R$-twisted volume, (simply said to be  volume-preserving if $R$ is understood) if one has  $(\tau_R \otimes  {\rm id} ) (\alpha_U(x))= \tau_R(x).1_\cls$ for all $x $ in $ \cle_D$, where $\cle_D$ and $\tau_R$ are as in Proposition \ref{5678}. We shall also call $(\cls, U)$ a quantum family of orientation-preserving isometries of the $R$-twisted spectral triple. 

If, furthermore, the  $C^*$-algebra $\cls$ has a coproduct $\Delta$ such that $(\cls,\Delta)$ is a CQG and $U$ is a unitary representation  of $(\cls, \Delta)$ on $\clh$, we say that $(\cls, \Delta)$ acts  by volume and orientation-preserving isometries  on the $R$-twisted spectral triple.

 We shall consider the categories ${\bf Q}_R \equiv {\bf Q}_R(D)$ and ${\bf Q}^\prime_R \equiv {\bf Q}^\prime_R(D)$ which are the full subcategories of ${\bf Q}$ and ${\bf Q}^\prime$ respectively, obtained by restricting the object-classes to the volume-preserving quantum families. 
\edfn

\brmrk

We shall not need the full strength of the definition of twisted spectral data here; in particular the third condition  in the definition \ref{qorient_R_twisted_data}. However, we shall continue to work with the usual definition of $R$-twisted spectral data, keeping in mind that all our results are valid without assuming the third condition.

\ermrk

Let us now fix a spectral triple $(\cla^\infty, \clh, D)$ which is of compact type. The $C^*$-algebra generated by $\cla^\infty$ in $\clb(\clh)$ will be denoted by $\cla$. Let $\lambda_0=0, \lambda_1, \lambda_2,\cdot \cdot \cdot  $ be the eigenvalues of $D$ with $V_i$ denoting the ( $d_i$-dimensional, $ 0 \leq d_i<\infty$ ) eigenspace for $\lambda_i$. Let $\{ e_{ij}, j=1,..., d_i \}$ be an orthonormal basis of $V_i$. 
We also assume that there is a positive  invertible $R$ on $\clh$ such that $(\cla^\infty, \clh, D, R)$ is 
 an $R$-twisted   spectral data. The operator $R$ must have the form $R|_{V_i}=R_i,$ say, with $R_i$ positive invertible $d_i \times d_i$ matrix. Let us denote the CQG $A_{u,d_i} (R^T_i)$ by $\clu_i$, with its canonical unitary representation $\beta_i$ on $V_i \cong \IC^{d_i}$, given by $\beta_i(e_{ij})=\sum_k e_{ik} \ot u^{R^T_i}_{kj}$. Let $\clu$ be the free product of $\clu_i, i=1,2,...$ and $\beta=\ast_i \beta_i$ be the corresponding free product representation of $\clu$ on $\clh$. We shall also consider the  corresponding  unitary element $\tilde{\beta}$ in $\clm(\clk(\clh) \ot \clu)$.
 
 \blmma
 
\label{lem2} 

Consider the $R$-twisted spectral triple $(\cla^\infty,\clh,D)$  and let $(\cls,U)$ be a quantum family of volume and orientation preserving  isometries of the given spectral triple.  Moreover, assume that the map  $U$ is faithful in the sense that there is no
proper  $C^*$-subalgebra $\cls_1$ of $\cls$ such that
$\widetilde{U} $ belongs to $  \clm(\clk(\clh ) \ot \cls_1)$. 
Then  we
can find a $C^*$-isomorphism  $\phi : \clu/
\cli \raro \cls$ between $\cls$ and a quotient of $\clu$ by a
 $C^*$-ideal
 $\cli$ of $\clu$, such that $ U= ({\rm id}\ot \phi) \circ ({\rm id} \ot \Pi_\cli) \circ
 \beta$, where $\Pi_\cli$ denotes the quotient map from $\clu$ to
 $\clu/\cli$.

 If, furthermore, there is a compact quantum group structure on $\cls$ given by a coproduct  $\Delta$ such that $(\cls,\Delta, U)$ is an object in ${\bf Q}^\prime_R$, the ideal $\cli$ is a Woronowicz $C^*$-ideal and the $C^*$-isomorphism $\phi : \clu/ \cli \raro \cls$ is a morphism of compact quantum groups.
 
\elmma

{\it Proof :} It is clear that $U$  maps $V_i $ into $V_i \otimes
\cls$ for each $i$. Let $v^{(i)}_{kj}$ ($j,k=1,...,d_i$) be the
elements of $\cls$ such that $U(e_{ij})=\sum_k e_{ik} \otimes
v^{(i)}_{kj}$.  Note that $v_i:=((v^{(i)}_{kj} ))$ is a unitary in
$M_{d_i}(\IC) \otimes \cls$. Moreover, the $\ast$-subalgebra
generated by all $ \{ v^{(i)}_{kj}, i \geq 0, ,j,k \geq 1 \}$ must be
dense in $\cls$ by the assumption of faithfulness.

Consider the $\ast$-homomorphism $\alpha_i $ from the finite dimensional $C^*$ algebra $\cla_i \cong M_{d_i}(\IC)$ generated by the rank one operators $\{ |e_{ij}><e_{ik}|, j,k=1,..., d_i \}$ to $\cla_i \ot \cls$ given by $\alpha_i(y)=\widetilde{U}(y \ot 1){\widetilde{U}}^*|_{V_i}$. Clearly,  the restriction of the  functional $\tau_R$ on $\cla_i$ is nothing but the  functional given by ${\rm Tr}( R_i~ \cdot)$, where ${\rm Tr}$ denotes the usual trace of matrices. Since $\alpha_i$ preserves this functional by assumption, we get,  by the universality of $\clu_i$,  a
$C^*$-homomorphism from $\clu_i$ to $\cls$ sending $u^{(i)}_{kj} \equiv u_{kj}^{R^T_i} $
to $v^{(i)}_{kj}$, and by definition of the free product, this
induces a $C^*$-homomorphism, say $\Pi$, from $\clu$ onto $\cls$,
 so that $\clu/\cli
\cong \cls$, where  $\cli:={\rm Ker}(\Pi)$.

In case $\cls$ has a coproduct $\Delta$ making it into a compact quantum group and $U$  is a quantum group representation, it is easy to see that  the subalgebra of $\cls$ generated by $\{ v^{(i)}_{kj},~i \geq 0, j,k=1,...,d_i \}$ is a Hopf
algebra, with $\Delta(v^{(i)}_{kj})=\sum_l v^{(i)}_{kl} \ot v^{(i)}_{lj}$. From this, it follows that  $\Pi$ is Hopf-algebra morphism, hence $\cli$ is a Woronowicz $C^*$-ideal.
\qed 

\vspace{2mm}

 \bthm
\label{main} 
For any  $R$-twisted  spectral triple of compact type $(\cla^\infty, \clh, D)$, the category ${\bf Q}_R$ of quantum families of volume and orientation preserving isometries has a universal (initial) object, say  $(\widetilde{\clg}, U_0)$. Moreover, $\widetilde{\clg}$ has a coproduct $\Delta_0$ such that $(\widetilde{\clg},\Delta_0)$ is a compact quantum group and $(\widetilde{\clg},\Delta_0,U_0)$ is a universal object in the category ${\bf Q}^\prime_R$.  The representation  $U_0$ is faithful.
 \ethm
 {\it Proof :}
Recall the $C^*$-algebra $\clu$ considered before, along with  the representation $\beta $ and the corresponding  unitary  $\widetilde{\beta} $ belonging to $ \clm(\clk(\clh) \ot \clu ) $. For any $C^*$-ideal $\cli$ of $\clu$, 
we shall denote by $\Pi_\cli$ the canonical quotient map from $\clu$ onto $\clu/\cli$, and let $\Gamma_\cli=({\rm id} \ot \Pi_\cli) \circ \beta$. Clearly,   $\widetilde{\Gamma_\cli}=({\rm id} \ot \pi_\cli)\circ \widetilde{\beta}$ is a unitary element of  $\clm(\clk(\clh) \ot \clu/\cli )$. 
Let $\clf$ be the collection of all those  $C^*$-ideals
$\cli$ of $\clu$ such that  $({\rm id} \ot \omega) \circ \alpha_{\Gamma_\cli } \equiv ({\rm id} \ot \omega) \circ {\rm ad}_{\widetilde{\Gamma_\cli} }$ maps $\cla^\infty$ into $\cla^{\prime \prime}$ for every state $\omega$ (equivalently, every bounded linear functional) on $\clu/\cli$.  This
collection is nonempty, since the trivial one-dimensional $C^*$-algebra $\IC$ gives an object in ${\bf Q}_R$ and by Lemma \ref{lem2} we do get a member of $\clf$. Now,
let $\cli_0$ be the intersection of all ideals in $\clf$. We claim
that $\cli_0$ is again a member of $\clf$.  
Indeed, in the notation of Lemma \ref{lim}, it is clear that for $a $ in $ \cla^\infty$, $({\rm id} \ot \phi)\circ \widetilde{\Gamma}_{\cli_0}(a) $ belongs to $ \cla^{\prime \prime}$ for all $\phi $ in the convex hull of $ \Omega \bigcup (-\Omega)$. Now, for any state $\omega$ on $\clu/\cli_0$, we can find, by Lemma \ref{lim}, a net $\omega_j $ in the above convex hull  (so in particular $\| \omega_j \| \leq 1 $ for all $j$) such that $\omega(x+\cli_0)=\lim_j \omega_j(x+ \cli_0)$ for all $x $ in $ \clu/\cli_0$. 

 It follows from  Lemma \ref{classical_case_family_ofmaps2} that $({\rm id} \ot \omega_j)(X) \raro ({\rm id} \ot \omega)(X)$ (in the strong operator topology)  for all $X $ belonging to $ \clm(\clk(\clh) \ot \clu/\cli_0)$. 
Thus, for $a $ in $ \cla^\infty$, $ ({\rm id} \ot \omega)\circ {\alpha}_{ \widetilde{\Gamma}_{\cli_0}}(a)$ is the limit of $({\rm id} \ot \omega_i)\circ {\alpha}_{ \widetilde{\Gamma}_{\cli_0}}(a)$ in the strong operator topology, hence belongs to $\cla^{\prime \prime}$.

We now show   that $(\widetilde{\clg}:=\clu/\cli_0, \Gamma_{\cli_0})$ is a 
universal object in ${\bf Q}_R$. To see this, consider any object $(\cls, U)$ of ${\bf Q}_R$.  Without loss
of generality we can assume $U$ to be faithful, since
otherwise we can replace $\cls$ by the  $C^*$-subalgebra
generated by the elements $\{ v^{(i)}_{kj} \}$ appearing in the proof of Lemma \ref{lem2}.  But by
Lemma \ref{lem2} we can further assume that $\cls$ is isomorphic
with $\clu/\cli$ for some $\cli $ belonging to $ \clf$. Since $\cli_0 \subseteq
\cli$, we have a $C^*$-homomorphism from
$\clu/\cli_0$ onto $\clu/\cli$, sending $x+\cli_0$ to $x +\cli$, which is clearly a morphism in the category ${\bf Q}_R$. This is indeed the unique such morphism, since it is uniquely determined on the dense subalgebra generated by $\{ u^{(i)}_{kj}+ \cli_0,~i\geq 0,~j,k \geq 1 \}$ of $\widetilde{\clg}$.

To construct the coproduct on $\widetilde{\clg}=\clu/\cli_0$, we first consider $U^{(2)}: \clh \raro \clh \ot \widetilde{\clg} \ot \widetilde{\clg}$ given by   $$U^{(2)}=(\Gamma_{\cli_0})_{(12)}(\Gamma_{\cli_0})_{(13)},$$ where $U_{ij}$ is the usual `leg-numbering notation'. It is easy to see that $(\widetilde{\clg} \ot \widetilde{\clg}, U^{(2)})$ is an object in the category ${\bf Q}_R$, so by the universality of $(\widetilde{\clg},\Gamma_{\cli_0})$, we have a unique unital $C^*$-homomorphism $\Delta_0 : \widetilde{\clg} \raro \widetilde{\clg} \ot \widetilde{\clg}$ satisfying $$ ({\rm id} \ot \Delta_0) (\Gamma_{\cli_0}) =U^{(2)}.$$ 
Letting both sides act on  $e_{ij}$, we get $$ \sum_l e_{il} \ot (\pi_{\cli_0} \ot \pi_{\cli_0}) \left( \sum_k u^{(i)}_{lk} \ot u^{(i)}_{kj} \right)=\sum_l e_{il} \ot \Delta_0(\pi_{\cli_0}(u^{(i)}_{lj})).$$ Comparing coefficients of $e_{il}$, and recalling that $\tilde{\Delta}(u^{(i)}_{lj})=\sum_k u^{(i)}_{lk} \ot u^{(i)}_{kj}$ (where $\tilde{\Delta}$ denotes the coproduct on $\clu$), we have \be \label{coprod}(\pi_{\cli_0} \ot \pi_{\cli_0}) \circ \tilde{\Delta}=\Delta_0 \circ \pi_{\cli_0} \ee on the linear span of $\{ u^{(i)}_{jk}, i \geq 0,~j,k \geq 1 \}$, and hence on the whole of $\clu$. 
This implies that $\widetilde{\Delta}$ maps $\cli_0={\rm Ker}(\pi_{\cli_0})$ into ${\rm Ker}(\pi_{\cli_0} \ot \pi_{\cli_0})=(\cli_0 \ot 1 + 1 \ot \cli_0) \subset \clu \ot \clu$. In other words, $\cli_0$ is a Hopf $C^*$-ideal, and hence $\widetilde{\clg}=\clu/\cli_0$ has the canonical compact quantum group structure as a quantum subgroup of $\clu$. It is clear from the relation (\ref{coprod}) that  $\Delta_0$ coincides with the canonical coproduct of the quantum subgroup $\clu/\cli_0$ inherited from that of $\clu$. It is also easy to see that the object $(\widetilde{\clg}, \Delta_0, \Gamma_{\cli_0})$ is universal in the category ${\bf Q}^\prime_R$, using the fact that (by Lemma \ref{lem2}) any compact quantum group $(\cls, U)$ acting by volume and orientation preserving isometries on the given spectral triple is isomorphic with a quantum subgroup $\clu/\cli$, for some Hopf $C^*$-ideal $\cli$ of $\clu$.

Finally, the faithfulness of $U_0$ follows from the universality by standard arguments which we briefly sketch. If $\clg_1 \subset \widetilde{\clg}$ is a $\ast$-subalgebra of $\widetilde{\clg}$ such that $\widetilde{U_0} $ is an element of $ \clm(\clk(\clh) \ot \clg_1)$, it is easy to see that $(\clg_1, U_0)$ is also an object in $ {\bf \clq_{R}},$ and by definition of universality of $\widetilde{\clg}$ it follows that there is a unique morphism, say $j$, from $\widetilde{\clg}$ to $\clg_1$. But the map $i \circ j$ is a morphism from $\widetilde{\clg}$ to itself, where $i : \clg_1 \raro \widetilde{\clg}$ is the inclusion. Again by universality, we have that $i \circ j ={\rm id}_{\widetilde{\clg}}$, so in particular, $i$ is onto, that is,. $\clg_1=\widetilde{\clg}$.  
\qed

\vspace{4mm}

Consider the  $\ast$-homomorphism $\alpha_0:=\alpha_{U_0}$, where  $(\widetilde{\clg},U_0)$ is the universal object obtained in the previous theorem. For every state $\phi$ on $\widetilde{\clg}$, $({\rm id} \ot \phi) \circ \alpha_0$ maps $\cla$ into $\cla^{\prime \prime}$. However, in general $\alpha_0$ may not be faithful  even if $U_0$ is so, and let $\clg$ denote the $C^*$-subalgebra of $\widetilde{\clg}$ generated by the elements $\{ (f \ot {\rm id})\circ \alpha_0(a):  f \in  \cla^*, a \in \cla \}$. 

\brmrk

If the spectral triple is even, then all the proofs above go through with obvious modifications.

\ermrk

\bdfn

We shall call $\clg$ the quantum  group of orientation-preserving isometries of  R-twisted spectral triple $(\cla^\infty, \clh, D, R)$  and denote it by  $ QISO^{+}_R (\cla^{\infty}, \clh, D, R) $ or even simply as $QISO^{+}_{R}(D)$. The quantum group $\widetilde{\clg}$ is denoted by $\widetilde{QISO^{+}_{R}}(D)$. 

If the spectral triple is even, then we will denote 
$ \clg $ and $ \widetilde{\clg} $ by $ QISO^{+}_{R} ( D, \gamma ) $ and  $ \widetilde{QISO^{+}_{R}} ( D, \gamma ) $ respectively.

\edfn

\subsection{Stability and topological action}

\label{qorient_subsection_stability_top_action}

In this subsection, we are going to use the notations as in subsection \ref{qorient_subsection_proof_of_main_theorem}, in particular $ \widetilde{\clg}, \clg, U_0, \alpha_0.$
It is not clear from the definition and construction of $QISO_R^+(D)$ whether the $C^*$ algebra $\cla$ generated by $\cla^\infty$ is stable under $\alpha_0$   in the sense that $({\rm id} \ot \phi) \circ \alpha_0$ maps $\cla$ into $\cla$ for every $\phi$. Moreover, even if $\cla$ is stable, the question remains whether $\alpha_0$ is a $C^*$-action of the CQG $QISO^+_R(D)$. In chapter \ref{quantumsphere},  subsection \ref{quantumsphere_subsection_chakpal_computation} we have given an example of a spectral triple for which the $ \ast $-homomorphism $ \alpha_0 $ is not a $ C^* $ action. However, one can prove that $ \alpha_0 $ is a  $ C^* $ action for a rather large class of spectral triples, including the cases mentioned below.

(i) For any spectral triple for which there is a `reasonable' Laplacian in the sense of \cite{goswami}. This includes all classical spectral triples as well as their Rieffel deformation (with $R=I$).

(ii) Under the assumption that there is an eigenvalue of $D$ with a one-dimensional eigenspace spanned by a cyclic separating vector $\xi$ such that any 
eigenvector of $D$ belongs to the span of $\cla^\infty \xi$ and $\{ a \in \cla^\infty:~a \xi ~{\rm is}~{\rm an}~{\rm eigenvector}~{\rm of}~D \}$ is norm-dense in $\cla$ ( to be proved in subsection \ref{qorient_subsection_univ_in_Q_Q_prime} ).

\vspace{2mm}

Now we prove the sufficiency of the condition (i).

We begin with a sufficient condition for stability of $\cla^\infty$ under $\alpha_0$.
 Let $ ( \cla^{\infty}, \clh, D ) $ be a (compact type) spectral triple such that \vspace{1mm}\\
({\bf 1}) ${\cla^\infty}$ and $\{ [D,a], ~a \in {\cla^\infty} \}$ are 
contained in the domains of all powers of the derivations $[ D, \cdot ] $ and $ [|D|, \cdot] .$ \\

We will denote by $ \widetilde{T_t} $ the one parameter group of $ \ast$-automorphisms on $ \clb ( \clh ) $ given by $ \widetilde{T_t} ( S ) =
 e^{ it D} S e^{-it D} $  for all $ S $  in $ \clb ( \clh )$ which is clearly continuous in SOT. We will denote the generator of this group by $ \delta.$  
For $ X $ such that $ [ D, X ] $ is bounded, we have $ \delta ( X ) = i[ D, X ]  $ and hence
$${\left\| \widetilde{T_t} ( X ) - X \right\|} = {\left\| \int^{t}_{0} \widetilde{T_s} ([ D, X ])ds \right\|} \leq ~ t {\left\| [ D, X ] \right\|} .$$

Let us say that the spectral triple satisfies the {\it Sobolev condition} if $$\cla^\infty=\cla^{\prime \prime} \bigcap_{n \geq 1} {\rm Dom}(\delta^n).$$  Then we have the following result, which is a natural generalization of the classical situation, where a measurable isometric action automatically becomes topological (in fact smooth).
\bthm
\label{stable}
(i) For every state $\phi$ on $\clg$,  $$({\rm id} \ot \phi)\circ \alpha_0(\cla^\infty) \text{ belongs to } \cla^{\prime \prime} \bigcap_{n \geq 1} {\rm Dom}(\delta^{n}).$$ 

(ii) If the spectral triple satisfies the  Sobolev condition then  
  $\cla^\infty$ ( and hence $\cla$ ) is stable under $\alpha_0$.   
\ethm
{\it Proof:} Since $U_0$ commutes with $D \ot I$, it is clear that the automorphism group $\widetilde{T_t}$ commutes with $\alpha_0^\phi \equiv ({\rm id} \ot \phi) \circ \alpha_0$, and thus by the continuity of $\alpha_0$ in the strong operator topology it is easy to see that, for $a $ in $ {\rm Dom}(\delta),$  \bean \lefteqn {\lim_{t \raro 0+} \frac{\widetilde{T_t}(\alpha_0^\phi(a))-\alpha_0^\phi(a)}{t}}\\
&=&   \lim_{t \raro 0+} \alpha_0^\phi(\frac{\widetilde{T_t}(a)-a}{t})\\
&=& \alpha_0^\phi(\delta(a)).\eean
Thus, $\alpha_0^\phi$ leaves ${\rm Dom}(\delta)$ invariant and commutes with $\delta$. Proceeding similarly, we prove (i). The assertion (ii)  is a trivial consequence of (i) and the  Sobolev condition. \qed

\vspace{4mm}

Let us now assume 

({\bf 2}) The spectral triple is $\Theta$-summable, that is, for every $t>0$, $e^{-tD^2}$ is trace-class and the functional  $ \tau ( X ) = 
{\rm Lim}_{t \rightarrow 0} \frac{ {\rm Tr} ( X e^{- t D^2} )}{ {\rm Tr} ( e^{- t D^2} )}$ ( where ${\rm Lim}$ is as in subsection \ref{preliminaries_subsection_spaceofforms_NCG} ), is a positive faithful trace on the $ \ast $ algebra, say $\cls^\infty$, generated by $ \{ \widetilde{T}_s ( \cla^{\infty} ),  \widetilde{T}_s ( \cla^{\infty} ) ( [ D, a ] ) : a \in \cla^{\infty}\}. $

The functional $\tau$ is to be interpreted as the volume form ( we refer to \cite{fro}, \cite{goswami} for the details ). The completion  of $\cls^\infty$ in the norm of $\clb(\clh)$ is denoted by $\cls,$ and we shall denote  by $\| a \|_2$ and $\| \cdot \|_\infty $ the $L^2$-norm $\tau(a^*a)^{\frac{1}{2}}$  and the operator norm of $\clb(\clh)$ respectively.

From the definition of $ \tau ,$ it is also clear that $ \widetilde{T_t} $ preserves $ \tau ,$ so extends to a group of unitaries on $\cln:=L^2( \cls^\infty, \tau)$. 
 Moreover,  for $X$ such that $[D,X] $ is in $ \clb(\clh)$, in particular for $X $ in $ \cls^{\infty} ,$ we  have 
\bean \lefteqn{\left\| \widetilde{T_s} ( X ) - X  \right\|^2_2 }\\
& =& \tau ( {\widetilde{T_s} ( X )}^{*} ( \widetilde{T_s} ( X ) - X ) ) + \tau ( X^* ( X - \widetilde{T_s} ( X ) )  ) \\
  &  \leq &   2  \left\| X - \widetilde{T_s} ( X )  \right\|_{\infty} \left\| X \right\|_2 \\
& \leq & 2 s \| [D,X] \|_\infty \| X \|_2,\\
\eean
which clearly shows that $s \mapsto \widetilde{T_s}(X)$ is $L^2$-continuous for $X $ belonging to $ \cls^\infty$, hence (by unitarity of $\widetilde{T_s}$) on the whole of $\cln$,
 that is, it is a strongly continuous one-parameter group of unitaries. Let us denote its generator by $\tilde{\delta}$, which is a skew adjoint map, that is, $i \tilde{\delta}$ is self adjoint, and $\widetilde{T_t}={\rm exp}(t \tilde{\delta})$. Clearly, $\tilde{\delta}=\delta=[D, \cdot]$ on $\cls^\infty.$ 

We will denote $ L^{2}  ( \cla^\infty, \tau ) \subset \cln$ by $ \clh^{0}_{D}$ and the restriction of $\tilde{\delta}$ to $\clh^0_D$ (which is a closable map from $\clh^0_D$ to $\cln$) by $d_D$. Thus, $d_D$ is closable too. 

We now recall the assumptions made in chapter \ref{qisol}, subsection \ref{qisol_subsection_definition_existence_of_quantum_isometry_group} for defining  the `Laplacian' and the corresponding  quantum isometry group of a spectral triple $(\cla^\infty,\clh,D).$ 

The following conditions will also be assumed throughout the rest of this subsection:\vspace{1mm}

({\bf 3}) $ \cla^{\infty} \subseteq {\rm Dom} ( \cll ) $ where $ \cll \equiv \cll_D := - d^{*}_{D} d_{D} .$\\
({\bf 4}) $ \cll $ has compact resolvent.\\
({\bf 5}) Each eigenvector of $ \cll $ ( which has a discrete spectrum , hence a complete set of eigenvectors ) belongs to $ \cla^{\infty}.$\\
({\bf 6}) The complex linear span of the eigenvectors of $ \cll $, say $ \cla^{\infty}_{0} $ ( which is a subspace of $ \cla^{\infty} $ by assumption ({\bf 5})  ), is norm dense in $ \cla^{\infty}.$\vspace{2mm}\\

It is clear that   $\cll $ maps $ (\cla^\infty_0) $ into itself. The $\ast$-subalgebra  of $\cla^\infty$ generated by $\cla^\infty_0$ is denoted by $\cla_0$.
  We also note that $\cll=P_0\tilde{\cll} P_0$, where $\tilde{\cll}:= (i \tilde{\delta})^2$ (which is a self adjoint operator on $\cln$) and $P_0$ denotes the orthogonal  projection in $\cln$ whose range is  the subspace $\clh^0_D$.

\bthm

\label{QISO_I_D_<_Q_L_D}
 
 Let $( \cla^{\infty}, \clh, D )$ be a spectral triple satisfying the assumptions $({\bf 1})-({\bf 6})$ made above. In addition, assume that  at least one of conditions (a) and (b) mentioned below is satisfied:\\ 
 
(a) $ {\cla}^{\prime \prime} \subseteq \clh^{0}_{D} .$\\
(b) $ \alpha_0^{\phi} ( \cla^{\infty} ) \subseteq \cla^{\infty} $ for every state $\phi$ on $\clg=QISO^+_I(D)$.\\
 
 Then  $ \alpha_0 $ is a $C^*$-action of $ QISO^{+}_{I} ( D ) $  on $ \cla .$ 

\ethm

{\it Proof :}
Under either of the conditions (a) and (b), for any fixed $\phi$, the map $\alpha^\phi_0$ maps $\cla^\infty$ into the subset $\clh^0_D $ of $ \cln$. Since $\alpha^\phi_0$ also commutes with $[D, \cdot]$ on $\cla^\infty$, it is clear that $\alpha^\phi_0 $ maps $\cls^\infty$ into $\cln.$ In fact, 
 using the complete positivity of the map $ \alpha_0^{\phi} $ 
  and the $\alpha_0$-invariance of $\tau$, we see that
 $$ {\tau ( \alpha_0^{\phi} ( a )}^{*} \alpha_0^{\phi} ( a ) ) \leq \tau ( \alpha_0^{\phi} ( a^* a ) ) = ( {\rm id} \otimes \phi ) ( ( \tau \otimes {\rm id} ) \alpha_0 ( a^* a)) = \tau ( a^* a ).1 $$
  which implies that $ \alpha_0^{\phi} $ extends to a bounded operator from $ \cln$ to itself.  Since $U_0$ commutes with $D$, it is clear that $\alpha_0^\phi$ (viewed as a bounded operator on $\cln$) 
 will commute with the group of unitaries $\widetilde{T_t}$, hence with its generator $\tilde{\delta}$ and also with the self adjoint operator $\tilde{\cll}=(i \tilde{\delta})^2$. 

On the other hand, it follows from the  definition of $  \clg=QISO^{+}_{I}( D )$ that $(\tau \ot {\rm id})(\alpha_0(X))=\tau(X) 1_\clg $ for all 
$X $ in $ \clb(\clh)$, in particular for $X $ belonging to $ \cls^\infty$, and thus the map $\cls^\infty \ot_{\rm alg} \clg \ni (a \ot q) \mapsto \alpha_0(a)(1 \ot q)$ 
extends to a $\clg$-linear unitary, denoted by $W$ (say), on the Hilbert $\clg$-module $\cln \ot \clg$. 
Note that here we have used the fact( which that for any $\phi$, $({\rm id} \ot \phi)(W) (\cls^\infty \ot_{\rm alg} \clg) \subseteq \cln$, 
since $\alpha^\phi_0(\cls^\infty) \subseteq \cln$. The commutativity of $\alpha^\phi_0$ with $\widetilde{T_t}$ for every $\phi$ clearly implies that 
$W$ and $\widetilde{T_t} \ot {\rm id}_\clg$ commute on $\cln \ot \clg$. Moreover, $\alpha^\phi_0$ maps $\clh^0_D$ into itself, so $W$ maps 
$\clh^0_D \ot \clg$ into itself, and  hence (by unitarity of $W$) it commutes with the projection $P_0 \ot 1$. 
It follows that $\alpha^\phi_0$ commutes with $P_0$, and (since it also commutes with $\tilde{\cll}$) commutes with $\cll=P_0 \tilde{\cll} P_0$ 
as well.

Thus, $ \alpha^\phi_0 $ preserves each of  the (finite dimensional)  eigenspaces of the Laplacian $\cll$,  and so is a Hopf algebraic action on  the subalgebra $\cla_0$ spanned algebraically by these eigenvectors. 
Moreover, the   $\clg$-linear unitary $W$ clearly restricts to a unitary representation on each of the above eigenspaces. 
If we denote by  $ (( q_{ij} ))_{( i,j )} $ the $\clg$-valued unitary matrix corresponding to one  such particular eigenspace, then by Proposition \ref{preliminaries_CQG_Peter_Weyl}, $q_{ij}$ must belong to $\clg_0$ and we must have $ \epsilon ( q_{ij} ) = \delta_{ij} $ (Kronecker delta).  This  implies $ ( {\rm id} \otimes \epsilon ) \circ \alpha_0  = {\rm id} $ on each of the eigenspaces, hence on the norm-dense subalgebra $\cla_0$ of $\cla$,  completing the proof of the fact that $ \alpha_0 $ extends to a $ C^* $ action on $ \cla. $ \qed

\vspace{4mm}

 Combining the above theorem with Theorem \ref{stable}, we get the following immediate corollary.
\bcrlre
If the spectral triple satisfies the Sobolev condition mentioned before, in addition to the assumptions ${\bf 1}-{\bf 6}$, then $QISO^+_I(D)$ has a $C^*$-action.  In particular, for a classical spectral triple, $QISO^+_I(D)$ has $C^*$-action.
\ecrlre
\brmrk 
Let us remark here that in case the restriction of $\tau$ on $\cla^\infty$ is normal, that is, continuous with respect to the weak operator topology inherited from $\clb(\clh)$, 
then $\clh^0_D$ will contain $\cla^{\prime \prime}$, which is the closure of $\cla^\infty$ in the weak operator topology of $\clb(\clh)$, so the condition (a) of Theorem \ref{QISO_I_D_<_Q_L_D} (and hence its conclusion) holds.  
\ermrk

\brmrk
In a private communication to us, Shuzhou Wang has kindly pointed out  that a possible alternative approach  to the formulation of quantum group of isometries may involve the category of CQG which has a $C^*$-action on the underlying $C^*$ algebra and a unitary representation with respect to which the Dirac operator is equivariant. However, we see from Corollary \ref{sphere_chak_pal_non_existence_wang_univ} of chapter \ref{quantumsphere} that the category proposed by Wang does not admit a universal object in general. 
\ermrk

\subsection{ Universal object in the  categories  ${\bf Q} $ or ${\bf Q}^{\prime}$}

\label{qorient_subsection_univ_in_Q_Q_prime}

We shall now investigate   further conditions  on the spectral triple which will ensure the existence of a  universal object in the category ${\bf Q}$ or  ${\bf Q}^\prime$.
 Whenever such a universal object exists we shall denote it by $\widetilde{QISO^+}(D)$, and denote by $QISO^+(D)$ its  largest Woronowicz subalgebra  for which $\alpha_U$ on $\cla^\infty$  (where $U$ is the unitary representation of $\widetilde{QISO^+}(D)$ on $\clh$) is faithful.
 \brmrk
If ${\widetilde{QISO^+}}(D)$ exists, then by Proposition \ref{5678}, there will exist some $R$ such that ${\widetilde{QISO^+}}(D)$ is an object in the category ${\bf Q}^\prime_R(D)$. Since the universal object in this category, that is, $\widetilde{QISO^+_{R}}(D)$, is clearly a sub-object of $\widetilde{QISO^+}(D)$, we have $\widetilde{QISO^+}(D) \cong \widetilde{QISO^+_{R}}(D)$ for this choice of $R$.
\ermrk

  Let us state and prove a result below, which gives some sufficient conditions for the existence of $\widetilde{QISO^+}(D)$.
\bthm

\label{unrestricted}

Let $(\cla^\infty, \clh, D)$ be a spectral triple of compact type as before and assume that $D$ has an one-dimensional eigenspace spanned by a unit vector $\xi$, which is cyclic and separating for the algebra $\cla^\infty$. Moreover, assume that each eigenvector of $D$ belongs to the dense subspace $\cla^\infty \xi$ of $\clh$. Then there is a universal object, 
$(\widetilde{\clg}, U_0)$. Moreover, $\widetilde{\clg}$ has a coproduct $\Delta_0$ such that $(\widetilde{\clg},\Delta_0)$ is a compact quantum group and $(\widetilde{\clg},\Delta_0,U_0)$ is a universal object in the category ${\bf Q}^\prime.$ 
 
 If we denote by $ \clg $ the Woronowicz $ C^{\ast} $ subalgebra of $ \widetilde{\clg} $ generated by elements of the form $ \left\langle \alpha_{U_{0}} ( a ) ( \eta \otimes 1 ) , ~ \eta^{\prime} \otimes 1 \right\rangle_{\widetilde{\clg}} $ where $ \eta, \eta^{\prime} $ are in $ \clh,~ a $ in $ \cla^{\infty} $ and $ \left\langle . ~,~ . \right\rangle_{\widetilde{\clg}} $ denotes the $ \widetilde{\clg} $ valued inner product of $ \clh \otimes \widetilde{\clg} ,$ we have $ \widetilde{\clg} \cong \clg \ast C ( \IT ).$
\ethm
{\it Proof:} Let $V_i$, $\{ e_{ij} \}$ be as before, and by assumption we have $e_{ij}=x_{ij} \xi$ for a unique $x_{ij} $ in $ \cla^\infty$. Clearly, since $\xi$ is separating, the vectors $\{ \overline{e_{ij}}=x^*_{ij} \xi, j=1,..., d_i \}$ are linearly independent, so the matrix $Q_i=(( \lgl \overline{e_{ij}}, \overline{e_{ik}} \rgl ))_{j,k=1}^{d_i}$ is positive and invertible. Now, given a quantum family of orientation-preserving  isometries $(\cls, U)$, we must have $\widetilde{U} ( \xi \ot 1)= \xi \ot q$, say, for some $q$ in $\cls$, and from the unitarity of $\widetilde{U}$ it follows that $q$ is a unitary element. Moreover, $U$ leaves $V_i$ invariant, so let $\widetilde{U}(e_{ij} \ot 1)=\sum_k e_{ik} \ot v^{(i)}_{kj}$. But this can be rewritten as $$ \alpha_U(x_{ij}) (\xi \ot q)=\sum_k x_{ik} \xi \ot v^{(i)}_{kj}.$$ Since $\xi$ is separating and $q$ is unitary, this implies $ \alpha_U(x_{ij})=\sum_k x_{ik} \ot v^{(i)}_{kj}q^*,$  and thus we have $$ \widetilde{U}(\overline{e_{ij}}\ot 1)=\alpha_U(x_{ij})^*(\xi \ot q)=\sum_k x^*_{ik} \xi \ot q (v^{(i)}_{kj})^* q =\sum_k \overline{e_{ik}} \ot q (v^{(i)}_{kj})^* q.$$ 
  Taking the $\cls$-valued inner product $\lgl \cdot, \cdot \rgl_\cls$ on both sides of the above expression,   and using the fact that $U$ preserves this $\cls$-valued inner product,  we obtain  
  \bean \lefteqn{ \left\langle \overline{e_{ij}} ~ , ~ \overline{e_{i j^{\prime} }} \right\rangle.1}\\
  &=& \left\langle \widetilde{U} ( \overline{e_{ij}} \otimes 1 ) ~, ~ \widetilde{U} ( \overline{e_{i j^{\prime}}} \otimes 1 ) \right\rangle\\
  &=& \left\langle  \sum_k \overline{e_{ik}} \otimes q {(v^{(i)}_{kj})}^* q, ~ \sum_{k^{\prime}} \overline{e_{ik^{\prime}}} \otimes q {(v^{(i)}_{k^{\prime}j^{\prime}})}^* q   \right\rangle\\
  &=& \sum_{k,k^{\prime}} \left\langle \overline{e_{ik}} ~,~ \overline{e_{ik^{\prime}}}  \right\rangle q^* v^{(i)}_{kj} q^* q {(v^{(i)}_{k^{\prime}j^{\prime}})}^* q. \eean
  
  This implies, 
    \bean \lefteqn{q {(Q_i)}_{j j^{\prime}} q^* .1 }\\
                &=& \sum_{k,k^{\prime}} v^{(i)}_{kj} \left\langle  \overline{e_{ik}} ~, ~ \overline{e_{ik^{\prime}}} \right\rangle {(v^{(i)}_{k^{\prime}j^{\prime}})}^* .\eean
                
   Thus, $$ {( Q_i )}_{j j^{\prime}} = \sum_{k, k^{\prime}} v^{(i)}_{kj} {( Q_i )}_{k k^{\prime}} {(v^{(i)}_{k^{\prime}j^{\prime}})}^*.$$ 
    Hence, we have  $Q_i=v_i^\prime Q_i \overline{v_i}$ (where $v_i=(( v^{(i)}_{kj}))$). Thus, $Q_i^{-1}v_i^\prime Q_i$ must be the (both-sided) inverse of $\overline{v_i}$. Thus, we get a canonical surjective morphism from $ A_{u, d_i}(Q_i) $ to the $C^{\ast}$ algebra generated by $\{ v^{(i)}_{kj}: j,k = 1,2,... d_i \}$. This induces a surjective morphism from the free product of $A_{u,d_i}(Q_i)$, $i=1,2,...$ onto $\cls$. The rest of the arguments for showing the existence of $\widetilde{\clg}$  will be quite similar to the arguments used in the proof of Theorem \ref{main}, hence omitted.
 
 Now we come to the proof of the last part of the theorem. For $ a $ in $ \cla^{\infty}, ~ \widetilde{U} ( a \xi \otimes 1 ) = \alpha_{U} ( a ) \widetilde{U} ( \xi \otimes 1 ) = \alpha_{U} ( a ) ( \xi \otimes q ).  $ Now, recalling that $ {\rm Span} \{ a \xi : a \in \cla^{\infty} \} $ is dense in $ \clh, $ it is clear that $\widetilde{\clg} =\clg \ast C^*(q) \cong \clg \ast C(\IT)$.
\qed

\brmrk

Some of the examples considered in section \ref{qorient_section_examples_computations} will show that the conditions of the above theorem are not actually necessary; $\widetilde{QISO}^+(D)$ may exist without the existence of a single cyclic separating eigenvector as above. 
 
\ermrk

Let $ ( {\cla}^{\infty}, \clh, D ) $ be a spectral triple of compact type satisfying the conditions of the above theorem. Let the faithful vector state corresponding to the cyclic separating vector $ \xi $ be denoted by $ \tau.$
Let $ \cla_{00} = {\rm span} \{ a \in \cla^{\infty} : a \xi $ is an eigenvector of $ D  \} .$

Moreover, assume that $ \cla_{00} $ is norm dense in $ \cla^{\infty} .$

Let $ \hat{D} : \cla_{00} \rightarrow  \cla_{00} $ be defined by :

$ \hat{D} ( a ) \xi = D ( a \xi ). $

This is well defined since $ \xi $ is cyclic and separating.

\bdfn
Let $ \cla $ be a $ C^{\ast} $ algebra and $ {\cla}^{\infty} $ be a dense $ \ast$-subalgebra.
Let $ ( {\cla}^{\infty} , \clh, D ) $ be a spectral triple of compact type as above.

Let $ \widehat{\bf C} $ be the category with objects $ ( \clq, \alpha ) $ such that $ \clq $ is a compact quantum group with a $ C^* $ action $\alpha$  on $ \cla $ such that:

1. $ \alpha $ is $ \tau $ preserving ( where $ \tau $ is as above ), that is, $ ( \tau \otimes {\rm  id} ) ( \alpha ( a ) ) = \tau ( a ).1 .$

2. $ \alpha $ maps $ \cla_{00} $ inside $  \cla_{00} \otimes_{alg} \clq .$

3. $ \alpha \widehat{D} = ( \widehat{D} \otimes I ) \alpha  .$

\edfn

\bcrlre

\label{unrestrictedcorollary}
 
 There exists a universal object $ \widehat{\clq} $ of the category $ \widehat{\clc} $ and it is isomorphic to the Woronowicz $ {C}^{*} $ subalgebra $\clg={QISO}^+(D)$ of $ \widetilde{\clg} $  obtained in Theorem \ref{unrestricted}.

\ecrlre

{\it Proof :} The proof of the  existence of the universal object follows verbatim from the proof of Theorem \ref{qisol_main_existence_theorem} replacing $ \cll $ by $ \widehat{D} $ and noting that $ D $ has compact resolvent. We denote by $\widehat{\alpha} $ the action of $ \widehat{\clq} $ on $ \cla.$ 

Now, we prove that $ \widehat{\clq} $ is isomorphic to  $ \clg .$

 Each eigenvector of $ D$ is in  $  \cla^{\infty} $ by assumption. It is easily observed from the proof of Theorem \ref{unrestricted} that  $ \alpha_{U_{0}}$ maps the norm-dense $\ast$-subalgebra $  \cla_{00} $ into $ \cla_{00} \otimes_{alg} \clg_{0} ,$
  and $({\rm id} \ot \epsilon)\circ \alpha_{U_0}={\rm id}$, so that $\alpha_{U_0}$ is indeed a  $ C^* $ action of the CQG $\clg$.   Moreover, it can be easily seen that $ \tau $ preserves  $  \alpha_{U_{0}} $
   and that $ \alpha_{U_{0}} $ commutes with $ \widehat{D}.$
    Therefore, $ ( \clg, \alpha_{U_{0}} ) $ is an element of $ {\rm Obj} ( \widehat{\bf C} ) $  and  hence  $ \clg $ is a quantum subgroup of $ \widehat{\clq} $ by the universality of $\widehat{\clq}$.

For the converse, we start by showing that $ \widehat{\alpha} $ induces a unitary representation $W$ of $ \widehat{\clq} \ast C( \IT ) $ on $ \clh $ which commutes with $ D$, and the corresponding conjugated action $\alpha_W$ coincides with $\widehat{\alpha}$. 

Define $ W( a \xi ) = \widehat{\alpha}( a ) ( \xi ) ( 1 \otimes q^{*} ) $  for all $ a $ in $ \cla^{\infty}_{0} $ where $ q  $ is a generator of $ C( \IT ).$

Since we have $ ( \tau \otimes {\rm  id} ) ( \alpha ( a ) ) = \tau ( a ).1,$ it follows that $ \widetilde{W}  $ is a ($\widehat{\clq} \ast C(\IT)$-linear)  isometry on the dense subspace $ \cla_{00} \xi \ot_{\rm alg} \widehat{\clq} $ and thus extends to $ \clh \ot ( \widehat{\clq}\ast C(\IT) )$ as an isometry. Moreover, since $ \widehat{\alpha}( \cla ) ( 1 \otimes \widehat{\clq} )$ is norm dense in $ \cla \otimes \widehat{\clq}$ (by the definition of a CQG action)  it is clear that the range of $\widetilde{W}$ is dense, so $\widetilde{W}$ is indeed a unitary.  It is quite obvious that it is a unitary representation of $\widehat{\clq} \ast C(\IT)$.

We also have, \bean 
\lefteqn{W D ( a \xi ) }\\
& =&  W( \widehat{D} ( a ) \xi ) =
   \widehat{\alpha} ( \widehat{D} ( a ) ) ( \xi ) ( 1 \otimes q^{*} ) \\
& =&  ( D \otimes I ) ( \widehat{\alpha} ( a ) \xi ) ( 1 \otimes q^{*} ) 
 = ( D \otimes I ) W ( a \xi ) , \eean

that is,   $ W $ commutes with $ D .$

Moreover, it is easy to observe that $ \alpha_{W} = \widehat{\alpha}.$
 This gives a surjective CQG morphism from $ \widetilde{\clg} = \clg \ast C( \IT ) $ to $ \widehat{\clq} \ast C( \IT ),$ sending $ \clg $ onto $ \widehat{\clq},$ which completes the proof. \qed

\section{Comparison with the approach of \cite{goswami} based on Laplacian}

\label{qorient_section_comparison_with_QISO_L}

Throughout this section, we shall assume the set-up of  subsection \ref{qorient_subsection_stability_top_action}  for the existence of a `Laplacian',  including assumptions ${\bf 1}-{\bf 6}$. Let us also use the notation of that subsection.

 We recall from chapter \ref{qisol} that a CQG $(\cls,\Delta)$ which has an action $\alpha$ on $\cla$ is said to act smoothly and isometrically on the noncommutative manifold ( of compact type ) $(\cla^\infty, \clh, D)$ if $({\rm id} \ot \phi) \circ \alpha(\cla^\infty_0) \subseteq \cla^\infty_0$ for every state $\phi$ on $\cls$, and also $({\rm id} \ot \phi) \circ \alpha$ commutes with the Laplacian $\cll \equiv \cll_D$ on $\cla^\infty_{0}$ ( where $ \cla^{\infty}_0 $ is the complex linear span of the eigenvectors of $ \cll $ ). One can consider the category $ {\bf Q^{\prime}_{\cll_D}} $ of all compact quantum groups acting smoothly and isometrically on $\cla$, where the morphisms are quantum group morphisms which intertwin the actions on $\cla$. 
We make the following additional assumption throughout the present section:
\vspace{2mm}

 ({\bf 7}) There exists a universal object in $ {\bf Q^{\prime}_{\cll_D}}  $ ( the quantum isometry group for the Laplacian $\cll \equiv \cll_D$ in the sense of \cite{goswami}), and it is denoted by  $ {QISO}^\cll \equiv {QISO}^{\cll_{D}} $ 
 
 \vspace{2mm}

The following result now follows immediately from Theorem \ref{QISO_I_D_<_Q_L_D} of subsection 2.3. 

\bcrlre

\label{QISO_I_D_<_Q_L}

If $ ( \cla^{\infty}, \clh, D ) $ is a spectral triple ( of compact type   ) satisfying any of the two conditions ( a ) or ( b ) of Theorem \ref{QISO_I_D_<_Q_L_D}, then $ QISO^+_{I}( D ) $ is a sub-object of $ {QISO}^{\cll_{D}}  $ in the category $ {\bf Q^{\prime}_{\cll_D}}.$

\ecrlre
{\it Proof:} The proof is a consequence of the fact that $QISO^+_I(D)$ has the $C^*$-action $\alpha_0$ on $\cla$, and the observation already made in the proof of  the Theorem \ref{QISO_I_D_<_Q_L_D} that this action commutes with the Laplacian $\cll_D$. 
\qed 

\vspace{4mm}

Now, we will need the Hilbert space of forms $ \clh_{d + d^*} $ corresponding to a $\Theta$-summable spectral triple $( \cla^{\infty}, \clh, D )$ as discussed in subsection \ref{preliminaries_subsection_spaceofforms_NCG}. We recall that one obtains an associated spectral triple $( \cla^{\infty}, \clh_{d + d^*}, d + d^* ) .$ We assume that this spectral triple is of compact type, that is, $d+d^*$ has compact resolvents.

We will denote the inner product on the space of $ k $ forms coming from the spectral triples $ ( \cla^{\infty}, ~ \clh, ~ D ) $ and   
$( \cla^{\infty}, ~ \clh_{d + d^*},~ d + d^* ) $ by $ {\left\langle ~ , ~ \right\rangle}_{\clh^{k}_{D}} $ and 
$ {\left\langle ~ , ~ \right\rangle}_{\clh^{k}_{d + d^*}} $ respectively, $ k = 0,1.$ 

We will denote by $ \pi_{D} , ~ \pi_{d + d^*} $ the representations of $ \cla^{\infty} $ in $ \clh $ and $ \clh_{d + d^*} $ respectively.

Let $ U_{d + d^*} $ be the canonical unitary representation of 
$ QISO^{+}_{I}( d + d^* )$ on $ \clh_{d + d^*}.$ 

$ \clh_{d + d^*} $ breaks up into finite dimensional orthogonal subspaces corresponding to the distinct eigenvalues of $ \Delta:=(d + d^*)^2=d^* d + d d^*.$ 
It is easy to see that $\Delta$ leaves each of the subspaces $\clh^i_D$ invariant, and we will denote by $ V_{\lambda,i} $ the subspace of $ \clh^{i}_{d+d^*} $ spanned by eigenvectors of $ \Delta $ corresponding to the eigenvalue 
$ \lambda.$ Let $ \{ e_{j,\lambda,i} \}_{j} $ be an orthonormal basis of $ V_{\lambda,i}.$ Note that $\cll_D$ is the restriction of $\Delta$ to $\clh^0_D$.

Now we recall Proposition \ref{qisol_unitary_rep_of_qisol}. It was shown there that $ {QISO}^{\cll_{D}} $ has a unitary representation $ U \equiv U_\cll $ on  $ \clh^{d + d^*} $ such that $ U $ commutes with $ d + d^* .$   Thus, $ ( \cla^{\infty}, ~  \clh_{d + d^*}, ~ d + d^* ) $ is a $ {QISO}^{\cll_{D}} $ 
equivariant spectral triple. Moreover, by Remark \ref{qiso_laplacian_qiao_tracial_haar_state}, $ {QISO}^{\cll_{D}} $   has tracial Haar state, which implies, by Proposition \ref{5678} and Remark \ref{qorient_haarstate_tracial_implies_R=I} that $ \alpha_U $ keeps the functional $ \tau_I $ invariant. Summarizing, we have the following result:

\bppsn

\label{Q_L_leq_QISO(d + d*)}

The quantum isometry group $ ( {QISO}^{\cll_{D}}, ~ U_\cll ) $ is a sub-object of $ ( QISO^{+}_{I} ( d + d^* ), ~ U_{d + d^*} ) $ in the category $ {\bf Q}_{I} ( d + d^* ),$ so in particular,
  $ {QISO}^{\cll_{D}} $ is  isomorphic to a quotient of $ QISO^{+}_{I}( d + d^* ) $ by a Woronowicz $ C^* $ ideal.

 \eppsn

We shall  give (under mild conditions) a concrete description of the above Woronowicz ideal. 

Let $ \cli $ be the $ C^* $ ideal of $QISO^{+}_{I}( d + d^* )$ generated by

 $$ \cup_{\lambda \in \sigma( \Delta )} \{ \left\langle ( P^{\bot}_{0} \otimes {\rm  id} ) U_{d + d^*} ( e_{j \lambda 0} ), e_{j \lambda i^{\prime}} 
\otimes 1 \right\rangle : j, i^{\prime} \geq 1 \},$$
where $ P_0 $ is the projection onto $ \clh^{0}_D, ~  \left\langle . ~,~ . \right\rangle $ denotes the $ QISO^{+}_{I}( d + d^* ) $ 
valued inner product and $ \sigma ( \Delta ) $ denotes the spectrum of $ \Delta. $

Since $ U_{d + d^*} $ keeps the eigenspaces of $ \Delta=(d + d^*)^2 $ invariant, we can write 
$$ U_{d + d^*} ( e_{j \lambda 0} ) = \sum_{k} e_{k \lambda 0} \otimes ~ q_{k j \lambda 0} + \sum_{i^{\prime}\neq 0, k^{\prime} } e_{k^{\prime} \lambda i^{\prime}} 
\otimes ~ q_{k^{\prime} j \lambda i^{\prime}}, $$ for some $ q_{k j \lambda 0}, ~ q_{k^{\prime} j \lambda i^{\prime}} $ in $ QISO^{+}_{I}( d + d^* ).$

We note that $ q_{k^{\prime} j \lambda i^{\prime}} $ is in $ \cli $ if $ i^{\prime} \neq 0.$

\blmma

$ \cli $ is a co-ideal of $ QISO^{+}_{I} ( d + d^* ).$

\elmma

{\it Proof :} It is enough to prove the relation $ \Delta ( X ) \in   \cli \otimes QISO^{+}_{I}( d + d^* ) + QISO^{+}_{I}( d + d^* ) \otimes \cli $ for the elements $ X $ in $ \cli $ of the form $ \left\langle ( P^{\bot}_{0} \otimes {\rm  id} ) U_{d + d^*} ( e_{j \lambda  0} ), e_{j \lambda i_0} \otimes 1 \right\rangle .$
We have:
\bean \Delta ( \left\langle ( P^{\bot}_{0} \otimes {\rm  id} ) U_{d + d^*} ( e_{m \lambda 0} ), e_{j \lambda i_0} \otimes 1 \right\rangle ) \eean
 \bean = \left\langle ( P^{\bot}_{0} \otimes {\rm  id} ) ( {\rm  id} \otimes \Delta ) U_{d + d^*} ( e_{m \lambda 0} ), e_{j \lambda i_0} \otimes 1 \ot 1 \right\rangle ) \eean 
 \bean = \left\langle ( P^{\bot}_{0} \otimes {\rm  id} ) U_{(12)} U_{(13)} ( e_{m \lambda 0} ), e_{j \lambda i_0} \otimes 1 \ot 1 \right\rangle  \eean 
 \bean = \left\langle ( P^{\bot}_{0} \otimes {\rm  id} ) U_{(12)} ( \sum_{k} e_{k \lambda 0} \otimes 1 \otimes q_{k m \lambda 0} ) ~~, ~~ e_{j \lambda i_{0}} \otimes 1 \ot 1 \right\rangle  \eean  
\bean +  \sum_{i^{\prime} \neq 0,l} \left\langle ( P^{\bot}_{0} \otimes {\rm  id} ) U_{(12)} ( e_{l \lambda i^{\prime}} \otimes 1 \otimes q_{l m \lambda i^{\prime}} ) ~~, ~~ e_{j \lambda i_{0}} \otimes 1 \ot 1 \right\rangle \eean 
 \bean = \sum_{k, k^{\prime}}  \left\langle ( P^{\bot}_{0} \otimes {\rm  id} ) ( e_{k^{\prime} \lambda 0} \otimes q_{k^{\prime} k \lambda 0} \otimes q_{k  m \lambda 0} ) ~~, ~~ e_{j \lambda i_0} \otimes 1 \ot 1 \right\rangle \eean  
\bean  + \sum_{i^{\prime }  \neq 0,~ k,~ k^{\prime \prime}} \left\langle ( P^{\bot}_{0} \otimes {\rm  id} ) (   e_{k^{\prime \prime} k \lambda i^{\prime}} \otimes q_{k^{\prime \prime},k,\lambda,i^{\prime}} \otimes q_{k m \lambda 0} ) ~~ , ~~ e_{j \lambda i_0} \otimes 1 \ot 1 \right\rangle   \eean 
 \bean + \sum_{i^{\prime} \neq 0,~ l,~l^{\prime}} \left\langle ( P^{\bot}_{0} \otimes {\rm  id} ) (    e_{l^{\prime} \lambda i^{\prime}} \otimes q_{l^{\prime} l \lambda i^{\prime}} \otimes q_{l m \lambda i^{\prime}} ) ~~ , ~~ e_{j \lambda i_0} \otimes 1 \ot 1 \right\rangle \eean 
 \bean + \sum_{i^{\prime} ~ \neq 0,~ i^{\prime \prime} \neq i^{\prime},~ l,~ l^{\prime \prime}} \left\langle ( P^{\bot}_{0} \otimes {\rm  id} ) (   e_{l^{\prime \prime} \lambda i^{\prime \prime}} \otimes q_{l^{\prime \prime} l \lambda i^{\prime \prime}} \otimes q_{l m \lambda i^{\prime}} ) ~~ , ~~ e_{j \lambda i_0} \otimes 1 \ot 1 \right\rangle \eean 
 \bean = \sum_{i^{\prime} \neq 0,~ k^{\prime},~ k^{\prime \prime}} \left\langle   e_{k^{\prime \prime} \lambda i^{\prime}} \otimes q_{k^{\prime \prime} k \lambda i^{\prime}} \otimes q_{k m \lambda 0} ~~,~~ e_{j \lambda i_0} \otimes 1  \ot 1 \right\rangle   \eean  
 \bean + \sum_{i^{\prime} \neq 0,~ l,~ l^{\prime}} \left\langle  e_{l^{\prime} \lambda i^{\prime}} \otimes q_{l^{\prime} l \lambda i^{\prime}} \otimes q_{l m \lambda i^{\prime}}  ~~,~~ e_{j \lambda i_0} \otimes 1  \ot 1 \right\rangle   \eean 
 \bean + \sum_{i^{\prime} \neq 0,~ i^{\prime \prime} \neq i^{\prime},~ i^{\prime \prime} \neq 0, ~ l,~ l^{\prime \prime} }  \left\langle e_{l^{\prime \prime} \lambda i^{\prime \prime}} \otimes q_{l^{\prime \prime} l \lambda i^{\prime \prime}} \otimes q_{l m \lambda i^{\prime}} ~~,~~  e_{j \lambda i_0} \otimes 1  \ot 1 \right\rangle, \eean 
which is clearly in $ \cli \otimes QISO^{+}_{I}( d + d^* ) + QISO^{+}_{I}( d + d^* ) \otimes \cli ,$ as $ q_{k j \lambda i^{\prime}} $ is an element of $ \cli $ for $ i^{\prime} \neq 0  .$ \qed

%
%
%
%
%
%
%
%

\bthm

If $ \alpha_{U_{d + d^{*}}} $ is a $ C^* $ action on $ \cla ,$ then we have $ {QISO}^{\cll_{D}} \cong QISO^{+}_{I} ( d + d^* ) / \cli .$
\ethm

{\it Proof :}  By   Proposition \ref{Q_L_leq_QISO(d + d*)}, 
 we conclude that there exists a surjective CQG morphism $ \pi : QISO^{+}_{I} ( d + d^* ) \rightarrow QISO^{\cll_{D}}.$ 
By construction ( as in Proposition \ref{qisol_unitary_rep_of_qisol} ), the unitary representation  $ U_\cll $ of $ QISO^{\cll_{D}} $ preserves each of the $ \clh^{i}_{D},$ in particular
 $\clh^0_D$. It is then clear from the definition of $\cli$ that  
 $ \pi $ induces a surjective CQG morphism (in fact, a morphism in the category ${\bf Q}_I^\prime(d+d^*)$) $ \pi^{\prime} : QISO^{+}_{I} ( d + d^* ) / \cli \rightarrow QISO^{\cll_{D}}.$

Conversely,  if $ V=({\rm id} \ot \rho_\cli)\circ U_{d+d^*}$ is the representation of $QISO^{+}_{I} ( d + d^* )/\cli $ on $ \clh_{d + d^*} $ 
induced by $U_{d+d^*}$ (where  $\rho_\cli: QISO^+_I(d+d^*)\raro QISO^+_I(d+d^*)/\cli$ denotes the quotient map), then $V$ preserves $\clh^0_D$ (by definition of 
$\cli$), so commutes with $P_0$. Since $V$ also commutes with $(d+d^*)^2$, it follows that $V$ must commute with $(d+d^*)^2P_0=\cll$, that is,
 $$ \widetilde{V} ( d^*d ~ P_0 \otimes 1 ) = ( d^*d P_0 \otimes 1 ) \widetilde{V}. $$  It 
 is easy to show from the above that $\alpha_V$ (which is a $C^*$ action 
   on $\cla$ since $\alpha_{U_{d+d^*}}$ is so by assumption) is a smooth isometric action of $QISO^+_I(d+d^*)/\cli$ in the sense of \cite{goswami}, with respect to
 the Laplacian $\cll$. This implies that $QISO^+_I(d+d^*)/\cli$ is a sub-object of ${QISO}^{\cll_{D}}$ in the category ${\bf Q^{\prime}_{\cll_D}} $, and completes the proof.
\qed
\vspace{2mm} 

%

Now we prove that under some further assumptions which are valid for classical manifolds as well as their Rieffel deformation,
  one even has the isomorphism $ {QISO}^{\cll_{D}} \cong QISO^{+}_{I}( d + d^* ).$

We assume the following:

\vspace{2mm}

({\bf A})  Both the spectral triples $ ( \cla^{\infty}, ~ \clh, ~ D ) $ and $ ( \cla^{\infty}, ~ \clh_{d + d^{*}}, ~ d + d^{*} ) $ satisfy the assumptions $({\bf 1})-({\bf 7})$, so in particular
 both  ${QISO}^{\cll_D}$ and ${QISO}^{\cll_{D^\prime}}$ exist (here $D^\prime=d+d^*$).
 
 \vspace{1mm}
 
({\bf B}) For all $a, b$  in $\cla^{\infty},$ we have $$ {\left\langle a,~ b \right\rangle}_{\clh^{0}_{ D }} = {\left\langle a,~ b \right\rangle}_{\clh^{0}_{ D^{\prime} }} ~ , ~  {\left\langle d_{D}a,~ d_{D}b \right\rangle}_{\clh^{1}_{ D }} = {\left\langle d_{D^{\prime}}a,~ d_{D^{\prime}}b 
\right\rangle}_{\clh^{1}_{ D^{\prime} }} .$$

\brmrk

For classical compact spin  manifolds these assumptions can be verified by comparing the local expressions of $D^2$ and the `Hodge Laplacian' $(D^\prime)^2$ in suitable coordinate charts. In fact, in this case, both these operators turn out to be essentially same, upto a `first order term', which is relatively compact with respect to $D^2$  or $(D^\prime)^2$.


\ermrk

 By assumption $({\bf  B})$, we observe that the identity map on $ \cla^{\infty} $ extends to a unitary, say $ \Sigma $, from $ \clh^{0}_{D} $ to 
$ \clh^{0}_{D^{\prime}}.$ Moreover,  we have  $$ \cll_{D} = \Sigma^{*} \cll_{D^{\prime}} \Sigma,$$ from which we conclude the following:

\bppsn

\label{Q_L_D_isom_Q_L_(d + d*)}

Under the above assumptions, $ {QISO}^{\cll_{D}} \cong {QISO}^{\cll_{D^{\prime}}} .$

\eppsn


%
%
We conclude this section with the following result, which identifies the quantum isometry group ${QISO}^{\cll_D}$ of \cite{goswami} as the $QISO^+_I$ of a spectral triple, and thus, in some sense, accommodates the construction of \cite{goswami} in the framework of the present article. 
\bthm

\label{QISO( D ) vs  QISO( L )}

If in addition to the assumptions already made, the spectral triple ( of compact type ) $(\cla^\infty, \clh_{D^\prime}, D^\prime)$ also satisfies the conditions of  Theorem \ref{QISO_I_D_<_Q_L_D}, so that $QISO^+_I(D^\prime)$ has a $C^*$-action, then we have the following isomorphism of CQG s: $$ {QISO}^{\cll_{D}} \cong QISO^{+}_{I}( D^{\prime} ) \cong {QISO}^{\cll_{D^{\prime}}}.$$

\ethm

{\it Proof :} 
By  Proposition \ref{Q_L_leq_QISO(d + d*)} we have that $ {QISO}^{\cll_{D}} $ is a sub-object of $ QISO^{+}_{I}( D^{\prime} ) $ 
in the category ${\bf Q}^\prime_{I}( D^{\prime} ) $. 
On the other hand, by Theorem \ref{QISO_I_D_<_Q_L_D} we have $ QISO^+_{I}( D^{\prime} ) $ as a sub-object of $ {QISO}^{\cll_{D^{\prime}}} $ in the category $ {\bf Q^{\prime}_{\cll_{D^{\prime}}}} .$
Combining these facts with the conclusion of Proposition \ref{Q_L_D_isom_Q_L_(d + d*)}, we get the required isomorphism.\qed




 

\brmrk

The assumptions, and hence the conclusions, of this section are valid also for spectral triples obtained by Rieffel deformation of a classical spectral triple, to be discussed in details in chapter \ref{deformation}.

\ermrk

\section{Examples and computations}

\label{qorient_section_examples_computations}

 In this section we compute the quantum group of orientation preserving isometries for spectral triples on $ SU_{\mu} ( 2 ) $ and $ C ( \IT^2 ).$ The computations for the Podles' spheres and Rieffel deformed manifolds are given in chapter \ref{quantumsphere} and Chapter \ref{deformation} respectively. 

\subsection{ Equivariant spectral triple on $SU_\mu(2)$}

 We recall from subsection \ref{preliminaries_subsection_SU_mu_2} that by $ t^n_{i,j} $ s, we will denote the $ ( i,j ) $ th matrix element of the $ ( 2n + 1 ) $ dimensional representation of $ SU_{\mu} ( 2 ) $ and $ e^n_{i,j} $ s will denote the normalized ( with respect to the Haar state $h$ ) $ t^n_{i,j} $s. We consider the  spectral triple on $ SU_{\mu}( 2 ) $ constructed by Chakraborty and Pal ( \cite{partha} ) and also discussed thoroughly in \cite{con2} which is defined by $ ( \cla^{\infty}, \clh, D ) $ where $\cla^\infty$ is the linear span of $t^n_{ij}$ s, $ \clh = L^{2} ( SU_{\mu} ( 2 ) ) $  and $ D $ is defined by : 
  \bean \lefteqn{ D ( e^{n}_{ij} ) }\\
     &=& ( 2n + 1 )e^{n}_{ij},~   n \neq i \\
     &=&  -( 2n + 1 ) e^{n}_{ij},~ n = i. \eean
      
 Here, we have a cyclic separating vector $ 1_{SU_{\mu}( 2 )}$, and the corresponding faithful state is  the Haar state $h$. Thus, we are in the set up of the subsection \ref{qorient_subsection_univ_in_Q_Q_prime}, and as  $ \xi = 1, ~ \cla_{00} = \cla^{\infty} $ in this case. Therefore, an operator commuting with $ D $ ( equivalently with $ \widehat{D} $ ) must keep $ V^{l}_{i}: = {\rm Span} \{ t^{l}_{ij} : j = -l,........l \} $ invariant for all  fixed $ l $ and $ i $ where $ \widehat{D} $ is the operator as in subsection \ref{qorient_subsection_univ_in_Q_Q_prime}.

 In the notation of Corollary \ref{unrestrictedcorollary}
, we have $ \cla_{00} = {\rm Span} \{ t^{l}_{i,j} : l = 0, 1/2,......... \}.=\cla^\infty$ in this case. 
All the conditions of Theorem \ref{unrestricted} and Corollary \ref{unrestrictedcorollary} are satisfied. Thus,  the universal object of the category $ \widehat{\bf C} $ exists ( notation as in Corollary \ref{unrestrictedcorollary} ) and we denote it by $ \widehat{\clq}.$

Before proving the next result, we note the following fact. We recall the fundamental unitary of $ SU_{\mu} ( 2 ) $ given by $ \left ( \begin {array} {cccc}
   \alpha & - \mu \gamma^* \\ \gamma & \alpha^*  \end {array} \right ) $ which is the matrix corresponding to the coproduct $ \Delta $ on $ {\rm span} ~ \{ \alpha. - \mu {\gamma}^* \} $ as given in subsection \ref{preliminaries_subsection_SU_mu_2}. This implies that $ V^{\frac{1}{2}}_{- \frac{1}{2}} = {\rm span} ~ \{ \alpha, \gamma^* \}  $ and $ V^{\frac{1}{2}}_{ \frac{1}{2}} = {\rm span} ~ \{ \alpha^*, \gamma \}.$

\blmma

\label{equivalent}

Given a CQG $\clq$ with a $ C^* $ action $\Phi$ on $\cla$, the following are equivalent :

1.$ ( \clq, \Phi ) $ is an element of $ Obj ( \widehat{{\bf C}} ) .$

2. The action is linear, in the sense that $ V^{1/2}_{-1/2} $ ( equivalently $,V^{1/2}_{1/2} $ ) is invariant under $ \Phi $ and the  representation obtained by restricting  $\Phi$  to $ V^{1/2}_{1/2} $ is a unitary representation.

3. $ \Phi $ is linear and Haar state preserving.

4. $ \Phi $ keeps $ V^{l}_{i} $ invariant for each fixed $ l $ and $ i .$

\elmma

{\it Proof :}  $ 1. \Rightarrow 2. $ Since  $ \Phi $ commutes with $ \widehat{D} $, $ \Phi $ keeps each of the eigenspaces of $ \widehat{D} $ invariant and so in particular  preserves $ V^{1/2}_{-1/2} $ , that is $ \Phi $ is linear. The condition $ ( h \otimes {\rm  id} ) \Phi = h( \cdot).1 $ implies the unitarity of the corresponding representation.

\vspace{4mm}

$ 2 \Rightarrow 3. $ By linearity, write  $ \Phi ( \alpha ) = \alpha \otimes X + {\gamma}^{*} \otimes Y $ and $ \Phi ( {\gamma}^{*} ) = \alpha \otimes Z + {\gamma}^{*} \otimes W.$

Firstly, $ \Phi $-invariance of ${\rm Span} \{t^{k}_{i,j}\} $ for $ k = 0 $ and $k= \frac{1}{2} $ follow from  the linearity and the fact that $ \Phi( 1 ) = 1.$ 

Next, we show that $ \Phi $ keeps $ {\rm Span} \{ t^{1}_{ij} : i,j = -1,0,1 \} $ invariant.

We recall the explicit form of the matrix (($ t^{1}_{ij} $)) from \cite{podles}: 

$ \left ( \begin {array} {cccc}
   {{\alpha}^{*}}^{2} & - ( {\mu}^{2} + 1) \alpha^{*} \gamma & - \mu {\gamma}^{2}  \\ {\gamma}^{*} {\alpha}^{*} &  1 - ( {\mu}^{2} + 1) {\gamma}^{*} \gamma & \alpha \gamma \\ - \mu {{\gamma}^{*}}^{2} & - ( {\mu}^{2} + 1) {\gamma}^{*} \alpha & {\alpha}^{2}  \end {array} \right ) .$ 

By inspection, we see that $ \Phi( V^{1}_{i} ) \subseteq V^{1}_{i} \otimes \clq $ for $ i = -1, 1.$

Hence, it is enough to check the $\Phi$-invariance   for $ \alpha \gamma $ and $ 1 - ( {\mu}^{2} + 1 ) {\gamma}^{*} \gamma .$

We have \bean \lefteqn{\Phi ( \alpha \gamma )}\\
        &=& ( \alpha \otimes X + \gamma^* \otimes Y ) ( \alpha^* \otimes Z^* + \gamma \otimes W^* )\\
        &=& \alpha \alpha^* \otimes X Z^* + \gamma^* \gamma \otimes Y W^* + \alpha \gamma \otimes X W^* + \gamma^* \alpha^* \otimes Y Z^*\eean
        \bean &=& \alpha \gamma \otimes X W^* + \gamma^* \alpha^* \otimes Y Z^* + 1 \otimes X Z^* + ( 1 - ( 1 + \mu^2 ) \gamma^* \gamma  ) \otimes \frac{\mu^2 X Z^* - Y W^*}{1 + \mu^2}\\
        && + 1 \otimes \frac{Y W^* - \mu^2 X Z^*}{1 + \mu^2}\\
        &=& \alpha \gamma \otimes X W^* + \gamma^* \alpha^* \otimes Y Z^* + 1 \otimes ( XZ^* + \frac{YW^* - \mu^2 X Z^* }{1 + \mu^2} ) + ( 1 - ( 1 + \mu^2 ) \gamma^* \gamma  )\\
        && \otimes \frac{\mu^2 X Z^* - Y W^*}{1 + \mu^2}. \eean

Thus, comparing coefficient of $ 1 $ in $ \Phi ( \alpha \gamma ),$ we can see that it belongs to $V^1_0$  if and only if $ X{Z}^{*} + Y{W}^{*} = 0.$ 

In the case of $ 1 - ( 1 + \mu^2 ) \gamma^* \gamma ,$ 
  \bean \lefteqn{\Phi ( 1 - ( 1 + \mu^2 ) \gamma^* \gamma ) }\\
          &=& 1 \otimes 1 - ( 1 + \mu^2 ) ( \alpha \alpha^* \otimes Z Z^* + \alpha \gamma \otimes Z W^* + \gamma^* \alpha^* \otimes W Z^* + \gamma^* \gamma \otimes W W^* )\\
          &=& 1 \otimes 1 - ( 1 + \mu^2 ) ( 1 - \mu^2 \gamma^* \gamma ) \otimes Z Z^* - \alpha \gamma \otimes ( 1 + \mu^2 ) Z W^* - \mu \alpha^* \gamma^* \otimes ( 1 + \mu^2 )\\ 
         && W Z^* - ( 1 + \mu^2 ) \gamma^* \gamma \otimes W W^*\\
          &=& 1 \otimes 1 - ( 1 + \mu^2 ).1 \otimes Z Z^* + ( - 1 + 1 - ( 1 + \mu^2 ) \gamma^* \gamma  ) \otimes - \mu^2 Z Z^* - \alpha \gamma \otimes ( 1 + \mu^2 )\\
          && Z W^* - \mu \alpha^* \gamma^* \otimes ( 1 + \mu^2 ) W Z^* + ( - 1 + 1 - ( 1 + \mu^2 ) \gamma^* \gamma ) \otimes W W^*\\
          &=& 1 \otimes ( 1 - ( 1 + \mu^2 ) Z Z^* + \mu^2 Z Z^* - W W^* ) + ( 1 - ( 1 + \mu^2 ) \gamma^* \gamma ) \otimes (  - \mu^2 Z Z^* +\\
          &&  W W^* ) - \alpha \gamma \otimes ( 1 + \mu^2 ) Z W^* - \mu \alpha^* \gamma^* \otimes ( 1 + \mu^2 ) W Z^*. \eean 

Comparing the coefficient of $ 1 $ in this case, we have the condition $ 1 - ( 1 + \mu^2 ) Z Z^* + \mu^2 Z Z^* - W W^* = 0,$ that is, $ Z{Z}^{*} + W{W}^{*} = 1.$

But these conditions follow from the unitarity of the matrix $ \left(  \begin {array} {cccc}
     X^{*}   &  Z^{*}  \\ Y^{*} & W^{*} \end {array} \right ) ,$ 
        which is nothing but the matrix corresponding  to the restriction of $ \Phi $ to $V^{1/2}_{1/2} .$ 
 Thus, $ \Phi $ keeps $ {\rm Span} \{ t^{1}_{ij} : i,j = -1,0,1 \} $ invariant.

Moreover, we claim that by using the recursive relations ( \ref{rec_1} ), ( \ref{rec_2} ) and the multiplication rule ( \ref{qorient_SUmu2_multiplication_rule_klimyk} ),
we obtain that for all $ l \geq 3/2, ~ \Phi ( V^{l + 1/2}_{i} )  \subseteq V^{l - 1/2}_{i} \oplus V^{l + 1/2}_{i}.$ We prove this for $ t^{l + \frac{1}{2}}_{i,j}, ~ - l + \frac{1}{2} \leq i \leq l - \frac{1}{2}, ~ j \leq l  $ only, as the proofs of the others are exactly similar. We have
\bean \lefteqn{ \Phi ( t^{l + \frac{1}{2}}_{i,j} )}\\
 &=& c ( l,i,j ) \Phi ( \alpha ) \Phi ( t^l_{i + \frac{1}{2},j + \frac{1}{2}} ) + c^{\prime} ( l,i,j ) \Phi ( \gamma ) \Phi ( t^l_{i - \frac{1}{2}, j + \frac{1}{2}} )\\
 &=& c ( l,i,j ) \Phi ( t^{\frac{1}{2}}_{- \frac{1}{2}, - \frac{1}{2}} ) \Phi ( t^l_{i + \frac{1}{2},j + \frac{1}{2}} ) + c^{\prime} ( l,i,j ) \Phi ( t^{\frac{1}{2}}_{\frac{1}{2}, - \frac{1}{2}} ) \Phi ( t^l_{i - \frac{1}{2}, j + \frac{1}{2}} )\\
 &\in& {\rm Span}\{ t^{\frac{1}{2}}_{- \frac{1}{2}, k}. t^{l}_{i + \frac{1}{2},m}, ~ t^{\frac{1}{2}}_{ \frac{1}{2}, k}. t^{l}_{i - \frac{1}{2},m} : k = \pm \frac{1}{2}, ~ m = - l,...l  \} \otimes \clq \eean
 $$\subseteq V^{l - \frac{1}{2}}_{i} \otimes \clq +  V^{l + \frac{1}{2}}_{i} \otimes \clq.$$  

Using these observations, we conclude that $ \Phi $ maps  ${\rm  Span} \{ t^{l}_{ij} : l \geq 1/2 \} $ into itself. 

So, in particular, $ {\rm Ker}( h) = {\rm Span} \{ V^{l}_{i} : i = -l,,,,l, l \geq 1/2 \} $ is invariant under $ \Phi $ which ( along with $ \Phi ( 1 ) = 1 $ ) implies that $ \Phi $ preserves $ h.$
 
\vspace{4mm}

$3. \Rightarrow 4.$  

 We proceed by induction.
 The induction hypothesis holds for $ l =  \frac{1}{2}$ since  linearity means that span $ \{ \alpha, {\gamma}^{*} \} $ is invariant under $ \Phi $ and hence Span $ \{ {\alpha}^{*}, \gamma  \} $ is also invariant. The case for $ l = 1 $ can be checked by inspection as in the proof of $2 \Rightarrow 3$.
 Consider the induction hypothesis that $\Phi$ keeps $V^k_i$ invariant  for all $k,i$ with $k \leq l$. From the proof of $ 2 \Rightarrow 3$ we also have for all $ l \geq \frac{3}{2},$ $ \Phi (  V^{l + 1/2}_{i} )  \subseteq V^{l - 1/2}_{i} \oplus V^{l + 1/2}_{i},$ by using  linearity only. 
  Thus, $\widetilde{\Phi}$ leaves invariant the Hilbert $\clq$ module $(V^{l-\frac{1}{2}}_i \oplus V^{l+ \frac{1}{2}}_i) \ot \clq$, and is a unitary there since $\Phi$ is Haar-state preserving. Since $\widetilde{\Phi}$ leaves invariant $V^{l-\frac{1}{2}}_i \ot \clq$ by the induction hypothesis, it must keep its orthocomplement,  $V^{l +\frac{1}{2}}_i$ invariant as well. 
  
  \vspace{4mm}

$ 4. \Rightarrow 3.$

 The fact that $ \Phi $ keeps $ V^{l}_{i} $ invariant for $ l = 1/2 $ will imply that $ \Phi $ is linear. The proof of Haar state preservation is exactly the same as in $ 2 \Rightarrow 3.$
 
 \vspace{4mm}

$ 4 \Rightarrow 1.$ 

That  $ \Phi $ preserves the Haar state follows from arguments used in the proof of the implication $ 2 \Rightarrow 3. $.  Since  $ \cla_{00} = 
 {\rm Span} \{ t^{l}_{ij} : l \geq 0, i,j = -l,.......l \},$ and  $ \Phi $ keeps each $ V^{l}_{i} $ invariant, it is obvious that $\Phi(\cla_{00}) \subseteq \cla_{00} \ot_{\rm alg} \clq_0$ and   $ \Phi \widehat{D} = ( \widehat{D} \otimes {\rm  id} ) \Phi .$

 \qed

 By  Lemma \ref{equivalent}, we have identified the category $ \widehat{\bf C} $ with the category  of CQG having $ C^* $ actions on $ SU_{\mu} ( 2 ) $  satisfying condition 3. of Lemma \ref{equivalent}. Let the universal object of this category be denoted by $ ( \widehat{\clq}, \Gamma ).$

Then by linearity we can write:

$$ \Gamma ( \alpha ) = \alpha \otimes A + {\gamma}^{*} \otimes B, $$

$$ \Gamma ( {\gamma}^{*} ) = \alpha \otimes C + {\gamma}^{*} \otimes D. $$

Now we shall exploit the fact that $\Gamma$ is a $ \ast $-homomorphism to get relations satisfied by  $ A, B, C, D $ where $ \widehat{\clq} $ is generated as a $ C^* $ algebra by the elements $A,B,C,D.$  
 
\blmma

\label{SU_mu_2_alpha* alpha + gamma* gamma = 1}

\be \label {qorient_sumu2_1} {A}^{*}A + C{C}^{*} = 1, \ee
\be \label {qorient_sumu2_2} {A}^{*}A + {\mu}^{2}C{C}^{*} = {B}^{*}B + D{D}^{*}, \ee
\be \label {qorient_sumu2_3} {A}^{*}B = - \mu D{C}^{*}, \ee
\be \label {qorient_sumu2_4} {B}^{*}A = - \mu C{D}^{*}. \ee

\elmma

{\it Proof :} The proof follows from the relation ( \ref{su2def1} ) by comparing coefficients of $ 1, {\gamma}^{*}\gamma ,{\alpha}^{*}{\gamma}^{*} $ and $ \alpha \gamma $ respectively. \qed

\blmma

\label{SU_mu_2_alpha alpha* + mu2 gamma gamma* = 1}

\be \label {qorient_sumu2_5} A {A}^{*} + {\mu}^{2} C{C}^{*} = 1, \ee
\be \label {qorient_sumu2_6} B {B}^{*} + {\mu}^{2}D{D}^{*} = {\mu}^{2}.1, \ee
\be \label {qorient_sumu2_7} B {A}^{*} = - {\mu}^{2} D{C}^{*}. \ee

\elmma

{\it Proof :} From the equation ( \ref{su2def2} ) by equating coefficients of $ 1 $ and $ {\alpha}^{*} {\gamma}^{*} ,$ we get respectively ( \ref{qorient_sumu2_5} ) and ( \ref{qorient_sumu2_7} ) whereas ( \ref{qorient_sumu2_6} ) is obtained by equating coefficients of $ {\gamma}^{*} \gamma $ and using ( \ref{qorient_sumu2_5} ).

  \qed

\blmma

\label{SU_mu_2_gamma* gamma = gamma gamma*}

\be \label {qorient_sumu2_8}  {C}^{*} C = C {C}^{*}, \ee
\be \label {qorient_sumu2_9} ( 1 - {\mu}^{2} ) {C}^{*}C = {D}^{*}D - D{D}^{*}, \ee
\be \label {qorient_sumu2_10}  {C}^{*}D =  \mu D{C}^{*}. \ee

\elmma

{\it Proof :} The proof follows from the equation ( \ref{su2def3} ) by comparing the coefficients of $ 1,{\gamma}^{*}\gamma, {\alpha}^{*} {\gamma}^{*}, $ respectively. \qed

\blmma

\label{SU_mu_2_alpha gamma = mu gamma alpha}

\be \label {qorient_sumu2_11}  -{\mu}^{2} A{C}^{*} + B{D}^{*} - \mu D^{*}B + \mu C^{*} A = 0, \ee
\be \label {qorient_sumu2_12}  A{C}^{*} = \mu {C}^{*}A,  \ee
\be \label {qorient_sumu2_13}  B{C}^{*} =  {C}^{*}B, \ee 
\be \label {qorient_sumu2_14}  A{D}^{*} = {D}^{*} A. \ee

\elmma

{\it Proof :} The proof follows from the equation ( \ref{su2def4} ) comparing the coefficients of $ {\gamma}^{*}\gamma, 1, {\alpha}^{*} {\gamma}^{*} $ and $ \alpha \gamma $ respectively. \qed

\blmma

\label{SU_mu_2_alpha gamma* = mu gamma* alpha}

\be \label {qorient_sumu2_15}  AC = \mu CA,  \ee
\be \label {qorient_sumu2_16}  BD =  \mu DB, \ee 
\be \label {qorient_sumu2_17}  AD -\mu CB = DA -  {\mu}^{-1} BC.   \ee

\elmma

{\it Proof :} The proof follows from ( \ref{su2def5}  ) from the coefficients of $ {\alpha}^{2}, {\gamma}^{*2}, {\gamma}^{*}\alpha $ respectively. \qed

\vspace{8mm}

Now we consider the antipode, say $ \kappa .$

From the condition $ ( h \otimes {\rm id} )\Gamma( a ) = h( a ).1,$ we see that $ \Gamma $ induces a unitary representation of the compact quantum group via $ \widetilde{\Gamma} ( a \otimes q ) = \Gamma( a )( 1 \otimes q ) .$

Now, the restriction of this unitary representation to the orthonormal set 

$ \{ \sqrt{\frac{1 + \mu^2}{\mu^2}}  \alpha,  \sqrt{1 + \mu^2} {\gamma}^{*} \} $ is given by the matrix : 
$ \left ( \begin {array} {cccc}
   A & \mu C  \\ {\mu}^{-1} B & D \end {array} \right ) .$ 
                        
 Similarly, with respect to the orthonormal set $ \{ \sqrt{ 1 + \mu^2 } {\alpha}^{*}, \sqrt{ 1 + \mu^2 }  \gamma \}  ,$ this representation is given by the matrix: 
    $ \left ( \begin {array} {cccc}
   A^{*} & C^{*}  \\ B^{*} & D^{*} \end {array} \right ) .$

Thus, we have:
 
 $ \kappa( A ) = {A}^{*},~ \kappa( D ) = {D}^{*}, ~ \kappa( C ) = {\mu}^{-2} B^{*}, ~ \kappa( B ) = {\mu}^{2} C^{*}, ~ \kappa( A^{*} ) = A, ~ \kappa(C^{*} ) = B, ~ \kappa(B^{*}) = C, ~ \kappa( D^{*} ) = D .$

 \blmma
 
 \be \label {qorient_sumu2_27} AB = \mu BA, \ee
 \be \label {qorient_sumu2_28} CD = \mu DC, \ee
\be \label {qorient_sumu2_31}  B{C}^{*} = {C}^{*}B. \ee

  \elmma

 {\it Proof :} The relations ( \ref{qorient_sumu2_27} ), ( \ref{qorient_sumu2_28} ), ( \ref{qorient_sumu2_31} ) follow by applying $ \kappa $ to the equations ( \ref{qorient_sumu2_15} ), ( \ref{qorient_sumu2_16} ) and ( \ref{qorient_sumu2_13} ) respectively.  
 
 
 \blmma
 
 \label{oneway}
 
 There exists a  $\ast$-homomorphism $ \phi : U_{\mu}( 2 ) \rightarrow \widehat{\clq} $ defined by $ \phi ( u_{11} ) = A, ~ \phi( u_{12} ) = \mu C, ~ \phi( u_{21} ) = {\mu}^{-1} B, ~ \phi( u_{22} ) = D .$            

 \elmma
 
 {\it Proof :} It is enough to check that the defining relations of  $ U_{\mu}( 2 ) $ are satisfied.
 
 1. $  \phi( u_{11} u_{12} ) = \phi( \mu u_{12} u_{11} ) \Leftrightarrow  \phi ( u_{11} ) \phi( u_{12} ) = \mu \phi( u_{12} ) \phi( u_{11} ) \Leftrightarrow  A ( \mu C ) = \mu ( \mu C )A \Leftrightarrow  AC = \mu CA $ which is the same as ( \ref{qorient_sumu2_15} ).
 
 2. $  \phi( u_{11} u_{21} ) = \phi( \mu u_{21} u_{11} ) \Leftrightarrow A( {\mu}^{-1}B ) = \mu ( {\mu}^{-1} B ) A \Leftrightarrow AB = \mu BA $ which is the same as equation ( \ref{qorient_sumu2_27} ).
 
 3. $  \phi ( u_{12} u_{22} ) = \phi( \mu u_{22} u_{12}  ) \Leftrightarrow \mu CD = \mu D ( \mu C ) \Leftrightarrow CD = \mu DC $ which is the same as equation ( \ref{qorient_sumu2_28} ).
 
 4. $  \phi(u_{21} u_{22} ) = \phi( \mu u_{22} u_{21}  ) \Leftrightarrow {\mu}^{-1} BD = \mu D {\mu}^{-1} B \Leftrightarrow BD = \mu DB $ which is the same as equation ( \ref{qorient_sumu2_16} ).

 5. $ \phi( u_{12}u_{21}  ) = \phi( u_{21}u_{12}    ) \Leftrightarrow \mu C {\mu}^{-1} B = {\mu}^{-1} B \mu C \Leftrightarrow CB = BC . $ 
  
  Now, $ B{C}^{*} = {C}^{*}B $ follows from equation ( \ref{qorient_sumu2_31} ). But by ( \ref{qorient_sumu2_8} ), $ C $ is normal, which implies $ BC = CB .$

 6.  $  \phi( u_{11} u_{22} - u_{22} u_{11}) = ( \mu - {\mu}^{-1} )\phi( u_{12} u_{21}) \Leftrightarrow AD - DA = ( \mu - {\mu}^{-1} ) \mu C {\mu}^{-1} B.$
 
 From ( \ref{qorient_sumu2_17} ), we have $ AD - DA = \mu CB - {\mu}^{-1} BC = ( \mu - \mu^{- 1} ) CB,$ using $ BC = CB.$
 
 
 \qed
 
 \blmma
 
 \label{qorient_U_mu_2_long_set_of relations}
 
 The equations ( \ref{qorient_sumu2_1} ) - ( \ref{qorient_sumu2_17} ) are true when $ A,B,C,D $ are replaced by $ u_{11}, \mu u_{21}, {\mu}^{-1}u_{12} $ and $ u_{22} $ respectively. 
 
  \elmma
 
 {\it Proof :} We check some of the relations ( \ref{qorient_sumu2_1} ) - ( \ref{qorient_sumu2_17}  by using the facts that $ D_{\mu} $ is a central element of $ U_{\mu} ( 2 ) , ~ \kappa ( u_{ij} ) = u^*_{ji} $ ( ( \ref{qorient_U_mu_2_CQG_structure} ) ), Proposition \ref{qorient_U_mu_2_k(u_ij)}, the equations ( \ref{qorient_Umu2_def_1}  ) - (  \ref{qorient_Umu2_def_6}  )  and (  \ref{qorient_U_mu_2_quantum_determinant}  ). The proofs of the others are exactly similar.
 
 Proof for ( \ref{qorient_sumu2_1} ) that is, $ {u_{11}}^{*}u_{11} + ( {\mu}^{-1}u_{12} ) { ({\mu}^{-1} u_{12} ) }^{*} = 1.$     
  \bean \lefteqn{ u^*_{11} u_{11} + \mu^{- 2} u_{12} u^*_{12} ~ = ~  u_{22} D^{-1}_{\mu} u_{11} + {\mu}^{-2} u_{12} ( - \mu u _{21} D^{-1}_{\mu} )}\\
    &=& ( u_{22} u_{11} - {\mu}^{-1} u_{12} u_{21} )D^{-1}_{\mu}  =  D_{\mu} D^{-1}_{\mu}  = 1.\eean

Proof for ( \ref{qorient_sumu2_2} ) that is, $ u^*_{11} u_{11} + \mu^2 ( \mu^{- 1} u_{12} ) {( \mu^{- 1} u_{12} )}^* - ( {( \mu u_{21} )}^* \mu u_{21} + u_{22} u^*_{22} ) = 0 .$
\bean \lefteqn{ u^*_{11} u_{11} + \mu^2 ( \mu^{- 1} u_{12} ) {( \mu^{- 1} u_{12} )}^* - ( {( \mu u_{21} )}^* \mu u_{21} + u_{22} u^*_{22} )}\\
&=& \kappa ( u_{11} ) u_{11} + u_{12} \kappa ( u_{21} ) - ( \mu^2 \kappa ( u_{12} ) u_{21} + u_{22} \kappa ( u_{22} ) )\\
&=& ( u_{22} u_{11} - \mu u_{12} u_{21} ) D^{- 1}_{\mu} - ( - \mu u_{12} u_{21} + u_{22} u_{11} ) D^{- 1}_{\mu}  \\
&=& 0. \eean



  
  Proof for ( \ref{qorient_sumu2_6} ) that is, $ \mu^2 u_{21} u^*_{21} + \mu^2 u_{22} u^*_{22} - \mu^2.1 = 0.$  
  \bean \lefteqn{\mu^2 u_{21} u^*_{21} + \mu^2 u_{22} u^*_{22} - \mu^2.1}\\
   &=& \mu^2 ( u_{21} \kappa ( u_{12} ) + u_{22} \kappa ( u_{22} ) - 1)\\
   &=& \mu^2 ( u_{21} ( - \mu^{- 1} u_{12} D^{- 1}_{\mu} ) + u_{22} u_{11} D^{- 1}_{\mu} - 1 )\\
   &=& \mu^2 ( ( u_{22} u_{11} - \mu^{- 1} u_{21} u_{12}  ) D^{- 1}_{\mu} - 1 )\\
   &=& \mu^2 ( D_{\mu} D^{- 1}_{\mu} - 1 )\\
   &=& 0. \eean
   
   Proof for ( \ref{qorient_sumu2_7} ) that is, $ u_{21} u^*_{11} + u_{22} u^*_{12} = 0.$  
   \bean \lefteqn{u_{21} u^*_{11} + u_{22} u^*_{12}}\\
   &=& u_{21} \kappa ( u_{11} ) + u_{22} \kappa ( u_{21} )\\
   &=& u_{21} u_{22} D^{- 1}_{\mu} - \mu u_{22} u_{21} D^{- 1}_{\mu}\\
   &=& ( u_{21} u_{22} - \mu u_{22} u_{21} ) D^{- 1}_{\mu}\\
   &=& 0. \eean
   
   
   Proof for ( \ref{qorient_sumu2_9} ) that is, $ ( 1 - \mu^2 ) u^*_{12} u_{12} - \mu^2 ( u^*_{22} u_{22} - u_{22} u^*_{22} ) = 0.$   
   \bean \lefteqn{( 1 - \mu^2 ) u^*_{12} u_{12} - \mu^2 ( u^*_{22} u_{22} - u_{22} u^*_{22} )}\\
   &=& ( 1 - \mu^2 ) \kappa ( u_{21} ) u_{12} - \mu^2 (  \kappa ( u_{22} ) u_{22} - u_{22} \kappa ( u_{22} )   )\\
   &=& - \mu ( 1 - \mu^2 ) ( u_{21} u_{12} D^{- 1}_{\mu} ) - \mu^2 ( u_{11} u_{22} D^{- 1}_{\mu} - u_{22} u_{11} D^{- 1}_{\mu} )\\
   &=& - \mu [ ( 1 - \mu^2 ) ( u_{21} u_{12} D^{- 1}_{\mu} ) - \mu ( u_{22} u_{11} D^{- 1}_{\mu} - u_{11} u_{22} D^{- 1}_{\mu} )]\\
   &=& - \mu [ ( 1 - \mu^2 ) u_{21} u_{12} D^{- 1}_{\mu} - \mu ( \mu^{- 1} - \mu ) u_{12} u_{21} D^{- 1}_{\mu}]\\
   &=& - \mu ( 1 - \mu^2 ) ( u_{12} u_{21} - u_{12} u_{21} ) D^{ - 1}_{\mu}\\
   &=& 0. \eean
    
   Proof for ( \ref{qorient_sumu2_10} ) that is, $ u^*_{12} u_{22} - \mu u_{22} u^*_{12} = 0.$   
   \bean \lefteqn{u^*_{12} u_{22} - \mu u_{22} u^*_{12}}\\
   &=& \kappa ( u_{21} ) u_{22} - \mu u_{22} \kappa ( u_{21} )\\
   &=& \mu^2 u_{22} u_{21} D^{- 1}_{\mu} - \mu u_{21} u_{22} D^{- 1}_{\mu}\\
   &=& \mu (  \mu u_{22} u_{21} - u_{21} u_{22} ) D^{- 1}_{\mu}\\
   &=& 0. \eean
   
   Proof for ( \ref{qorient_sumu2_11} ) that is, $ - \mu u_{11} u^*_{12} + \mu u_{21} u^*_{22} - \mu^2 u^*_{22} u_{21} + u^*_{12} u_{11} = 0.$   
   \bean \lefteqn{- \mu u_{11} u^*_{12} + \mu u_{21} u^*_{22} - \mu^2 u^*_{22} u_{21} + u^*_{12} u_{11} }\\
    &=& - \mu u_{11} \kappa ( u_{21} ) + \mu u_{21} \kappa ( u_{22} ) - \mu^2 \kappa ( u_{22} ) u_{21} + \kappa ( u_{21} )u_{11}\\
    &=& - \mu u_{11} ( - \mu u_{21} D^{- 1}_{\mu} ) + \mu u_{21} ( u_{11} D^{- 1}_{\mu} ) - \mu^2  u_{11} u_{21} D^{- 1}_{\mu} - \mu u_{21} u_{11} D^{- 1}_{\mu}\\
    &=& \mu^2 ( u_{11} u_{21} - u_{11} u_{21} ) D^{- 1}_{\mu} + \mu ( u_{21} u_{11} - u_{21} u_{11} ) D^{- 1}_{\mu}\\
    &=& 0. \eean     
     Proof for ( \ref{qorient_sumu2_17} ) that is, $  u_{11} u_{22} - \mu ( \mu^{- 1} u_{12} ) \mu u_{21} = u_{22} u_{11} - \mu^{- 1} ( \mu u_{21} ) ( \mu^{- 1} u_{12} ). $      
     \bean \lefteqn{u_{11} u_{22} - \mu ( \mu^{- 1} u_{12} ) \mu u_{21} - u_{22} u_{11} + \mu^{- 1} ( \mu u_{21} ) ( \mu^{- 1} u_{12} )}\\
     &=& u_{11} u_{22} - \mu u_{12} u_{21} - u_{22} u_{11} + \mu^{- 1} u_{21} u_{12}\\
     &=& 0. \eean  \qed

\blmma

 There is a $ C^{*} $ action $ \Psi $ of $ U_{\mu} ( 2 ) $ on $ SU_{\mu} ( 2 ) $ such that $ ( U_{\mu}( 2 ), \Psi ) $ is an object of $ Obj( \widehat{\bf C} ) $ and $ \Psi $ is given by :
 
 $$ \Psi( \alpha ) = \alpha \otimes u_{11} + {\gamma}^{*} \otimes \mu u_{21}, $$
     
 $$  \Psi( {\gamma}^{*} ) = \alpha \otimes {\mu}^{-1}u_{12} + {\gamma}^{*} \otimes  u_{22} .$$

\elmma

{\it Proof :}  The homomorphism conditions are exactly the conditions ( \ref{qorient_sumu2_1} ) - ( \ref{qorient_sumu2_17} )  with $ A,B,C,D $ replaced by $ u_{11}, \mu u_{21}, {\mu}^{-1}u_{12} $ and $ u_{22} $ respectively which are true by Lemma \ref{qorient_U_mu_2_long_set_of relations}.

 Clearly, $ \Psi $ keeps  $ V^{1/2}_{-1/2} $ invariant and the corresponding representation is a unitary.
    
    It follows from Lemma \ref{equivalent} that  $ ( U_{\mu}( 2 ), \Psi ) $ is an object of $ \widehat{\bf C} .$
    
    \qed
    
    \bcrlre
    
    \label{qorient_sumu2_corollaryabc}
    
    There exists a surjective CQG morphism from $ \widehat{\clq} $ to $ U_{\mu}( 2 ) $ sending $ A,~ \mu C, ~ {\mu}^{-1} B, $ and $ D $ to $ u_{11}, u_{12}, u_{21} $ and $ u_{22} $ respectively.  
    
    \ecrlre
    
    \bthm

 We have $  \widehat{\clq} \cong U_{\mu}( 2 ) $ and hence $ \widetilde{{QISO}^{+}}( D ) \cong U_{\mu}( 2 ) \ast C( \IT ) .$
    
    \ethm
   
   {\it Proof :} The first part follows from  Lemma \ref{oneway} and  Corollary \ref{qorient_sumu2_corollaryabc} and the second part follows from Theorem \ref{unrestricted}.
    \qed

\subsection{A commutative example : spectral triple on $\IT^2$}

We consider the spectral triple $ ( \cla^{\infty}, \clh, D ) $ on $ \IT^{2} $ given by $ \cla^{\infty} = C^{\infty} ( \IT^{2} ),~ \clh = L^{2} ( \IT^{2} ) \oplus L^{2} ( \IT^{2} )$ and $ D =  
  \left ( \begin {array} {cccc}
    0 &  d_{1} + i d_{2}   \\  d_{1} - i d_{2} &  0  \end {array} \right ) ,$ 
    
    where we view $ C( \IT^{2} ) $ as the universal $ C^{*} $ algebra generated by two commuting unitaries $ U  $ and $ V,$ and $ d_{1} $ and $ d_{2} $ are derivations on $ \cla^{\infty} $ defined by :    
     \be d_{1} ( U ) = U, ~ d_{1} ( V ) = 0, ~ d_{2} ( U ) = 0,~ d_{2} ( V ) = V .\ee    
The vectors $ e_{1} = ( 1, ~ 0 ) $ and $ e_{2} = ( 0, ~ 1 ) $ form an orthonormal basis of the eigenspace corresponding to the eigenvalue zero.
  
  The Laplacian in the sense of chapter \ref{qisol} exists in this case, and is given by $ \cll ( U^{m}V^{n} ) = - ( m^{2} + n^{2} ) U^{m}V^{n}  .$  We recall that we denote the quantum isometry group from the Laplacian $ \cll $ in the sense of \cite{goswami} by $ {QISO}^{\cll_{D}} .$
  
  \blmma
  
  \label{qorient_torus_along with 2_equations}
  
  Let $ ( \widetilde{\clq},  W ) $ be an object of $ {\bf Q^{\prime}} ( D ) $. Then the $ \ast  $-homomorphism $ \alpha = \alpha_{W} $  must be of the following form:
  
  \be \label{torus1} \alpha ( U ) = U \otimes z_{1}, \ee  
  \be  \label{torus2} \alpha ( V ) = V \otimes z_{2}, \ee   
   where $ z_{1}, z_{2} $ are two commuting unitaries.
  
  \elmma
 
 {\it Proof:} We denote the the maximal Woronowicz $ C^{*} $ subalgebra of $ \widetilde{\clq} $ which acts on $ C ( \IT^{2} ) $ faithfully by $ \clq .$
 
 We observe that $ D^2 ( a e_i ) = \cll ( a ) e_i $ for $i = 1,2.$ 
 Now, $ W $ commutes with $ D $ implies that $ W $ commutes with $ {D}^{2} $ as well. Using this, we can show that $ ( \cll \otimes {\rm id} ) \alpha ( a ) e_i  = \alpha \cll ( a ) e_i, ~ i = 1,2.$ As the pair $ \{ e_1, e_2 \} $ is together separating for $ C ( \IT^2 ), $ we conclude that $ \alpha $ commutes with the Laplacian $ \cll.$
    Therefore, $ \clq $ is a quantum subgroup of $ {QISO}^{\cll_{D}}.$ From Theorem \ref{QISO_Laplacian_computations_Torus_finaltheorem}, we conclude that $ {QISO}^{\cll_{D}} =  C ( {\IT}^{2} >\!\!\! \lhd ( \IZ^{2}_2 >\!\!\! \lhd \IZ_2 ) ) .$ Thus $ \clq $ must be of the form $C(G)$ for  a classical subgroup $ G $ of the orientation preserving isometry group of  $ \IT^{2} $, which is $ \IT^{2} $ itself  and whose (co )action is given by $ U \mapsto U \otimes U $ and  $ V \mapsto V \otimes V .$ \qed
    
    \bthm
    
    \label{qorient_torus_computation_QISO_final theorem}
    
 The universal CQG $  {\widetilde{QISO}}^+( C^\infty( \IT^{2} ), \clh, D ) $ exists and is isomorphic with $ C( \IT^{2} ) \ast C ( \IT ) \cong C^*({\IZ}^2 \ast \IZ)$ (as a CQG). Moreover, $QISO^+$ of this spectral triple is isomorphic with $C(\IT^2)$.
    
    \ethm
    
 {\it Proof:}  Let $ ( \widetilde{\clq}, W ) $ be an object in $ {\bf Q^{\prime}} ( D ) $ as in Lemma \ref{qorient_torus_along with 2_equations}. Since $ \{ e_1, ~ e_2 \} $ is an orthonormal basis for an eigenspace of $ D ,$ we must have  
    \be \label{torus3} W ( e_{1} ) = e_{1} \otimes q_{11} + e_{2} \otimes q_{12}, \ee    
    \be \label{torus4} W ( e_{2} ) = e_{1} \otimes q_{21} + e_{2} \otimes q_{22}, \ee
  for some $ q_{ij} $ in $ \widetilde{\clq}.$
  
 We now make use of the equation $ ( D \otimes {\rm id} ) \widetilde{W} ( U e_{1} \otimes 1 ) =  \widetilde{W} ( D \otimes {\rm id} ) ( U e_{1} \otimes 1 ).$ Let $ z_1, z_2 $ are as in Lemma \ref{qorient_torus_along with 2_equations}.  We compute 
    \bean \lefteqn{( D \otimes {\rm id} ) \widetilde{W} ( U e_1 \otimes 1 ) }\\
        &=& ( D \otimes {\rm id} ) ( \alpha ( U ) \widetilde{W} ( e_1 \otimes 1 ) )\\
        &=& ( D \otimes {\rm id} ) ( U \otimes z_1 ) ( e_1 \otimes q_{11} + e_2 \otimes q_{12} )\\
        &=& ( D \otimes {\rm id} ) ( U e_1 \otimes z_1 q_{11} + U e_2 \otimes z_1 q_{12} )\\
        &=& U e_2 \otimes z_1 q_{11} + U e_1 \otimes z_1 q_{12}. \eean
        
        On the other hand, \bean \lefteqn{\widetilde{W} ( D \otimes {\rm id} ) ( U e_1 \otimes 1 )}\\
                           &=& \widetilde{W} ( U e_2 \otimes 1 )\\
                           &=& \widetilde{W} ( U \otimes {\rm id} ) {\widetilde{W}}^* \widetilde{W} ( e_2 \otimes 1 )\\
                           &=& \alpha ( U ) \widetilde{W} ( e_2 \otimes 1 )\\
                           &=& ( U \otimes z_1 ) ( e_1 \otimes q_{21} + e_2 \otimes q_{22} )\\
                           &=& U e_1 \otimes z_1 q_{21} + U e_2 \otimes z_1 q_{22}.\eean
                           
        By comparing coefficients of $ U e_{1} $ and $ U e_{2}$ in the both sides of the equality $ ( D \otimes {\rm id} ) W ( U e_{1} ) =  W D U e_{1} ,$ we have,    
   \be z_{1} q_{12} = z_{1} q_{21} \ee   
 and  \be z_{1} q_{11} = z_{1} q_{22}. \ee   
  Since $ z_{1} $ is a unitary, we have $ q_{11} = q_{22} $ and $ q_{12} = q_{21}.$
   
   Similarly, from the relation $ ( D \otimes I ) W ( V e_{1} ) = W D V e_{1},$ we have $ q_{12} = - q_{21}, q_{22} = q_{11}.$ 
    
   By the above two sets of relations, we obtain :
   
  $ q_{12} = q_{21} = 0,~ q_{11} = q_{22} = q $ ( say ).
  
  But the matrix $ \left ( \begin {array} {cccc}
    q_{11} &  q_{12}   \\  q_{21} &  q_{22}  \end {array} \right ) $ is a unitary in $  M_{2} ( \widetilde{\clq} ) ,$ so  $ q  $ is a unitary.
    
    Moreover, we note that $ W ( a e_{i} ) = \alpha( a ) W ( e_{i} ) $  for all $ a $ in $ C^{\infty}( \IT^{2} ) .$
    Using Lemma \ref{qorient_torus_along with 2_equations} and the above observations, we deduce that any CQG which has a unitary representation commuting with the Dirac operator is a quantum subgroup of $ C( \IT^{2} ) \ast C( \IT ).$
    
   On the other hand, $ C( \IT^{2} ) \ast C( \IT ) $ has a unitary representation commuting with $ D ,$ given by the formulae ( \ref{torus1} ) - ( \ref{torus4} ) taking $ q_{12} = q_{21} = 0,~ q_{11} = q_{22} = q^{\prime} $ where $ q^{\prime} $ is the generator of $ C( \IT ) $ and  $ z_{1}, z_{2} $ to be the generator of $ C( \IT^{2} ).$ 
  This completes the proof. \qed
  
  \brmrk 
  
  The canonical grading on $ C ( \IT^2 ) $ is given by the operator $ ( {\rm  id} \otimes \gamma ) $ on $ L^{2} ( \IT^2 \otimes \IC^2 ) $ where $ \gamma $ is the matrix  $ \left ( \begin {array} {cccc}
    0 &  1  \\  - 1 &  0  \end {array} \right ) .$
    
    The representation of $ C( \IT^{2} ) \ast C( \IT ) $ clearly commutes with the grading operator and hence is isomorphic with  $ \widetilde{QISO} ( C ( \IT^2 ),  L^{2} ( \IT^2 \otimes \IC^2 ), D, \gamma ) .$\\
    
    \ermrk

 \brmrk
 
 This example shows that  the conditions of Theorem \ref{unrestricted} are not necessary for the existence of $  \widetilde{{QISO}^{+}}.$
 
 \ermrk

\section{$ QISO^{+} $ for zero dimensional manifolds }

\subsection{Inductive limit construction for quantum isometry groups}

In this section we use the limiting construction for an inductive system of compact quantum groups ( Lemma \ref{basic} ) and give an application for quantum isometry groups which is fundamental for the results of the next section.

The next theorem connects the inductive construction done in Lemma \ref{basic} with some specific quantum isometry groups.

\bthm
 
\label{induc}

Suppose that $\cla$ is a $C^*$-algebra acting on a Hilbert space $\clh$ and that  $D$ is a (densely defined) self adjoint operator
on $\clh$ with compact resolvent, such that $D$ has a one-dimensional eigenspace spanned by a vector $\xi$ which is cyclic and
separating for $\cla$. Let $(\cla^{\infty}_n)_{n \in \IN}$ be an increasing net of a unital $^*$-subalgebras of $\cla$ and put $\cla^{\infty} =
\bigcup_{n \in \IN} \cla^{\infty}_n$. Suppose that $\cla^{\infty}$ is dense in $\cla$ and that for each $a \in \cla^{\infty}$ the commutator $[D,a]$ is
densely defined and bounded. Additionally put $\clh_n = \ol{\cla^{\infty}_n \xi}$, let $P_n$ denote the orthogonal projection on $\clh_n$
and assume that each $P_n$ commutes with $D$. Then each $(\cla^{\infty}_n, \clh_n, D|_{\clh_n})$ is a spectral triple satisfying the
conditions of Theorem \ref{main}, there exist natural compatible CQG\ morphisms $\pi_{m,n}:\widetilde{{QISO}^{+}}(\cla^{\infty}_m, \clh_m,
D|_{\clh_m}) \to \widetilde{{QISO}^{+}}(\cla^{\infty}_n, \clh_n, D|_{\clh_n})$ ($n, m \in \IN, m \leq n)$ and
\[ \widetilde{{QISO}^{+}}(\cla^{\infty}, \clh, D) = \lim_{n \in \IN} \widetilde{{QISO}^{+}}(\cla^{\infty}_n, \clh_n, D|_{\clh_n}). \]
Similar conclusions hold if we replace everywhere above $\widetilde{{QISO}^{+}}$ by $QISO^{+}$.
\ethm
{\it Proof:} We prove the assertion corresponding to $\widetilde{{QISO}^{+}}$ only, since the proof for $QISO^{+}$  follows by very similar arguments.
 Let us denote $ \widetilde{{QISO}^{+}} ( \cla^{\infty}_{n}, \clh_{n}, D_{n} ) $ by $ \cls_{n} $ and the corresponding unitary representation (in $\clh_n$) by  $ U_{n}.$
Let us denote the category of compact quantum groups acting by orientation preserving isometries on
$ (\cla^{\infty}_{n}, \clh_n, D|_{\clh_n}) $ and $ (\cla^{\infty}, \clh, D) $ respectively by $ \bf Q_{n} $ and $ \bf
Q .$

Since $ U_{n} $ is a unitary which commutes with $ D_n \equiv D|_{\clh_n} $ and hence preserves the
eigenspaces of $ D_n$,  it restricts to a unitary representation of $ S_{n} $ on each $ H_{m} $ for
$ m \leq n.$ In other words, $(\cls_n, U_n|_{\clh_m}) \in {\rm Obj}( {\bf Q_m} )$, and by
 the universality of $ \cls_{m}$ there exists a compact quantum group morphism, say, $ \pi_{m,n} : \cls_{m} \rightarrow \cls_{n} $ such that $ ( {\rm id} \otimes \pi_{m,n} ) U_{m}|_{\clh_{m}} = U_{n}|_{\clh_{m}}.$

Let $ p \le m \le n .$ Then  we have $ ( {\rm id}  \otimes \pi_{m,n} \pi_{p,m} )U_{p}|_{\clh_{p}} =
U_{n}|_{\clh_{p}} .$ It follows  by the uniqueness of the map $ \pi_{p,n}$ that $ \pi_{p,n} =
 \pi_{m,n}\pi_{p,m} $, that is \   $(\cls_n)_{n \in \IN}$ forms an inductive system of compact quantum
groups satisfying the assumptions of  Lemma \ref{basic}.  Denote by $\cls_\infty$ the inductive
limit CQG\ obtained in that lemma, with $\pi_{n,\infty} : \cls_n \to \cls$ denoting the
corresponding CQG\ morphisms. The family of formulas $U|_{\clh_n}:=({\rm id} \ot \pi_{n,\infty})
\circ U_n$ combine to define a unitary representation $U$ of $\cls_\infty$ on $\clh$.  It is also
easy to see from the construction that $U$ commutes with $D$. This means that $(\cls_\infty, U) \in
{\rm Obj}({\bf Q})$, hence there exists a unique surjective CQG\ morphism from $\cls:=\widetilde{{QISO}^{+}}(\cla^{\infty},
\clh, D)$ to $\cls_\infty$ identifying $\cls_\infty$ as a quantum subgroup of $\cls$.

The proof will now be complete if we can show that there is a surjective CQG\ morphism in the
reverse direction, identifying $\cls$ as a quantum subgroup of $\cls_\infty$. This can be deduced
from Lemma \ref{basic} by using the universality property of the inductive limit. Indeed,
 for each $n \in \IN$ the unitary representation, say $V_n$, of $\widetilde{{QISO}^{+}}(\cla^{\infty}, \clh, D)$ restricts to $\clh_n$ and commutes with $D$ on that subspace, thus inducing a CQG\ morphism $\rho_n$  from $\cls_n=\widetilde{{QISO}^{+}}(\cla^{\infty}_n,\clh_n, D_n)$ into $\cls$. The family of morphisms $(\rho_n)_{n \in \IN}$ satisfies the compatibility conditions required in Lemma \ref{basic}. It remains to show  that the induced CQG\ morphism $\rho_{\infty}$  from $\cls_\infty$ into $\cls$ is surjective.
By the faithfulness of the representation $V$ of $\widetilde{{QISO}^{+}}(\cla^{\infty}, \clh, D)$, we know that the span of matrix elements corresponding to all $V_n$  forms a norm-dense subset of $\cls$.  As the range of $\rho_n$ contains the matrix elements corresponding to $V_n=V|_{\clh_n}$, the proof of surjectivity of $\rho_{\infty}$ is finished. \qed

\vspace{4mm}

The assumptions of the theorem might seem very restrictive. In the next section however we will describe a natural family of spectral triples on $AF$-algebras, constructed in \cite{Chrivan}, for which we have exactly the situation as above.

\subsection{Quantum isometry groups for spectral triples on AF algebras}

We first recall the construction of natural spectral triples on $AF$ algebras due to
E.\,Christensen and C.\,Ivan ( \cite{Chrivan} ). Let $\cla$ be a unital $AF$ $C^*$-algebra, the norm
closure of an increasing sequence $(\cla_n)_{n \in \IN}$ of finite dimensional $C^*$-algebras. We
always put $\cla_0 = \IC 1_{\cla}$, $\cla^{\infty} = \bigcup_{n=1}^{\infty} \cla_n$  and assume that the
unit in each $\cla_n$ is the unit of $\cla$.
 Suppose that $\cla$ is acting on a Hilbert space $\clh$ and that $\xi \in \clh$ is a separating and cyclic unit vector for $\cla$.
Let $P_n$ denote the orthogonal projection onto the subspace $\clh_n:=\cla_n \xi$ of $\clh$ and
write $Q_0=P_0=P_{\IC \xi}$, $Q_n=P_n - P_{n-1}$ for $n \in \IN$. There exists a (strictly
increasing) sequence of real numbers $(\alpha_n)_{n=1}^{\infty}$ such that the self adjoint operator
$D=\sum_{n \in \IN} \alpha_n Q_n$ yields a spectral triple $(\cla^{\infty}, \clh,D).$  Due to the existence of a cyclic and separating vector the quantum group of orientation
preserving isometries exists by Theorem \ref{unrestricted}.

In \cite{Chrivan}, the following fact was also observed:

\bppsn

 If $\cla$ is infinite-dimensional and $p>0$ then one can choose $(\alpha_n)_{n=1}^{\infty}$ in such a way that the spectral triple is $p$-summable. For this reasons, the spectral triple should be thought of as
{\bf $0$-dimensional noncommutative manifolds}.

\eppsn

Note that for each $n\in\IN$ by restricting we obtain a (finite-dimensional) spectral triple
$(\cla_n, \clh_n, D|_{\clh_n})$. As we are precisely in the framework of Theorem \ref{induc}, to
compute $QISO^{+}(\cla^{\infty},\clh,D)$  we need to understand the quantum isometry groups $QISO^{+}(\cla^{\infty}_n, \clh_n,
D|_{\clh_n})$ and embeddings relating them. To simplify the notation we will write
$\cls_n:=QISO^{+}(\cla_n, \clh_n, D|_{\clh_n})$.

We begin with some general observations.

\blmma

Let  ${\clq \clu}_{\cla_n, \omega_{\xi}}$ denote the universal quantum group acting on $\cla_n$ and
preserving the (faithful) state on $\cla_n$ given by vector $\xi$ (see \cite{wang}). There exists
a  CQG\ morphism from ${\clq \clu}_{\cla_n, \omega_{\xi}}$ to $\cls_n$.

\elmma

{\it Proof:} The proof is based on considering the spectral triple given by $(\cla_n, \clh_n, D'_n)$, where $D'_n = P_n - P_0$. It is then
easy to see that $QISO^{+}(\cla_n, \clh_n, D'_n)$ is isomorphic to the universal compact quantum group acting on $\cla_n$ and
preserving $\omega_{\xi}$. On the other hand universality assures the existence of the   CQG\ morphism from
$QISO^{+}(\cla_n, \clh_n, D'_n)$ to $\cls_n$. \qed


\blmma 

\label{univ}

Assume that each $\cla_n$ is commutative, $\cla_n = \IC^{k_n}$, $n \in \IN$. There exists a
CQG\ morphism from ${\clq \clu}_{k_n}$ to $\cls_{k_n}$, where ${\clq \clu}_{k_n}$ denotes  the universal
quantum group acting on $k_n$ points (\cite{wang}).

\elmma

{\it Proof:} We observe that for any measure $\mu$ on the set $\{1,\ldots,k_n\}$ which has full support there is a natural CQG morphism from ${\clq \clu}_{k_n}$ to ${\clq \clu}_{\IC^{k_n}, \mu}$. In case when $\mu$ is uniformly distributed, we simply have ${\clq \clu}_{\IC^{k_n}, \mu} = {\clq \clu}_{k_n} $, as follows from Lemma \ref{fact}. \qed
 
 \vspace{4mm}

Let $\alpha_n: \cla_n \to \cla_n \ot \cls_n$ denote the universal action (on the $n$-th level). Then we have the following
important property, being the direct consequence of the Theorem \ref{induc}. We have
\begin{equation} \label{invar}\alpha_{n+1} (\cla_n)
\subset \cla_n \ot \cls_{n+1}\end{equation} (where we identified $\cla_n$ with a subalgebra of $\cla_{n+1}$) and $\cls_n$ is
generated exactly by these coefficients of $\cls_{n+1}$ which appear in the image of $\cla_n$ under $\alpha_{n+1}$. This in
conjunction with the previous lemma suggests the strategy for computing relevant quantum isometry groups inductively. Suppose
that we have determined the generators of $\cls_n$. Then $\cls_{n+1}$ is generated by generators of $\cls_n$ and these of the
${\clq \clu}_{\cla_n, \omega_{\xi}}$, with the only additional relations provided by the equation ( \ref{invar} ).

This will be used below to determine the concrete form of relations determining $\cls_n$ for the commutative AF algebras.

Before stating the next result,  we fix some notations. Let $ \cla_n $ be a sequence of commutative finite dimensional $ C^* $ algebras as above. Let $ \cla_n = C ( X_n ) $ where $ X_n = \{ x_1, x_2,..., x_m \}.$ Dualizing  the embedding from $ \cla_n $ to $\cla_{n + 1},$ there is a surjective map, say $ f_{n + 1,n} $  from $ X_{n + 1} $ to $ X_n. $ Let $ l_i  $ denote the cardinality of the set $ \{ x \in X_{n + 1} : f_{n + 1,n} ( x ) = x_i \}.$ Thus the embedding of $ \cla_n $ into $ \cla_{n + 1} $ is determined by the sequence $ \{ l_i : i = 1,2,..., m \}.$ We note that the cardinality of $ X_{n + 1}  $ equals $ \sum^m_{i = 1} l_i.$ Moreover, a basis of  $ \cla_{n + 1} $ is given by $ \{ e_{i, r_i}: r_i \in \{ 1,2,... l_i \} \}$  where  $ e_{i, r_i} $ is the indicator function of an element $ y $ of $ X_{n + 1} $ such that  $ f_{n + 1, n} ( y )  = x_i $ and $ y $ is the $ r_i $ th element in  $  X_{n + 1} $ belonging to $ f^{- 1}_{n + 1, n} \{ x_i \}.$   
 
\blmma 

\label{embed}

Let $\cla$ be a commutative AF algebra. Suppose that $\cla_n$ is isomorphic to $\IC^{m}$ and the embedding of $\cla_{n}$ into
$\cla_{n+1}$ is given by a sequence $(l_i)_{i=1}^{m}$. Let $m'= \sum_{i=1}^{m} l_i$. Suppose that the `copy' of ${\clq \clu}_{m}$ in
$\cls_n$ is given by the family of projections $a_{i,j}$ ( $i,j\in\{1, \ldots m\}$ ) and that the `copy' of ${\clq \clu}_{m'}$ in
$\cls_{n+1}$ is given by the family of projections $a_{(i,r_i),(j,s_j)}$ ($i,j\in \{1,\ldots,m\},$ $r_i\in\{1, \ldots, l_i\} $,
$s_j \in \{1, \ldots,l_j\}$). Then the formula ( \ref{invar} ) is equivalent to the following system of equalities:
\begin{equation} \label{finite} a_{i,j} = \sum_{r_{i}=1}^{l_i} a_{(i,r_i), (j,s_j)}\end{equation}
for each $i,j \in \{1, \ldots,m\},s_{j}\in \{1,\ldots, l_j\}$.
\elmma

{\it Proof:} We have (for the universal action $\alpha: \cla_n \to \cla_n \ot \cls_n$)
\[\alpha(\widetilde{e_i} ) = \sum_{j=1}^m \widetilde{e_j} \ot a_{i,j},\]
where by $\widetilde{e_i}$ we denote the image of the basis vector $e_i\in \cla_n$ in $\cla_{n+1}.$
  As $\widetilde{e_j} = \sum_{r_j=1}^{l_j}
e_{(j,s_j)}$,
\[\alpha(\widetilde{e_i}) = \sum_{r_i=1}^{l_i} \alpha(e_{i,r_i}) = \sum_{r_i=1}^{l_i} \sum_{j=1}^m \sum_{s_j=1}^{l_j}
e_{(j, s_j)}  \ot a_{(i,r_i), (j,s_j)}.\] On the other hand we have
\[\alpha(\widetilde{e_i} ) = \sum_{j=1}^m \sum_{s_j=1}^{l_j}
e_{(j,s_j)}  \ot a_{i,j},\] and the comparison of the formulas above yields exactly ( \ref{finite} ).\qed

\vspace{4mm}

One can deduce from the above lemma the exact structure of generators and relations between them for each $\cls_n$ associated with
a commutative AF algebra. To be precise, if $\cla_n= \IC^{k_n}$ for some $k_n \in \IN$, then the quantum isometry group $\cls_n$
is generated as a unital $C^*$-algebra by the family of self adjoint projections $\bigcup_{i=1}^{n} \{a_{\alpha_i,
\beta_i}:\alpha_i, \beta_i = 1,\cdots, k_i\}$ such that for each fixed $i=1,\ldots,n$ the family $ \{a_{(\alpha_i,
\beta_i)}:\alpha_i, \beta_i = 1,\cdots, k_i\}$ satisfies the relations of ${\clq \clu}_{k_n}$ and the additional relations between
$a_{(\alpha_i, \beta_i)}$ and $a_{(\alpha_{i+1}, \beta_{i+1})}$ for $i \in \{1,\ldots,n-1\}$ are given by the formulas
( \ref{finite} ), after suitable reinterpretation of indices according to the multiplicities in the embedding of $\IC^{k_i}$ into
$\IC^{k_{i+1}}$.

\cleardoublepage

\chapter{Quantum isometry groups for Rieffel deformed manifolds}

\label{deformation}

    In this chapter, we give a general scheme for computing $ QISO^{\cll} $ and $  \widetilde{QISO^{+}_{R}}  $  by proving that $ QISO^{\cll} ,$ ( respectively $ \widetilde{{QISO}^{+}_{R}} $ ) of a deformed noncommutative manifold coincides with (under reasonable assumptions) a similar deformation of  $ QISO^{\cll}, $ ( respectively $  \widetilde{{QISO}^{+}_{R}} $ ) of the original manifold. 
    
    \section{Deformation of spectral triple} 

     \label{deformation_section_deformation_of_spectral_triple}

 We recall from Chapter \ref{preliminaries} the generalities of CQG s and Hopf algebras, in particular, the dense unital Hopf $\ast$-subalgebra $ \cls_0 $ of a CQG $\cls$ generated by the matrix elements of the irreducible unitary representations, the Sweedler convention for CQG action, as well as the convolutions $f \triangleleft c,  ~ c \triangleright f $ and $ f \diamond g $ for  functionals $ f, g $ on $ \cls $ and $ c $ in $ \cls.$  Moreover,  given an action $\gamma : \clb \raro \clb \ot \cls$ of the compact quantum group $(\cls, \Delta)$  on a unital $C^*$-algebra $\cla$, the dense, unital  $\ast$-subalgebra of $\cla$ on which $\gamma$ becomes an action by the Hopf $\ast$-algebra $\cls_{0}$ is going to be denoted by $ \cla_0.$ 
 
 {\it A word of caution:}~The algebra $ \cla_0 $ should not be confused with the Rieffel deformed $ C^* $ algebra $ \cla_{J} $ in the case $ J = 0 ,$ for which we simply write $ \cla. $ 
 
 Let $ (\cls, \Delta_\cls) $ be a compact quantum group. We also adopt the convention of calling a vector space $M$ an $ \cls $ co-module if it is an algebraic $\cls_0$ co-module in the sense of definition \ref{preliminaries_CQG_comodule}.

\vspace{4mm}

 Before introducing the set up in which we are going to work, we prove the following well known fact for the sake of completeness.
 
 \bppsn
 
 Let $ E $ be a Banach space and $ G $ a second countable Lie group with a strongly continuous action $ \alpha $ on $ E $ such that $ \| \alpha_g ( x ) \| = \| x \| $ for all $g$ in $G$ and for all $x$ in $E.$ Then  $ E^{\infty} = \{ e \in E : g \rightarrow \alpha_g ( e ) ~ {\rm is} ~ C^{\infty} \} $ is norm dense in $ \cla. $  
  
 \eppsn
 
 {\it Proof:} For a compactly supported continuous function $ f $ on $ G $ and $ a $ in $ E $, we will denote by $ \alpha ( f ) ( a )  $ the norm convergent integral $ \int_{G} f ( h ) \alpha_{h} ( a ) dh $ where $ dh $ denotes a left invariant Haar measure on $ G.$ Then, it can be seen that  $ \alpha_g ( \alpha ( f ) a ) = \int_{G} f ( g^{- 1} h ) \alpha_h ( a ) dh. $ Thus, for $ f $ in $ C^{\infty}_c ( G ), ~ \alpha ( f ) ( a )  $ is in $ E^{\infty}.$ Now, for any $ \epsilon > 0,   $ we choose a small enough neighbourhood $ U $ of identity of $ G, $ such that $ \left\| \alpha_g ( a ) - a \right\| \leq \epsilon $ for all $ g $ in $ U .$ Next, we choose $ f  $ in $ C^{\infty}_c ( G ) $ with $ f \geq 0, ~ \int_{G} f dh = 1 $ and $ {\rm supp} ( f ) \subseteq U.$ Then, 
 \bean \lefteqn{ \left\| \alpha ( f ) ( a ) - a \right\| }\\
       &=& \left\| \int_{G} f ( g ) \alpha_g ( a ) dg - a \int_{G} f ( g ) dg  \right\|\\
       &=& \left\| \int_{G} f ( g ) (  \alpha_g ( a ) - a  ) dg  \right\|\\
       &\leq& \int_G f ( g ) \left\| \alpha_g ( a ) - a \right\| dg\\
       &\leq& \epsilon .\eean 
       
       This shows that $ E^{\infty} $ is dense in $ E.$ \qed
       
 
 
 
 \blmma
 
 \label{preliminaries_Rieffel_positivity}
 
Let $\cla$ be a $C^*$ algebra with a strongly continuous action $\alpha$ of $G$ as above. Then $ \cla^{\infty} $ is closed under holomorphic functional calculus. Let $ \phi $ be a positive linear map from $ \cla^{\infty} $ to another $C^*$ algebra  $ \clb.$ Then, for any self adjoint element $ x $ in $ \cla^{\infty}, ~ \left\| \phi ( x ) \right\| \leq \left\| x \right\| \phi ( 1 ).$
 
 \elmma
 
 {\it Proof:} The first fact is quite well known. We refer to \cite{goswamibook} for a proof. For the second part, let $ x $ be a self adjoint element of $ \cla^{\infty}.$ Then, $ y = (1 + \epsilon ) \left\| x \right\|  - x $ is a positive and invertible element ( since its spectrum does not contain zero ) of $ \cla^{\infty},$ which being closed under holomorphic functional calculus, is closed under taking square root of an invertible element. Thus, $ y^{\frac{1}{2}} $ belongs to $ \cla^{\infty} $ and therefore $ \phi ( y ) = \phi ( ( y^{\frac{1}{2}} )^*  y^{\frac{1}{2}} ) \geq 0 .$ This proves the Lemma. \qed           

\vspace{4mm}

Let $ ( \cla, \IT^{n}, \beta ) $ be a $ C^{*} $ dynamical system, that is, $\cla$ is endowed with a strongly continuous action of $\IT^n$ by $\ast$ automorphisms. Moreover, $ \pi_{0} : \cla \rightarrow \clb( \clh ) $ be a faithful representation, where $ \clh $ is a separable Hilbert space.
  
Let $ \cla^{\infty}  $ be the smooth algebra corresponding to the $ \IT^{n} $ action $ \beta.$ 
  
  Assume now that we are given  a spectral triple $ ( \cla^{\infty}, \pi_{0}, \clh , D ) $ of compact type. Suppose that $ D $ has eigenvalues $ \{ \lambda_{0}, \lambda_{1},......... \} $ and $ V_{i} $  denotes the (finite dimensional) eigenspace of $ \lambda_{i} $ and let $ \cls_{00} $ denote the linear span of $ \{ V_{i}: ~ i = 0, 1,2,.. \}.$  
  
 Suppose, furthermore,  that there exists a compact abelian Lie group $ \widetilde{\IT^{n}} $, with a  covering map $ \gamma :  \widetilde{\IT^{n}} \rightarrow \IT^{n}  .$ The Lie algebra of both $ \IT^{n} $ and  $ \widetilde{\IT^{n}} $ are isomorphic with $ \IR^{n} $ and we denote by $ e $ and $ \widetilde{e} $ respectively the corresponding exponential maps, so that $ e( u ) = e( 2 \pi iu ) , u \in \IR^{n} $ and $ \gamma ( \widetilde{e}( u ) ) = e ( u ).$ By a slight abuse of notation we shall denote the $\IR^n$-action $\beta_{e(u)}$ by $\beta_u$.

   {\bf Assumption}:\\
  There exists a strongly continuous unitary representation $   V_{\tilde{g}}$ , $\tilde{g} \in \widetilde{\IT^{n}} $ of $  \widetilde{\IT^{n}} $ on $\clh$ such that \\
     (a) $  V_{\tilde{g}} D = D V_{\tilde{g}}$ for all $ \tilde{g},$\\
   (b) $ V_{\tilde{g}} \pi_{0} ( a ) {V_{\tilde{g}}}^{-1} = \pi_{0} ( \beta_{g} ( a ) )$, where $a $ belongs to $\cla, ~ \tilde{g} $ belongs to $ \widetilde{\IT^{n}},$ and $ g = \gamma ( \tilde{g} )  .$\\
   
 We shall now show that we can `deform' the given spectral triple along the lines of \cite{connes_etal}.   For each $J$, the map $ \pi_{J}: \cla^\infty \raro {\rm Lin}(\clh^\infty)$ (where $ \clh^{\infty} $ is the smooth subspace corresponding to the representation $ V $ and ${\rm Lin}(\clv)$ denotes the space of linear maps on a vector space $\clv$) defined by $$ \pi_{J} ( a ) s \equiv  a \times_{J} s :=\int \int \beta_{Ju} ( a ) \widetilde{\beta}_{v} ( s ) e( u.v ) du dv$$  extends to a $\ast$-representation of the $C^*$-algebra $ \cla^{\infty} $ in $\clb( \clh )$  where   $ \widetilde{\beta_{v}} = V_{\widetilde{e}( v )} $ ( which clearly maps $ \clh^{\infty}$ into $\clh^{\infty}        ).$  
    
    We can extend the action of $ \IT^n $ on the $ C^* $ subalgebra $ \cla_1 $ of $ \clb ( \clh ) $  generated by $ \pi_{0} ( \cla ), ~ \{ e^{itD} : t \in \IR \} $ and elements of the form  $ \{ [ D, a ]: a \in \cla^{\infty} \} $  by $  \beta_g ( X ) = V_{\widetilde{g}} X {V_{\widetilde{g}}}^{- 1} $  for all $ X $ in  $ \cla_1 $ where by an abuse of notation, we denote the action by the same symbol $ \beta .$ Let $ \cla^{\infty}_1 $  denote the smooth vectors of $ \cla_1 $ with respect to this action. We note that  for all $ a $ in $ \cla^{\infty}_{1}, ~ [ D, a ]  $ belongs to  $ {\cla}^{\infty}_{1} .$
    
    \blmma
    
    \label{qorient_deformation_extension}
    
    $ \beta $ is a strongly continuous action (in the $C^*$-sense)  of $ \IT^n $ on $ \cla_1 $ and hence    for all $ X $ in $ { \cla_1 }^{\infty}, ~ \pi_{J} ( X ) $ defined by 
    $$ \pi_{J} ( X ) s = \int \int \beta_{Ju} ( X ) \widetilde{\beta}_{v} ( s ) e( u.v ) du dv $$ is a bounded operator.
    
    \elmma
    
  {\it Proof:} We note that $ \beta $ is already strongly continuous on the $ C^* $ algebra generated by $ \pi_{0} ( \cla ), ~ \{ e^{itD} : t \in \IR \} .$ Thus it suffices to check the statement for elements of the form  $ [ D, a ] $ where $ a $ belongs to $ \cla^{\infty} .$
  
  To this end, fix any one parameter subgroup $g_t$ of $\IT^n$ such that $g_t$ goes to the identity of $\IT^n$ as $t \raro 0$. 
Let $ T^{\prime}_{t}, ~ \widetilde{T_{t}} $ denote the group of normal $\ast$-automorphisms  on $ \clb(\clh) $ defined by $ T^{\prime}_{t} ( X ) = V_{g_{t}} X V_{g^{- 1}_{t}} $ and $ \widetilde{T_t} ( X ) = e^{i t D} X e^{- i t D} .$  As  $ V_{g_{t}} $ and $ D $ commute, so do 
the  generators of $ T^{\prime}_{t} $ and $ \widetilde{T_t} $ . In particular, each of these generators leave the domain of the other invariant. Note also that $\cla^\infty$ is in the domain of the both the generators, and the generator  of $\widetilde{T_t}$ is given by $[D, \cdot]$ there. Thus, for 
 $a $ in $ \cla^\infty,$ we have $a, ~ [D, a] $ belong to  $ {\rm Dom}(\Xi)$ (where $\Xi$ is the generator of $T^\prime_t$),  and $\Xi([D,a])=[D,\Xi(a)] $ belongs to $ \clb(\clh).$

Using this, we obtain 
  $$ \left\| T^{\prime}_{t} ([ D, a ])  - [ D, a ] \right\|  = \int^{t}_{0}  T^{\prime}_{s} (\Xi( [ D, a ])) ds  \leq ~ t  \left\| \Xi([ D, a ]) \right\|  .$$ 
 The required strong continuity follows from this. Then applying Proposition \ref{deformation_content_of_theorem4.6_rieffel} to the $ C^{*} $ algebra $ \cla_1 $ and the action $\beta, $ we deduce that $ \pi_{J} ( X ) $ is a bounded operator. \qed
   
  \blmma
  
  \label{deformation_deformed representation_is_a_spectral triple}
   
  For each $J$, $ ( \cla^{\infty}_{J}, \pi_{J}, \clh, D ) $ is a spectral triple, that is, $ [ D, \pi_{J}( a ) ] $ belongs to $ \clb ( \clh ) $ for all $ a $ in $ \cla^\infty_J .$
   
  \elmma 
  
  {\it Proof:} $ [ D, \pi_{J} ( a ) ] ( s )
   = D \int \int \beta_{Ju} ( a ) \widetilde{\beta_{v}} ( s) e ( u. v ) du dv - \int \int \beta_{Ju} ( a ) \widetilde{\beta_{v}} ( D s ) e ( u. v ) du dv .$
   
   Using ( \ref{preliminaries_oscillatory_definition} ) and closability of $ D ,$ we have
   
    $ D \int \int \beta_{Ju} ( a ) \widetilde{\beta_{v}} ( s) e ( u. v ) du dv = \int \int D ( \beta_{Ju} ( a ) \widetilde{\beta_{v}} ( s) ) e ( u. v ) du dv .$
   
 As $ D $ commutes with $ V, $  the above expression equals
  \bean \int \int D ( \beta_{Ju} ( a ) \widetilde{\beta_{v}} ( s) ) e ( u. v ) du dv - \int \int \beta_{Ju} ( a ) D  \widetilde{\beta_{v}} ( s ) e ( u. v ) du dv. \eean


   So we have   
  \bean  [D, \pi_J(a)](s)
   &=& \int \int [ D, \beta_{Ju} ( a) ] \widetilde{\beta_{v}} ( s) e ( u. v ) du dv\\
   &=& \int \int V_{\widetilde{Ju}} [ D, a ] {V_{\widetilde{Ju}}}^{- 1} \widetilde{\beta_{v}} ( s) e ( u. v ) du dv\\
   &=& \pi_{J} ( [ D, a ] ), \eean
   
   which is a bounded operator by Lemma \ref{qorient_deformation_extension}. \qed
   
   \vspace{8mm}
   
   \section{Some preparatory results}
   
   \label{deformation_preparatory results}
   
   In this section, we prove some preparatory results which will be needed in the next two sections. Let $ \IT^n, ~ \widetilde{\IT^n}, \beta, ~ \widetilde{\beta}, \gamma $ be as in the previous subsection. By abuse of notation, we will use the symbols $ \beta $ and $ \widetilde{\beta} $ for the corresponding comodule maps also. Let $ \gamma^*, \gamma_{*} $ be the canonical maps induced by $ \gamma $ from $ C ( \IT^n ) \rightarrow C ( \widetilde{\IT^n} ) $ and $ {\rm Lie} ( \widetilde{\IT^n} ) \rightarrow {\rm Lie} ( \IT^n ) $ respectively. Moreover, from now on, we will identify $ \cla^{\infty}_{J} $ with $ \pi_{J}( \cla^{\infty} )$ and often write $\pi_0(a)$ simply as $a$. 
   
  {\bf Assumption}
  
  {\bf 2a.} Let $ ( \widetilde{\clq},\Delta ) $ be a CQG and $ \clq $ a Woronowicz $ C^* $ subalgebra of $ \widetilde{\clq}.$
  Let there exist unital $ \ast $-subalgebra $ \cla_0 \subseteq \cla, $ which is norm dense in every $ \cla_J, $ such that $ \alpha $ is an action : $ \alpha ( \pi_0 ( \cla_0 ) ) \subseteq \pi_0 ( \cla_0 ) \otimes \clq_0.$ Let $ \cls_0 $ be a vector subspace of $ \clh $ ( not necessarily closed ) such that there is a map $ \widetilde{\alpha} : \cls_0 \rightarrow \cls_0 \otimes_{alg} \widetilde{\clq_0}$ making it into an algebraic $ \widetilde{\clq_0} $ co module.
  Moreover, $ ( {\rm id} \otimes \pi_{\widetilde{\clq}} ) \widetilde{\alpha} = \widetilde{\beta}.$  
  
    {\bf 2b.} $ C ( \widetilde{\IT^n} ) $ is a quantum subgroup of $ \widetilde{\clq},$ the quotient map being denoted by $ \pi_{\widetilde{\clq}}.$ 
    
    {\bf 2c.} $ \widetilde{\alpha} ( as ) = \alpha ( a ) \widetilde{\alpha} ( s ) $ for $ a  $ in $ \cla_0, ~ s $ in $ \cls_0.$    
    
    
    
    
    
    
    \vspace{4mm}
    
We recall that we shall denote by $ \eta $  the canonical homomorphism from $\IR^n$ to $\IT^n$ given by $ \eta (x_1,x_2,......,x_n ) = ( e (x_1),e (x_2),.....,e (x_n) ).$ Moreover, we define $ \Omega ( u ) := {\rm ev}_{e (u)} \circ \pi_\clq$, $\widetilde{\Omega}(u):={\rm ev}_{\tilde{e}(u)} \circ \pi_{\widetilde{\clq}}$, for $u $ in $ R^n,$ where  ${\rm ev}_x$ (respectively ${\rm ev}_{\tilde{x}}$ )  denotes  the state on $C(\IT^n)$ ( respectively, on $C(\widetilde{\IT}^n)$ ) obtained by evaluation of a function at the point $x$ (respectively $\tilde{x}$).
 
 
  
  \vspace{2mm} 
  
We now make some observations. 

\blmma

\label{deformation_observation_from_assumptions}

1. From assumption {\bf 2c.}, it follows that $ {\rm ad}_{ \widetilde{\alpha}}   = \alpha .$

 2. $ ( {\rm id} \otimes \pi_{\widetilde{\clq}} ) {\rm ad}_{ ( \widetilde{\alpha} )} = {\rm ad}_{ \widetilde{\beta}}.$

3. $ \widetilde{\beta}_x =  ( {\rm id} \otimes \widetilde{\Omega} ( x ) ) \widetilde{\alpha}.$

4. $ \beta_x =  ( {\rm id} \otimes \Omega ( x ) ) \alpha.$

5. $ (\gamma^*)^{-1} \circ \pi_{\widetilde{\clq}}  $ is a surjective $ C^* $ homomorphism from $ \clq $ to $ C( \IT^n ) $ identifying $ C( \IT^n ) $ as a quantum subgroup of $ \clq.$
\elmma

 {\it Proof :} By using ( \ref{preliminaries_comodule_inverse} ), we have \bean \lefteqn{{\rm ad}_{ \widetilde{\alpha}} ( a ) s  }\\
                   &=& \widetilde{\alpha} ( a \otimes {\rm id} ) {\widetilde{\alpha}}^{- 1} ( s )\\
                   &=& \widetilde{\alpha} ( a \otimes {\rm id} ) ( s_{(1)} \otimes \kappa ( s_{(2)} )\\
                   &=& \widetilde{\alpha} ( a s_{(1)} ) \otimes \kappa ( s_{(2)} )\\
                   &=& \alpha ( a ) \widetilde{\alpha}( s_{(1)} ) \otimes \kappa ( s_{(2)} )\\
                   &=&\alpha ( a ) ( \widetilde{\alpha} \otimes {\rm id} ) ( {\rm id} \otimes \kappa ) \widetilde{\alpha} ( s )\\
                   &=& \alpha ( a ) ( \widetilde{\alpha} \otimes {\rm id} ) ( {\widetilde{\alpha}}^{- 1} \otimes {\rm id} ) ( s )\\
                   &=& \alpha ( a ) s, \eean where we have used Sweedler notations. This proves 1.          
                   
  2. follows from 1. and the fact that $ \pi_{\widetilde{\clq}} $ is a homomorphism.
    \bean \lefteqn{ \widetilde{\beta}_x ( h )}\\
    &=& \widetilde{\beta} ( h ) ( \widetilde{e} ( x ) ) = ( {\rm id} \otimes \widetilde{\pi_{\clq}} ) \widetilde{\alpha} ( h ) ( \widetilde{e} ( x ) ) = ( {\rm id} \otimes {\rm ev}_{\widetilde{e}( x )} \widetilde{\pi}_{\clq} ) \widetilde{\alpha} ( h )\\
    &=& (  {\rm id} \otimes \widetilde{\Omega} ( x )  ) \widetilde{\alpha} ( h ) .\eean
    Therefore, $ \widetilde{\beta}_x =  ( {\rm id} \otimes \widetilde{\Omega} ( x ) ) \widetilde{\alpha}.$ Similarly, 4. follows from 2.          
                   
 We now prove 5. Let us denote by $\gamma^*$ the dual map of $\gamma$, so that $\gamma^* : C(\IT^n) \raro C(\widetilde{\IT}^n)$ is an injective $C^*$-homomorphism. It is quite clear that $({\rm id} \ot \pi_{\widetilde{\clq}})  \circ \alpha (\cla_0) \subseteq {\rm Im}({\rm id} \ot \gamma^*)$, hence it follows that $\pi_{\widetilde{\clq}}(\clq_0) \subseteq {\rm Im}(\gamma^*)$. Thus, $\pi_\clq:=(\gamma^*)^{-1} \circ \pi_{\widetilde{\clq}} $ is a surjective CQG morphism from $\clq$ to $C(\IT^n)$, which identifies $C(\IT^n)$ as a quantum subgroup of $\clq$. \qed
 
 \vspace{4mm}

    
    For a fixed $J$, we shall work with several multiplications on the vector space ${\cla_0} \ot_{\rm alg} \clq_0$. We shall denote the counit and antipode of $\widetilde{\clq_0}$ by $\epsilon$ and $\kappa$ respectively. Let us define the following operation :    
     $$ x \odot y = \int_{\IR^{4n}} e( -u.v )e( w.s )(\widetilde{\Omega} ( -Ju ) \triangleleft x \triangleright (\widetilde{\Omega} ( Jw ) ) ) (\widetilde{\Omega} ( -v ) \triangleleft y \triangleright \widetilde{\Omega}( s ))  du dv dw ds ,$$
     where $x,y $ belong to $ \widetilde{\clq_0}.$ Then $ \odot $ is a  bilinear maps, and will be seen to be  associative multiplication later on. 
     
     We note that when $ x $ is in $ \clq_0, f $ in $C ( \IT^n ), ~ \gamma^* ( f ) ( \widetilde{e} ( u ) = f ( \gamma ( \widetilde{e} ( u ) ) = f ( e ( \gamma_* ( u ) ) ) = f ( e ( u ) )$ ( as $ \gamma $ is a covering map ). Using this, we have 
      \bean \lefteqn{ ( \widetilde{\Omega} ( u ) \otimes {\rm id} ) \Delta ( x )}\\
      &=& ( ev_{\widetilde{e}( u )} \pi_{\widetilde{\clq}} \otimes {\rm id} ) \Delta ( x )\\
      &=& ( ev_{\widetilde{e}( u )} \gamma^* \pi_{\clq} \otimes {\rm id} ) \Delta ( x )\\
      && ( {\rm as} ~ \gamma^* \pi_{\clq} = \pi_{\widetilde{\clq}} )\\ 
      &=& ( ev_{e( u )} \pi_{\clq} \otimes {\rm id} ) \Delta ( x )\\
      &=& ( \Omega ( u ) \otimes {\rm id} ) \Delta ( x ) \eean
      and thus when $x$ belongs to $\clq_0,$ \be \label{deformation_widetilde_omega=omega_on_Q} ( \widetilde{\Omega} ( u ) \otimes {\rm id} ) \Delta ( x ) = ( \Omega ( u ) \otimes {\rm id} ) \Delta ( x ). \ee

    Moreover, we define bilinear maps $\bullet, ~ \bullet_J, $ by setting  $(a \ot x) \bullet (b \ot y):=ab \ot x \odot y, ~ (a \ot x) \bullet_J (b \ot y):=(a \times_J b) \ot (x \odot y),$ for $a,b $ in $ {\cla_0}$, $x,y $ in $ \widetilde{\clq}_0.$


\blmma

\label{deformation_general_Omega}

For $ x $ in $ \widetilde{\clq_0}, $ we have 
$$  \widetilde{\Omega}(u)  \triangleleft   ( \widetilde{\Omega}(v) \triangleleft  x ) = ( \widetilde{\Omega}(u) \diamond \widetilde{\Omega}(v) ) \triangleleft  x.  $$
For $ x $ in $ \clq_0, $ we have
$$ \Omega(u) \triangleleft   ( \Omega(v) \triangleleft x ) = ( \Omega(u) \diamond \Omega(v) ) \triangleleft x.$$
        
\elmma

 {\it Proof :} 
          We will denote by $ \Delta_{\widetilde{\IT^{n}}} $ the coproduct on $ C ( \widetilde{\IT^n} ) ,$ hence, we have     
         \be \label{deformation_general_pi_delta_0} ( \pi_{\widetilde{\clq}} \otimes \pi_{\widetilde{\clq}} ) \Delta = \Delta_{\widetilde{\IT^{n}}} \pi_{\widetilde{\clq}} \ee 
          
          Moreover, we note that as $ \IT^{n} $ is a commutative group, $ f \diamond g = g \diamond f $ for any two functionals $f $ and $ g $ on  $ C ( \IT^{n} ).$  
        \bean \lefteqn{ \widetilde{\Omega}(u) \triangleleft   ( \widetilde{\Omega}(v) \triangleleft  x )}\\
              &=& ( \widetilde{\Omega} ( u ) \otimes {\rm id} ) \Delta ( \widetilde{\Omega} ( v ) \triangleleft x )\\
              &=& ( \widetilde{\Omega} ( u ) \otimes {\rm id} ) \Delta ( \widetilde{\Omega} ( v ) ( x_{(1)} ). x_{(2)}  )\\
              &=& ( \widetilde{\Omega} ( u ) \otimes {\rm id} ) \Delta ( x_{(2)} ) \widetilde{\Omega} ( v ) ( x_{(1)} )\\
              &=& ( \widetilde{\Omega} ( v ) \otimes \widetilde{\Omega} ( u ) \otimes {\rm id} )  ( x_{(1)} \otimes x_{(2)(1)} \otimes x_{(2)(2)} )\\
              &=& ( \widetilde{\Omega} ( v ) \otimes \widetilde{\Omega} ( u ) \otimes {\rm id} ) ( ( {\rm id} \otimes \Delta ) \Delta ( x ) )\\
              &=& ( \widetilde{\Omega} ( v ) \otimes \widetilde{\Omega} ( u ) \otimes {\rm id} ) ( ( \Delta \otimes {\rm id} ) \Delta ( x ) )\\
              &=& ( \widetilde{\Omega} ( v ) \otimes \widetilde{\Omega} ( u ) ) \Delta ( x_{(1)} ) \otimes x_{(2)}\\
              &=& ( ev_{\eta( v )} \otimes ev_{\eta(u)}  ) ( \pi_{\widetilde{\clq}} \otimes \pi_{\widetilde{\clq}} ) \Delta ( x_{(1)} ) \otimes x_{(2)}, \eean
              which by ( \ref{deformation_general_pi_delta_0} ) equals              
              \bean \lefteqn{( ev_{\eta( v )} \otimes ev_{\eta(u)}  ) \Delta_{\IT^{n}} \pi_{\widetilde{\clq}} ( x_{(1)} ) \otimes x_{(2)} }\\
              &=&  ( ev_{\eta( v )}  \diamond ev_{\eta(u)}  ) ( \pi_{\widetilde{\clq}} ( x_{(1)} ) \otimes x_{(2)} )\\
              &=&  ( ev_{\eta( u )}  \diamond ev_{\eta(v)}  ) ( \pi_{\widetilde{\clq}} ( x_{(1)} ) \otimes x_{(2)} )\\
              &=&  ( \widetilde{\Omega} ( u ) \otimes \widetilde{\Omega} ( v ) ) \Delta ( x_{(1)} ) \otimes x_{(2)} \\
              &=&  ( \widetilde{\Omega} ( u ) \diamond \widetilde{\Omega} ( v ) ) ( x_{(1)} ) \otimes x_{(2)}\\
              &=&  ( ( \widetilde{\Omega} ( u ) \diamond \widetilde{\Omega} ( v ) ) \otimes {\rm id} ) \Delta ( x ) \\
              &=&  ( \widetilde{\Omega} ( u ) \diamond \widetilde{\Omega} ( v ) ) \triangleleft x. \eean
              The second part follows from this and using ( \ref{deformation_widetilde_omega=omega_on_Q} ). \qed

    \blmma
     
    \label{qiso_deformation_Lemma1}
       
  The map $\odot$  satisfies 
   $$ \int_{\IR^{2n}} ( \widetilde{\Omega}( J u ) \triangleleft  x ) \odot ( \widetilde{\Omega}( v ) \triangleleft  y )e(u.v ) du dv = \int_{\IR^{2n}} (x \triangleright  ( \widetilde{\Omega}( J u )) )( y \triangleright  \widetilde{\Omega}( v ) )e( u.v ) du dv ,$$
   for $x,y $ in $ \widetilde{\clq_0}.$ When $ x,y $ are in $ \clq_0, $ we have
   $$ \int_{\IR^{2n}} ( \Omega( J u ) \triangleleft  x ) \odot ( \Omega( v ) \triangleleft  y )e(u.v ) du dv = \int_{\IR^{2n}} (x \triangleright  ( \Omega( J u )) )( y \triangleright  \Omega( v ) )e( u.v ) du dv .$$
   
   \elmma
    {\it Proof :}  The expression in the left hand side equals
 \bean \lefteqn{\int(    \widetilde{\Omega}(Ju^{\prime}) \triangleleft  x ) \odot (\widetilde{\Omega}( v^{{\prime}} ) \triangleleft  y )e( u^{\prime}.v^{\prime} ) du^{\prime}dv^{\prime}}\\
 & =&  \int_{\IR^{2n}}\{ \int_{\IR^{4n}} e(-u.v )e(w.s )(\widetilde{\Omega}(-Ju )  \triangleleft  (\widetilde{\Omega}( Ju^{\prime}) \triangleleft  x ) \triangleright \widetilde{\Omega}( Jw )) \\
 && ( \widetilde{\Omega}( -v ) \triangleleft (\widetilde{\Omega}( v^{\prime} ) \triangleleft  y ) \triangleright \widetilde{\Omega}( s ) ) ~    dudvdwds \}e(u^{\prime}.v^{\prime})du^{\prime}dv^{\prime} \\
  &  =&  \int_{\IR^{6n}} ( \widetilde{\Omega}( J( u^{\prime}-u ) ) \triangleleft  x ) \triangleright  \widetilde{\Omega}( Jw )) (\widetilde{\Omega}( v^{{\prime}}- v )  \triangleleft  y  \triangleright  \widetilde{\Omega}( s ) ) \\
  && e( u^{{\prime}}.v^{{\prime}}) e(-u.v )e ( w.s ) du dv dw ds du^{{\prime}} dv^{{\prime}}\\
  &  = & \int_{\IR^{2n}} e( w.s )dw ds \{ \int_{\IR^{4n}} e ( u^{{\prime}}.v^{{\prime}} )e ( -u.v ) du dv du^{{\prime}}dv^{{\prime}} \\
  && ( \widetilde{\Omega}( J ( u^{{\prime}} -u )) \triangleleft  x_w )(  \widetilde{\Omega}( v^{{\prime}} - v ) \triangleleft  y_s ) \},\eean
       where $ x_w = x \triangleright \widetilde{\Omega}( Jw ) , y_s =y \triangleright \widetilde{\Omega}( s ) $.

  The proof of the lemma will be complete if we show 
  $$ \int_{\IR^{4n}} e ( u^{{\prime}}.v^{{\prime} } )e ( -u.v ) (\widetilde{\Omega}( J( u^{{\prime}} -u )) \triangleleft  x_w ) ( \widetilde{\Omega}( v^{{\prime} } -v ) \triangleleft  y_s ) du dv du^{{\prime}} dv^{{\prime}}  = x_w.y_s. $$
  
  By changing variable in the above integral, with  $z = u^{{\prime}} - u,  
       t = v^{\prime} - v,$ it becomes 
     
  $ \int_{\IR^{4n}} e( -u.v ) e ( ( u + z ).( v + t ) ) \phi ( z,t ) du dv dz dt $
  
 $ = \int_{\IR^{4n}} \phi ( z,t ) e ( u.t + z.v ) e ( z.t ) du dv dz dt, $
 where  $$ \phi ( z,t ) = (\widetilde{\Omega}( J ( z ) ) \triangleleft  x_w )( \widetilde{\Omega}( t ) \triangleleft  y_s ) .$$
  By taking 
  $ ( z,t ) = X, ( v,u ) = Y,$
      and 
   $     F ( X ) = \phi ( z,t ) e( z.t )$,  the integral can be written as 
       \bean \lefteqn{ \int\int F ( X ) e ( X.Y )dX dY }\\
                    & = & F ( 0 )   ~({\rm ~by~ Proposition ~ \ref{preliminiaries_rieffel_1.12}})\\
                     &  =& ( \widetilde{\Omega}( J ( 0 ) ) \triangleleft  x_w ) ( \widetilde{\Omega}( 0 ) \triangleleft  y_s ) \\
                                        & =& x_w.y_s, \eean since 
                       $$  \widetilde{\Omega}( J(0) ) \triangleleft  x_{w} =  (ev_{\eta(0)} \pi_{\widetilde{\clq}} \otimes {\rm id} ) \Delta (x_{w})   =  (\epsilon_{\widetilde{\IT^n}} \circ \pi_{\widetilde{\clq}} \otimes {\rm id} ) \Delta( x_{w} ) = ( \epsilon \otimes {\rm id} ) \Delta( x_{w} )= x_{w} $$ and similarly  $ \widetilde{\Omega}(0) \triangleleft  y_{s} =y_s,$ where $ \epsilon_{\widetilde{\IT^{n}}} $ denotes the counit of the quantum group   $ C( \widetilde{\IT^{n}} ) .$ 
                                        
                      This proves the claim and hence the first part of the Lemma. The second part follows from this and ( \ref{deformation_widetilde_omega=omega_on_Q} ).\qed
   
   \blmma
   
  \label{qiso_deformation_Lemma2}
  
   We have for $ a $ in $ {\cla_0},~ s $ in $ \cls_0,$
   
   \be \label{qiso_deformation_Lemma2.1} \widetilde{\alpha} ( \widetilde{\beta}_{u} ( s ) ) = s_{(1)} \otimes ( {\rm id} \otimes \widetilde{\Omega} ( u )   ) ( \Delta ( s_{(2)} ) ), \ee
    
   \be \label{qiso_deformation_Lemma2.2} \alpha ( \beta _u ( a ) ) = a_{(1)} \otimes ( {\rm id} \otimes \Omega( u ) ) ( \Delta ( a_{( 2 )} ) ). \ee
   
   \elmma
   
    {\it Proof :} $ \widetilde{\beta}_{u} = ( {\rm id} \otimes ev_u \circ \widetilde{\pi}  ) \widetilde{\alpha}.$
   We have \bean \lefteqn{   \widetilde{\beta} _u ( s )}\\
   & =&  ( {\rm id} \otimes \widetilde{\Omega} ( u ) ) \widetilde{\alpha} ( s )\\
                     & =&  ( {\rm id} \otimes \widetilde{\Omega} ( u ) ) ( s_{( 1 )} \otimes s_{( 2 )} )\\
                    &= & s_{( 1 )} ( \widetilde{\Omega} ( u )) ( s_{( 2 )} ).\eean
                   
    This gives, \bean \lefteqn{ \widetilde{\alpha} ( \widetilde{\beta}_{ u } ( s ) )}\\
           &=&  \widetilde{\alpha} ( s_{ ( 1 ) } ) \widetilde{\Omega} ( u ) ( s_{ ( 2 ) } )\\
       & = & ( {\rm id} \otimes {\rm id} \otimes \widetilde{\Omega} ( u ) ) ( \widetilde{\alpha} ( s_{(1)} ) \otimes s_{ (2) }  )\\
       &= & ( {\rm id} \otimes {\rm id} \otimes \widetilde{\Omega} ( u ) ) ( ( \widetilde{\alpha} \otimes {\rm id} ) \widetilde{\alpha} ( s ) )\\
       &= & ( {\rm id} \otimes {\rm id} \otimes \widetilde{\Omega} ( u ) ) ( ( {\rm id} \otimes \Delta ) \widetilde{\alpha} ( s ) )\\
       &= & s_{(1)} \otimes ( {\rm id} \otimes \widetilde{\Omega} ( u ) ) \Delta ( s_{(2)}).\eean 
       
    Proceeding in a similar way, we obtain $ \beta_u (a) = a_{(1)} (\Omega(u) ) (a_{(2)}) $ for all $ a $ in $ \cla_0 $ and hence $ \alpha ( \beta _u ( a ) ) = a_{(1)} \otimes ( {\rm id} \otimes \Omega( u ) ) ( \Delta ( a_{( 2 )} ) ) $ for all $ a $ in $ \cla_0 .$  
        \qed
        
      \blmma
      
       \label{qiso_deformation_Lemma3}
       
    For all $ s $ in $ \cls_0, a $ in $ \cla_0, $ we have    
   \be  \label{qiso_deformation_Lemma3.1} \widetilde{\alpha} ( a \times_J s ) ) = a_{(1)} s_{(1)} \otimes (  \int \int ( a_{(2)} \triangleright  \Omega ( Ju ) ) (s_{(2)} \triangleright  \widetilde{\Omega} ( v ) ) e ( u.v ) du dv ). \ee   
       For $a,b $ in $ {\cla_0}$, we have      
     \be  \label{qiso_deformation_Lemma3.2} \alpha ( a \times _{J} b ) = a_{(1)}b_{(1)} \otimes ( \int\int ( a_{(2)} \triangleright \Omega( J u ) )  ( b_{(2)} \triangleright  \Omega( v )  ) e ( u.v ) du dv  ). \ee
     
     \elmma
     {\it Proof :} Using the notations and definitions in section \ref{preliminaries_section_Rieffel_deformation}, we note that for any $ f: \IR^{2} \rightarrow \IC $ belonging to  $ \IB( \IR^{2} ) $ and fixed $ x $ in $ E $( where $ E $ is a Banach algebra ), the function $ F( u, v ) = xf(u,v )$ belongs to  $ \IB^{E} ( \IR^{2} ) $ and we have  \bean \lefteqn{ x ~( \int \int f( u,v ) e( u.v ) du dv )}\\
      & =&  x ~( \lim_L \sum_{p \in L} \int \int ( f \phi_{p})( u,v ) e( u.v ) du dv )\\
    &=& \lim_L \sum_{p \in L} \int \int x~( f \phi_{p})( u,v ) e( u.v ) du dv )\\
      &= & \int\int x~f(u,v ) e(u.v) du dv.\eean
      Then,
        \bean
     \lefteqn{\widetilde{\alpha} ( a \times _{J} s )}\\
      &=&  \widetilde{\alpha} ( \int\int \beta_{Ju} ( a ) \widetilde{\beta} _{v} ( s ) e( u.v ) du dv )\\
      &=& \widetilde{\alpha} ( \int\int a_{(1)} ( \Omega(Ju) )( a_{(2)} ) s_{(1)} ( \widetilde{\Omega}(v) )(s_{(2)}) e(u.v) du dv ) \\
      &=& \widetilde{\alpha} ( (a_{(1)}  s_{(1)} )  \int\int  ( \Omega(Ju) )( a_{(2)} )  ( \widetilde{\Omega}(v) )(s_{(2)}) e(u.v) du dv ) \\
      &=& \alpha( a_{(1)} ) \widetilde{\alpha}( s_{(1)} ) \int\int  ( \Omega(Ju) )( a_{(2)} )  ( \widetilde{\Omega}(v) )(s_{(2)}) e(u.v) du dv \\
      && ( ~ {\rm by} ~ {\rm assumption} ~ 2.c ~ )\\      
     &=&  \int \int \alpha( a_{(1)} ) ( \Omega(Ju) )(a_{(2)}) \widetilde{\alpha}( s_{(1)} ) ( \widetilde{\Omega}(v) )( s_{(2)} )   e(u.v) du dv  \\
     &=& \int \int \alpha ( a_{(1)} \Omega ( Ju ) ( a_{(2)}  ) ) \widetilde{\alpha} ( s_{(1)} \widetilde{\Omega} ( v ) ( s_{(2)} ) ) e ( u.v ) du dv \\
     &=& \int\int \alpha ( \beta_{Ju}(a) ) \widetilde{\alpha} ( \widetilde{\beta} _{v} ( s ) ) e( u.v ) du dv \\
      &= & \int\int ( a_{(1)} \otimes ( {\rm id} \otimes \Omega(Ju) ) ( \Delta ( a_{( 2 )} )) ) ( s_{(1)} \otimes ( {\rm id} \otimes \widetilde{\Omega} ( v ) ) ) ( \Delta ( s_{( 2 )}) ) \\
      && e ( u.v ) du dv\\& & ( {\rm ~using ~Lemma ~\ref{qiso_deformation_Lemma2} )}\\
     &=&  a_{(1)}s_{(1)} \otimes \int\int ( a _{(2)} \triangleright  \Omega( Ju ) ) ( s_{(2)} \triangleright  \widetilde{\Omega}( v ) ) e ( u.v ) du dv .\eean
     \qed  
        
\blmma 
    
    \label{qiso_deformation_Lemma4}
    
  For $ s $ in $ \cls_0, ~ a $ in $ \cla_0,$
  
  \be \label{qiso_deformation_Lemma4.1} \alpha ( a ) \bullet_J \widetilde{\alpha} ( s ) = a_{(1)} s_{(1)} \otimes ( \int \int ( \Omega ( Ju ) \triangleleft  a_{(2)} ) \odot  ( \widetilde{\Omega} ( v ) \triangleleft  s_{(2)} ) e ( u.v ) du dv.  \ee
   For $a,b $ in $ {\cla_0},$         
    \be \label{qiso_deformation_Lemma4.2}   \alpha ( a ) \bullet_J \alpha ( b ) = a_{(1)}b_{(1)} \otimes \{ \int\int ( \Omega( Ju )  \triangleleft  a_{(2)} ) \odot (  \Omega( v ) \triangleleft   b_{(2)}) e( u.v ) du dv \}. \ee
    
     \elmma
     {\it Proof :} We have 
     \bean
     \lefteqn{\alpha ( a )\bullet_J \widetilde{\alpha} ( s )}\\
         & = & ( a_{(1)} \otimes a_{(2)} ) \bullet_{J} ( s_{(1)} \otimes s_{(2)} )\\
        &=&  a_{(1)} \times _{J} s_{(1)} \otimes ( a_{(2)} \odot s_{(2)} )\\
        &=&  \int\int \beta_{Ju} ( a_{(1)} ) \widetilde{\beta} _{v} ( s_{(1)} ) e( u.v ) du dv \otimes ( a_{(2)} \odot s_{(2)} ). \eean
      Let $\epsilon $ be the counit of $ \widetilde{\clq} .$
     So we have $ ( {\rm id} \otimes \epsilon ) \alpha = {\rm id} $ and $ ( {\rm id} \otimes \epsilon ) \widetilde{\alpha} = {\rm id}.$ This gives, 
      \bean
      \lefteqn{\alpha ( a )\bullet_J \widetilde{\alpha} ( s )}\\
     &     =&  \int\int ( {\rm id} \otimes \epsilon ) \alpha ( \beta_{Ju} ( a_{(1)} ) ) ( {\rm id} \otimes \epsilon ) \widetilde{\alpha} ( \widetilde{\beta} _{v} ( s_{(1)} ) e( u.v ) du dv \otimes ( a_{(2)} \odot s_{2} ).\eean
     Note that by Lemma \ref{qiso_deformation_Lemma2}, 
 $ \int\int ( {\rm id} \otimes \epsilon ) ( \alpha ( \beta_{Ju} ( a_{(1)} ) ) ( {\rm id} \otimes \epsilon ) ( \widetilde{\alpha} ( \widetilde{\beta} _{v} ( s_{(1)} ) ) )$ 
 
 $e(u.v ) du dv $
     
   $  = \int\int ( {\rm id} \otimes \epsilon ) ( a_{(1)(1)} \otimes ( {\rm id} \otimes \Omega ( Ju) ) ( \Delta (a_{(1)(2)})) ) ({\rm id} \otimes \epsilon ) ( s_{(1)(1)} \otimes ( {\rm id} \otimes $
   
   $ \widetilde{\Omega} ( v ) ) ( \Delta ( s_{(1)(2)} )) ) e( u.v ) du dv $

   $  = \int\int ( {\rm id} \otimes \epsilon ) ( a_{(1)(1)} \otimes ( a_{(1)(2)} \triangleright  \Omega( Ju )  ) ( {\rm id} \otimes \epsilon ) ( s_{(1)(1)} \otimes ( s_{(1)(2)} \triangleright  \Omega( v )  ) ) e( u.v ) $
   
   $ du dv $
     
   $  = \int\int a_{(1)(1)} s_{(1)(1)} \epsilon (  a_{(1)(2)} \triangleright \Omega( Ju )  ) \epsilon ( s_{(1)(2)} \triangleright \widetilde{\Omega}( v )  ) e( u.v ) du dv. $

     Using the fact that  $f \diamond \epsilon=\epsilon \diamond f=f$ for any functional on $\clq_0$, one has 
       $ \epsilon ( a_{(1)(2)} \triangleright  \Omega( Ju )  )
                =   \Omega( Ju )   ( a_{(1)(2)} )$ and 
    $ \epsilon ( s_{(1)(2)} \triangleright  \widetilde{\Omega}( v ) )
                        =  \widetilde{\Omega}( v )  ( s_{(1)(2)} )  $, from which it follows that  
                 
        \bean \lefteqn{ \alpha ( a )\bullet_J \widetilde{\alpha} ( s )}\\
        &=&  a_{(1)(1)}s_{(1)(1)} \int\int \Omega ( Ju ) ( a_{(1)(2)} ) \widetilde{\Omega} ( v ) ( s_{(1)(2)} ) e(u.v ) du dv                    \otimes ( a_{(2)} \odot s_{(2)} )\\
       &=&  \int\int ( {\rm id} \otimes \Omega ( Ju ) \otimes {\rm id} ) ( a_{(1)(1)} \otimes a_{(1)(2)} \otimes a_{(2)} )            \bullet ( {\rm id} \otimes \widetilde{\Omega} ( v ) \otimes {\rm id} ) ( s_{(1)(1)} \otimes \\ 
       & & s_{(1)(2)} \otimes s_{(2)} ) e( u.v ) du dv\eean
  \bean \lefteqn{=  \int\int  ( {\rm id} \otimes \Omega ( Ju ) \otimes {\rm id} ) ( a_{(1)} \otimes \Delta ( a_{(2)} ) ) \bullet ( {\rm id} \otimes ( \widetilde{\Omega} ( v ) \otimes {\rm id} ) ) ( s_{(1)} \otimes \Delta ( s_{(2)} ) )}\\  &&e( u.v ) du dv \\
   &   =&  \int\int \{  a_{(1)} \otimes ( \Omega ( Ju ) \otimes {\rm id} ) \Delta ( a_{(2)} ) \} \bullet \{  s_{(1)} \otimes ( \widetilde{\Omega} ( v ) \otimes {\rm id} ) \Delta ( s_{2} ) \} e( u.v ) du dv\\
      &=&  a_{(1)}s_{(1)} \otimes \int\int ( ( \widetilde{\Omega} ( Ju ) \otimes {\rm id} ) \Delta ( a_{2} ) ) \odot ( \widetilde{\Omega} ( v ) \otimes {\rm id} ) ) \Delta ( s_{(2)} ) e( u.v ) du dv\\
      && ( ~ by ~ ( \ref{deformation_widetilde_omega=omega_on_Q}   ) )\\
     &=& a_{(1)}s_{(1)} \otimes \int\int (  \widetilde{\Omega}( Ju ) \triangleleft  a_{(2)} ) \odot (  \widetilde{\Omega}( v )   \triangleleft  s_{(2)} ) e( u.v ) du dv,\eean
              where we have used the relation $(\alpha \ot {\rm id}) \alpha=({\rm id} \ot \Delta) \alpha$ to get  $a_{(1)(1)} \otimes a_{(1)(2)} \otimes a_{(2)}
   =a_{(1)} \otimes \Delta ( a_{(2)} )$ 
                            and similarly $ s_{(1)(1)} \otimes s_{(1)(2)} \otimes s_{(2)} = s_{(1)} \otimes \Delta ( s_{(2)} )   .$ \qed 
                            
   \vspace{4mm}

 Combining Lemma \ref{qiso_deformation_Lemma1}, Lemma \ref{qiso_deformation_Lemma3} and Lemma \ref{qiso_deformation_Lemma4}, we conclude the following.  
     
   \blmma 
             
   \label{qiso_deformation_Lemma5}
         
    For $ a $ in $ \cla_0, ~ s $ in $ \cls_0,$ we have
    \be \label{qiso_deformation_Lemma5.1} \alpha(a) \bullet_J \widetilde{\alpha}(s)= \widetilde{\alpha}(a \times_J s) . \ee    
For $a,b $ in $ {\cla_0},$ we have   
   \be  \label{qiso_deformation_Lemma5.2} \alpha(a) \bullet_J \alpha(b)=\alpha(a \times_J b) . \ee

  \elmma

  We shall now identify $\odot$ with the multiplication of a Rieffel-type deformation of $\clq$( $ \widetilde{\clq} $ ). We discuss the case for $ \widetilde{\clq} ,$ that of $ \clq $ being similar. Since $\widetilde{\clq}$ has a quantum subgroup isomorphic with $C ( \widetilde{\IT^n} ),$ we can consider the following canonical action $\chi$ of $\IR^{2n}$ on $\widetilde{\clq}$ ( as in ( \ref{preliminiaries_deformation_wang_formula} ) ) given by          
      $$  \chi_{( s,u )} = ( \widetilde{\Omega}( -s ) \otimes {\rm id} ) \Delta ( {\rm id} \otimes \widetilde{\Omega} ( u ) ) \Delta.$$  
          Now, let 
       $ \widetilde J := -J \oplus J $, which is a skew-symmetric ${2n} \times 2n$ real matrix, so one can deform $\widetilde{\clq}$ by defining the product of $x$ and $y$ ($x,y $ belonging to $ \widetilde{\clq_0},$ say) to be the following: $$ \int\int \chi_{\widetilde J ( u,w )}( x ) \chi_{v,s} (  y ) e ( ( u,w ).( v,s ) ) d ( u,w ) d ( v,s ).$$ We claim that this is nothing but $\odot$ introduced before.
         
         \blmma
         
         \label{qiso_deformation_Lemma6}
              
         \bean  x \odot y = x \times_{\widetilde{J}} y ~ {\rm ~ for ~ all} ~  x,y \in \widetilde{\clq_{0}}. \eean
           
        \elmma 
             
          {\it Proof :} Let us first observe that 
          \bean
          \lefteqn{ \chi _{ \widetilde J ( u,w ) } ( x )}\\
           &=& ( \widetilde{\Omega} ( Ju ) \otimes {\rm id} ) \Delta ( {\rm id} \otimes \widetilde{\Omega} ( Jw ) ) \Delta ( x )\\
          &=& \widetilde{\Omega}( Ju ) \triangleleft  x \triangleright \widetilde{\Omega}( Jw ), \eean
          and similarly
           $ \chi_{(v,s)}(y)= \widetilde{\Omega}( -v ) \triangleleft  y \triangleright  \widetilde{\Omega}( s ).$ 
                     
     Thus, we have                
     \bean     
          \lefteqn{ x \odot y}\\
          & =&  \int_{\IR^{4n}} (  \widetilde{\Omega}( -Ju )  \triangleleft x \triangleright  \widetilde{\Omega}( Jw ) ) (  \widetilde{\Omega}( -v ) \triangleleft  y \triangleright \widetilde{\Omega}( s ) ) e ( -u.v ) e ( w.s ) du dv dw ds\\
      &=&  \int_{\IR^{4n}}  ( \widetilde{\Omega}( J u^{{\prime}} ) \triangleleft  x \triangleright \widetilde{\Omega}( J w ) ) ( \widetilde{\Omega}( -v )               \triangleleft  y  \triangleright \widetilde{\Omega}( s )) e ( u^{{\prime}}.v ) e ( w.s ) du^{{\prime}} dv dw ds\\
      &=& 
     \int_{\IR^{2n} } \int_{\IR^{2n}} \chi _{\widetilde J ( u,w )} ( x )\chi _{( v,s )} ( y ) e ( ( u,w ).( v,s )            ) d( u,w ) d ( v,s ),\eean which proves the claim. \qed 

\vspace{2mm}

     Let us denote by $\widetilde{\clq}_{\widetilde{J}}$ ( $ \clq_{\widetilde{J}} $ ) the $C^*$ algebra obtained from $\widetilde{\clq}$ ( $ \clq $ ) by the  Rieffel deformation w.r.t. the matrix $\widetilde{J}$ described above. We recall from  subsection \ref{preliminaries_subsection_Rieffel_deformation_CQG}  that the coproduct $\Delta$ on $\widetilde{\clq_0}$ ( $ \clq_0 $ ) extends to a coproduct for the deformed algebra  as well and $ (\widetilde{\clq}_{\widetilde{J}}, \Delta) $ ( $(\clq_{\widetilde{J}}, \Delta)$ ) is a compact quantum group.

   \vspace{4mm}

 \section{ $ \widetilde{QISO^{+}_{R}} $ of a Rieffel deformed noncommutative manifold}
 
 \subsection{Derivation of the result}
 
 \label{qorient_subsection_Derivation of the result}
 
In this subsection, our set up is as in section \ref{deformation_section_deformation_of_spectral_triple} so that we have spectral triples on $ \cla^{\infty}_{J} $ for each $ J. $
 
 \blmma
   
   \label{1111}
   
   Suppose that $ ( \widetilde{\clq}, U ) $ belongs to $ {\rm Obj} ( {\bf Q}( \cla, \clh, D ) ) , $ and  
       there exists a unital $ \ast $-subalgebra $ \cla_{0} \subseteq \cla $ which is norm dense in every $ \cla_{J}$ such that 

      $ \alpha_{U} ( \pi_{0} ( \cla_{0} )  ) \subseteq \pi_{0} ( \cla_{0} ) \otimes_{alg} \clq_{0}, $ where $\clq \subseteq \widetilde{\clq} $ is the smallest Woronowicz $ C^{*} $ subalgebra such that $ \alpha_{U}( \cla_{0} ) \subseteq \pi_{0} ( \cla_{0} ) \otimes \clq,$ and $ \clq_{0} $ is the Hopf $\ast$-algebra obtained by matrix coefficients of irreducible unitary (co)-representations of $\clq.$ Also, let $ S_0 = {\rm span} \{ a s : a \in \cla_0, s \in S_{00} \} ,$
Then we have the following:

(a) $U(\cls_0) \subseteq \cls_0 \ot_{\rm alg} \widetilde{\clq}_0$.\\
(b) $\widetilde{\alpha}:=U|_{\cls_0}: \cls_0 \raro \cls_0 \ot_{\rm alg} \widetilde{\clq}_0$  makes $\cls_0$ an algebraic $\widetilde{\clq}_0$ co-module, satisfying $$ \widetilde{\alpha}(\pi_0(a) s)=\alpha_{U}(a) \widetilde{\alpha}(s) ~{\rm for ~ all} ~   a ~ {\rm ~ in} ~ \cla_0, ~ s ~ {\rm in} ~ \cls_0.$$ 

Moreover, if  $ ( C ( \widetilde{\IT^n} ), V ) $ is a sub object of $ \widetilde{\clq} $ in $ {\bf Q}( \cla, \clh, D ) $ such that $ V ( . \otimes {\rm id} ) V^* = \beta,$ then $ C ( \IT^n ) $ is a quantum subgroup of $ \clq .$ 
\elmma
 {\it Proof:}  $ U $ commutes with $ D $ and hence preserves the eigenspaces of $ D $ which shows that $ U $ preserves $ S_{00}.$ Then, $ \widetilde{U} ( a s \otimes 1 ) = \alpha ( a ) \widetilde{U} ( s \otimes 1 ) \subseteq ( \cla_0 \otimes \clq_0 ) ( S_{00} \otimes \widetilde{\clq_0} ) \subseteq S_0 \otimes \widetilde{\clq_0}.$ Thus, the first assertion follows.
 
  The second assertion follows from the definition of $ \widetilde{\alpha} $ and $ \alpha_{U}.$ The third assertion follows as in Lemma \ref{deformation_observation_from_assumptions}.\qed

\brmrk

 From the definitions of $ \cla_0 $ and $ S_0 ,$ it follows that\\
   (i)  $ \pi_{0} ( \cla_{0} ) \cls_{0} \subseteq \cls_{0} $,\\
   (ii) $ \beta_{g}( \cla_{0} ) \subseteq \cla_{0} $  for all $ g.$\\
   
 \ermrk  

\vspace{4mm}

Let us now fix the object $(\widetilde{\clq}, U)$ as in the statement of Lemma \ref{1111}. We recall that using the identification of $\clq_0$ as a common vector-subspace of all $\clq_{\widetilde{J}},$  we shall sometimes denote this identification map from $\clq_0 $ to $\clq_{\widetilde{J}}$ by $\rho_J$.

\vspace{4mm}

Let us consider the finite dimensional unitary representations $U^{(i)}:=U|_{V_i}$, where $V_i$ is the eigenspace of $D$ corresponding to the eigenvalue $\lambda_i$. By Corollary \ref{repdef}, we can view $U^{(i)}$ as a unitary representation of $\clq_{\widetilde{J}}$ as well, and let us denote it by $U^{(i)}_J$. In this way, we obtain a unitary representation $U_J$ on the Hilbert space $\clh$, which is the closed linear span of all the $V_i$'s. It is obvious from the construction (and the fact that the linear span of $V_i$'s, that is $\cls_0$, is a core for $D$) that $U_J D=(D \ot I)U_J$.  Let $\alpha_J:=\alpha_{U_J}$.  
With this, we have the following:

  \blmma
For $ a $ in $ \cla_0,$ we have $\alpha_J(a)=(\alpha(a))_J \equiv (\pi_J \ot \rho_J)(\alpha(a))$, hence in particular, for every state $\phi$ on $\clq_{\widetilde{J}}$, $({\rm id} \ot \phi) \circ \alpha_J(\cla_J) \subseteq \cla_J^{\prime \prime}$.

  \elmma  
  
  Using the equation ( \ref{qiso_deformation_Lemma5.1} ), we have  for all $ s  $ in $ \cls_0, a $ in $ \cla_0,$
   \bean
  \lefteqn{\alpha_{J} ( a ) U_{J} ( s )}\\
       &=& {U_J} ( \pi_J(a)s)\\
       &=& \widetilde{\alpha}(a \times_J s)\\
&=& \alpha(a) \bullet_{J} \widetilde{\alpha}(s)\\
&=& ( \alpha(a) )_J U_J(s),
\eean
 from which we conclude by the density of $\cls_0$ in $\clh$ that $\alpha_J(a) =(\alpha(a))_J $ belongs to $ \pi_J(\cla_0) \ot \clq_{\widetilde{J}}$. The lemma now follows  using the norm-density of $\cla_0$ in $\cla_J$. 
\qed 

\bcrlre

\label{deformation_qorient_subobject}
   
 $(\widetilde{\clq}_{\widetilde{J}}, U_J)$ is an orientation preserving isometric action on the spectral triple $(\cla^\infty_J, \clh, D)$.
 
 \ecrlre
    
    We shall now show that if we fix a `volume-form' in terms of an $R$-twisted structure, then the `deformed' action $\alpha_J$ preserves it.
    
  \blmma
       \label{volumepreserving}
    Suppose, in addition to the set-up already assumed, that  there is an invertible positive operator $R$ on $\clh$ such that $(\cla^\infty, \clh, D, R)$ is an $R$-twisted $ \Theta $-summable spectral triple, and let $\tau_R$ be the corresponding `volume form'. Assume that  $\alpha_U$  preserves 
 the functional $\tau_R$. Then the action $\alpha_{U_J}$  preserves $\tau_R$ too.

    \elmma
      
  {\it Proof :} Let the (finite dimensional) eigenspace corresponding to the eigenvalue $ \lambda_n $ of $D$ be $ V_n .$   As $ U $ commutes with $ D ,$ there exists subspaces 
$ V_{n,k} $ of $ V_n $ and an orthonormal basis $ {\{e^{n,k}_{j}\}}_{j} $ for $ V_{n,k} $  such that the restriction of $U$ to $V_{n,k}$  is irreducible. Write $\widetilde{U}( e^{n,k}_{j} \otimes 1 ) = \sum_{i} e^{n,k}_{i} \otimes t^{n}_{i,j}.$ Then, $ {\widetilde{U}}^{*} ( e^{n,k}_{j} ) = \sum_{i} e^{n,k}_{i} \otimes t^{n *}_{j,i} .$
             
 Then $ \clh $ will be decomposed as $ \clh = \oplus_{n \geq 1, ~ k} V_{n,k} .$
             
     Let $ R ( e^{n,i}_{j} ) = \sum_{s,t} F_n ( i,j,s,t ) e^{n,s}_{t} .$
     
    
 By hypothesis, $ \widetilde{U} ( . \otimes {\rm id} ) {\widetilde{U}}^{*} $ preserves the functional $ \tau_{R} ( \cdot) = Tr ( R~ \cdot ) $ on $ \cle_{D} $ where $ \cle_{D} $ is as in Proposition \ref{5678}, that is the weakly dense $ \ast $ subalgebra of $ \clb ( \clh ) $ generated by the rank one operators $ | \xi > < \eta | $ where $ \xi, \eta $ are eigenvectors of $ D .$      
   Thus, $( \tau_{R} \otimes {\rm id} )( \widetilde{U} ( X \otimes {\rm id} ){\widetilde{U}}^{\ast} ) = \tau_{R}( X ).1_{Q} $ for all $ X $ in $ \cle_{D}.$

   Then, for $a $ in $ \cle_D,$ we have: \bean
         \lefteqn{ ( \tau_{R} \otimes h )( \widetilde{U}_{J} ( a \otimes 1 ) {\widetilde{U}_{J}}^{*} ) } \\ 
         &=& \sum_{n,i,j} \left\langle e^{n,i}_{j} \otimes 1,\widetilde{U}_{J} ( a \otimes 1 ) {\widetilde{U}_{J}}^{*} ( R   e^{n,i}_{j} \otimes 1 ) \right\rangle  \\
         & = &  \sum_{n,i,j,s,t}  \left\langle  {\widetilde{U}_{J}}^{*} ( e^{n,i}_{j} \otimes 1 ), ( a \otimes 1 ){\widetilde{U}_{J}}^{*} (  F_{n}( i,j,s,t ) e^{n,s}_{t} \otimes 1 ) \right\rangle \\
         & = & \sum_{n,i,j,s,t,k,l} F_{n}(i,j,s,t ) \left\langle e^{n,i}_{k} \otimes {(t^{n}_{j,k})}^{*}, ( a \otimes 1 )( e^{n,s}_{l}  \otimes {(t^{n}_{t,l})}^{*} \right\rangle  \\
          &=& \sum_{n,i,j,s,t,k,l} F_{n}( i,j,s,t ) \left\langle e^{n,i}_{k} , a e^{n,s}_{l} \right\rangle h_{J} ( ( t^{n}_{j,k}) \times_{J} {(t^{n}_{t,l})}^{*}) \\
        & = &  \sum_{n,i,j,s,t,k,l} F_{n}( i,j,s,t ) \left\langle e^{n,i}_{k} , a e^{n,s}_{l} \right\rangle h_{0} ( t^{n}_{j,k} {t^{n}_{t,l}}^{*})  \\
   &=& ( \tau_{R} \otimes h )( \widetilde{U} ( a \otimes 1 ){ \widetilde{U} }^{*} ) \\ 
   &=&  \tau_{R} ( a ).1 \eean 
  
  where $ h_{J} ( ( t^{n}_{j,k}) \times_{J} {(t^{n}_{t,l})}^{*}) = h_{0} ( t^{n}_{j,k} {t^{n}_{t,l}}^{*}) $ as  deduced by using Lemma \ref{deformationLemma7}.
  
  Thus $ ( \tau_R \otimes h ) \widetilde{U_{J}} ( a \otimes {\rm id} ) { \widetilde{U_{J}} }^{*}  = \tau_R ( a ).1 $
  
  Let $ ( \tau_R \otimes h ) \widetilde{U_{J}} ( X \otimes {\rm id} ) { \widetilde{U_{J}} }^{*}  = ( \tau_R \ast h ) ( X ) .$  As $ \widetilde{U_{J}} ( \cdot  \otimes {\rm id} ) { \widetilde{U_{J}} }^{*} $ keeps $ \cle_{D} $ invariant, we can use Sweedler notation: $\widetilde{U_J}(a \ot 1) \widetilde{U_J}^*=
 a_{(1)} \ot a_{(2)}$, with $a, a_{(1)} $ belonging to $ \cle_D,~a_{(2)} $ belonging to $ \widetilde{\clq_{\tilde{J}}},$   to have
 
 \bean \lefteqn{( \tau_{R} \ast h \otimes {\rm id} )\widetilde{U_{J}} ( a \otimes {\rm id} ) { \widetilde{U_{J}} }^{*}  }\\
 &=& ( \tau_R \otimes h ) ( \widetilde{U_{J}} ( a_{(1)} \otimes {\rm id} ) { \widetilde{U_{J}} }^{*} ) \otimes a_{(2)}\\
 &=& ( \tau_{R} \otimes h \otimes {\rm id} ) ( \alpha_U \otimes {\rm id} ) \alpha_U ( a )\\
 &=& ( \tau_R \otimes h \otimes {\rm id} ) ( {\rm id} \otimes \Delta ) \alpha_U ( a )\\
 &=& ( \tau_R \otimes {\rm id} ) ( {\rm id} \otimes ( h \otimes {\rm id} ) \Delta ) \alpha_U ( a )\\
 &=& ( \tau_R \otimes {\rm id} ) ( {\rm id} \otimes h ( . ).1 ) \alpha_U ( a )\\
 &=& ( {\rm id} \otimes h ( . ).1 ) ( \tau_R \otimes {\rm id} ) \alpha_U ( a )\\
 &=& ( {\rm id} \otimes h ( . ).1 ) \tau ( a ).1\\
 &=& \tau ( a ).1.\eean
  
 Thus, 
\bean \lefteqn{ ( \tau_{R} \otimes {\rm id} ) (\widetilde{U_{J}} ( a \otimes 1 ) { \widetilde{U_{J}} }^{*})}\\
& =&
 ( \tau_{R} \ast h \otimes {\rm  id} )(\widetilde{U_{J}} ( a \otimes 1 ) { \widetilde{U_{J}} }^{*})  = ( \tau_{R} \ast h ) ( a_{(1)})  a_{(2)}\\
& =&  ( \tau_{R} \otimes h  \ot {\rm id})  ( a_{(1)(1)} \ot a_{(1)(2)} \otimes a_{(2)}) = 
 ( \tau_{R} \otimes h  \otimes {\rm id} )( {\rm id } \otimes \Delta_{\widetilde{J}} )\\
 && (\widetilde{U_J}( a \ot 1 )\widetilde{U_J}^*)\\
 &= & \tau_{R}( a_{(1)} )  ( h \otimes {\rm  id} )\circ \Delta_{\widetilde{J}} ( a_{(2)} )
 = \tau_{R}( a_{(1)} )  h( a_{(2)} ).1_{\clq_{\widetilde{J}}} \\
&=&  ( \tau_{R} \otimes h ) ( a_{(1)} \ot a_{(2)} ) = ( \tau_{R} \ast h )( a ).1_{\clq_{\widetilde{J}}} = \tau_{R}( a ).1_{\clq_{\widetilde{J}}}.\eean

   \qed 
 
 \brmrk
   
 If  $ QISO^{+}_{R}( \cla^{\infty}, \clh, D ) $ (  $ QISO^{+}( \cla^{\infty}, \clh, D ),$ if it exists )  has a $ C^* $ action, then from the definition of a $ C^{*} $ action, we get a subalgebra $ \cla_0  $ as in Lemma \ref{1111}. Thus, the assumptions of section \ref{deformation_preparatory results} are satisfied so that the conclusions in that subsection hold for $ QISO^{+}_{R}( \cla^{\infty}, \clh, D ) $ (  $ QISO^{+}( \cla^{\infty}, \clh, D ) ) .$. Similarly, the conclusions of Lemma \ref{1111} and the subsequent Lemmas hold for $ QISO^{+}_{R}( \cla^{\infty}, \clh, D ) $ (  $ QISO^{+}( \cla^{\infty}, \clh, D ) ) .$
 
 \ermrk

  For any two compact quantum groups $ ( \cls_1, U^{\cls_1} ) $ and $ ( \cls_2, U^{\cls_2} ) $ in ${\bf Q^{\prime}}( \cla_{J}, \clh, D )$, we write $ \cls_1 < \cls_2 $ if $ \cls_1 $ is a sub object of $ \cls_2 $ in the category ${\bf Q^{\prime}} ( \cla_{J}, \clh, D ).$

       \blmma
       
       \label{deformation_qorient_<}
        
        If $G_1, G_2 $ be two CQG s such that $ G_1 < G_2 $ in the category ${\bf Q^{\prime}}( \cla_{J}, \clh, D ),$ If $ ( G_1 )_{\widetilde{J}} $ and $ ( G_2 )_{\widetilde{J}} $ make sense, then  $ ( G_1 )_{J} < ( G_2 )_{J} $ in the category ${\bf Q^{\prime}}( \cla_{J}, \clh, D ).$
         
        \elmma
        
        {\it Proof :} From Corollary \ref{deformation_qorient_subobject}, we see that $ ( G_i )_{\widetilde{J}} $
is an object in the category ${\bf Q^{\prime}}( \cla_{J}, \clh, D ).$ Let $ \pi_2 $ be the morphism from  $ G_2 $ to $ G_1 $ in the category $ {\bf \clq}^{\prime} $ and $ \pi_1 $ be the morphism from $ G_1 $ to $ \IT^{n} $ in the same category. Let $ \Delta^{i}, ~ \times^{i}_{J}, ~ \chi^{i} $ denote respectively  the coproducts, products and $ \IR^{2n} $ action  on  $ ( G_i )_{J} ,~ i = 1,2.$ 
        
        As the quantum group structure is not altered under Rieffel deformation, to prove the Lemma,  it is enough to show that $ \pi_2 $ is a homomorphism from $ ( G_2 )_{J} $ to $ ( G_1 )_{J}.$
        
       In any CQG $( Q, ~ \Delta ), ~ f,g $ linear functionals on $ Q $ and  for all $ a $ in $ Q_0,$
              $ ( f \otimes {\rm id} )\Delta ( {\rm id} \otimes g ) \Delta ( a ) = ( f \otimes {\rm id} ) \Delta ( a_{(1)} ) g ( a_{(2)} ) = ( f \otimes {\rm id} ) ( a_{(1)(1)} \otimes a_{(1)(2)} ) g ( a_{(2)} ) = f ( a_{(1)(1)} ) g ( a_{(2)} ) a_{(1)(2)} = ( f \otimes {\rm id} \otimes g ) ( a_{(1)(1)} \otimes a_{(1)(2)} \otimes a_{2} ) = ( f \otimes {\rm id} \otimes g ) ( \Delta ( a_{(1)} ) \otimes a_{(2)} ) = ( f \otimes {\rm id} \otimes g ) ( \Delta \otimes {\rm id} ) \Delta ( a ).$
              
              Hence,   
        \be \label{QISO_Laplacian_deformation_lemma} ( f \otimes {\rm id} ) \Delta ( {\rm id} \otimes g ) \Delta = ( f \otimes {\rm id} \otimes g ) ( \Delta \otimes {\rm id} ) \Delta .\ee 
         
         Moreover, we will also need the equation         
         \be \label{QISO_Laplacian_deformation_subgroup} ( \pi_2 \otimes \pi_2 ) \Delta^{2} = \Delta^{1} \pi_{2}, \ee        which holds as $ \pi_2 $ is a morphism of  CQG s, $ G_2 \rightarrow G_1.$
            
         Let $\lambda, ~ \rho $ be as in ( \ref{preliminiaries_deformation_wang_lambda} ) and ( \ref{preliminiaries_deformation_wang_rho} ) and $ a $ belongs to $ ( G_2 )_{J}.$
         
         Then, \bean \lefteqn{ \pi_{2} \chi^{2}_{( s, u )} ( a ) }\\
         &=& \pi_2 ( \lambda_{\eta( - s )} \rho_{\eta( u )} ) ( a )\\
         &=& \pi_2 ( {\rm ev}_{\eta( - s )} ( \pi_1 \circ \pi_2 ) \otimes {\rm id} ) \Delta^{2} ( {\rm id} \otimes {\rm ev}_{\eta ( u )} ( \pi_1 \circ \pi_2 ) ) \Delta^{2} ( a ) )\\    
      &=& \pi_2 ( {\rm ev}_{\eta( - s )} ( \pi_1 \circ \pi_2 ) \otimes {\rm id} \otimes {\rm ev}_{\eta ( u )} \circ ( \pi_1 \circ \pi_2 ) ) ( \Delta^{2} \otimes {\rm id} ) \Delta^{2} ( a ) \\
      && ( ~ {\rm by} ~ ( ~ \ref{QISO_Laplacian_deformation_lemma} ~ ) ~ )\\ 
         &=& ( {\rm ev}_{\eta( - s )} ( \pi_1 \circ \pi_2 ) \otimes \pi_2 \otimes {\rm ev}_{\eta ( u )} \circ ( \pi_1 \circ \pi_2 ) ) ( \Delta^{2} \otimes {\rm id} ) \Delta^{2} ( a ) \\
         &=& ( {\rm ev}_{\eta( - s )} \pi_1 \otimes {\rm id} \otimes {\rm ev}_{\eta ( u )} \pi_1 ) ( \pi_2 \otimes \pi_2 \otimes \pi_2 ) ( \Delta^{2} \otimes {\rm id} ) \Delta^{2} ( a ) \\
         &=& ( {\rm ev}_{\eta( - s )} \pi_1 \otimes {\rm id} \otimes {\rm ev}_{\eta( u )} \pi_1 ) ( ( \pi_2 \otimes \pi_2 ) \Delta^{2} \otimes \pi_2 ) \Delta^{2} ( a )\\ 
     &=& ( {\rm ev}_{\eta( - s )} \pi_1 \otimes {\rm id} \otimes {\rm ev}_{\eta( u )} \pi_1 ) ( \Delta^{1} \pi_2 \otimes \pi_2 ) \Delta^{2} ( a ) \\
         && ( ~ {\rm using} ~ ( \ref{QISO_Laplacian_deformation_subgroup} ) ~ )\\
         &=& ( {\rm ev}_{\eta( - s )} \pi_1 \otimes {\rm id} \otimes {\rm ev}_{\eta( u )} \pi_1 ) ( \Delta^{1} \otimes {\rm id} ) ( \pi_2 \otimes \pi_2 ) \Delta^{(2)} ( a )\\
         &=& ( {\rm ev}_{\eta( - s )} \pi_1 \otimes {\rm id} \otimes {\rm ev}_{\eta( u )} \pi_1 ) ( \Delta^{1} \otimes {\rm id} ) \Delta^{1} ( \pi_{2} ( a ) ) \\
         &=& ({\rm ev}_{\eta( - s )} \pi_1 \otimes {\rm id} ) \Delta^{1} ( {\rm id} \otimes {\rm ev}_{\eta( u )} \pi_1 ) \Delta^{1} ( \pi_{2} ( a ) )\\
         &=& ( \lambda_{\eta( - s )} \rho_{\eta( u )} ) \pi_2 ( a ) \eean
         $$= \chi^{1}_{(s,u)} \pi_2 ( a ).$$
         
       Thus,  for all $ ( s, u ) $ in $ \IR^{2n}, ~ \pi_2 \chi^{2}_{(s,u)} = \chi^{1}_{(s,u)} \pi_2.$ 
         
    Therefore,  for all $ a, ~ b $ in $  ( G_2 )_{J},$    
     \bean \lefteqn{ \pi_2 ( a  \times^{2}_{J} b )}\\
          &=& \pi_2 ( \int \int \chi^{2}_{Ju} ( a ) \chi^{2}_{v} ( b ) e ( u.v ) du dv )\\
          &=& \int \int \pi_2 ( \chi^{2}_{Ju} ( a ) ) \pi_2 ( \chi^{2}_{v} ( b ) ) e ( u . v ) du dv \\
          &=& \int \int \chi^{1}_{Ju} ( \pi_2 ( a ) ) \chi^{1}_{v} ( \pi_2 ( b ) ) e ( u.v ) du dv \\
          &=& \pi_2 ( a ) \times^{1}_{J} \pi_2 ( b ). \eean 
                   
          where the third step is permissible by Proposition \ref{preliminiaries_rieffel_1.14}.
          
          This proves that $ \pi_2 $ is indeed a homomorphism. \qed
 
\vspace{4mm}
     
\bthm
\label{main_def}

1. If $ QISO^{+}_{R}(\cla^\infty_J, \clh, D) $ and $ (QISO^{+}_{R}(\cla^\infty, \clh, D))_{\widetilde{J}} $ have $ C^* $ actions on $ \cla $ and $ \cla_J $ respectively, we have

$$ \widetilde{QISO^{+}_{R}} ( \cla^{\infty}_{J}, \clh, D ) \cong {( \widetilde{QISO^{+}_{R}} ( \cla^{\infty}, \clh, D ) )}_{\widetilde{J}}  $$

$$ QISO^{+}_{R}(\cla^\infty_J, \clh, D) \cong (QISO^{+}_{R}(\cla^\infty, \clh, D))_{\widetilde{J}}.$$

2. If moreover, $\widetilde{QISO}^+(\cla^\infty, \clh, D)$ and $ \widetilde{{QISO}^{+}} ( \cla^{\infty}_{J}, \clh, D ) $ both exist and have $ C^* $ actions on $ \cla $ and  $ \cla_{J} $ respectively, then 
  $${\widetilde{QISO}}^+(\cla^\infty_J, \clh, D) \cong \left({\widetilde{QISO}}^+(\cla^\infty, \clh, D)\right)_{\widetilde{J}},$$  $$    QISO^+(\cla^\infty_J, \clh, D) \cong (QISO^+(\cla^\infty, \clh, D))_{\widetilde{J}}.$$
  
\ethm

{\it Proof :} We prove 1 only. From Corollary \ref{deformation_qorient_subobject} and Lemma  \ref{volumepreserving} we see that $  \widetilde{QISO^+_R} ( \cla, \clh, D )_{\widetilde{J}} $ is an object of $ {\bf Q^{\prime}}_R ( \cla_J, \clh, D ).$ Thus, 

$$ ( \widetilde{QISO^+_R}( \cla, \clh, D ) )_{\widetilde{J}} < \widetilde{QISO^+_R} ( \cla_{J}, \clh, D ) ~~{\rm in}~ {\bf Q^{\prime}}_R ( \cla_{J}, \clh, D ).$$
           
           So, by Lemma \ref{deformation_qorient_<}, $ ( ( \widetilde{QISO^+_R}( \cla, \clh, D ) )_{\widetilde{J}} )_{- \widetilde J} < ( \widetilde{QISO^+_R} ( \cla_{J}, \clh, D ) )_{- \widetilde J} $ 
           
           in $ {\bf Q^{\prime}}_R ( \cla, \clh, D ) $, hence $ \widetilde{QISO^+_R}( \cla, \clh, D ) <  ( \widetilde{QISO^+_R} ( \cla_{J}, \clh, D ) )_{- \widetilde J}.$
           
           Replacing $ \cla $ by $ \cla_{-J} $, we have 
           \bean \lefteqn{\widetilde{QISO^+_R} ( \cla_{- J}, \clh, D )}\\
           & <&  \widetilde{QISO^+_R} ( ( \cla_{- J} )_{J}, \clh, D )_{- \widetilde{J}} ~~( {\rm  in}~~   {\bf Q^{\prime}}_R ( \cla_{- J}, \clh, D )  )\\
            & \cong &   {\widetilde{QISO^+_R}( \cla, \clh, D )}_{- \widetilde J} ~~({\rm  in}~~ {\bf Q^{\prime}}_R ( \cla_{- J}, \clh, D )  )~
                           \cong ( \widetilde{QISO^+_R}( \cla, \clh, D ) )_{\widetilde {-J}}. \eean                 
           Thus,  $ \widetilde{QISO^+_R}( \cla_J, \clh, D ) < ( \widetilde{QISO^+_R}( \cla, \clh, D ) )_{\widetilde J} $ in $ {\bf Q^{\prime}}_R ( \cla_{J}, \clh, D ) $ which implies 
            $ \widetilde{QISO^+_R}( \cla_J, \clh, D )  \cong ( \widetilde{QISO^+_R}( \cla, \clh, D )  )_{\widetilde J} $ in $ {\bf Q^{\prime}}_R ( \cla_{ J}, \clh, D ). $   \qed

\vspace{8mm}

\subsection{Computations}

Fix a real number $\theta$, and then we recall from subsection \ref{preliminaries_subsection_C*algebras} that the $ C^* $ algebra $ \cla_{\theta} $ is the universal $ C^{*} $ algebra generated by two unitaries $ U $ and $ V $ such that $ U V = \lambda V U $, where $\lambda:=e^{2 \pi i \theta}.$  We also  recall from section \ref{preliminaries_section_Rieffel_deformation}  that $ \cla_{\theta} $ is a Rieffel type deformation of $ C( \IT^2 ) $ by using the canonical action of $ \IR^2 $ on $ \IT^2 $ and the skew symmetric matrix $J =  \left(  \begin {array} {cccc}
       0 & - \frac{\theta}{2}\\   \frac{\theta}{2} & 0 \end {array} \right ) .$ 
     It is well-known (see \cite{con}) that the set $ \{ U^{m}V^{n} : m,n \in \IZ \}  $ is  an orthonormal basis for  $ L^{2} ( \cla_{ \theta } , \tau ),     $ where $ \tau$ denotes the unique faithful normalized trace on $\cla_{\theta}$ given by, $\tau ( \sum a_{m n} U^{m} V^{n} ) = a_{0 0} $.     
      We will denote the GNS space $ L^2 (\cla_\theta,\tau) $ by $ \clh_0.$ 
      Let $\cla_\theta^{\rm fin}$ be the unital  $\ast$-subalgebra generated by finite complex linear combinations of $U^mV^n$, where $m,n $ belong to $ \IZ,$ and $d_1,d_2$ be the maps on $\cla^{\rm fin}_\theta$ defined by $d_1(U^mV^n)=mU^mV^n$, $d_2(U^mV^n)=nU^mV^n.$ 
      
      
      
      We consider the spectral triple obtained from the classical spectral triple on $  \IT^2 $ as described in section \ref{deformation_section_deformation_of_spectral_triple}.
   
   \bthm
    
     $\widetilde{QISO}^+(\cla^\infty_\theta, \clh, D ) = \widetilde{QISO}^+(C^\infty(\IT^2))=C(\IT^2) \ast C(\IT)$, and $QISO^+(\cla^\infty_\theta)=QISO^+(C^\infty(\IT^2))=C(\IT^2).$
     
     \ethm
    
{\it Proof:}   We use Theorem \ref{main_def} and recall that $ QISO^+(C^\infty(\IT^2))=C(\IT^2)$ ( Theorem \ref{qorient_torus_computation_QISO_final theorem} ) which is generated by $ z_{1} $ and $ z_{2} ,$ say. $ QISO^+(C^\infty(\IT^2)) $ contains $ C ( \IT^2 ) $ itself as a quantum subgroup. Hence, by Theorem \ref{main_def}, $ QISO^+(\cla^\infty_\theta) $ is the CQG obtained from the Rieffel deformation of $ C ( \IT^2 ) $ via the action of $ \IR^4 $ and the skew symmetric matrix $ \widetilde{J} = J \oplus - J = \left ( \begin {array} {cccc}
    0 & - \frac{\theta}{2} & 0 & 0   \\  \frac{\theta}{2} &  0 & 0 & 0 \\ 0 & 0 & 0 & \frac{\theta}{2} \\ 0 & 0 & - \frac{\theta}{2} & 0  \end {array} \right ) $ so that $ \widetilde{J} ( r_1, r_2, r_3, r_4 ) = ( - \frac{\theta}{2} r_2, \frac{\theta}{2} r_1, \frac{\theta}{2} r_4, - \frac{\theta}{2} r_3 ).$ For $ f_1, ~ f_2  $ in $ C^{\infty} ( \IT^2 ),~ r = ( r_1, r_2, r_3, r_4 ) $ in $ \IR^4 , ~ r^{\prime} = ( r^{\prime}_1, r^{\prime}_2, r^{\prime}_3, r^{\prime}_4 ) $ in $ \IR^4,~ ( t_1, t_2 ) $ in $ \IT^2, $ the deformed product is given by
$$ ( f_1 \times_{\widetilde{J}} f_2 ) ( t_1, t_2 )  = \int \int \chi_{\widetilde{J}(r_1, r_2, r_3, r_4)} ( f_1 ) ( t_1, t_2 ) \chi_{(r^{\prime}_1, r^{\prime}_2, r^{\prime}_3, r^{\prime}_4)} ( f_2 ) ( t_1, t_2 ) e ( r.r^{\prime} ) dr dr^{\prime}. $$

Here, for  $ f $ in $ C^{\infty} ( \IT^2  ),$
\bean \lefteqn{ \chi_{(r_1, r_2, r_3, r_4)} f ( t_1, t_2 ) }\\
&=& ( ev_{\eta( - ( r_1, r_2  ) )} \otimes {\rm id} ) \Delta (  {\rm id} \otimes ev_{\eta ( r_3, r_4 )} ) \Delta ( f ) ( t_1, t_2 )\\
&=& f ( \eta (  - ( r_1, r_2 ) ) ( t_1, t_2 ) \eta (  r_3, r_4 )  )\\
&=& f ( ( e ( - r_1 ), e ( - r_2  )   )  ( t_1, t_2 ) ( e ( r_3 ), e ( r_4 ) )  )\\
&=& f ( e ( r_3 - r_1 ) t_1, ~ e ( r_4 - r_2 ) t_2   ). \eean

Therefore, \bean \lefteqn{z_1 \times_{\widetilde{J}} z_2}\\
          &=& \int \int z_1 ( e ( \frac{\theta}{2} r_4 + \frac{\theta}{2} r_2 ) t_1, ~ e ( - \frac{\theta}{2} r_1 - \frac{\theta}{2} r_3 ) t_2  ) z_2 (  e ( r^{\prime}_3 - r^{\prime}_1 ) t_1, e ( - r^{\prime}_2 + r^{\prime}_4 ) t_2 )\\
          && e ( r. r^{\prime} ) dr dr^{\prime} \\
          &=& \int \int e ( \frac{\theta}{2} r_4 + \frac{\theta}{2} r_2 ) t_1 e ( - r^{\prime}_2 + r^{\prime}_4 ) t_2 e ( r_1. r^{\prime}_1 ) e ( r_2. r^{\prime}_2 ). e ( r_3. r^{\prime}_3 ) e ( r_4. r^{\prime}_4 ) dr dr^{\prime}\\
          &=& t_1 t_2 \int \int e ( \frac{\theta}{2} r_2 ) e ( - r^{\prime}_2 ) e ( r_2. r^{\prime}_2 ) dr_2 dr^{\prime}_2 \int \int e ( \frac{\theta}{2} r_4 ) e ( r^{\prime}_4 ) e ( r_4. r^{\prime}_4 ) dr_4 dr^{\prime}_4 \\
          & \int \int & e ( r_1. r^{\prime}_1 ) dr_1 dr^{\prime}_1 \int \int e ( r_3. r^{\prime}_3 ) dr_3 dr^{\prime}_3\\
          &=& t_1 t_2 \int \int e ( - \frac{\theta}{2} r_2 ) e ( r^{\prime}_2 ) e ( ( - r_2 ). ( - r^{\prime}_2 )  ) dr_2 dr^{\prime}_2 \int \int e ( - \frac{\theta}{2} r_4 ) e ( - r^{\prime}_4 ) e ( - r_4. - r^{\prime}_4 )\\
          && dr_4 dr^{\prime}_4.1.1. \eean
( by Proposition \ref{preliminiaries_rieffel_1.12} )          
          
 Similarly, \bean \lefteqn{z_2 \times_{\widetilde{J}} z_1 }\\
                   &=& \int \int e ( - \frac{\theta}{2} r_1 - \frac{\theta}{2} r_3 ) t_2 e ( r^{\prime}_3 - r^{\prime}_1 ) t_1 e ( r_1. r^{\prime}_1 ) e ( r_3. r^{\prime}_3 ) dr_1  dr^{\prime}_1 dr_3 dr^{\prime}_3.1.1\\
                   &=& t_1 t_2 \int \int e ( - \frac{\theta}{2} r_1 ) e ( - r^{\prime}_1 ) e ( r_1. r^{\prime}_1 ) dr_1 dr^{\prime}_1 \int \int e ( - \frac{\theta}{2} r_3 ) e ( r^{\prime}_3 ) e ( r_3. r^{\prime}_3 ) dr_3 dr^{\prime}_3\\
                   &=& z_1 \times_{\widetilde{J}} z_2 .\eean
                             
  This proves the theorem.

\qed


\section{ $ QISO^{\cll} $ of a Rieffel deformed noncommutative manifold}

\subsection{Derivation of the main result}

In this subsection, our set up is as in section \ref{deformation_section_deformation_of_spectral_triple} so that we have spectral triples on $ \cla^{\infty}_{J} $ for each $ J $ such that the spectral triple on $ \cla^{\infty} $ satisfies all the assumptions mentioned in chapter \ref{qisol} for ensuring the existence of $ QISO^{\cll}$ where $ \cll $ is the Laplacian coming from the spectral triple. As $ QISO^{\cll} $ has a $ C^* $ action on $ \cla, $ we get a subalgebra $ \cla_0 $ as in Lemma \ref{1111}. Thus, the assumptions of section \ref{deformation_preparatory results} are satisfied so that the conclusions in that subsection hold for $ QISO^{\cll}( \cla^{\infty}, \clh, D ).$  
   

    We recall from Chapter \ref{qisol} that $\cla_0$ will denote the $\ast$-algebra generated by complex linear (algebraic, not closed) span $\cla^\infty_0$ of the eigenvectors of $\cll$ which has a countable discrete set of eigenvalues each with finite multiplicities, by assumptions for the existence of $QISO^{\cll}.$ By the same assumptions, $\cla^\infty_0$ is a subset of $\cla^\infty$ and is norm-dense in $\cla$. Here, we make the following additional assumptions.
     
     {\bf Assumptions}
     {\bf(i)} $\cla_0$ is dense in $\cla^\infty$ w.r.t. the Frechet topology coming from the action of $\IT^n.$
     
     {\bf( ii )} $\bigcap_{n \geq 1 } {\rm Dom}  ({\cll}^n)=\cla^\infty.$
     
     {\bf (iii)} $\cll$ commutes with the $\IT^n$-action $\beta$, hence $C(\IT^n)$ can be identified as a sub object  of $QISO^{\cll}$ in the category $ {\bf Q^{\prime}_{\cll}} $.
     
      Let $\pi$ denote the surjective map from $QISO^{\cll} $ to its quantum subgroup $C(\IT^n)$, which is a morphism of compact quantum groups. We denote by $\alpha : \cla \raro \cla \ot QISO^{\cll} $ the action of $QISO^{\cll} $ on $\cla$, and note that on $\cla_0$, this action is algebraic, that is, it is an action of the Hopf $\ast$-algebra $\clq_0$ consisting of matrix elements of finite dimensional unitary representations of $\clq$. We have $ ( {\rm id} \otimes \pi ) \circ \alpha=\beta$. 
 
   We recall from Proposition \ref{preliminiaries_rieffel_4.10_corollary} that $\cla^\infty=\cla^\infty_J$ as topological spaces, that is they coincide as sets and the corresponding Frechet topologies are also equivalent. In view of this, we shall denote this space simply by $\cla^\infty$, unless one needs to consider it as Frechet algebra, in which case the suffix $J$ will be used. 
  

\blmma

\label{deformation_trace_and_integral_interchanged}

For $ F $ in $ \clb^{\clb ( \clh )} ( \IR^2 ) $ ( notation as in section \ref{preliminaries_section_Rieffel_deformation}   ) and  a trace class operator  $ W ,$ 
$$ Tr ( \int \int F ( u, v ) W e ( u.v ) du dv   ) = \int \int Tr (  F ( u, v ) W ) e ( u . v ) du dv. $$

\elmma

{\it Proof :} From the definition of oscillatory integral, we have 
$ \pi_J ( F ( u, v ) ) = \sum_{p \in L} \int ( F ( u, v )  \phi_p ) ( u, v ) e ( u.v ) du dv  $ ( notations as in section \ref{preliminaries_section_Rieffel_deformation}    ). Let $ \{ L_n \}_{n} $ be a sequence of subsets of $ L $ such that it increases to $ L.$ Then, as the above sum is strongly convergent, $ lim_n (  \sum_{p \in L_n} \int ( F ( u, v ) \phi_p ) ( u, v ) e ( u.v ) du dv ) $ converges in SOT. We deduce the result by using Lemma \ref{preliminaries_SOT_convergence_Trace}. \qed

\bppsn

Let $ \cll_{J} $ denote the Laplacian from the spectral triple $ ( \cla^{\infty}_{J}, \clh, D ).$ Then $ \cll_{J} $ coincides with $ \cll $ on $ \cla^{\infty} \subseteq \cla_{J}.$   

\eppsn

{\it Proof :} We recall that from the proof of Lemma \ref{deformation_deformed representation_is_a_spectral triple}, we have $[ D, \pi_J ( a ) ] = \pi_J ( [ D, a ] ).$  Denoting the inner product on $ \Omega^{1}_{J} ( \cla^{\infty}  )$ by $ \left\langle . , . \right\rangle_{J}$ and letting $a, b $ in $ \cla^{\infty}, ~  {\rm Lim} $ as in subsection \ref{preliminaries_subsection_spaceofforms_NCG}, we have by using Lemma \ref{deformation_trace_and_integral_interchanged}
\bean \lefteqn{ \left\langle \cll_{J} a, b \right\rangle_{J}  }\\
&=&{\rm Lim}_{t \rightarrow 0^+} \frac{ Tr ( {[ D, \pi_J ( a )]}^* [ D, \pi_J ( b ) ] e^{- t D^2} )}{Tr ( e^{- t D^2} ) }\\
&=& {\rm Lim}_{t \rightarrow 0^+} \frac{ Tr ( \pi_J ( {[ D, a ]}^* [ D, b ] ) e^{- t D^2} )}{Tr ( e^{- t D^2} ) }\\
&=& {\rm Lim}_{t \rightarrow 0^+} \frac{ Tr ( \int \int \beta_{Ju} ( {[ D, a ]}^* [ D, b ]  ) \widetilde{\beta_v} e ( u.v ) du dv e^{- t D^2}   )}{Tr ( e^{- t D^2} ) }\\
&=& {\rm Lim}_{t \rightarrow 0^+} \frac{  \int \int Tr ( \beta_{Ju} ( {[ D, a ]}^* [ D, b ]  ) \widetilde{\beta_v} e^{- t D^2}  ) e ( u.v ) du dv}{Tr ( e^{- t D^2} ) }   \\
&& (  {\rm ~by ~ Lemma ~ \ref{deformation_trace_and_integral_interchanged}} )\\
&=& {\rm Lim}_{t \rightarrow 0^+} \frac{ \int \int Tr ( V_{\widetilde{Ju}} ( {[ D, a ]}^* [ D, b ] ) V_{{\widetilde{Ju}}^{- 1}} \beta_v e^{- t D^2} ) e ( u.v ) ) du dv}{Tr ( e^{- t D^2} ) }\\
&=& {\rm Lim}_{t \rightarrow 0^+} \frac{ \int \int Tr ( V_{\widetilde{Ju}} ( {[ D, a ]}^* [ D, b ] ) \beta_v e^{- t D^2}  V_{{\widetilde{Ju}}^{- 1}} ) e ( u.v )  du dv}{Tr ( e^{- t D^2} ) }\\ 
&=& {\rm Lim}_{t \rightarrow 0^+} \frac{  \int \int Tr (( {[ D, a ]}^* [ D, b ] )  \widetilde{\beta_v} e^{- t D^2} ) e ( u.v ) ) du dv}{Tr ( e^{- t D^2} ) }\\
&=&{\rm Lim}_{t \rightarrow 0^+} \frac{ Tr (  \int \int {[ D, a ]}^* [ D, b ] \widetilde{\beta_v} e^{- t D^2} e ( u.v ) du dv )}{Tr ( e^{- t D^2} ) }\\
& &
\mbox{which by Proposition \ref{preliminiaries_rieffel_1.12}, equals}\\
&=&{\rm Lim}_{t \rightarrow 0^+} \frac{ Tr ( {[ D, a ]}^* [ D, b ]  e^{- t D^2} )}{Tr ( e^{- t D^2} ) }\\ 
&=& \left\langle  \cll a, b \right\rangle. \eean.  \qed

Thus, the quantum isometry group  $QISO^{\cll} ( \cla_J )$ is the universal compact quantum group acting on $\cla_J$, with the action keeping each of the eigenspaces of $\cll$ invariant. Note that the algebraic span of eigenvectors of $\cll_J$ coincides with that of $\cll$, that is $\cla^\infty_0$, which is already assumed to be Frechet-dense in $\cla^\infty = \cla^\infty_J$, hence in particular norm-dense in $\cla_J.$ Moreover, all the results of section \ref{deformation_preparatory results} hold for $QISO^{\cll} ( \cla_J ).$

  \vspace{4mm}
  
   We now state and prove a criterion, to be used later, for extending positive maps defined on $\cla_0.$ 
  
    \blmma
    
    \label{deformation_laplacian_extension_positive}
    
    Let $\clb$ be another unital $C^*$-algebra equipped with a $\IT^n$-action, so that we can consider the $C^*$-algebras $\clb_J$ for any skew symmetric $n \times n$ matrix $J$. 
    Let $\phi : \cla^\infty \raro \clb^\infty$ be a linear map, satisfying the following :\\
        (a)  $\phi$  is positive w.r.t. the deformed products $\times_J$ on $\cla_0$ and $\clb^\infty$, that is $\phi(a^* \times_J a) \geq 0 $ ( in $\clb^\infty_J \subset \clb_J$ )  for all $a $ in $ \cla_0,$ and \\
    (b) $\phi$ extends to a norm-bounded map (say $\phi_0$) from $\cla$ to $\clb$.\\
    Then $\phi$ also have an extension $\phi_J$   as a $\|~\|_J$-bounded positive map from $\cla_J$ to $\clb_J$ satisfying $\| \phi_J \| =\| \phi(1) \|_J.$
    
    \elmma
    
    {\it Proof :} We can view $\phi$ as a map between the Frechet spaces $\cla^\infty$ and $\clb^\infty$,  which is clearly closable, since it is continuous w.r.t. the norm-topologies of $\cla$ and $\clb$, which are weaker than the corresponding Frechet topologies.  By the Closed Graph Theorem, we conclude that $\phi$ is continuous in the Frechet topology. Since  $\cla^\infty=\cla^\infty_J$ and $\clb^\infty=\clb^\infty_J$ as Frechet spaces, consider $\phi$ as a continuous map from $\cla^\infty_J$ to $\clb^\infty_J$, and  it  follows  by the Frechet-continuity of $\times_J$ and $\ast$ and the Frechet-density of $\cla_0$ in $\cla^\infty_J$ that  the positivity (w.r.t. $\times_J$) of the restriction of $\phi$ to $\cla_0 \subset \cla^\infty_J$ is inherited by the extension on $\cla^\infty=\cla^\infty_J$. Indeed,  given $a $ in $ \cla^\infty_J=\cla^\infty$, choose a sequence $a_n $ belonging to $  \cla_0$ such that $a_n \raro a$ in the Frechet topology. We have $\phi(a^* \times_J a)=\lim_n \phi(a_n^* \times_J a_n)$ in the Frechet topology, so in particular, $\phi(a_n^* \times_J a_n) \raro \phi(a^* \times_J a)$ in the norm of $\clb_J$, which implies that $\phi(a^* \times_J a)$ is a positive element of $\clb_J$ since $\phi(a_n^* \times_J a_n)$ is so for each $n$. 
        Next, by Lemma \ref{preliminaries_Rieffel_positivity}, we note that $\cla^\infty$ is closed under holomorphic functional calculus as a unital $\ast$-subalgebra of $\cla_J$ (the identity of $\cla^\infty_J$ is same as that of $\cla$), and hence, by Lemma \ref{preliminaries_Rieffel_positivity}, for any self adjoint element $ x $ in $ \cla^{\infty}, ~ \left\| \phi ( x ) \right\| \leq \left\| x \right\| \phi ( 1 ).$  Thus, for any $ y $ in $ \cla^{\infty}, ~ \left\| \phi ( y ) \right\| = \left\| \phi ( \frac{y + y^*}{2} + i \frac{y - y^*}{2i} ) \right\| \leq  \left\| \phi ( \frac{y + y^*}{2} ) \right\| +  \left\| \phi ( \frac{y - y^*}{2i} ) \right\| \leq (  \left\| \frac{y + y^*}{2} \right\| + \left\| \frac{y - y^*}{2i} \right\| ) \phi ( 1 ) \leq 2 \left\| y \right\| \phi ( 1 ).$ Thus, $ \phi  $ is bounded on $ \cla^{\infty} $ and hence admits a bounded extension (say $\phi_J$) on $\cla_J$, which will still be a positive map, so in particular the norm of $\phi_J$ is same as $\| \phi_J(1)\|$. \qed 
        
        \vspace{4mm}
    
Now we note that due to the assumption {\bf (iii)}, $ {( QISO^{\cll} )}_{\widetilde{J}}$ makes sense. 
    
 \vspace{4mm}
 
    \brmrk
    
    \label{deformation_laplacian_subobject}
     
     By Lemma \ref{qiso_deformation_Lemma5} along with Lemma \ref{qiso_deformation_Lemma6}, $ \alpha $ is a homomorphism from $\cla_J$ to $\cla_J \otimes ( QISO^{\cll}(\cla) )_{\widetilde{J}}$ and hence $ ( QISO^{\cll} ( \cla ) )_{\widetilde{J}} $ is an object in the category $ {\bf \clq^{\cll}}( \cla_{J} ).$ 
    
    \ermrk

     \bthm
     \label{premain}
   If the Haar state is faithful on $\clq$, then $\alpha : {\cla_0} \raro {\cla_0} \ot \clq_0$ extends to an action of the compact quantum group $\clq_{\widetilde{J}}$ on $\cla_J$, which is isometric, smooth and faithful. 
     \ethm
        {\it Proof :}
        We have already seen in ( \ref{qiso_deformation_Lemma5.2} ) that $\alpha$ is an algebra homomorphism 
     from $\cla_0$ to $\cla_0 \ot_{\rm alg} \clq_0$ (w.r.t. the deformed products), and it is also $\ast$-homomorphism since it is so for the undeformed case and the involution $\ast$ is the same for the deformed and undeformed algebras.  
     It now suffices to show that $\alpha$ extends to $\cla_J$ as a $C^*$-homomorphism. 
     Let us fix any faithful imbedding $\cla_J \subseteq \clb(\clh_0)$ (where $\clh_0$ is a Hilbert space) and consider the imbedding $\clq_{\widetilde{J}} \subseteq \clb(L^2(h_J))$ ( by virtue of Lemma \ref{deformation_laplacian_haar state faithful_1} ). By definition, the norm on $\cla_J \ot \clq_{\widetilde{J}}$ is the minimal (injective) $C^*$-norm, so it is equal to the norm inherited from the imbedding $\cla_J \ot_{\rm alg} \clq_{\widetilde{J}} \subseteq \clb(\clh_0 \ot L^2(h_J))$. Let us consider the dense subspace $\cld \subset \clh_0 \ot L^2(h_J)$ consisting of vectors which are finite linear combinations of the form $\sum_i u_i \ot x_i,$ with $u_i $ belonging to $ \clh_0$, $x_i $ belonging to $ \clq_0 \subset L^2(h_J)$. Fix such a vector $\xi=\sum_{i=1}^k u_i \ot x_i$ and consider $\clb:=\cla \ot M_k(\IC)$, with the $\IT^n$-action $\beta  \otimes \it {\rm id}$ on $\clb$. Let $\phi: \cla^\infty \raro \clb^\infty$ be the map given by 
     $$ \phi(a):= \left( \left( ({\rm id} \ot \phi_{(x_i,x_j)})(\alpha(a)) \right) \right)_{i,j=1}^k,$$ where $\phi_{(x,y)}(z):=h(x^* \times_{\widetilde{J}} z \times_{\widetilde{J}} y)$ for $x,y,z $ in $ \clq_0.$ By Remark \ref{haarrem}, $\phi_{(x_i,x_j)}$ extends to $\clq$ as  a bounded linear  functional. 
     Note that the range of $\phi$ is in $\clb^\infty=\cla^\infty \ot M_k(\IC)$ since we have  $ ( {\rm id} \otimes \phi_{(x,y)} ) \alpha (\cla^\infty) \subseteq \cla^{\infty}$ by Proposition \ref{QISO_Laplacian_formulation_remark_2.16}, using our assumption (ii) that $\bigcap_{n \geq 1}{\rm Dom}(\cll^n)=\cla^\infty$.
     
      Since  $\alpha$ maps $\cla_0$ into $\cla_0 \ot_{\rm alg} \clq_0$ and $h=h_J$ on $\clq_0$, it is easy to see that   for $a $ in $ \cla_0,$ $\phi(a^* \times_J a)$ is positive in $\clb_J$. As $\phi_{(x_i,x_j)}$ extends to $\clq$ as  a bounded linear  functional, $\phi$ extends to a bounded linear (but not necessarily positive) map from $\cla$ to $\clb$.  
      Thus, the hypotheses of  Lemma \ref{deformation_laplacian_extension_positive} are satisfied and we conclude that $\phi$ admits a positive extension, say $\phi_J$, from $\cla_J$ to $\clb_J=\cla_J \ot M_k(\IC)$. Thus, we have for $a $ in $ \cla_0,$
       \begin{eqnarray*} 
       &&{\sum_{i,j=1}^k \lgl u_i, \phi(a^* \times_J a)u_j \rgl}\\
      & \leq & \| a \|_J^2 \sum_{ij} \lgl u_i, \phi(1)u_j \rgl 
      = \| a\|_J^2 \sum_{ij} \lgl u_i,u_j \rgl h(x_i^* \times_{\widetilde{J}} x_j)\\
      &=& \| a \|_J^2 \sum_{ij} \lgl u_i \ot x_i, u_j \ot x_j \rgl =\|a\|_J^2 \| \sum_{i=1}^k u_i \ot x_i \|^2.\end{eqnarray*}
       This implies 
      $$ \| \alpha(a) \xi \|^2=\lgl \xi, \alpha(a^* \times_J a) \xi \rgl = \sum_{i,j=1}^k \lgl u_i, \phi(a^* \times_J a)u_j \rgl \leq \| a \|^2_J \| \xi \|^2$$
       for all $\xi $ in $ \cld$ and $a $ in $ \cla_0,$  hence $\alpha$ admits a bounded extension which is clearly a $C^*$-homomorphism.    
                 \qed \\

           
        
 
 For any two objects $ \cls_1 $ and $ \cls_2 $ in 
$ {\bf Q^{\prime}_{\cll}}),$ we write $ \cls_1 <_{\cll} \cls_2 $ if $ \cls_1 $ is a sub-object of $ \cls_2 $ in the category $ {\bf Q^{\prime}_{\cll}}.$

 \blmma
 
 \label{deformation_laplacian_<}
        
        If $G_1, ~ G_2 $ be two CQG s such that $ G_1 <_{\cll} G_2 $ in the category $ {\bf Q^{\prime}_{\cll}}( \cla, \clh, D ).$ If $ ( G_1 )_{\widetilde{J}} $ and $ ( G_2 )_{\widetilde{J}} $ make sense, then  $ ( G_1 )_{\widetilde{J}} <_{\cll} ( G_2 )_{\widetilde{J}} $ in the category $ {\bf Q^{\prime}_{\cll}}( \cla_J, \clh, D ).$
         
        \elmma
        
        {\it Proof :} By Remark \ref{deformation_laplacian_subobject}, $ ( G_i )_{\widetilde{J}} $ are objects in the category $ {\bf Q^{\prime}_{\cll}} ( A_{J} ), ~ i = 1,2.$ The rest of the proof is the same as that of Lemma \ref{deformation_qorient_<}  and hence we omit it.  \qed

           \bthm          
           \label{abcd}
           If the Haar state on $QISO^{\cll}( \cla )$ is faithful, we have the isomorphism of compact quantum groups:
           $$( QISO^{\cll} ( A ) )_{ \widetilde J } \cong  QISO^{\cll} ( A_{J} ). $$ 
           
           \ethm
           
        {\it Proof :}   By Theorem \ref{premain},  we have seen that $ ( QISO^{\cll} ( \cla ) )_{\widetilde J} $ also acts faithfully, smoothly and isometrically on $ \cla_{J} $, which implies,           
             $$ ( QISO^{\cll} ( \cla, \clh, D ))_{\widetilde  J} <  QISO^{\cll} ( \cla_{J}, \clh, D ) ~~{\rm in}~ {\bf Q^{\prime}_{\cll}} ( \cla_J, \clh, D )  .$$
             The rest of the proof is same as that of Theorem \ref{main_def} by using Lemma \ref{deformation_laplacian_<} and hence we omit it. \qed
           
           

    \subsection{Computations}

     {\bf The case of the noncommutative tori}
     
     \vspace{2mm} 

We recall ( chapter \ref{preliminaries} ) that the noncommutative $n$-tori $ \IT^n_{\theta} $ is the universal $ C^* $ algebra generated by $n$ unitaries $ U_1, U_2,..., U_n$ such that $ U_i U_j = e ( \theta_{ij} ) U_j U_i, ~ i,j, = 1,2, ...n $ where $ \theta \equiv (( \theta_{ij} )) $ is a $ n \times n $ skew symmetric matrix. We also recall that  $ \IT^n_{\theta} $ is a Rieffel deformation of $ C ( \IT^n )$ via the $ \IR^n $ action on $ C ( \IT^n ) $ ( section \ref{preliminaries_section_Rieffel_deformation} ) and the matrix $ J = \frac{\theta}{2}.$     
       
      
      We consider the isospectral deformation of the classical spectral triple on $ C ( \IT^n ) $ which will give a spectral triple on $ \IT^n_{\theta}.$
   so that the corresponding Laplacian $\cll$ is given by $\cll(U^{m_1}_1 ... U^{m_n}_n)= - (m^2_1 + ...m^2_n) U^{m_1}_1 ... U^{m_n}_n,$ and it is also easy to see that all the assumptions in chapter \ref{qisol} required for defining $ QISO^{\cll} ( \IT^n_{\theta} ) $ are satisfied.
      
  \bthm 
  
  \label{deformation_laplacian_all_NC_tori}      
       
     Using Theorem \ref{abcd}, we conclude that $ QISO^{\cll} ( \IT^n_{\theta} ) $ is a Rieffel deformation of $ QISO^{\cll} ( C ( \IT^{n} ) ).$

     \ethm 
     
     \vspace{4mm} 
     
     Next, we will use  Theorem \ref{abcd} to compute the exact structure of $ QISO^{\cll} ( \cla_{\theta} ).$ Our target is to show that $ QISO^{\cll} ( \cla_{\frac{1}{2}} )$ is a commutative $ C^* $ algebra which will give an example of a non commutative $ C^* $ algebra with $ QISO^{\cll} $ a commutative CQG.       
          
  We have seen ( Theorem \ref{QISO_Laplacian_computations_Torus_finaltheorem} ) that $ QISO^{\cll}( C ( \IT^{n} ) ) \cong C( {\IT}^{n} >\!\!\! \lhd  ( {\IZ}^{n}_{2} >\!\!\! \lhd S_n )  ).$
  
   In particular, $ QISO^{\cll} ( C ( \IT^2  )  ) \cong C ( {\IT}^{2} >\!\! \lhd ( {\IZ}^{2}_{2} >\!\!\! \lhd S_2  )) $ where the definition of the two semi-direct products are described below.
  Let the generators of the first, second and the third copy of $ \IZ_2 $ in $ \IZ^{2}_{2}  >\!\! \lhd \IZ_2 $ be denoted by $( 1, 0, 0 ), ( 0, 1, 0), ( 0, 0, 1 ) .$ Therefore, as a set $  ( \IZ^{2}_{2}  >\!\! \lhd \IZ_2   )   $ equals $ \{ ( \gamma_1, ~ \gamma_2,~ \gamma_3 ) : ~ \gamma_i \in \{0,~ 1\}, i = 1,2,3  \} $ and $ \IT^2 ~ \!\!\! \lhd ( \IZ^{2}_{2}  >\!\! \lhd \IZ_2   ) = \{ ( z_1, ~ z_2, \gamma_1, ~ \gamma_2,~ \gamma_3 ) : z_1, z_2 \in S^{1},  ~ \gamma_i \in \{0,~ 1\}, i = 1,2,3  \} .$  The  action of $ \IZ_2 $ on $ \IZ^2_2 $ is given by $ ( 0, 0, 1) ( x, y ) = ( y, x ) $ (  ( $ x,y $ )  belongs to $ \IT^2 $ ) and the action of $  \IZ^{2}_{2}  >\!\! \lhd \IZ_2  $ on $ \IT^2 $ ( denoted by $ \circ $ ) is given by  
 \bean ( 1, 0, 0 ) \circ ( z_1, z_2 ) = ( \overline{z_1},  z_2 ),   ( 0, 1, 0 ) \circ ( z_1, z_2 ) = ( z_1, \overline{z_2} ),  ( 0,0,1 ) \circ ( z_1, z_2 ) = ( z_2, z_1 ).  \eean

     \blmma
     
     Let $ ( z_1, z_2 ),~ ( z^{\prime}_1, ~ z^{\prime}_2 ), ~ ( z^{\prime \prime}_1, ~ z^{\prime \prime}_2 ) $ belong to $ \IT^{2}.$ Then we have the following group multiplication formulae in $ \IT^2 \!\!\! \lhd ( \IZ^{2}_{2}  >\!\! \lhd \IZ_2   ) : $

     \be \label{deformation_laplacian_torus_101} ( z_1, z_2, 0, 0, 0 ) ( z^{\prime}_1, z^{\prime}_2, 1, 0, 1 ) ( z^{\prime \prime}_1, z^{\prime \prime}_2,0,0,0 ) = ( z_1 z^{\prime}_1 \overline{z^{\prime \prime}_2}, ~ z_2 z^{\prime}_2 z^{\prime \prime}_1, 1,0,1 ) ,  \ee     
     \be \label{deformation_laplacian_torus_001} ( z_1, z_2, 0, 0, 0 ) ( z^{\prime}_1, z^{\prime}_2, 0, 0, 1 ) ( z^{\prime \prime}_1, z^{\prime \prime}_2,0,0,0 ) = ( z_1 z^{\prime}_1 z^{\prime \prime}_2, ~ z_2 z^{\prime}_2 z^{\prime \prime}_1, 0,0,1 ) ,  \ee     
     \be \label{deformation_laplacian_torus_000} ( z_1, z_2, 0, 0, 0 ) ( z^{\prime}_1, z^{\prime}_2, 0, 0, 0 ) ( z^{\prime \prime}_1, z^{\prime \prime}_2,0,0,0 ) = ( z_1 z^{\prime}_1 z^{\prime \prime}_1, ~ z_2 z^{\prime}_2 z^{\prime \prime}_2, 0,0,0 )  . \ee

     \elmma 
     
       {\it Proof :} Using the definition of semi direct product , we have       
        \bean \lefteqn{ ( z_1, z_2, 0, 0, 0 ) ( z^{\prime}_1, z^{\prime}_2, 1, 0, 1 ) ( z^{\prime \prime}_1, z^{\prime \prime}_2,0,0,0 )}\\
       &=&( z_1 z^{\prime}_1,  z_2 z^{\prime}_2,1,0,1 ) ( z^{\prime \prime}_1, z^{\prime \prime}_2,0,0,0 )\\
       &=&( z_1 z^{\prime}_1,  z_2 z^{\prime}_2 ) ( ( 1,0,1 ) \circ ( z^{\prime \prime}_1, z^{\prime \prime}_2 ),1,0,1 )\\
  &=& ( z_1 z^{\prime}_1, z_2 z^{\prime}_2 ) (  ( ( 1, 0 ), 0 ) ( ( 0, 0 ), 1 ) \circ ( z^{\prime \prime}_1, z^{\prime \prime}_2 ), 1,0,1 )\\     
  &=& ( ( z_1 z^{\prime}_1,  z_2 z^{\prime}_2 ) ( \overline{z^{\prime \prime}_2}, z^{\prime \prime}_1 ), 1,0,1  )\\
  &=& ( z_1 z^{\prime}_1 \overline{z^{\prime \prime}_2}, z_2 z^{\prime}_2 z^{\prime \prime}_1, 1,0,1 ). \eean
       
       Thus, we have ( \ref{deformation_laplacian_torus_101} ).
       
       Similarly, 
       \bean \lefteqn{ ( z_1, z_2, 0, 0, 0 ) ( z^{\prime}_1, z^{\prime}_2, 0, 0, 1 ) ( z^{\prime \prime}_1, z^{\prime \prime}_2,0,0,0 )}\\
       &=& ( z_1 z^{\prime}_1,  z_2 z^{\prime}_2,0,0,1 ) ( z^{\prime \prime}_1, z^{\prime \prime}_2,0,0,0 )\\ 
       &=& ( ( z_1 z^{\prime}_1,  z_2 z^{\prime}_2 ) ( ( 0,0,1 ) \circ ( z^{\prime \prime}_1, z^{\prime \prime}_2 ) ),0,0,1  )\\
       &=& ( ( z_1 z^{\prime}_1,  z_2 z^{\prime}_2 )( z^{\prime \prime}_2, z^{\prime \prime}_1 ),0,0,1 )\\
       &=& ( z_1 z^{\prime}_1 z^{\prime \prime}_2, z_2 z^{\prime}_2  z^{\prime \prime}_1 ,0,0,1  ), \eean 
       
       which proves ( \ref{deformation_laplacian_torus_001} ).
       
       ( \ref{deformation_laplacian_torus_000} ) follows similarly and therefore we omit the proof. \qed

\vspace{4mm}
       
       $ QISO^{\cll} ( C ( \IT^{2} ) ) = C( {\IT}^{2} >\!\!\! \lhd ( \IZ^{2}_{2}  >\!\!\! \lhd \IZ_2  ) )  $ which as a $ C^* $ algebra is isomorphic to $ \oplus_{i = 1,2,...8} C ( \IT^{2} ).$ 
   We will denote the generators of $ QISO^{\cll} ( \IT^{2} )$ by $\{ A_{ ( \gamma_1, \gamma_2, \gamma_3 )}, ~ B_{ ( \gamma_1, \gamma_2, \gamma_3 )} : ( \gamma_1, \gamma_2, \gamma_3 ) \in  \IZ^{3}_{2} \},$ more precisely, we fix the following convention:
   $$ A_{ ( \gamma_1, \gamma_2, \gamma_3 )} ( z_1, z_2,\gamma_1, \gamma_2, \gamma_3   ) = z_1, ~ B_{ ( \gamma_1, \gamma_2, \gamma_3 )} ( z_1, z_2, \gamma_1, \gamma_2, \gamma_3  ) = z_2.$$
   
   Now, $ \IT^2 $ sits as a subgroup of $ {\IT}^{2} >\!\!\! \lhd Z^{3}_{2} $ as $ \{ ( z_1, z_2, 0,0,0 ) : z_i \in S^{1} \}.$
   
   Hence, the action of $ \IT^2 $ on $ {\IT}^{2} >\!\!\! \lhd Z^{3}_{2} $ with which we are concerned with is given by the group multiplication in $ {\IT}^{2} >\!\!\! \lhd Z^{3}_{2} ,$ that is, the action of $ ( z_1, z_2 ) ~ ( (z_i \in  S^1) ) $ on $ ( z^{\prime}_1, ~ z^{\prime}_2, \gamma_1, \gamma_2, \gamma_3  ) ( \in {\IT}^{2} >\!\!\! \lhd Z^{3}_{2}  )  $ is given by $ ( z_1, z_2, 0,0,0 ) ( z^{\prime}_1, ~ z^{\prime}_2, \gamma_1, \gamma_2, \gamma_3  ). $
  
    The action of $ \IR^{4} $ on $ QISO^{\cll} ( C ( \IT^2 ) ) $ as prescribed by Wang ( \ref{preliminiaries_deformation_wang_formula} ) is given by    
    $$  \alpha_{( s,u )} = ( \Omega( -s ) \otimes {\rm id} ) \Delta ( {\rm id} \otimes \Omega ( u ) ) \Delta $$    
    where $s, u $ belong to $ \IR^2.$
    
    If $ s = ( s_1, s_2 ), ~ u = ( u_1, u_2 ) $ belong to $\IR^2$ and $ z_1, z_2 $ belong to $ S^{1}, $ we have
    
     $  \alpha_{( s,u )} f ( z_1, z_2, \gamma_1, \gamma_2, \gamma_3  ) = f ( ( e ( - s_1 ), e ( - s_2 ), 0,0,0 ) ( z_1, z_2, \gamma_1, \gamma_2, \gamma_3 ) ( e ( u_1 ),  e ( u_2 ),$
     
    $  0,0,0 ) ) .  $
       
  Moreover,   $ \widetilde{J} = - J \oplus J =  \left(  \begin {array} {cccc}
       0 & \frac{\theta}{2} & 0 & 0  \\  - \frac{\theta}{2} & 0 & 0 & 0 \\ 0 & 0 & 0 & - \frac{\theta}{2} \\ 0 & 0 & \frac{\theta}{2} & 0 \end {array} \right ) .$ 
       
       Hence, $ \widetilde{J} ( s_1, s_2, t_1, t_2 )^{t} = ( \frac{\theta}{2} s_2, - \frac{\theta}{2} s_1, - \frac{\theta}{2} t_2, \frac{\theta}{2} t_1  )^{t} $ where $ t $ stands for the transpose of a matrix.  
       
     Let $ s = ( s_1, s_2 ), t = ( t_1, t_2 ), s^{\prime} = ( s^{\prime}_1, s^{\prime}_2 ), t^{\prime} = ( t^{\prime}_1, t^{\prime}_2 ) $ belong to $ \IR^{2} ~ u = ( s, t ), ~ v = ( s^{\prime}, t^{\prime} ).$

     \blmma
     
     \label{deformation_torus_commutation_relation0}
     
     $ \{ A_{(\gamma_1, \gamma_2, \gamma_3)}, ~ B_{( \gamma_1, \gamma_2, \gamma_3)}  \}  $ are unitaries for the multiplication $ \times_{\widetilde{J}} $ and for $  ( \gamma_1, \gamma_2, \gamma_3 ) $ belonging to $ \IZ^{3}_{2}.$
     
     \elmma 
     
     {\it Proof :} We give the proof for  $ A_{(1,0,1)},$ the proof for the rest being exactly similar.     
           \bean \lefteqn{  A^{*}_{(1,0,1)} \times_{\widetilde{J}} A_{1,0,1} ( u_1, u_2, 1,0,1 ) }\\
             &=& \int \int \chi_{\widetilde{J} ( u ) } ( A^{*}_{(1,0,1)} ) ( u_1, u_2, 1,0,1 ) \chi_v (  A_{1,0,1} ) ( u_1, u_2, 1,0,1 ) \\
             && e ( u.v ) du dv\eean
            \bean &=& \int \int \chi_{(\frac{\theta}{2} s_2, - \frac{\theta}{2} s_1, - \frac{\theta}{2} t_2, \frac{\theta}{2} t_1 )} A^*_{(1,0,1)} ( u_1, u_2,1,0,1 ) \chi_{(s^{\prime}_1,s^{\prime}_2,t^{\prime}_1,t^{\prime}_2)} A_{(1,0,1)} ( u_1,\\
             && u_2, 1,0,1 ) e ( u.v ) du dv \\
             &=& \int \int A^*_{(1,0,1)} ( ( e ( - \frac{\theta}{2} s_2  ), e ( \frac{\theta}{2} s_1 ) ) ( u_1, u_2,1,0,1   ) ( e ( - \frac{\theta}{2} t_2 ), e ( \frac{\theta}{2} t_1 ) ) ) \\
         && A_{(1,0,1)} ( ( e ( - s^{\prime}_1  ) , e ( - s^{\prime}_2  ) ) ( u_1, u_2, 1,0,1 ) ( e ( t^{\prime}_1 ) , e ( t^{\prime}_2 ) ) ) e ( u.v ) du dv \eean
             
             which by ( \ref{deformation_laplacian_torus_101} ) equals             
             \bean& &  \int \int A^*_{(1,0,1)} (   e ( - \frac{\theta}{2} s_2 )  u_1 e (  - \frac{\theta}{2} t_1 ), ~ e ( \frac{\theta}{2} s_1 ) u_2   e ( - \frac{\theta}{2} t_2 ) ) A_{(1,0,1)} ( e( - s^{\prime}_1 ) u_1 e ( - t^{\prime}_2 ),\\
             && e ( - s^{\prime}_2 )  u_2 e ( t^{\prime}_1 ) ) e ( u.v ) du dv \\
             &=& \int ...\int e (  \frac{\theta}{2} s_2 )e (   \frac{\theta}{2} t_1 ) e( - s^{\prime}_1 )e ( - t^{\prime}_2 ) e ( s_1. s^{\prime}_1 + t_1 . t^{\prime}_1 + s_2 . s^{\prime}_2 +  t_2 . t^{\prime}_2  ) ds_1 ds^{\prime}_1 ds_2 ds^{\prime}_2 dt_1\\
             && dt^{\prime}_1 dt_2 dt^{\prime}_2 \\
             &=& \int \int e( - s^{\prime}_1 ) e ( s_1 . s^{\prime}_1 ) ds_1 ds^{\prime}_1 \int \int  e (  \frac{\theta}{2} s_2 ) e ( s_2 . s^{\prime}_2 )ds_2 ds^{\prime}_2 \\
             & \int \int & e (   \frac{\theta}{2} t_1 ) e ( t_1 . t^{\prime}_1 ) dt_1 dt^{\prime}_1 \int \int e ( - t^{\prime}_2 ) e (t_2 . t^{\prime}_2 ) dt_2 dt^{\prime}_2 \eean
             
             which by Proposition \ref{preliminiaries_rieffel_1.12}  equals \bean  1.1.1.1 = 1 \eean
             
      Similarly, \bean A_{(1,0,1)} \times_{\widetilde{J}} A^{*}_{1,0,1} ( u_1, u_2, 1,0,1 ) = 1 \eean    \qed       
          
 \brmrk
 
 It can be easily checked that 
 \bean A^*_{\gamma_1, \gamma_2,\gamma_3} \times_{\widetilde{J}} A_{\gamma_1, \gamma_2,\gamma_3} ( u_1, u_2, {\gamma}^{\prime}_1, {\gamma}^{\prime}_2,{\gamma}^{\prime}_3 ) = 0 \eean
 if $ ( \gamma_1, \gamma_2,\gamma_3 ) \neq ( \gamma^{\prime}_1, \gamma^{\prime}_2,\gamma^{\prime}_3 ).$
 
 Similar is the case with $ B_{\gamma_1, \gamma_2,\gamma_3}. $

 \ermrk
          
   \blmma
   
   \label{deformation_torus_commutation_relation1}
   
   $ A_{000} \times_{\widetilde{J}} B_{000} = B_{000} \times_{\widetilde{J}} A_{000}. $ 
   
   \elmma          
   {\it Proof :}  \bean \lefteqn{ ( A_{000} \times_{\widetilde{J}} B_{000} )( (u_1, u_2,0,0,0))}\\
    &=& \int \int \chi_{\widetilde{J} u} ( A_{000} ) ( (u_1, u_2,0,0,0) ) \chi_v ( B_{000} ) ( (u_1, u_2,0,0,0) ) e ( u.v ) du dv\\
   &=& \int.....\int \chi_{\widetilde{J}( s_1, s_2, t_1, t_2)} A_{000} ( (u_1, u_2,0,0,0) ) \chi_{(s^{\prime}_{1}, s^{\prime}_{2}, t^{\prime}_{1}, t^{\prime}_{2})} ( B_{000} ) ( ( u_1, u_2, 0, 0, 0)  ) \\
   && e ( u.v ) ds dt ds^{\prime} dt^{\prime} \\
   &=& \int... \int \chi_{( \frac{\theta}{2} s_2, - \frac{\theta}{2} s_1, - \frac{\theta}{2} t_2, \frac{\theta}{2} t_1  )} A_{000} ( ( u_1, u_2,0,0,0)  ) \chi_{(s^{\prime}_{1}, s^{\prime}_{2}, t^{\prime}_{1}, t^{\prime}_{2})} ( B_{000} ) ( ( u_1, u_2,0,0,0)  )\\
   && e ( u.v ) ds dt ds^{\prime} dt^{\prime} \\
   &=& \int...\int A_{000} [ ( e ( - \frac{ \theta}{2} s_2 ), e ( \frac{\theta}{2} s_1 ) ) ( u_1, u_2,0,0,0) ( e ( - \frac{\theta}{2} t_2 ), e ( \frac{\theta}{2} t_1 ) ) ] \\
   &B_{000}&[ ( e( - s^{\prime}_{1} ), e ( - s^{\prime}_{2} )  ) ( u_1, u_2,0,0,0) ( e ( t^{\prime}_1 ), e( t^{\prime}_2 )   )  ] e ( u. v ) ds dt ds^{\prime} dt^{\prime} \eean
   
   which by ( \ref{deformation_laplacian_torus_000} ) equals
   
   $ \int...\int A_{000} [  e ( - \frac{\theta}{2} s_2 ) u_1 e( - \frac{\theta}{2} t_2 ), e ( \frac{\theta}{2} s_1 )  u_2 e ( \frac{\theta}{2} t_1 ), 0, 0, 0 ] B_{000} [ (e(- s^{\prime}_1 ) u_1 e ( t^{\prime}_1 )  ), ( e (- s^{\prime}_2)$
   
   $u_2 e (t^{\prime}_2) ),0,0,0   ] e ( u. v ) ds dt ds^{\prime} dt^{\prime} $
  \bean \lefteqn{ = \int....\int e ( - \frac{\theta}{2} s_2 ) u_1 e ( - \frac{\theta}{2} t_2 ) e ( - s^{\prime}_{2} ) u_2 e ( t^{\prime}_{2} ) e ( s_1 s^{\prime}_1 ) e ( s_2 s^{\prime}_2 ) e ( t_1 t^{\prime}_1 ) e ( t_2 t^{\prime}_2 ) ds dt ds^{\prime} dt^{\prime}}\\
   &=& u_1 u_2 \int \int e ( - \frac{\theta}{2} s_2 ) e ( - s^{\prime}_2 ) e ( s_2 . s^{\prime}_2 ) ds_2 ds^{\prime}_2 \int \int e ( - \frac{\theta}{2} t_2 ) e ( t^{\prime}_2 ) e ( t_2 . t^{\prime}_2 ) dt_2 dt^{\prime}_2 \\
   && \int \int e ( s_1 . s^{\prime}_{1} ) ds_1 ds^{\prime}_{1} \int \int e ( t_1 . t^{\prime}_1 ) dt_1 dt^{\prime}_1 \eean
   
   which by Corollary \ref{preliminiaries_rieffel_2.20_corollary} and Proposition \ref{preliminiaries_rieffel_1.12}   
  \bean \lefteqn{ = u_1 u_2 e ( - \frac{\theta}{2} ) e ( \frac{\theta}{2} ). 1.1 }\\
   &=& u_1 u_2. \eean
   
   Similarly, \bean \lefteqn{ ( B_{000} \times_{\widetilde{J}} A_{000} ) ( ( u_1, u_2 ,0, 0, 0 )  ) }\\
              &=& \int.... \int B_{000} [ ( e ( - \frac{\theta}{2} s_2 ) u_1 e ( - \frac{\theta}{2} t_2 ) ), ~ (  e ( \frac{\theta}{2} s_1 ) u_2 e ( \frac{\theta}{2} t_1 )  ), 0,0,0)] \\
              && A_{000} [ ( e ( - s^{\prime}_1 ) u_1 e ( t^{\prime}_1 ) ), ( e ( - s^{\prime}_2 ) u_2 e ( t^{\prime}_2 )  )  0,0,0)] e ( u. v ) ds dt ds^{\prime} dt^{\prime} \\
            &=& \int... \int e ( \frac{\theta}{2} s_1 ) u_2 e ( \frac{\theta}{2} t_1  ) e ( - s^{\prime}_{1} ) u_1 e ( t^{\prime}_{1} ) e ( u. v ) ds dt ds^{\prime} dt^{\prime} \\
            &=& u_1 u_2 \int \int e ( \frac{\theta}{2} s_1 ) e ( - s^{\prime}_{1} ) e ( s_1. s^{\prime}_1 ) ds_1 ds^{\prime}_1 \int \int e ( \frac{\theta}{2} t_1 ) e ( t^{\prime}_1 ) e ( t_1. t^{\prime}_1 ) dt_1 dt^{\prime}_1 \\
            && \int \int e ( s_2 s^{\prime}_2 ) ds_2 ds^{\prime}_2 \int \int e ( t_2 t^{\prime}_2 ) dt_2 dt^{\prime}_2, \eean
 
 which by  Corollary \ref{preliminiaries_rieffel_2.20_corollary}  and Proposition \ref{preliminiaries_rieffel_1.12}       
     \bean = u_1 u_2. \eean
            
            Therefore, $ A_{000} \times_{\widetilde{J}} B_{000} = B_{000} \times_{\widetilde{J}} A_{000} .$ \qed


            \blmma
            
            \label{deformation_torus_commutation_relation2}
            
            $ A_{001} \times_{\widetilde{J}} B_{001} = e ( - 2 \theta ) B_{001} \times_{\widetilde{J}} A_{001}  $
            
            \elmma

   {\it Proof :}  Proceeding exactly as in the previous Lemma, we obtain   
                  \bean \lefteqn{ ( A_{001} \times_{\widetilde{J}} B_{001} )  ( u_1, u_2, 0,0,1  )}\\
                     &=& \int.... \int A_{001} [ ( e ( - \frac{\theta}{2} s_2 ), e ( \frac{\theta}{2} s_1 ) ) ( u_1, u_2 ,0,0,1) ( e ( - \frac{\theta}{2} t_2  ), e ( \frac{\theta}{2} t_1 ) ) ] \\
                     && B_{001} [ ( e ( - s^{\prime}_1 ), e ( - s^{\prime}_2 ) ) ( u_1, u_2,0,0,1) ( e ( t^{\prime}_1  ), e ( t^{\prime}_2 )  ) ] e ( u. v ) ds dt ds^{\prime} dt^{\prime}, \eean
                     
                     which by ( \ref{deformation_laplacian_torus_001} ) equals 
                     \bean \lefteqn{ \int....\int A_{001} [ e ( - \frac{\theta}{2} s_2 ) u_1 e ( \frac{\theta}{2} t_1 ), ~ e ( \frac{\theta}{2} s_1 ) u_2 e ( - \frac{\theta}{2} t_2 ), 0,0,1 ] }\\
                     &B_{001}& [ e ( - s^{\prime}_1 ) u_1 e ( t^{\prime}_2 ), e ( - s^{\prime}_2 ) u_2 e ( t^{\prime}_1 ), 0,0,1 ] e ( u. v ) ds dt ds^{\prime} dt^{\prime} \\
        &=& u_1 u_2 \int...\int e ( - \frac{\theta}{2} s_2 ) e ( \frac{\theta}{2} t_1 ) e ( - s^{\prime}_2 ) e ( t^{\prime}_1 )  e ( u. v ) ds dt ds^{\prime} dt^{\prime} \\
        &=& u_1 u_2 \int \int e ( - \frac{\theta}{2} s_2 ) e ( - s^{\prime}_2 ) e ( s_2. s^{\prime}_2 ) ds_2 ds^{\prime}_2 \int \int e ( \frac{\theta}{2} t_1 ) e ( t^{\prime}_1 ) e ( t_1. t^{\prime}_1 ) dt_1 dt^{\prime}_1 \\
        && \int \int e ( s_1 . s^{\prime}_1 ) ds_1 ds^{\prime}_1 \int \int e ( t_2 . t^{\prime}_2 ) dt_2 dt^{\prime}_2 \eean
        
        which by Corollary \ref{preliminiaries_rieffel_2.20_corollary} and Proposition \ref{preliminiaries_rieffel_1.12}   equals
        \bean \lefteqn{ u_1 u_2 e ( - \frac{\theta}{2} ) e ( - \frac{\theta}{2} ) }\\
        &=& u_1 u_2 e ( - \theta ). \eean
        
        Similarly, using Corollary \ref{preliminiaries_rieffel_2.20_corollary}  and Proposition \ref{preliminiaries_rieffel_1.12} 
         \bean \lefteqn{ ( B_{001} \times_{\widetilde{J}} A_{001} ) (  u_1, u_2,0,0,1  )}\\ 
                  &=& \int....\int B_{001} [ e ( - \frac{\theta}{2} s_2 ) u_1 e ( \frac{\theta}{2} t_1 ), ~ e ( \frac{\theta}{2} s_1 ) u_2 e ( - \frac{\theta}{2} t_2 ), 0,0,1 ] \\
                     &&A_{001}  [ e ( - s^{\prime}_1 ) u_1 e ( t^{\prime}_2 ), e ( - s^{\prime}_2 ) u_2 e ( t^{\prime}_1 ), 0,0,1 ] e ( u. v ) ds dt ds^{\prime} dt^{\prime}\\
                     &=&  u_1 u_2 \int....\int e ( \frac{\theta}{2} s_1 ) e ( - \frac{\theta}{2} t_2 ) e ( - s^{\prime}_1 ) e ( t^{\prime}_2 ) e ( s_1 . s^{\prime}_1 + s_2 . s^{\prime}_2 + t_1 . t^{\prime}_1 + t_2 . t^{\prime}_2 ) ds dt ds^{\prime} dt^{\prime} \\                     
      &=& u_1 u_2 e ( \frac{\theta}{2}  ) e ( \frac{\theta}{2}  )\\
            &=& u_1 u_2 e ( \theta ). \eean
            
      This proves the Lemma. \qed  
      
      \brmrk
      
      \label{deformation_torus_commutation_relation3}
      
      Proceeding similarly, one can prove $ A_{0,1,1} \times_{\widetilde{J}} B_{0,1,1} = B_{0,1,1}    \times_{\widetilde{J}} A_{0,1,1}, ~ A_{1,1,0} \times_{\widetilde{J}} B_{1,1,0} = B_{1,1,0}    \times_{\widetilde{J}} A_{1,1,0}, ~ A_{1,0,1} \times_{\widetilde{J}} B_{1,0,1} = B_{1,0,1}    \times_{\widetilde{J}} A_{1,0,1}  $
      
      and $A_{0,0,1} \times_{\widetilde{J}} B_{0,0,1} = e ( - 2 \theta ) B_{0,0,1} \times_{\widetilde{J}} A_{0,0,1}, ~A_{1,0,0} \times_{\widetilde{J}} B_{1,0,0} = e ( - 2 \theta ) B_{1,0,0}    \times_{\widetilde{J}} A_{1,0,0}, ~A_{1,1,1} \times_{\widetilde{J}} B_{1,1,1} = e ( - 2 \theta ) B_{1,1,1} \times_{\widetilde{J}} A_{1,1,1}. $
      
      \ermrk

   Let us now consider  a $C^*$ algebra $\clb$, which has
eight direct summands, four of which are isomorphic with the
commutative algebra $C(\IT^2)$, and the other four are irrational
rotation algebras. More precisely, we take $$\clb=\oplus_{k=1}^8
C^*(U_{k1},U_{k2}), $$   where for odd $k$, 
$U_{k1},U_{k2}$ are the two commuting unitary generators of
$C(\IT^2)$, and  for even $k$, $U_{k1}U_{k2}={\rm e}(2 \pi i
\theta)U_{k2}U_{k1}$, that is they generate $\cla_{2 \theta}$. 

       
     \vspace{4mm}          
               
  Now we are in a position to describe $QISO^{\cll}(\cla_\theta)$ explicitly.              
          
       \bthm
             $ QISO^{\cll} ( \cla_{ \theta } )$ is isomorphic
             with $\clb = C ( \IT^{2} ) \oplus \cla_{ 2 \theta } \oplus  C ( \IT^{2} ) \oplus \cla_{ 2 \theta } \oplus C ( \IT^{2} ) \oplus \cla_{ 2 \theta } \oplus  C ( \IT^{2} ) \oplus \cla_{ 2 \theta } .$
         \ethm   
            
     {\it Proof :} Define $ \phi : \clb \rightarrow QISO^{\cll} ( \cla_{ \theta } ) $ by 
     
     
      $ \phi ( U_{11} ) = A_{000}, ~ \phi ( U_{12} ) = B_{000}, ~  \phi ( U_{21} ) = A_{010}, ~  \phi ( U_{22} ) = A_{010}, ~  \phi ( U_{31} ) = A_{011}, ~  \phi ( U_{32} ) = B_{011}, ~  \phi ( U_{41} ) = A_{001}, ~  \phi ( U_{42} ) = B_{001}, ~  \phi ( U_{51} ) = A_{110}, ~  \phi ( U_{52} ) = B_{110}, ~  \phi ( U_{61} ) = A_{111}, ~  \phi ( U_{62} ) = B_{111}, ~ \phi ( U_{71} ) = A_{101}, ~  \phi ( U_{72} ) = B_{101}, ~  \phi ( U_{81} ) = A_{100}, ~  \phi ( U_{82} ) = B_{100} .$

     From Lemma \ref{deformation_torus_commutation_relation0} - Lemma \ref{deformation_torus_commutation_relation2} and  Remark \ref{deformation_torus_commutation_relation3}, it is clear that $ \phi $ is indeed a $ \ast $ isomorphism.  
    
                         \qed

 \brmrk
   In particular, we note that if $ \theta $ is taken to be $ 1/2 $, then we have a commutative compact quantum group as the quantum isometry group of a noncommutative $ C^{\ast} $ algebra. 
                                                 \ermrk

    {\bf The case of $ \theta $ deformed spheres} 
    
    \vspace{2mm}
    
   We can apply Theorem \ref{abcd} on $ S^{n}_{\theta} $ which are obtained by Rieffel deformation  of $ C ( S^n ) $ as described in section \ref{preliminaries_section_Rieffel_deformation}. We will consider the isospectral deformation of the classical spectral triple on $ C ( S^n ) .$ Since we have proved in Theorem \ref{QISO_Laplacian_computations_sphere} that $QISO^{\cll}(S^n) \cong C(O(n))$, it will follow that $QISO^{\cll}(S^n_\theta) \cong O_\theta(n)$, where $O_\theta(n)$ is the compact quantum group discussed in subsection \ref{preliminaries_subsection_Rieffel_deformation_CQG} as the $\theta$-deformation of $C(O(n)).$

\cleardoublepage

\chapter{Quantum isometry groups of the Podles spheres}

\label{quantumsphere}

In this chapter, we compute quantum group of orientation preserving isometries for spectral triples on the Podles spheres $S^2_{\mu,c}.$ 
We do it for two different families of spectral triples, one constructed by Dabrowski et al in \cite{{Dabrowski_et_al}} 
 and the other by Chakraborty and Pal in \cite{chak_pal} for $ c > 0 .$ 
We get completely different kinds of quantum group of orientation preserving isometries for the two families; for the former, it is  $SO_\mu(3),$  
whereas, for the latter it  is $ C^* ( \IZ_2 \ast \IZ^{\infty} ) $ where $ \IZ^{\infty} $ denotes countably infinite copies of the group of integers.

\section{The Podles Spheres}
 
  The Podles spheres were discovered in \cite{podles}. We will also need the equivalent descriptions of the Podles spheres as given in \cite{Dabrowski_et_al}, \cite{klimyk} and \cite{Schmudgen_wagner_crossproduct}. Hence, we give the descriptions one by one.
  
  \vspace{4mm}
  
  \subsection{ The original definition by Podles}
  
Let $ c $ belongs to $ \IR .$ The Podles' sphere $ S^{2}_{\mu, c} $ is the universal $ C^{*} $ algebra generated by elements $ e_{-1}, e_0, e_1 $ such that :
$$ e^{*}_i = e_{-i}, ~ i = -1,0,1, $$
$$ ( 1 + {\mu}^2 ) ( e_{-1} e_1 + {\mu}^{-2} e_1 e_{-1} ) + {e_{0}}^{2} = ( {( 1 + \mu^2 )}^{2} \mu^{- 2} c + 1 ) I, $$
$$ e_0 e_{-1} - {\mu}^2 e_{-1} e_0 = ( 1 - {\mu}^2 ) e_{-1}, $$
$$ ( 1 + {\mu}^2 ) ( e_{-1} e_1 - e_1 e_{-1} ) + ( 1 - {\mu}^2 ){e_0}^2 = ( 1 - {\mu}^2 ) e_0, $$
$$ e_1 e_0 - {\mu}^2 e_0 e_1 =  ( 1 - {\mu}^2 ) e_1. $$

Let \be \label{sphere_Podles_definition_A,B_in_term_of_e_s} B = e_1, ~ A = {( 1 + \mu^2 )}^{- 1} ( 1 - e_0 ).  \ee  

Then we have an alternate description of the Podles spheres, that is, the universal $ C^{*} $ algebra generated by elements $ A $ and $ B $ satisfying the relations:
$$ A^{*} = A,~  AB = \mu^{-2} BA,$$
$$ B^* B = A - A^2 + cI, ~ B B^* = \mu^2 A - \mu^4 A^2 + cI. $$

These are the relations which we are going to exploit for getting homomorphism conditions while computing the quantum    group of orientation preserving isometries.

\vspace{4mm}

{\bf Notation :} We will denote the co-ordinate $ \ast $ algebra of $ S^{2}_{\mu,c} $ by $ \clo ( S^{2}_{\mu,c} ).$   

\vspace{4mm}

\subsection{ The Podles spheres as in \cite{Dabrowski_et_al}}

We take $ \mu $ in $ ( 0 , 1 )$ and $ t $ in $ ( 0, 1 ].$
For $n$ belonging to $ \IR, $ let $ [ n ]_{\mu} = \frac{\mu^n - \mu^{-n}}{\mu - \mu^{-1}} . $

Then $ S^{2}_{\mu, c} $ be the universal $ C^{*} $ algebra generated by elements $ x_{-1}, x_{0}, x_{1} $ satisfying the relations:
\be \label{sphere_S2_mu_c_definition_1} x_{-1} ( x_0 - t ) = \mu^2 ( x_0 - t ) x_{-1}, \ee
\be \label{sphere_S2_mu_c_definition_2} x_1 ( x_0 - t ) = \mu^{-2} ( x_0 - t ) x_1,  \ee
\be \label{sphere_S2_mu_c_definition_3} - [ 2 ] x_{-1} x_1 + ( \mu^2 x_0 + t ) ( x_0 - t ) = {[ 2 ]}^{2} ( 1 - t ), \ee
\be \label{sphere_S2_mu_c_definition_4} - [ 2 ] x_1 x_{-1} + ( \mu^{-2} x_0 + t ) ( x_0 - t ) = {[ 2 ]}^{2} ( 1 - t ), \ee
where $ c = t^{-1} - t, t > 0. $

The involution on $ S^{2}_{\mu, c} $ is given by 
\be \label{sphere_Dabrowski_et_al_involution} x^{*}_{-1} = - {\mu}^{-1} x_1, ~~ x^{*}_{0} = x_0. \ee

Setting

\be \label{sphere_A,B_S2_mu_t} A = \frac{ 1 - t^{-1} x_0}{ 1 + \mu^2}, ~~ B = \mu ( 1 + \mu^2 )^{- \frac{1}{2}} t^{-1} x_{-1}, \ee
one obtains ( \cite{Dabrowski_et_al} ) that $ S^{2}_{\mu,c} $ is the same as the Podles' sphere as in \cite{podles}. 

\vspace{4mm}

\subsection{ The Podles spheres as in \cite{klimyk}}

 This description is taken from \cite{klimyk},  page 124.

Let $ \alpha^{\prime},~ \beta $ be elements of $ \IC $ such that $ ( \alpha^{\prime}, ~ \beta ) \neq ( 0, 0 ).$ We denote by $ \chi_{q, \alpha^{\prime}, \beta } $ the algebra with three generators $x_{- 1}, ~ x_0, ~ x_1 $ and the following defining relations:
\be \label{sphere_chi_alpha_beta_definition_1} x^{2}_0 - q x_1 x_{- 1} - q^{- 1} x_{- 1} x_1 = \beta.1, \ee
\be \label{sphere_chi_alpha_beta_definition_2} ( 1 - q^2 ) x^{2}_0 + q x_{- 1} x_1 - q x_1 x_{- 1} = ( 1 - q^2 ) \alpha^{\prime} x_0, \ee
\be \label{sphere_chi_alpha_beta_definition_3} x_{- 1} x_0 - q^2 x_0 x_{- 1} = ( 1 - q^2 ) \alpha^{\prime} x_{- 1}, \ee
\be \label{sphere_chi_alpha_beta_definition_4} x_0 x_1 - q^2 x_1 x_0 = ( 1 - q^2 ) \alpha^{\prime} x_1. \ee
Let $ \rho^2 = \alpha^2 {( \beta - \alpha^2 )}^{ - 1}. $

Then for $ q $ and $ \rho $ real, the involution is defined by $ x^{*}_{- 1} = - q^{- 1} x_1, ~ x^{*}_0 = x_0, ~ x^{*}_{1} = - q x_{- 1}.$

Moreover, from page 125 of \cite{klimyk}, it follows that for $ \rho $ belonging to $ \IC, ~  \chi_{q, \alpha^{\prime}, \beta }  $ can be realized as a $ \ast $ subalgebra of $ SU_{\mu}( 2 ) $ via:
\be \label{sphere_x_-1_in_termsof_su_mu20} x_{-1} = {( 1 + q^2 )}^{- \frac{1}{2}} a^2 + \rho {( 1 + q^{- 2} )}^{\frac{1}{2}} ac - q {( 1 + q^2 )}^{- \frac{1}{2}} c^2, \ee
\be \label{sphere_x_0_in_termsof_su_mu20} x_0 = ba + \rho ( 1 + ( q + q^{- 1} ) bc   ) - cd, \ee
\be \label{sphere_x_1_in_termsof_su_mu20} x_1 = {( 1 + q^2 )}^{- \frac{1}{2}} b^2 + \rho {( 1 + q^2 )}^{ \frac{1}{2}} db - q {( 1 + q^2 )}^{- \frac{1}{2}} d^2. \ee
where $ a,b,c,d $ are as in subsection \ref{preliminaries_subsection_U_mu_su2}.

\bppsn

\label{sphere_Dabrowski=klimyk}

 $ S^{2}_{\mu,c} $ as defined above is the same as $ \chi_{q, \alpha^{\prime}, \beta } $ with $ q = \mu, ~ \alpha^{\prime} = t, ~ \beta = t^2 + \mu^{-2} {( \mu^2 + 1 )}^{2} ( 1 - t ) .$

\eppsn

{\it Proof :} We note that $ [2]_{\mu} = \frac{\mu^2 - \mu^{- 2}}{\mu - \mu^{- 1}} = \frac{\mu^2 + 1}{\mu}.$

From ( \ref{sphere_S2_mu_c_definition_4} ), we have
\be \label{sphere_S2_mu_c_definition_6}  - [2] \mu^2 x_1 x_{- 1} + x^{2}_0 + (t \mu^2 - t) x_0 - \mu^2 t^2 = {[2]}^{2} \mu^2 (1 - t). \ee
  Adding ( \ref{sphere_S2_mu_c_definition_3} ) with ( \ref{sphere_S2_mu_c_definition_6} ), we obtain 

$- [ 2 ] x_{- 1} x_1 - [2] {\mu}^2 x_1 x_{- 1} + ( \mu^2 + 1 ) x^{2}_0 - ( 1 + \mu^2 ) t^2 = {[2]}^{2} ( 1 - t ) ( 1 + \mu^2 ).$

Using $ [2]_{\mu} = \frac{1 + \mu^2}{\mu},$ we get this to be

$ - \frac{1 + \mu^2}{\mu} x_{- 1} x_1 - \frac{1 + \mu^2}{\mu} \mu^2 x_1 x_{- 1} + ( 1 + \mu^2 ) x^{2}_{0} - ( 1 + \mu^2 ) t^2 = \frac{{( 1 + \mu^2 )}^{2}}{\mu^2} ( 1 - t ) ( 1 + \mu^2 ).$

Hence, $$ x^2_0 - \mu x_1 x_{- 1} - {\mu}^{- 1} x_{- 1} x_{1} = \mu^{- 2} {( 1 + \mu^2 )}^{2} ( 1 - t ) + t^2 = \beta.1. $$

Thus, we have derived ( \ref{sphere_chi_alpha_beta_definition_1} ).

Multiplying (\ref{sphere_S2_mu_c_definition_3}) by $\mu^2,$ we have
\be \label{sphere_S2_mu_c_definition_7} - [2] \mu^2 x_{- 1} x_1 + \mu^4 x^{2}_{0} + t (1 - \mu^2) \mu^2 x_0 - t^2 \mu^2 = {[2]}^{2} \mu^2 (1 - t). \ee


Subtracting (\ref{sphere_S2_mu_c_definition_6}) from (\ref{sphere_S2_mu_c_definition_7}) gives
$$ - \frac{\mu^2 + 1}{\mu} \mu^2 x_{- 1} x_1 + (\mu^4 - 1) x^{2}_0 + t (1 - \mu^2) \mu^2 x_0 + \frac{\mu^2 + 1}{\mu} \mu^2 x_1 x_{- 1} - t (\mu^2 - 1)x_0 = 0 ,$$
and hence, $$ (1 - \mu^2) x^2_{0} + \mu x_{- 1} x_1 - \mu x_1 x_{- 1} = t (1 - \mu^2) x_0.$$
which proves ( \ref{sphere_chi_alpha_beta_definition_2} ).

Next, (\ref{sphere_S2_mu_c_definition_1}) gives $$ x_{- 1} x_0 - \mu^2 x_0 x_{- 1} = (1 - \mu^2) t x_{- 1}  $$ which is (\ref{sphere_chi_alpha_beta_definition_3}).

Finally, (\ref{sphere_S2_mu_c_definition_2}) gives $ \mu^2 x_1 x_0 - \mu^2 t x_1 = x_0 x_1 - t x_1 $, that is, (\ref{sphere_chi_alpha_beta_definition_4}) is obtained.

Thus, we obtain the relations of $ \chi_{q, \alpha^{\prime}, \beta } $ from the relations of $ S^{2}_{\mu,c} $ for $ q = \mu, ~ \alpha^{\prime} = t, ~ \beta = t^2 + \mu^{-2} {( \mu^2 + 1 )}^{2} ( 1 - t ) .$ 

Similarly, ( \ref{sphere_chi_alpha_beta_definition_3} ) is same as ( \ref{sphere_S2_mu_c_definition_1} ), ( \ref{sphere_chi_alpha_beta_definition_4} ) is same as ( \ref{sphere_S2_mu_c_definition_2} ). Subtracting ( \ref{sphere_chi_alpha_beta_definition_2} ) from ( \ref{sphere_chi_alpha_beta_definition_1} ) gives ( \ref{sphere_S2_mu_c_definition_3} ) and adding ( \ref{sphere_chi_alpha_beta_definition_1} ) with $ \mu^{- 2} $ times ( \ref{sphere_chi_alpha_beta_definition_2} ) gives  ( \ref{sphere_S2_mu_c_definition_4} ). Thus, we get back the relations of $ S^{2}_{\mu,c} $ from the relations of $ \chi_{q, \alpha^{\prime}, \beta } .$ \qed

 \qed 
 
 \vspace{4mm}

Thus, combining ( \ref{sphere_x_-1_in_termsof_su_mu20} ) - ( \ref{sphere_x_1_in_termsof_su_mu20}  ) with the correspondence ( \ref{SU_mu_2_woronowicz_klimyk_correspondence} ) and using Proposition \ref{sphere_Dabrowski=klimyk}, we have expressions of $x_{- 1}, ~ x_0, ~ x_1 $ in terms of  $ SU_{\mu}( 2 ) $ elements:
\be \label{sphere_x_-1_in_termsof_su_mu2} x_{-1} = \frac{\mu {\alpha}^{2} + \rho ( 1 + \mu^2 ) \alpha \gamma - \mu^2 \gamma^2 }{ \mu ( 1 + \mu^2 )^{\frac{1}{2}} }, \ee
\be \label{sphere_x_0_in_termsof_su_mu2} x_0 = - \mu \gamma^* \alpha + \rho ( 1 - ( 1 + \mu^2 ) \gamma^* \gamma ) - \gamma \alpha^*, \ee
\be \label{sphere_x_1_in_termsof_su_mu2} x_1 = \frac{ \mu^2 {\gamma}^{*2} - \rho \mu ( 1 + \mu^2 ) \alpha^* \gamma^* - \mu \alpha^{*2} }{( 1 + \mu^2 )^{\frac{1}{2}} }, \ee
where $ \rho^2 = \frac{\mu^2 t^2}{ {( \mu^2 + 1 )}^{2} ( 1 - t )} .$

\vspace{4mm}

\subsection{ The description as in \cite{Schmudgen_wagner_crossproduct}} 

We need the quantum group $ \clu_{\mu}( su ( 2 ) ) $ for this description.

Let
 $$ X_{c} = \mu^{\frac{1}{2}} {( \mu^{-1} - \mu )}^{-1} c^{ - \frac{1}{2}} ( 1 - K^2 ) + E K + \mu F K ~ ~ {\rm ~ for ~ all} ~ c \in ( 0, ~ \infty ), $$
$$ X_0 = 1 - K^2. $$


One has  $ \Delta ( X_c ) = 1 \otimes X_c + X_c \otimes K^2 .$

Then ( as in Page 9, \cite{Schmudgen_wagner_crossproduct} )  we have the following :
$$ \clo ( S^{2}_{\mu, c} ) = \{ x \in \clo ( SU_{\mu} ( 2 ) ) :  x \triangleleft X_c = 0 \} $$
where $ \triangleleft $ is as in subsection \ref{preliminaries_subsection_U_mu_su2}.
The following is a basis of the vector space $\clo(S^2_{\mu, c})$:
$$ \{ A^k, A^k B^l, A^k{B^*}^m,~k \geq 0,~ l,m > 0 \}.$$  So, any element of $\clo(S^2_{\mu,c})$ can be written as a {\it finite} linear combination of elements of the form $A^k, A^kB^l, A^k{B^*}^l$. 

Let $ \psi $ be the densely defined linear map on $ L^{2} ( SU_{\mu} ( 2 ) ) $ defined by $ \psi ( x ) = x \triangleleft X_c.$ 

From \cite{wagner}, ( Page 5 ),  we recall that $ v^l_{j,k} \triangleleft E = \mu \alpha^{l}_{k - 1} v^l_{j,k - 1}, ~ v^{l}_{j,k} \triangleleft F = \frac{\alpha^{l}_{k}}{\mu} v^{l}_{j,k + 1},~ v^l_{j,k} \triangleleft K = \mu^k v^l_{j,k}. $ for some constants $ \alpha^l_{j}.$


\blmma

\label{sphere_subset_ker_psi}
The map $\psi$ is closable and we have $ \overline{S^{2}_{\mu,c}} \subseteq {\rm Ker}( \overline{\psi}) $ where $ \overline{\psi} $ is the closed extension of $ \psi $ and 
 $ \overline{S^{2}_{\mu,c}} $ denotes  the Hilbert subspace generated by $ S^{2}_{\mu,c} $ in $ L^{2} ( SU_{\mu} ( 2 ) ) $ ( the G.N.S space corresponding to the Haar state on $ SU_{\mu} ( 2 ) ).$ Moreover, $\clo(S^2_{\mu,c})=\clo(SU_\mu(2)) \bigcap {\rm Ker}(\overline{\psi})=\clo(SU_\mu(2)) \bigcap {\rm Ker}(\psi).$
\elmma 

{\it Proof :} We note that $ v^l_{j,k} \triangleleft ( 1 - K^2 ) = ( 1 - \mu^{2k} ) v^{l}_{j,k}, ~ v^l_{j,k} \triangleleft EK = \mu^k \alpha^l_{k - 1} v^l_{j,k - 1}, ~ v^l_{j,k} \triangleleft F K = \mu^k \alpha^l_k v^l_{j,k+1} .$ Hence,  $ v^l_{j,k} $ belongs to $ {\rm Dom} ( \psi ) $ with $  ( \psi^* ) ( v^l_{j,k} ) = \mu^{\frac{1}{2}} {( \mu^{- 1} - \mu )}^{- 1} c^{- \frac{1}{2}}. $

$ ( 1 - \mu^{2k} ) v^l_{j,k} + \mu^k \alpha^l_{k - 1} v^l_{j,k + 1} + \mu^{k + 1} \alpha^l_{k} v^l_{j,k - 1}.$ Hence,  $ \clo ( SU_{\mu} ( 2 ) ) \subseteq {\rm Dom} ( \psi^* ) $ implying that $ \psi $ is closable, hence   $ {\rm Ker}( \overline{\psi}) $ is closed. The lemma now follows from the observation that $ \clo ( S^{2}_{\mu,c} )= {\rm Ker}( \psi) \subseteq {\rm Ker}(\overline{\psi}).$  \qed

\vspace{4mm}

 \subsection{Haar functional on the Podles spheres}

We recall from \cite{Schmudgen_wagner_crossproduct}, page- 33 that  for all  bounded complex Borel function $ f $ on $ \sigma ( A ) ,$ 

\be \label{sphere_haar_state_crossproductpaper} h ( f( A ) ) = \gamma_{+} \sum_{n = 0}^{\infty} f ( \lambda_{+} \mu^{2n} ) \mu^{2n} + \gamma_{-} \sum_{n = 0}^{\infty} f ( \lambda_{-} \mu^{2n} ) \mu^{2n} . \ee
where $ \lambda_{+} = \frac{1}{2} + {( c + \frac{1}{4} )}^{\frac{1}{2}}, ~ \lambda_{-} = \frac{1}{2} - {( c + \frac{1}{4} )}^{\frac{1}{2}},~ \gamma_{+} = ( 1 - \mu^2 ) \lambda_{+} {( \lambda_{+} - \lambda_{-} )}^{- 1}, ~ \gamma_{-} = ( 1 - \mu^2 ) \lambda_{-} {( \lambda_{-} - \lambda_{+} )}^{- 1} .$

\blmma
 
\label{sphere_haar_functional_x-1_x0_x1_orthogonal} 

 $ \{ x_{-1}, x_0, x_{1} \} $ is a set of orthogonal vectors.
 
 \elmma
 
{\it Proof :} From ( \ref{sphere_x_-1_in_termsof_su_mu2} ) and (  \ref{sphere_x_0_in_termsof_su_mu2}  ), we note that $ x^{*}_{- 1} x_0 $ belongs to $ {\rm span} \{  \alpha^{*2} \gamma^* \alpha,$

$ \alpha^{*2},   \alpha^{*2} \gamma^* \gamma,  \alpha^{*2} \gamma \alpha^*,  \gamma^* \alpha^* \gamma^* \alpha,  \gamma^* \alpha^*,  \gamma^* \alpha^* \gamma^* \gamma,  \gamma^* \alpha^* \gamma \alpha^*,  \gamma^{*3} \alpha,  \gamma^{*2}, \gamma^{*3} \gamma,  \gamma^{*2} \gamma \alpha^* \}. $

Further, using ( \ref{su2def1} ) - ( \ref{su2def5} ), we note that this  span equals $ {\rm span} \{   \alpha^* \gamma^*, ~ \alpha^* \gamma^{*2} \gamma, $ 

$~ \alpha^{*2}, ~ \alpha^{*2} \gamma^* \gamma, ~  \alpha^{*3} \gamma,~ \gamma^{*2}, ~ \gamma^{*3} \gamma,  ~ \alpha^* \gamma^*, ~ \alpha \gamma^{*3}, ~ \gamma^{*3}   \}.$

Then $ h (  x^{*}_{- 1} x_0   ) = 0 $ follows by using ( \ref{sphere_haar_state_su_mu_2_original} ). Similarly, one can prove  $ h (  x^{*}_{0} x_1   ) = h (  x^{*}_{1} x_{- 1}   ) =  0 .$ \qed 

\blmma

\label{sphere_haar_functional_A_A^2}

$$ h ( A ) = \frac{1}{1 + \mu^2}, ~~ h ( A^2 ) = \frac{(1 - \mu^2)(\lambda^{3}_{+} - \lambda^{3}_{-} )}{(\lambda_{+} - \lambda_{-})(1 - \mu^6)}. $$

\elmma

{\it Proof :} Recalling ( \ref{sphere_haar_state_crossproductpaper} ), we have
\bean \lefteqn{ h ( A )}\\
 &=& \gamma_{+} \sum^{\infty}_{n = 0} \lambda_{+} \mu^{4n} + \gamma_{-} \sum^{\infty}_{n = 0} \lambda_{-} \mu^{4n}\\
 &=& \frac{(1 - \mu^2)(\lambda^{2}_{+} - \lambda^{2}_{-})}{(\lambda_{+} - \lambda_{-}) (1 - \mu^4)} \\
 &=& \frac{\lambda_{+} + \lambda_{-}}{1 + \mu^2}\\
 &=& \frac{1}{1 + \mu^2}. \eean 
 Similarly, putting $ f ( x ) = x^2$ in ( \ref{sphere_haar_state_crossproductpaper} ), we have 
 \bean \lefteqn{h ( A^2 )}\\
 &=& \gamma_{+} \sum^{\infty}_{n = 0} \lambda^{2}_{+} \mu^{6n} + \gamma_{-} \sum^{\infty}_{n = 0} \lambda^{2}_{-} \mu^{6n}\\
 &=& \frac{(1 - \mu^2)(\lambda^{3}_{+} - \lambda^{3}_{-} )}{(\lambda_{+} - \lambda_{-})(1 - \mu^6)}. \eean \qed 

\bppsn

\label{sphere_haar_functional_x-1_x0_x1}

 $ h ( x^{*}_{-1} x_{-1} ) = h ( x^{*}_{0} x_{0} ) = h ( x^{*}_{1} x_{1} ) = t^2 ( 1 - \mu^2 ) {( 1 - \mu^6 )}^{-1} [ \mu^2 + t^{-1} ( \mu^4 + 2 \mu^2 + 1 ) + t ( - \mu^4 - 2 \mu^2 - 1 ) ]. $
 
 \eppsn

{\it Proof :} 

From ( \ref{sphere_A,B_S2_mu_t} ) we have  $ x^{*}_{-1} x_{-1} = \frac{t^2 (1 + \mu^2)}{\mu^2} B^* B $ and hence, using  Lemma \ref{sphere_haar_functional_A_A^2}, we obtain
\bean \lefteqn{ h ( x^{*}_{-1} x_{-1} )}\\
 &=& \frac{t^2 (1 + \mu^2)}{\mu^2} [ h ( A ) - h ( A^2 ) + (t^{- 1} - t).1 ]\eean
 $$= \frac{1 - \mu^6 - (1 - \mu^4) ( t^{- 1} - t + 1 ) + (t^{- 1} - t)(1 + \mu^2)(1 - \mu^6)}{(1 + \mu^2)(1 - \mu^6)},$$ 
 from which we get the desired result.
 
 Similarly, the second equality is derived from $ h ( x^{*}_{0} x_0 ) = t^2 - 2 t^2 (1 + \mu^2) h ( A ) + {(1 + \mu^2)}^{2} t^2 h ( A^2 ) .$

From ( \ref{sphere_Dabrowski_et_al_involution} ), $ x_1 = - \mu x^{*}_{- 1} $ and hence, $ x^{*}_1 x_1 = t^2 ( 1 + \mu^2 ) ( \mu^2 A - \mu^4 A^2 + c.I ),$ from which $ h ( x^*_1 x_1 ) $ is obtained and can be shown to be equal to the same value as $ h ( x^{*}_{- 1} x_{- 1} ) = h ( x^*_{0} x_0 ).$

 \qed 

\vspace{4mm}

\section{ Spectral triples on the Podles spheres}

\subsection{The spectral triple as in \cite{Dabrowski_et_al}}
We now recall the spectral triples on $ S^{2}_{\mu c} $ discussed in \cite{Dabrowski_et_al}.

Let $ s = - c^{- \frac{1}{2}} \lambda_{-}, ~ \lambda_{\pm} = \frac{1}{2} \pm {( c + \frac{1}{4} )}^{\frac{1}{2}}. $

 For all $ j $ belonging to $ \frac{1}{2} \IN,$

 $ u_j = ( {\alpha}^{*} - s \gamma^{*} ) ( \alpha^{*} - \mu^{-1} s \gamma^* )......( \alpha^* - \mu^{- 2j + 1} s \gamma^* ),$
 
 $ w_j = ( \alpha - \mu s \gamma ) ( \alpha - \mu^2 s \gamma )........( \alpha - \mu^{2j} s \gamma ), $
 
 $ u_{- j} = E^{2j} \triangleright w_j, $
 
 $ u_0 = w_0 = 1, $
 
 $ y_1 = {( 1 + \mu^{- 2} )}^{\frac{1}{2}} ( c^{\frac{1}{2}} \mu^2 \gamma^{*2} - \mu \gamma^* \alpha^* - \mu c^{\frac{1}{2}} \alpha^{*2} ), $
 
 $ N^{l}_{kj} = {\left\| F^{l - k} \triangleright ( {y_1}^{l - \left| j \right|} u_j )  \right\|}^{- 1}. $
 
 Define \be \label{sphere_Dabrowski_et_al_spectral triple} v^{l}_{k,j} = N^{l}_{k,j} F^{l - k} \triangleright ( y^{ l - \left| j \right| }_{1} u_j ), ~ l \in \frac{1}{2} \IN_{0}, ~ j,k = - l, - l + 1,......l. \ee
 
 Let $ \clm_{N} $ be the Hilbert subspace of $L^2(SU_\mu(2))$  with the orthonormal basis $ \{ v^{l}_{m,N} : l = \left| N \right|, ~ \left| N \right| + 1,~ ........,~ m = - l,.......l \}.$
 
 Set 
  $$ \clh = \clm_{- \frac{1}{2}} \oplus \clm_{\frac{1}{2}} ,$$ and define a 
   representation $ \pi $ of $ S^{2}_{\mu,c} $ on $ \clh $  by   
 \be \label{sphere_Dabrowski_et_al_representation} \pi ( x_i ) v^{l}_{m,N} = \alpha^{-}_{i} ( l, m; N ) v^{l - 1}_{m + i, N} + \alpha^{0}_{i} ( l, m; N ) v^{l}_{m + i, N} + \alpha^{+}_{i} ( l,m; N ) v^{l + 1}_{m + i,N}, \ee 
where $ \alpha^{-}_{i},~ \alpha^{0}_{i}, ~ \alpha^{+}_{i} $ are some constants.

We will often identify $ \pi ( S^{2}_{\mu,c} ) $ with $ S^{2}_{\mu,c} .$

Finally by Proposition 7.2 of \cite{Dabrowski_et_al}, the following Dirac operator $ D $ gives a spectral triple $ ( \clo(S^{2}_{\mu,c}), \clh, D ) $ which we are going to work with :
\be \label{sphere_Dabrowski_et_al_Dirac_operator} D ( v^{l}_{m, \pm \frac{1}{2}} ) = ( c_{1} l + c_2 ) v^{l}_{m, \mp \frac{1}{2}},  \ee
where $ c_1, c_2 $ are elements of $ \IR, c_1 \neq 0 .$

\vspace{4mm}

\subsection{$ SU_{\mu} ( 2 ) $ equivariance of the spectral triple}


From \cite{Dabrowski_et_al}, we see that the vector spaces $ \nu^{l}_{\pm \frac{1}{2}} = \rm{span} \{ v^{l}_{m, \pm \frac{1}{2}} : m = -l,....l \} $ are $( 2l + 1 ) $ dimensional Hilbert spaces on which the $ SU_{\mu}( 2 ) $ representation is unitarily equivalent to the standard $ l $th unitary irreducible representation of $ SU_{\mu}( 2 ) , $ that is, if the representation is denoted by $ U_0, $ then $ U_0 ( v^{l}_{i,\pm \frac{1}{2}} ) = \sum v^{l}_{j, \pm \frac{1}{2}} \otimes t^{l}_{i,j} $ where $ t^{l}_{i,j} $ denotes the matrix elements in the $ l $th  unitary irreducible representation of $ SU_{\mu} ( 2 ).$

We now recall Theorem 3.5 of \cite{goswamisuq2}.

  \bppsn
  
  \label{sphere_dabrowski_goswami_suq2_R0}
  
  Let $ R_0 $ be an operator on $ \clh $ defined by $ R_0 ( v^n_{i, \pm \frac{1}{2}  } ) = \mu^{- 2i \mp 1 } v^n_{i, \pm \frac{1}{2}}.$ Then ${\rm Tr}(R_0e^{-tD^2})<\infty ~ ({\rm ~ for ~ all}~t >0)$ and one has 
$$ (\tau_{R_0} \ot {\rm id})(\widetilde{U_0} (x \ot 1){\widetilde{U_0}}^*)=\tau_{R_0}(x).1,$$
 for all $x $ in $ \clb(\clh)$, where $\tau_{R_0}(x) = {\rm Tr}(xR_0e^{-tD^2}).$ 
  
  \eppsn

 We define a positive, unbounded operator  $ R $ on $ \clh $ by $ R ( v^{n}_{i,\pm \frac{1}{2}} ) = \mu^{- 2i} v^{n}_{i,\pm \frac{1}{2}} .$

\bppsn

\label{sphere_tau_R_=h}

$ \alpha_{U_{0}} $ preserves the $R$-twisted volume. In particular, for $x $ in $ \pi(S^2_{\mu,c})$ and $t>0$, we have $h(x)=\frac{\tau_R(x)}{\tau_R(1)}$, where $\tau_R(x):={\rm Tr}(xRe^{-tD^2})$, and $h$ denotes the restriction of the Haar state of $SU_\mu(2)$ to the subalgebra
$S^2_{\mu,c}$, which is the unique $SU_\mu(2)$-invariant state on $S^2_{\mu,c}$. 
 \eppsn
 {\it Proof :}   
It is enough to prove that $\tau_R$ is $\alpha_{U_0}$-invariant. 
 Let us denote by $ P_{\frac{1}{2}},P_{- \frac{1}{2}} $  the projections onto the closed subspaces generated by $ \{ v^{l}_{i, \frac{1}{2}}  \} $ and
 $ \{ v^{l}_{i, - \frac{1}{2}} \} $ respectively. Moreover, let $ \tau_{\pm }  $ be the functionals defined by $ \tau_{\pm } ( x ) = {\rm Tr} ( x R_0 P_{\pm \frac{1}{2}} e^{- t D^2} ) .$ Now observing that 
  $R_0$, $e^{-tD^2}$ and $U_0$ commute with $P_{\pm \frac{1}{2}}$  and using Proposition \ref{sphere_dabrowski_goswami_suq2_R0}, we have, for $x $ belonging to $ \clb(\clh),$   
\bean \lefteqn{
 ( \tau_{\pm} \otimes {\rm id} ) (\alpha_{U_0} ( x )) }\\
 & = & ( {\rm Tr} \otimes {\rm id} )( \widetilde{U_0} ( x \otimes 1) {\widetilde{U_0}}^{*} ( R_0 P_{\pm \frac{1}{2}} e^{- t D^2} \otimes {\rm id} ) )\\
   &=& ( {\rm Tr} \otimes {\rm id} ) ( \widetilde{U_0} ( x P_{\pm \frac{1}{2}} \otimes 1 ) {\widetilde{U_0}}^{*} ( R_0 e^{- t D^2} \otimes { \rm id} ) )\\
   &=& ( \tau_{R_0} \otimes {\rm id} ) \alpha_{U_0} ( x P_{\frac{1}{2}} )\\
    &=& \tau_{R_0} ( x P_{\pm \frac{1}{2}} )\\
    &=& \tau_{\pm} ( x ).1, \eean
    
    that is $\tau_{\pm}$ are $\alpha_{U_0}$-invariant.
    
      Thus,  $x \mapsto {\rm Tr}(xR_0P_{\pm \frac{1}{2}}e^{-tD^2})$ is invariant under $\alpha_{U_0}$. 
 Moreover, since we have 
 $ R P_{\pm \frac{1}{2}} = \mu^{\pm} R_0 P_{\pm \frac{1}{2}},$ the functional $\tau_{R}$ coincides with $\mu^{- 1} \tau_+ + \mu \tau_-$, hence is $\alpha_{U_0}$-invariant.   
 \qed
 
 \bthm
 
 \label{sphere_SU_mu_2_object}
 
 $ ( SU_{\mu} ( 2 ), ~ U_0 ) $ is an object in ${\bf Q}^\prime_R(D).$
 
 \ethm
 
 {\it Proof :}   \begin{align*}  ( D \otimes {\rm id} ) U_0 ( v^{l}_{i,\pm \frac{1}{2}}) 
        &= ( D \otimes {\rm id} ) \sum v^{l}_{j, \pm \frac{1}{2}} \otimes t^{l}_{i,j} &\\
        &= ( c_1 l + c_2 ) \sum v^{l}_{j, \mp \frac{1}{2}} \otimes t^{l}_{i,j}& \\
        &=  ( c_1 l + c_2 ) U_0 ( v^{l}_{i, \mp \frac{1}{2}}  )&\\
        &= U_0 D ( v^{l}_{i, \pm \frac{1}{2}}  ).& \end{align*}
 
 Thus, the above spectral triple is  equivariant w.r.t. the representation $ U_0 $ and it preserves $\tau_R$ by  Proposition \ref{sphere_tau_R_=h}, which completes the proof. \qed

\subsection{ The CQG $ SO_{\mu}( 3 ) $ and its action on the Podles sphere }

Here we recall the CQG   $ SO_{\mu}( 3 ) $ as described in \cite{podles_subgroup}.

$ SO_{\mu}( 3 ) $ is the universal unital $ C^* $ algebra generated by elements $  M,N, G, C, L $ satisfying :
\be \label{sphere_SO_mu_3_description}
 \left\{
 \begin{array}{cccc}
 L^* L = ( I - N ) ( I - \mu^{-2} N ), L L^* = ( I - \mu^2 N ) ( I - \mu^4 N ), G^* G = G G^* = N^2,\\ M^{ *} M = N - N^2,~ M M^* = \mu^2 N - \mu^4 N^2,~ C^{*} C = N - N^{2},\\ C C^{*} = \mu^{2} N - \mu^{4} N^{2},~ L N = \mu^{4} N L,~ G N  = N G,\\ M N = \mu^2 N M,~ C N = \mu^{2} N C,~ L G = \mu^{4} G L,\\~ L M = \mu^2 M L,~ M G = \mu^{2} G M,~ C M = M C,\\ L G^{*} = \mu^{4} G^{*} L,~ M^2 = \mu^{-1} L G,~ M^* L = \mu^{-1} ( I - N ) C,~ N^{*} = N .\end{array}
 \right\} 
 \ee

This CQG can be identified with a Woronowicz subalgebra of $SU_\mu(2)$ by taking:\\
$$ N = \gamma^{*} \gamma,~ M = \alpha \gamma,~ C = \alpha {\gamma}^{*},~ G = {\gamma}^{2},~ L = {\alpha}^{2} ,$$
where $ \alpha, \gamma $ are as in subsection \ref{preliminaries_subsection_SU_mu_2}. 
 
 The canonical action of $SU_\mu(2)$ on $S^2_{\mu,c}$, that is the action obtained by restricting the coproduct of $SU_\mu(2)$ to the subalgebra 
$S^2_{\mu,c}$, is actually a faithful action of $SO_\mu(3).$ 

With respect to the ordered basis $ \{ x_{-1},~ x_0,~ x_1 \} ,$ this action on the subspace generated by them is given by the following $SO_\mu(3)$-valued $3 \times 3$-matrix( \cite{klimyk} ): 

$$ \left ( \begin {array} {cccc}
  a^2  &  \frac{{( 1 + q^2 )}^{\frac{1}{2}}}{q} ab  & b^2 \\ \frac{{( 1 + q^2 )}^{\frac{1}{2}}}{q} ac  & I + ( q + q^{- 1} ) bc & \frac{{( 1 + q^2 )}^{\frac{1}{2}}}{q} bd  \\ c^2 & \frac{{( 1 + q^2 )}^{\frac{1}{2}}}{q} cd  &  d^2 \end {array} \right ) $$
( where $a,b,c,d$ are as in subsection \ref{preliminaries_subsection_U_mu_su2} ).

Recalling the correspondence in ( \ref{SU_mu_2_woronowicz_klimyk_correspondence} ) and Proposition \ref{sphere_Dabrowski=klimyk}, the matrix in the $ \alpha, ~ \gamma $ notation is :
$$ \left ( \begin {array} {cccc}
  \alpha^2  & - \mu {( 1 + \mu^{- 2} )}^{\frac{1}{2}} \alpha {\gamma}^*  &  \mu^2 \gamma^{*2} \\ {( 1 + \mu^{-2} )}^{\frac{1}{2}} \alpha \gamma  & I - \mu ( \mu + \mu^{- 1} ) \gamma^* \gamma  & - \mu {( 1 + \mu^{- 2} )}^{\frac{1}{2}} \gamma^* {\alpha}^{*}   \\ \gamma^2 & {( 1 + \mu^{-2} )}^{\frac{1}{2}} \gamma \alpha^*   &  {\alpha}^{*2} \end {array} \right ) .$$
  
  Finally, in the symbols $  M,N, G, C, L ,$ the above matrix is :  
           \be \label{sphere_matrix_of_SO_mu_3} Z_{1} = \left ( \begin {array} {cccc}
  L  & - \mu {( 1 + \mu^{- 2} )}^{\frac{1}{2}} C  &  \mu^2 G^* \\ {( 1 + \mu^{-2} )}^{\frac{1}{2}} M  & I - \mu ( \mu + \mu^{- 1} ) N  & - \mu {( 1 + \mu^{- 2} )}^{\frac{1}{2}} M^*  \\ G & {( 1 + \mu^{-2} )}^{\frac{1}{2}} C^*  &  L^* \end {array} \right ). \ee 
 As $ x_{- 1}, x_0, x_1 $ generates $ S^{2}_{\mu,c} ,$ the Woronowicz subalgebra of $ SU_{\mu} ( 2 ) $ generated by the elements of the form $<\xi \ot 1, \Delta_U(a)(\eta \ot 1)>$ is $ SO_{\mu}( 3 ) $ ( where $\xi, \eta $ are elements of $ \clh, a $ belongs to $ \clo ( S^{2}_{\mu,c} )$ and $< \cdot, \cdot >$ is the $  SU_{\mu} ( 2 ) $-valued inner product of $\clh \ot  SU_{\mu} ( 2 )$ ).

\section{Quantum Isometry Groups of the Podles sphere}

Here we will compute $ \widetilde{QISO^{+}_{R}} ( S^{2}_{\mu,c} ) $ with respect to the spectral triple given in \cite{Dabrowski_et_al} and show that it is isomorphic with $SO_\mu(3)$. 

  Let  $ (\tilde{ \clq}, U ) $ be an object in the category ${\bf Q}^\prime_R(D)$ and $\clq$ be the Woronowicz $ C^{*} $ subalgebra of $ \tilde{\clq}$ generated by $<(\xi \ot 1), \alpha_U(a) (\eta \ot 1)>_{\tilde{\clq}}$, for $\xi, \eta $ belonging to $ \clh,$ $a $ belongs to $ S^2_{\mu,c}$ (where $< \cdot, \cdot>_{\tilde{\clq}}$ is the $\tilde{\clq}$ valued inner product of $\clh \ot \tilde{\clq}$).   We shall denote $\alpha_U$ by $\phi$ from now on.

The computation has the following steps:

{\bf Step 1.} We prove that $ \phi $ is `linear', in the sense that it keeps the span of $ \{ 1, A, B, B^{*} \} $ invariant.

{\bf Step 2.} We shall exploit the facts that $\phi$ is a $\ast$-homomorphism and preserves the canonical volume form on $S^2_{\mu,c}$, that is the restriction of the Haar state of $SU_\mu(2)$. 

{\bf Step 3.} We will compute the antipode of  $ \clq  $ and apply it to get some more relations.

{\bf Step 4.} We will use the above steps to identify $ \clq $ as a subobject of $ SO_{\mu} ( 3 )$ which will finish the proof.

 \vspace{4mm}
 
 \brmrk

  The first step does not make use of the fact that $\alpha$ preserves the $R$-twisted volume, so linearity of the action follows for any object in the bigger category ${\bf Q^\prime}(D)$. 

\ermrk

\vspace{4mm}

We now note down some useful  facts for later use. 

\blmma

\label{sphere_Dabrowski_eigenspace}
   
  We observe :

1. The eigenspace of $ D $ corresponding to $ ( c_{1} l + c_2 ) $ and $ - ( c_{1} l + c_2 ) $ are $ \rm{span}  \{ v^{l}_{m, \frac{1}{2}} + v^{l}_{m, - \frac{1}{2}} : - l \leq m \leq l \} $ and span $ \{ v^{l}_{m, \frac{1}{2}} - v^{l}_{m, - \frac{1}{2}} : - l \leq m \leq l \} $ respectively.

\vspace{2mm}

2. The  eigenspace of $ \left| D \right| $ corresponding to the eigenvalue $ ( c_{1}.\frac{1}{2} + c_2 ) $ is $  \rm{span} \{ \alpha, ~ \gamma, ~ \alpha^*,$

$ \gamma^* \}. $ 

\vspace{2mm}

  
 \elmma
 
  {\it Proof :} 1. follows from ( \ref{sphere_Dabrowski_et_al_Dirac_operator}  ). To prove 2., we note that by 1, it is sufficient to identify $ \rm{span}\{ v^{\frac{1}{2}}_{- \frac{1}{2}, \frac{1}{2} }, ~ v^{\frac{1}{2}}_{ \frac{1}{2}, \frac{1}{2} }, ~ v^{\frac{1}{2}}_{- \frac{1}{2}, - \frac{1}{2} }, ~ v^{\frac{1}{2}}_{ \frac{1}{2}, - \frac{1}{2} }  \}. $
  
  Using ( \ref{sphere_Dabrowski_et_al_spectral triple}  ) and ( \ref{sphere_Uq_su2_left_action_su2} ), we have:  
  \bean \lefteqn{v^{\frac{1}{2}}_{- \frac{1}{2}, \frac{1}{2}}}\\
        &=& N^{\frac{1}{2}}_{- \frac{1}{2}, \frac{1}{2}} F \triangleright ( y^0_1 u_{\frac{1}{2}} ) ~ = ~ N^{\frac{1}{2}}_{- \frac{1}{2}, \frac{1}{2}} F \triangleright ( \alpha^* - s \gamma^* )\\
        &=& N^{\frac{1}{2}}_{- \frac{1}{2}, \frac{1}{2}} ( \gamma + \mu^{- 1} s \alpha ) ~ = ~ N^{\frac{1}{2}}_{- \frac{1}{2}, \frac{1}{2}} \gamma + \mu^{- 1} s N^{\frac{1}{2}}_{- \frac{1}{2}, \frac{1}{2}} \alpha. \eean 
        
        \bean \lefteqn{v^{\frac{1}{2}}_{\frac{1}{2},\frac{1}{2}}}\\
        &=& N^{\frac{1}{2}}_{\frac{1}{2},\frac{1}{2}} F^0 \triangleright ( y^0_1 u_{\frac{1}{2}}  ) ~ = ~ N^{\frac{1}{2}}_{\frac{1}{2},\frac{1}{2}} ( \alpha^* - s \gamma^* )\\
        &=& N^{\frac{1}{2}}_{\frac{1}{2},\frac{1}{2}} \alpha^* - s N^{\frac{1}{2}}_{\frac{1}{2},\frac{1}{2}} \gamma^*. \eean
        
        \bean \lefteqn{v^{\frac{1}{2}}_{- \frac{1}{2}, - \frac{1}{2}}}\\
        &=& N^{\frac{1}{2}}_{-\frac{1}{2}, - \frac{1}{2}} F \triangleright ( y^0_1 u_{- \frac{1}{2}} ) ~ = ~  N^{\frac{1}{2}}_{-\frac{1}{2}, - \frac{1}{2}} F \triangleright ( E \triangleright w_{\frac{1}{2}} )\\
        &=& N^{\frac{1}{2}}_{-\frac{1}{2}, - \frac{1}{2}} F \triangleright ( E \triangleright ( \alpha - \mu s \gamma  )  ) ~ = ~ N^{\frac{1}{2}}_{-\frac{1}{2}, - \frac{1}{2}} F \triangleright ( - \mu \gamma^* - \mu s \alpha^* )\\
        &=& N^{\frac{1}{2}}_{-\frac{1}{2}, - \frac{1}{2}} ( \alpha - \mu s \gamma ) ~ = ~ N^{\frac{1}{2}}_{-\frac{1}{2}, - \frac{1}{2}} \alpha - \mu s N^{\frac{1}{2}}_{-\frac{1}{2}, - \frac{1}{2}} \gamma. \eean
        
        \bean \lefteqn{v^{\frac{1}{2}}_{\frac{1}{2}, - \frac{1}{2}}}\\
        &=& N^{\frac{1}{2}}_{\frac{1}{2}, - \frac{1}{2}} F^{0} \triangleright ( y^0_	1 u_{- \frac{1}{2}} ) ~ = ~  N^{\frac{1}{2}}_{\frac{1}{2}, - \frac{1}{2}} ( E \triangleright ( \alpha - \mu s \gamma )  )\\
        &=& N^{\frac{1}{2}}_{\frac{1}{2}, - \frac{1}{2}} ( - \mu \gamma^* - \mu s \alpha^* ) ~ = ~ - \mu N^{\frac{1}{2}}_{\frac{1}{2}, - \frac{1}{2}} \gamma^* - \mu s N^{\frac{1}{2}}_{\frac{1}{2}, - \frac{1}{2}} \alpha^*. \eean
     Combining these, we have the result. \qed
 
 \blmma
 
 \label{sphere_Dabrowski_pi_A,B}
 
 $ {\rm 1.} ~  \pi ( A ) v^{l}_{m,N} \in ~  \rm{span}  \{ v^{l - 1}_{m,N}, ~ v^{l}_{m,N}, ~ v^{l + 1}_{m,N} \} ,$
  
  $ ~ \pi ( B ) v^{l}_{m,N} \in ~  \rm{span}  \{ v^{l - 1}_{m - 1,N},~ v^{l}_{m - 1,N},~ v^{l + 1}_{m - 1,N} \} ,$
  
  $ ~ \pi ( B^* ) v^{l}_{m,N} \in ~ \rm{span}  \{ v^{l - 1}_{m + 1,N},~ v^{l}_{m + 1,N},~ v^{l + 1}_{m + 1,N} \}.$
  
  \vspace{2mm}
 
 $ {\rm 2.} ~ \pi ( A^k ) ( v^{l}_{m,N} ) \in ~ \rm{span}  \{ v^{l - k}_{m,N}, ~ v^{l - k + 1}_{m,N},.......,~ v^{l + k}_{m,N} \}. $
 
 \vspace{2mm}

$ {\rm 3.} ~ \pi ( A^{m^{\prime}} B^{n^{\prime}} ) ( v^{l}_{m,N} ) \in  \rm{span}  \{ v^{l - m^{\prime} - n^{\prime}}_{m - n^{\prime},N},~ v^{l - ( n^{\prime} + m^{\prime} - 1 )}_{m - n^{\prime},N},~.......,~ v^{l + n^{\prime} + m^{\prime}}_{m - n^{\prime}, N} \}. $

\vspace{2mm}

 $ {\rm 4.} ~ \pi ( A^{r} B^{*s} ) ( v^l_{m,N} ) \in ~ \rm{span}  \{ v^{l - s - r}_{m + s,N},~ v^{l - s - r + 1}_{m + s,N},.......v^{l + s + r}_{m + s,N} \}. $ 
 
 \vspace{2mm}
 
 \elmma
 
{\it Proof :} Using ( \ref{sphere_A,B_S2_mu_t} ) and ( \ref{sphere_Dabrowski_et_al_representation} ), we have
 \bean \lefteqn{ \pi ( A ) v^{l}_{m, N}}\\
           &=& \frac{1}{1 + \mu^2} v^{l}_{m, N} - \frac{t^{- 1}}{1 + \mu^2} [ \alpha^{-}_0 ( l,m; N ) v^{l - 1}_{m,N} + \alpha^{0}_0 ( l,m; N ) v^{l}_{m,N} + \alpha^{+}_0 ( l,m;N ) v^{l + 1}_{m,N} ]. \eean
           
  Thus, $ \pi ( A ) v^{l}_{m,N} $ belongs to $ \rm{span} \{ v^{l}_{m,N} , ~ v^{l - 1}_{m,N}, ~ v^{l + 1}_{m,N} \}.$
  
  Similarly, using the expressions for $ B $ and $ B^* $ from  ( \ref{sphere_A,B_S2_mu_t} ) and then using ( \ref{sphere_Dabrowski_et_al_representation} ) just as above, we get the required statements for $ \pi ( B )  $ and $ \pi ( B^* ).$ This proves 1.
  
  Repeated use of 1. now yields 2.,3. and 4. \qed

\subsection{Linearity of the action}

For a vector $v $ in $ \clh,$ 
we shall denote by $T_v$ the map from $\clb(\clh)$   to $L^2(SU_\mu(2))$ given by $T_v(x)=xv \in \clh \subset L^2(SU_\mu(2))$.  It is clearly a continuous map w.r.t. the SOT on $\clb(\clh)$ and the Hilbert space topology of $L^2(SU_\mu(2))$. 


  For an element $a $ in $ SU_\mu(2),$ we consider the right multiplication $R_a $ as a bounded linear map on $L^2(SU_\mu(2))$. Clearly the composition $R_a T_v$ is a continuous linear map from $\clb(\clh) $ (with SOT) to the Hilbert space $L^2(SU_\mu(2))$.  We now define  



$$ T=R_{\alpha^*}T_\alpha+\mu^2R_\gamma T_{\gamma^*} .$$ 


%






\blmma

\label{sphere_T_tilda}
For any state $\omega$ on $\tilde{\clq}$ and $x $ belonging to $ S^2_{\mu, c}$, we have $T(\phi_\omega(x))=\phi_\omega(x) \equiv R_1(\phi_\omega(x)) \in \overline{S^2_{\mu,c}} \subseteq L^2(SU_\mu(2)),$ where $\phi_\omega(x)=({\rm id} \ot \omega)(\phi(x)).$

\elmma

{\it Proof :} It is clear from the definition of $T$  (using $\alpha \alpha^*+\mu^2 \gamma \gamma^*=1$) that $T(x) =x \equiv R_1(x) $ for $x $ in $ S^2_{\mu,c} \subset \clb(\clh)$, where $x$ in the right hand side of the above denotes the identification of $x $ in $ S^2_{\mu,c}$ as a vector in $L^2(SU_\mu(2))$.  Now, the lemma follows by noting that for $ x $ belonging to $ S^{2}_{\mu,c},   ~ \phi_\omega(x)$ belongs to  $(S^2_{\mu,c})^{\prime \prime}$, which is  the SOT closure of $S^2_{\mu, c}$, and the SOT continuity of $T$ discussed before. 
 \qed

\vspace{4mm}
     
     Let     
      $$ \clv^{l} = {\rm Span} \{ v^{l^{\prime}}_{i, \pm \frac{1}{2}}, -l^{\prime} \leq i \leq l^{\prime},~ l^{\prime} \leq l \}.$$ As $ {\rm Span} \{ v^{l}_{i, \pm \frac{1}{2}}, -l \leq i \leq l \} ,$ is the eigenspace of $|D|$ corresponding to the eigenvalue $c_1l+c_2$, $\tilde{U}$ and $\tilde{U}^*$ keep  $\clv^l$ invariant  for all $ l.$

 \blmma
 
 \label{sphere_Dabrowski_linearity_1}
 
  There is some finite dimensional subspace $\clv$ of $\clo(SU_\mu(2))$  such that $ {R}_{\alpha^{*}} (\phi_\omega ( A )v^{\frac{1}{2}}_{j,\pm \frac{1}{2}}), {R}_{\gamma}(\phi_\omega ( A )v^{\frac{1}{2}}_{j,\pm \frac{1}{2}}) $  belong to $\clv$ for all states $\omega$ on $\tilde{\clq}$.
 
 The same holds when $ A $ is replaced by $ B$ or $B^{*}.$ 
 
 \elmma

{\it Proof :}
We prove the result for $A$ only, since a similar argument will work for $B$ and $B^*$.

We have
$ \phi ( A )( v^{\frac{1}{2}}_{j, \pm \frac{1}{2}} \otimes 1 ) = \widetilde{U} ( \pi ( A ) \otimes 1 ) {\widetilde{U}}^{*} ( v^{\frac{1}{2}}_{j, \pm \frac{1}{2}} \otimes 1 ) .$ 

Now,  $ {\widetilde{U}}^{*} ( v^{\frac{1}{2}}_{j, \pm \frac{1}{2}} \otimes 1 ) $ belong to $ \clv^{\frac{1}{2}} \ot \tilde{\clq}$, and then 
using the definition of $ \pi $ as well as the  Lemma \ref{sphere_Dabrowski_pi_A,B}, we get $ ( \pi ( A ) \otimes 1 ) {\widetilde{U}}^{*} ( v^{\frac{1}{2}}_{j, \pm \frac{1}{2}} \otimes 1 ) $ belong to $ {\rm Span}\{ v^{l^{\prime}}_{j, \pm \frac{1}{2}} : - l^{\prime} \leq  j \leq l^{\prime}, l^{\prime} \leq \frac{3}{2}   \} \ot \tilde{\clq}=\clv^{\frac{3}{2}} \ot \tilde{\clq}.$ 
Again, $\tilde{U}$ keeps $\clv^{\frac{3}{2}} \ot \tilde{\clq}$ invariant, so  $R_{\alpha^*}(\phi_\omega(A)( v^{\frac{1}{2}}_{j,\pm \frac{1}{2}})) $ belong to $ {\rm Span}\{ v \alpha^*, v \in \clv^{\frac{3}{2}} \}$. Similarly, $R_\gamma(\phi_\omega(A)(v^{\frac{1}{2}}_{j,\pm \frac{1}{2}})) $ belong to $ {\rm Span}\{ v \gamma : v \in \clv^{\frac{3}{2}} \}.$ 
So, the lemma follows for $A$ by taking $\clv ={\rm Span} \{ v\alpha^*, v \gamma : ~v \in \clv^{\frac{3}{2}} \} \subset \clo(SU_\mu(2)).$ \qed

 %
%

\vspace{4mm}

Since $\alpha, \gamma^* $ are in $ {\rm Span}\{ v^{\frac{1}{2}}_{j,\pm \frac{1}{2}} \},$ we have the following immediate corollary:
\bcrlre
\label{dab_lin}
There is a finite dimensional subspace $\clv$ of $\clo(SU_\mu(2))$ such that for every state (hence for every bounded linear functional) $\omega$ on $\tilde{\clq}$, we have $T(\phi_\omega(A)) $ belongs to $ \clv.$ Similar conclusion holds for $B$ and $B^*$ as well. 
\ecrlre

\bppsn

\label{dab_lin1}

$\phi(A),~ \phi(B),~ \phi(B^*)$ belong to $\clo(S^2_{\mu,c} )\otimes_{\rm alg} \clq$. 

\eppsn

{\it Proof :} 
We give the proof for $\phi(A)$ only, the proof for $B,B^*$ being similar. 
From the Corollary  \ref{dab_lin} and Lemma \ref{sphere_T_tilda} it follows that  for every bounded linear functional $\omega$ on $\tilde{\clq}$, we have $T(\phi_\omega(A)) $ belongs to $ \clv \bigcap \overline{S^2_{\mu,c}} \subset \clo(SU_\mu(2)) \bigcap {\rm Ker}(\psi) $ ( by Lemma \ref{sphere_subset_ker_psi}   ) and hence $ \clv \bigcap \overline{S^2_{\mu,c}} = \clv \bigcap \clo(S^2_{\mu,c})$, where $\clv$ is the finite dimensional subspace mentioned in Corollary \ref{dab_lin}.  Clearly, $\clv \bigcap \clo(S^2_{\mu,c})$ is a finite dimensional subspace of $\clo(S^2_{\mu,c})$ implying that there must be finite $m$, say, such that  for every $\omega$, $ T(\phi_\omega ( A )) $belongs to $ {\rm Span} \{ A^k,~ A^k B^l, ~ A^k B^{*l} : 0 \leq k,l \leq m  \}.$ Denote by ${\cal W}$ the (finite dimensional) subspace of $\clb(\clh)$ spanned by $ \{ A^k,~ A^k B^l, ~ A^k B^{*l} : 0 \leq k,l \leq m \} $.  
Since for every state (and hence for every bounded linear functional) $\omega$ on $\tilde{\clq}$, we have $T(\phi_\omega(A))=R_1(\phi_\omega(A)) \equiv \phi_\omega(A).1 $, it is clear that $\phi_\omega(A) $ is in $ {\cal W} $ for every $\omega $ in $ {\tilde {\clq}}^*$.  Now, let us fix any faithful state  $\omega$ on the separable unital $C^*$-algebra $\tilde{\clq}$ and embed $\tilde{\clq}$ in $\clb(L^2(\clq, \omega))\equiv \clb(\clk)$. Thus, we get a canonical embedding of $\cll(\clh \ot \tilde{\clq})$ in $\clb(\clh \ot \clk)$.  Let us thus identify $\phi(A)$ as an element of $\clb(\clh \ot \clk)$, and then by choosing a countable family of elements $\{ q_1, q_2,... \}$ of $\tilde{\clq}$ which is an orthonormal basis in $\clk=L^2(\omega)$, we can write $\phi(A)$ as a  weakly convergent series of the form $\sum_{i,j=1}^\infty \phi^{ij}(A) \ot |q_i><q_j|$. But $\phi^{ij}(A)=
({\rm id} \ot \omega_{ij})(\phi(A))$, where $\omega_{ij}(\cdot)=\omega(q_i^* \cdot q_j)$. So we have $\phi^{ij}(A) $ belongs to $ {\cal W}$ for all $i,j$, and hence the sequence $\sum_{i,j=1}^n \phi^{ij}(A) \ot |q_i><q_j| $ belongs to $ {\cal W} \ot \clb(\clk)$  converges weakly, and ${\cal W}$ being finite dimensional ( hence weakly closed ), the limit, that is $\phi(A)$, must belong to ${\cal W} \ot \clb(\clk)$. In other words, if $A_1,..., A_k$ denotes a basis of ${\cal W}$, we can write $\phi(A)=\sum_{i=1}^k A_i \ot B_i$ for some $B_i $ in $ \clb(\clk).$ 

We claim that each $B_i$ must belong to $\tilde{\clq}$. For any trace-class positive operator $\rho$ in $\clh$, say of the form $\rho=\sum_j \lambda_j |e_j><e_j|$, where $\{ e_1, e_2,,...\}$ is an orthonormal basis of $\clh$ and $\lambda_j \geq 0, \sum_j \lambda_j< \infty$,  let us denote by $\psi_\rho$ the normal functional on $\clb(\clh)$ given by $x \mapsto {\rm Tr}(\rho x)$, and it is easy to see that it has a canonical extension $\tilde{\psi}_\rho:=(\psi_\rho \ot {\rm id})$ on $\cll(\clh \ot \tilde{\clq})$ given by $\tilde{\psi}_\rho(X)=\sum _j \lambda_j < e_j \ot 1, X(e_j \ot 1)>_{\tilde{\clq}}$, where $X $ is in $ \cll(\clh \ot \tilde{\clq})$ and $<\cdot, \cdot, >_{\tilde{\clq}}$ denotes that $\tilde{\clq}$ valued inner product of $\clh \ot \tilde{\clq}$. Clearly, $\tilde{\psi}_\rho$ is a bounded linear map from $\cll(\clh \ot \tilde{\clq})$ to $\tilde{\clq}$. Now, since $A_1,..., A_k$ in the expression of $\phi(A)$ are linearly independent, we can choose trace class operators $ \rho_1,..., \rho_n $ such that $\psi_{\rho_i}(A_i)=1$ and $\psi_{\rho_i}(A_j)=0$ for $j \neq i$. Then, by applying $\tilde{\psi}_{\rho_i}$ on $\phi(A)$ we conclude that $B_i $ belongs to $ \tilde{\clq}$. But by definition of $\clq$ as the Woronowicz subalgebra of $\tilde{\clq}$ generated by $<\xi \ot 1, \phi(x)(\eta \ot 1)>_{\tilde{\clq}}$, with $\eta, \xi $ belonging to $ \clh$, we must have $B_i $ belongs to $ \clq.$ 
 \qed

\bppsn

\label{sphere_Dabrowski_linearity_2}

$ \phi $ keeps the span of $ 1, A, B, B^* $ invariant.

\eppsn 

{\it Proof :} We prove the result for $ \phi ( A ) .$ The proof for the others are exactly similar.

Using Proposition \ref{dab_lin1}, we can write $ \phi ( A ) $ as a finite sum of the form $ \sum_{k \geq 0} A^k \otimes Q_k + \sum_{m^{\prime}, n^{\prime}, n^{\prime} \neq 0 } A^{m^{\prime}} B^{n^{\prime}} \otimes R_{m^{\prime},n^{\prime}} + \sum_{r,s, s \neq 0} A^r B^{*s} \otimes R^{\prime}_{r,s} .$

Let $ \xi = v^l_{m_0,N_{0}}. $

We have $ U ( \xi ) $ belongs to $ {\rm Span}  \{ v^{l}_{m,N}, m = - l,......l, ~ N = \pm \frac{1}{2} \} $.
 Let us write $$ \widetilde{U}( \xi \otimes 1 ) = \sum_{m = -l,....l, N = \pm \frac{1}{2} } v^{l}_{m, N} \otimes q^{l}_{(m, N),(m_{0},N_{0})},$$ where $q^{l}_{(m, N),(m_{0},N_{0})} $ belong to $ \clq.$ 
 Since $\alpha_U$ preserves the $R$-twisted volume, we have : \be \label{rvol111} \sum_{m^{\prime},N^{\prime}} q^l_{(m,~N),~(m^{\prime},~N^{\prime})} q^{l*}_{(m,N),~(m^{\prime},~ N^{\prime})} = 1 .\ee
  It also follows that $ U ( A \xi ) $ belong to $ {\rm Span} \{ v^{l^{\prime}}_{m,~N}, m = - l^{\prime},........l^{\prime}, l^{\prime} = l - 1, l, l + 1, ~ N = \pm \frac{1}{2} \} .$

Recalling Lemma \ref{sphere_Dabrowski_pi_A,B}, we have $ \phi ( A ) \widetilde{U}( \xi \otimes 1 )
 = \sum_{k, ~ m = - l,....l, N =
 \pm \frac{1}{2} } A^k v^{l}_{m,~ N} \otimes Q_k q^{l}_{(m,~N),~(m_{0},N_{0})} + \sum_{m^{\prime},~ n^{\prime}, n^{\prime} \neq 0, ~ m = - l,....l, N = \pm \frac{1}{2} }
 A^{m^{\prime}} B^{n^{\prime}} v^{l}_{m, N} \otimes R_{m^{\prime},~ n^{\prime}}  q^{l}_{(m,N),~(m_{0},N_{0})} $ 
 
 $ + \sum_{r,s, ~ s \neq 0, ~ m = -l,...l,~ N = 
\pm \frac{1}{2}} A^r B^{*s} v^{l}_{m, N} \otimes R^{\prime}_{r,s} q^{l}_{(m,N),~(m_{0},N_{0})}.$

Let $ m^{\prime}_{0}$ denote the  largest integer $ m^{\prime}$ such that there is a nonzero coefficient of $ A^{m^{\prime}} B^{n^{\prime}}, n^\prime \geq 1 $ 
 in the expression of $ \phi ( A ).$
We claim that the coefficient of $ v^{l - m^{\prime}_{0} - n^{\prime}}_{m - n^{\prime}, N} $ in 
$ \phi ( A ) \widetilde{U} ( \xi \otimes 1 )$ is $R_{m^{\prime}_{0},n^{\prime}}q^l_{(m,N),~(m_{0},N_{0})}. $

Indeed, the term $ v^{l - m^{\prime}_{0} - n^{\prime}}_{m - n^{\prime},N} $  
can arise in three ways: it can come from a term of the form $ A^{m^{''}} B^{n^{''}}v^l_{m,N} $ or $ A^kv^l_{m,N} $ or $ A^{r} B^{*s} v^l_{m.N}$ 
for some  $ m^{''}, ~ n^{''},~ k,~ r,~ s .$

In the first case, using Lemma \ref{sphere_Dabrowski_pi_A,B}, we must have $ l - m^{\prime}_{0} - n^{\prime} = l - m^{''} - n^{''} + t, ~ 0 \leq t \leq 2 m^{''} $ and $ m - n^{\prime} = m - n^{''} $ 
implying $ m ^{''} = m^{\prime}_{0} + t $, and since  $ m^{\prime}_{0} $ is the largest integer such that such that $ A^{m^{\prime}_{0}} B^{n^{\prime}} $ 
appears in $ \phi ( A ) ,$ we only have the possibility $t=0$, that is $ v^{l - m^{\prime}_{0} - n^{\prime}}_{m - n^{\prime}, N} $ appears only in 
$ A^{m^{\prime}_{0}} B^{n^{\prime}}.$

In the second case, we have $ m - n^{\prime} = m $ implying $ n^{\prime} = 0 .$ - a contradiction.
In the last case, we have $ m - n^{\prime} = m + s $ so that $ - n^{\prime} = s $ which is only possible when $ n^{\prime} = s = 0 $ which is again a contradiction.

Now, coefficient of $ v^{l - m^{\prime}_0 - n}_{m - n^{\prime}, N}  $ in $ \widetilde{U} ( A \xi \otimes 1 ) $ is zero if $ m^{\prime}_0 \geq 1 $ ( as $ n^{\prime} \neq 0 $ ).  It now follows from the above claim, using Lemma \ref{sphere_Dabrowski_pi_A,B} and comparing coefficients in the equality 
$ \widetilde{U} ( A \xi \otimes 1 ) =  \phi ( A ) \widetilde{U} ( \xi \otimes 1 ) $, that  
 $ R_{m^{\prime}_{0}, n^{\prime}} q^{l}_{(m,N),~(m_{0},N_{0})} = 0 $ for all $ n^{\prime} \geq 1,$  for all $ m,N $ when $ m^{\prime}_{0} \geq 1 .$ Now varying $ ( m_{0},~ N_{0} ) ,$ we conclude that the above holds  for all $ ( m_{0},~ N_{0} ).$
  Using (\ref{rvol111}), we conclude that
  
   $ R_{m^{\prime}_{0}, n^{\prime}}  \sum_{m^{\prime},N^{\prime}} q^l_{(m,N), ~(m^{\prime},N^{\prime})} q^{l*}_{(m,~N),(m^{\prime},~N^{\prime})} = 0 $ for all $ n^{\prime} \geq 1 ,$
   
    that is, $ R_{m^{\prime}_{0},~ n^{\prime}} = 0 $  for all $ n^{\prime} \geq 1$ 
 if $ m^{\prime}_{0} \geq 1 .$   Proceeding by induction on $ m^{\prime}_{0},$ we deduce  $ R_{m^{\prime},~n^{\prime}} = 0 $ for all $ m^{\prime} \geq 1,~ n^{\prime} \geq 1 .$

Similarly, we have $ Q_k = 0 $ for all $ k \geq 2 $ and $ R^{\prime}_{r,~s} = 0 $ for all $ r \geq 1,~ s \geq 1.$ 

Thus,  $ \phi ( A ) $ belongs to  $ {\rm span} \{ 1, A, B, B^{*}, B^{2},.......B^{n},B^{*2},.....B^{*m} \} .$
But the coefficient of $ v^{l - n^{\prime}}_{m - n^{\prime},~N} $ in $ \phi ( A ) \widetilde{U} ( \xi \otimes 1 )$ is $ R_{0, n^{\prime}} .$ 
Arguing as before,  we conclude that $ R_{0,~n^{\prime}} = 0 $  for all $ n^{\prime} \geq 2 .$ In a similar way, we can prove $ R^{\prime}_{0,~n^{\prime}} = 0 $ for all $ n^{\prime} \geq 2 .$
\qed

 In view of the above, let us write:  
 \be \label{alphaA} \phi( A ) = 1 \otimes T_{1} + A \otimes T_{2} + B \otimes T_{3} + B^{*} \otimes T_{4}, \ee
 \be \label{alphaB}  \phi( B ) = 1 \otimes S_{1} + A \otimes S_{2} + B \otimes S_{3} + B^{*} \otimes S_{4}, \ee
 for some $ T_i, S_i $ in $ \clq.$

   \subsection{Homomorphism conditions}

In this subsection, we shall use the facts  that $ \phi $ is a $\ast$-homomorphism and it preserves the $R$-twisted volume to derive relations among $T_i, S_i$ in ( \ref{alphaA} ), ( \ref{alphaB} ).

 \blmma
 
 \label{h_invariance}
 
 $$ T_1 = \frac{1 - T_2}{1 + \mu^2}, $$ 
 $$ S_1 = \frac{- S_2}{1 + \mu^2}. $$

 \elmma

{\it Proof :} We have the expressions of $ A $ and $ B $ in terms of the $ SU_{\mu} ( 2 ) $ elements from the equations ( \ref{sphere_x_-1_in_termsof_su_mu2} ), ( \ref{sphere_x_0_in_termsof_su_mu2} ) and  ( \ref{sphere_x_1_in_termsof_su_mu2} ). From these, we note that $ h ( A ) = {( 1 + \mu^2 )}^{-1} $ and $ h ( B ) = 0.$ By recalling Proposition \ref{sphere_tau_R_=h}, we use $ ( h \otimes {\rm id} ) \phi ( A ) = h ( A ).1 $ and $ ( h \otimes {\rm id} ) \phi ( B ) = h ( B ).1 $ to have the above two equations. \qed

 \blmma 
 
 \label{spherehomomorphism A* = A}
 
$ T^{*}_{1} = T_{1},~ T^{*}_{2} = T_{2},~ T^{*}_{4} = T_{3}. $

\elmma

{\it Proof :} Follows by comparing the coefficients of $ 1, A $ and $ B $ respectively in the equation $ \phi( A^{*} ) = \phi( A ) .$ \qed

\vspace{4mm}

We shall now assume that $ \mu \neq 1 .$ The case $ \mu = 1 $ will be discussed separately. 

\blmma

\label{spherehomomorphism B*B = A - A2}

\bean  S^{*}_{2}S_{2} + c ( 1 + \mu^2 )^{2} S^{*}_{3} S_{3} + c {( 1 + \mu^2 )}^{2} S^{*}_{4} S_4  \eean 
\be \label{spherehomomorphism B*B = A - A2 1} = ( 1 - T_2 ) ( \mu^2 + T_2 )  - c ( 1 + \mu^2 )^{2} T_3 T^{*}_3 - c {( 1 + \mu^2 )}^{2} T^{*}_{3} T_3 + c ( 1 + \mu^2 )^{2}.1, \ee
\be \label{spherehomomorphism B*B = A - A2 2} - 2 S^{*}_{2}S_{2} + ( 1 + \mu^2 ) S^{*}_{3}S_{3} + {\mu}^{2} ( 1 + \mu^2 ) S^{*}_{4}S_{4} = ( \mu^2 + 2 T_2 - 1 ) T_{2} - \mu^2 ( 1 + \mu^2 ) T_3 T^*_3 - ( 1 + \mu^2 ) T^*_3 T_3, \ee
\be \label{spherehomomorphism B*B = A - A2 3} S^{*}_{2}S_{2} - S^{*}_{3}S_{3} - {\mu}^{4} S^{*}_{4}S_{4} = - T^{2}_{2} + {\mu}^{4}T_{3}T^*_{3} + T^*_{3}T_{3}, \ee
\be \label{spherehomomorphism B*B = A - A2 5} S^{*}_{2}S_{4} + S^{*}_{3}S_{2} = - ( \mu^2 + T_2 )T^*_{3} + T^*_3 ( 1 - T_2 ), \ee
\be \label{spherehomomorphism B*B = A - A2 6} S^{*}_{2}S_{3} + {\mu}^{2}S^{*}_{4}S_{2} = - T_{2}T_{3} - {\mu}^{2}T_{3}T_{2},  \ee
\be \label{spherehomomorphism B*B = A - A2 7}  S^{*}_{4}S_{3} = - T^{2}_{3}. \ee

\elmma

{\it Proof :} It follows by comparing the coefficients of $ 1, A, A^{2}, B^{*}, AB $ and $ B^{2} $   in the equation $ \phi( B^{*} B ) = \phi( A ) - \phi( A^{2} ) + cI $ and then substituting $ S_1, ~ T_1, ~ T^*_2, T_4 $ by $ \frac{- S_2}{1 + \mu^2}, ~ \frac{1 - T_2}{1 + \mu^2}, ~ T_2, ~ T^*_3 $ respectively by  using the relations in  Lemma \ref{h_invariance}   and Lemma \ref{spherehomomorphism A* = A}. \qed


\blmma

\label{spherehomomorphism BA = mu2 AB}

 \bean - S_{2}( 1 - T_{2} ) + c {( 1 + \mu^2 )}^{2} S_3 T^{*}_3 + c {( 1 + \mu^2 )}^{2} S_4 T_3 \eean 
 \be \label{spherehomomorphism BA = mu2 AB 1} =  - {\mu}^{2} ( 1 -  T_{2} ) S_{2} + c \mu^2 {( 1 + \mu^2 )}^{2} T_3 S_4 + c \mu^2 ( 1 + \mu^2 )^{2} T^{*}_{3} S_3, \ee 
 \be \label{spherehomomorphism BA = mu2 AB A} S_{2} - 2 S_2 T_{2} + ( 1 + \mu^2 ) ( \mu^2 S_3 T^*_3 + S_4 T_3 ) = \mu^2 S_2 - 2\mu^2 T_2 S_2 + \mu^4 ( 1 + \mu^2 ) T_3 S_4 + \mu^2 ( 1 + \mu^2 ) T^*_3 S_3,  \ee 
 \be \label{spherehomomorphism BA = mu2 AB B} - S_{2}T_{3} + S_{3}( 1 - T_{2} ) = - {\mu}^{2} T_{3}S_{2} + \mu^2 ( 1 - T_2 )S_3, \ee 
 \be \label{spherehomomorphism BA = mu2 AB B*} - S_{2}T^*_{3} + S_{4}( 1 - T_{2} ) =  {\mu}^{2} ( 1 - T_2 )S_4  - \mu^2 T^*_{3}S_{2},  \ee 
 \be \label{spherehomomorphism BA = mu2 AB AB} S_{2}T_{3} + {\mu}^{2} S_{3}T_{2} = {\mu}^{2}( T_{2}S_{3} + {\mu}^{2} T_{3}S_{2} ), \ee 
 \be \label{spherehomomorphism BA = mu2 AB B2}  S_{3}T_{3} = {\mu}^{2} T_{3} S_{3}, \ee 
 \be \label{spherehomomorphism BA = mu2 AB B*2} S_{4}T^*_{3} = {\mu}^{2} T^*_{3} S_{4}. \ee
 
 \elmma
 
 {\it Proof :} It follows by equating the coefficients of $ 1, A, B, B^{*}, AB, B^{2} $ and $ {B^{*}}^{2}  $ in the equation $ \phi( BA ) = {\mu}^{2} \phi( AB ) $ and then using Lemma \ref{h_invariance}   and Lemma \ref{spherehomomorphism A* = A}. \qed 
 
 \blmma
 
 \label{spherehomomorphism BB* = mu2 A - mu4 A2}
 
 \bean  S_{2}S^{*}_{2} + c ( 1 + \mu^2 )^{2} S_3 S^{*}_3 + c ( 1 + \mu^2 )^{2} S_4 S^{*}_{4} \eean 
 \be \label{spherehomomorphism BB* = mu2 A - mu4 A2 1} = {\mu}^{2} ( 1 -  T_2 ) ( 1 + \mu^2 T_2 ) + c ( 1 + \mu^2 )^{2} T_3 T^{*}_3 + c ( 1 + \mu^2 )^{2} T^{*}_3 T_3 + c ( 1 + \mu^2 )^{2}.1, \ee 
 \bean - 2 S_{2}S^{*}_{2} + \mu^2 ( 1 + \mu^2 ) S_{3}S^{*}_{3} + ( 1 + {\mu}^{2} ) S_{4}S^{*}_{4} \eean 
 \be \label{spherehomomorphism BB* = mu2 A - mu4 A2 A} = {\mu}^{2} ( 1 + \mu^2 ) T_{2} - 2 {\mu}^{4} ( 1 - T_{2} )T_{2} - {\mu}^{6} ( 1 + \mu^2 ) T_3 T^*_3 - \mu^4 ( 1 + \mu^2 ) T^*_{3}T_{3}, \ee 
 \be \label{spherehomomorphism BB* = mu2 A - mu4 A2 B}  - S_{2}S^{*}_{4} - S_{3}S^{*}_{2} = {\mu}^{2} ( 1 + \mu^2 ) T_{3} - {\mu}^{4} ( 1 - T_{2} ) T_{3} - {\mu}^{4} T_{3} ( 1 - T_{2} ), \ee 
 \be \label{spherehomomorphism BB* = mu2 A - mu4 A2 A2} S_{2}S^{*}_{2} - {\mu}^{4} S_{3}S^{*}_{3} - S_{4}S^{*}_{4}  = - {\mu}^{4} T^{2}_{2} + {\mu}^{8}T_{3}T^*_{3} + {\mu}^{4} T^*_{3}T_{3},  \ee 
 \be \label{spherehomomorphism BB* = mu2 A - mu4 A2 AB} S_{2}S^{*}_{4} + {\mu}^{2} S_{3}S^{*}_{2} = - {\mu}^{4} T_{2}T_{3} - {\mu}^{6} T_{3}T_{2}, \ee 
 \be \label{spherehomomorphism BB* = mu2 A - mu4 A2 B2} S_{3}S^{*}_{4} = - {\mu}^{4} T^{2}_{3}. \ee

 \elmma
 
 {\it Proof :} The Lemma is proved by equating the coefficient of $1,~A,~ B,~A^2,~ AB,~  B^{2} $  in the equation $ \phi( BB^{*} ) = {\mu}^{2}\phi( A ) - {\mu}^{4}\phi( A^{2} ) + c.1 $ and then using Lemma \ref{h_invariance}   and Lemma \ref{spherehomomorphism A* = A}. \qed


\vspace{8mm}

\subsection{Relations from the antipode}

 Now, we compute the antipode, say $ \kappa $ of $ \widetilde{\clq} .$

To begin with, we recall from Lemma \ref{sphere_haar_functional_x-1_x0_x1_orthogonal} and Proposition \ref{sphere_haar_functional_x-1_x0_x1} that $ \{ x_{-1}, x_0, x_{1} \} $ is a set of orthogonal vectors with same norm.

\blmma

If $ x^{\prime}_{-1},~ x^{\prime}_{0},~ x^{\prime}_{1} $ is the normalized basis corresponding to $ \{ x_{-1},~ x_0,~ x_{1} \} $, then from ( \ref{alphaA} ) and ( \ref{alphaB} ) we obtain

$ \phi ( x^{\prime}_{-1} ) = x^{\prime}_{-1} \otimes S_3 + x^{\prime}_{0} \otimes - {\mu}^{-1} {( 1 + {\mu}^{2} )}^{- \frac{1}{2}} S_2 + x^{\prime}_{1} \otimes - {\mu}^{-1} S_4,  $

$ \phi ( x^{\prime}_{0}  ) = x^{\prime}_{-1} \otimes - \mu {( 1 + \mu^2 )}^{\frac{1}{2}} T_3   + x^{\prime}_{0} \otimes T_2  + x^{\prime}_{1} \otimes {( 1 + \mu^2 )}^{\frac{1}{2}} T_4 , $

$ \phi ( x^{\prime}_{1} )  = x^{\prime}_{-1} \otimes - \mu S^{*}_{4} + x^{\prime}_{0} \otimes {( 1 + \mu^2 )}^{- \frac{1}{2}} S^{*}_{2}  + x^{\prime}_{1} \otimes S^{*}_{3}.  $ 

\elmma

{\it Proof :} As $  x_{-1},~ x_0,~ x_{1}  $ have same norm, it follows that $ x^{\prime}_{i} = K x_{i},~  $ where $K = {\left\| x_i \right\|}^{- 1}, ~ i = \{ - 1, 0, 1 \}.$
 
Now, using ( \ref{sphere_A,B_S2_mu_t} ) and ( \ref{alphaB} ), we have 
 \bean \lefteqn{ \phi ( x^{\prime}_{- 1} )}\\
 &=& \frac{K t {( 1 + \mu^2 )}^{\frac{1}{2}}}{\mu} \phi ( B )\\
 &=& \frac{K t {( 1 + \mu^2 )}^{\frac{1}{2}}}{\mu} [ 1 \otimes S_1 + \frac{1 - t^{ - 1} x_0}{1 + \mu^2} \otimes S_2 + \frac{\mu x_{- 1}}{t {(1 + \mu^2)}^{\frac{1}{2}}} \otimes S_3 + \frac{\mu ( - \mu^{ - 1} x_1 )}{t {(1 + \mu^2)}^{\frac{1}{2}} } \otimes S_4 \\
 &=& K x_{ - 1} \otimes S_3 + K x_0 \otimes - \frac{S_2}{\mu {(1 + \mu^2)}^{\frac{1}{2}}} + K x_1 \otimes - \frac{S_4}{\mu}\eean 
 ( by Lemma \ref{h_invariance} ) 
 $$ = x^{\prime}_{-1} \otimes S_3 + x^{\prime}_{0} \otimes - {\mu}^{-1} {( 1 + {\mu}^{2} )}^{- \frac{1}{2}} S_2 + x^{\prime}_{1} \otimes - {\mu}^{-1} S_4 .$$ 
 By similar calculations, we get the second and the third equations. \qed
 
\vspace{4mm}

  Hence, $ \phi $ keeps the span of the orthonormal set $ \{ x^{\prime}_{-1},~ x^{\prime}_{0},~ x^{\prime}_{1} \} $ invariant. Moreover, $ \phi $ is kept invariant by the Haar state  $ h $  of $ SU_{\mu}( 2 ) .$ Therefore, we have a unitary representation of the CQG $ \widetilde{\clq} $ on span $ \{ x^{\prime}_{-1},~ x^{\prime}_{0},~ x^{\prime}_{1} \}.$

Using $ T_4 = T^{*}_3 $ from Lemma \ref{spherehomomorphism A* = A}, the unitary matrix, say $ Z $ corresponding to $ \phi $ and the ordered basis $ \{ x^{\prime}_{-1},~ x^{\prime}_{0},~ x^{\prime}_{1} \} $ is given by : 

\vspace{4mm} 


\be \label{sphere_QISO_matrix_Z} Z = \left ( \begin {array} {cccc}
  S_{3}  & - \mu \sqrt{1 + \mu^2} T_3  & - \mu S^{*}_4 \\ \frac{ - S_{2} }{ \mu  \sqrt{1 + {\mu}^2} }  &  T_2   &   \frac{ S^{*}_{2} }{ \sqrt{1 + {\mu}^2} } \\ - {\mu}^{- 1} S_{4}   & \sqrt{1 + \mu^2} T^{*}_{3}   &  S^{*}_{3} \end {array} \right ). \ee

\vspace{4mm}
  
  Recall that ( cf \cite{vandaelenotes} ), the antipode $\kappa$ on the matrix elements of a finite-dimensional   unitary representation $U^\alpha \equiv ( u_{pq}^\alpha)$ is given by $\kappa (u_{pq}^\alpha ) =( u_{qp}^\alpha )^* .$ Hence, the antipode is given by :
  
 $$ \kappa ( T_2 ) = T_2,~ \kappa ( T_3 ) = \frac{S^{*}_{2}}{{\mu}^{2}( 1 + \mu^2 )},~ \kappa ( S_2 ) = \mu^2 ( 1 + \mu^2 ) T^{*}_{3}, $$
 $$~ \kappa ( S_3 ) = S^{*}_{3},~ \kappa ( S_4 ) = {\mu}^2 S_4,~ \kappa ( T^{*}_{3} ) = \frac{S_2}{1 + \mu^2}, $$
 $$ ~ \kappa ( S^{*}_2 ) = ( 1 + \mu^2 ) T_3,~ \kappa ( S^{*}_{3} ) = S_3,~ \kappa ( S^{*}_{4} ) = {\mu}^{-2} S^{*}_{4}. $$

Now we derive some more equations by applying $ \kappa $ on the equations obtained by homomorphism condition.

\blmma

\label{k on B*B = A - A2}

 \bean  \mu^4 {( 1 + \mu^2 )}^2 T^{*}_3 T_3 + c \mu^2 ( 1 + \mu^2 )^{2} S^{*}_3 S_3 + c \mu^2 ( 1 + \mu^2 )^{2} S_4 S^{*}_4 \eean
  \be \label{k on B*B = A - A2 1} = \mu^2 ( 1 - T_2 ) ( \mu^2 + T_2 ) - c S_2 S^{*}_{2} - c S^{*}_2 S_2  + c \mu^2 ( 1 + \mu^2 )^{2}.1, \ee
 \bean  -2 \mu^4 {( 1 + \mu^2 )}^3 T^{*}_3 T_3 + \mu^2 {( 1 + \mu^2 )}^2 S^*_3 S_3 + \mu^4 {( 1 + \mu^2 )}^2 S_4 S^*_4 \eean
  \be \label{k on B*B = A - A2 2} = \mu^2 ( 1 + \mu^2 ) T_2 ( \mu^2 + 2T_2 - 1 ) - \mu^2 S_2 S^{*}_2 - S^{*}_2 S_2, \ee
 \bean  \mu^4 {( 1 + \mu^2 )}^4 T^{*}_3 T_3 - \mu^2 {( 1 + \mu^2 )}^2 S^*_3 S_3 - \mu^6 {( 1 + \mu^2 )}^2 S_4 S^*_4 \eean
 \be \label{k on B*B = A - A2 3}  = - \mu^2 {( 1 + \mu^2 )}^2 T^2_2 + \mu^4 S_2 S^*_2 + S^*_2 S_2, \ee

\be \label{k on B*B = A - A2 5} \mu^2 {( 1 + \mu^2 )}^2 S_4 T_3 + \mu^2 {( 1 + \mu^2 )}^2 T^*_3 S_3 = - S_2 ( \mu^2 + T_2 ) + ( 1 - T_2 )S_2, \ee
\be \label{k on B*B = A - A2 7} S_4 S_3 = - \frac{- S^2_2}{\mu^2 {( 1 + \mu^2 )}^2}. \ee

\elmma

{\it Proof :} The relations follow by applying $ \kappa $ on ( \ref{spherehomomorphism B*B = A - A2 1} ), ( \ref{spherehomomorphism B*B = A - A2 2} ), ( \ref{spherehomomorphism B*B = A - A2 3}), ( \ref{spherehomomorphism B*B = A - A2 5} ) and ( \ref{spherehomomorphism B*B = A - A2 7} ) respectively. \qed
\blmma

\label{k on BA = mu2 AB}

\be \label{k on BA = mu2 AB 1} - \mu^2 ( 1 - T_2 ) T^{*}_3 + c S_2 S^{*}_3 + c S^{*}_2 S_4 = - \mu^4 T^{*}_{3}( 1 - T_2 ) + c \mu^2 S_4 S^{*}_2 + c \mu^2 S^{*}_{3} S_2, \ee
\be \label{k on BA = mu2 AB B2} S_3 S_2 = \mu^2 S_2 S_3, \ee
\be \label{k on BA = mu2 AB B*2} S_2 S_4 = \mu^2 S_4 S_2, \ee
\be \label{k on BA = mu2 AB B} - S^*_2 T^*_3 + ( 1 - T_2 )S^*_3 = - \mu^2 T^*_3 S^*_2 + \mu^2 S^*_3 ( 1 - T_2 ), \ee
\be \label{k on BA = mu2 AB B*} - S_2 T^*_3 + ( 1 - T_2 ) S_4 = \mu^2 S_4 ( 1 - T_2 ) - \mu^2 T^*_3 S_2, \ee
  \be \label{k on BA = mu2 AB AB} T_3 S_2 + \mu^2 S_3 T_2 = \mu^2 ( T_2 S_3 + \mu^2 S_2 T_3 ). \ee

\elmma

{\it Proof :} The relations follow by applying $ \kappa $ on (  \ref{spherehomomorphism BA = mu2 AB 1} ), ( \ref{spherehomomorphism BA = mu2 AB B2} ), ( \ref{spherehomomorphism BA = mu2 AB B*2} ), ( \ref{spherehomomorphism BA = mu2 AB B} ), ( \ref{spherehomomorphism BA = mu2 AB B*} ) and ( \ref{spherehomomorphism BA = mu2 AB AB}  ) respectively. \qed


\blmma

\label{k on BB* = mu2 A - mu4 A2}

 \bean  \mu^4 {( 1 + \mu^2 )}^2 T_3 T^{*}_3 + c \mu^2 ( 1 + \mu^2 )^{2} S_3 S^{*}_{3} + c \mu^2 ( 1 + \mu^2 )^{2} S^{*}_{4} S_4 \eean
 \be \label{k on BB* = mu2 A - mu4 A2 1}  = \mu^4 ( 1 - T_2 ) ( 1 + \mu^2 T_2 ) + c S_2 S^{*}_2 + c S^{*}_{2} S_2 + c \mu^2 ( 1 + \mu^2 )^{2}.1, \ee
\be \label{k on BB* = mu2 A - mu4 A2 B2} S_3 S_4 = - \frac{{\mu}^2 S^2_2}{{( 1 + \mu^2 )}^2}, \ee
\be \label{k on BB* = mu2 A - mu4 A2 B} - \mu^2 {( 1 + \mu^2 )}^2 S^*_4 T^*_3 - \mu^2 {( 1 + \mu^2 )}^2 T_3 S^*_3 = \mu^2 ( 1 + \mu^2 ) S^*_2 - \mu^4 S^*_2 ( 1 - T_2 ) - \mu^4 ( 1 - T_2 )S^*_2, \ee
\be \label{k on BB* = mu2 A - mu4 A2 AB} {( 1 + \mu^2 )}^2 S^*_4 T^*_3 + \mu^2 {( 1 + \mu^2 )}^2 T_3 S^*_3 = - \mu^2 S^*_2 T_2 - \mu^4 T_2 S^*_2. \ee

\elmma

{\it Proof :} The relations follow by applying $ \kappa $ on  (\ref{spherehomomorphism BB* = mu2 A - mu4 A2 1}), (\ref{spherehomomorphism BB* = mu2 A - mu4 A2 B2}), ( \ref{spherehomomorphism BB* = mu2 A - mu4 A2 B} ) and (\ref{spherehomomorphism BB* = mu2 A - mu4 A2 AB}) respectively. \qed

\brmrk

\label{sphere_remark_k}

It follows from ( \ref{k on B*B = A - A2 7} ) and ( \ref{k on BB* = mu2 A - mu4 A2 B2} ) that $ \mu^4 S_4 S_3 = S_3 S_4 .$

\ermrk
























\subsection{Identification of $ SO_{\mu} ( 3 ) $ as the quantum isometry group}

Motivated by  (  \ref{sphere_matrix_of_SO_mu_3}  )  and  (  \ref{sphere_QISO_matrix_Z}  )  we are led to state and prove the following  statement : 

The map $ SO_{\mu}( 3 ) \rightarrow \clq $ sending $ M, L, G, N, C $ to $ - ( 1 + \mu^2 )^{-1} S_2,~ S_3,~ - \mu^{-1} S_4,~ ( 1 +\mu^2 )^{-1} ( 1 - T_2 ),~ \mu T_3 $ respectively is a $ \ast $ homomorphism ( See Proposition \ref{sphere_identification} ). 

To prove this, it is enough to show that all the relations of $ SO_{\mu} ( 3 ) $( as in ( \ref{sphere_SO_mu_3_description} ) ) when translated to relations of $ QISO^{+}_{R} ( \clo(S^{2}_{\mu,c}), ~ \clh, ~ D ) $ via the above map are satisfied. Hence, we list the relations one by one. 

$ L^* L = ( I - N  ) ( I - \mu^{- 2} N ) $ gives

$ S^*_3 S_3 = ( 1 - \frac{1 - T_2}{1 + \mu^2}  ) ( 1 - {\mu}^{- 2} \frac{1 - T_2}{1 + \mu^2} ) = \frac{\mu^2 + T_2}{1 + \mu^2} ( \frac{\mu^2 ( 1 + \mu^2 ) - ( 1 - T_2 )}{\mu^2 ( 1 + \mu^2 )}  ), $

which implies $$ \mu^2 {( 1 + \mu^2 )}^{2} S^*_3 S_3 = ( \mu^2 + T_2 ) ( \mu^2 ( 1 + \mu^2 ) - ( 1 - T_2 ) ). $$
$ L L^* = ( 1 - \mu^2 N ) ( 1 - \mu^4 N  ) $ gives 
$S_3 S^*_3 = ( 1 - \mu^2 ( \frac{1 - T_2}{1 + \mu^2}  ) ) ( 1 - \mu^4 ( \frac{1 - T_2}{1 + \mu^2} )  ) $ 
implying $$ {( 1 + \mu^2 )}^{2} S_3 S^*_3 = ( 1 + \mu^2 T_2  ) ( ( 1 + \mu^2 ) - \mu^4 ( 1 - T_2  )  )  .$$
$ G^* G = G G^* = N^2 $ gives 
$ - \frac{S^*_4}{\mu} ( - \frac{S_4}{\mu}   ) = ( - \frac{S_4}{\mu}  ) ( - \frac{S^*_4}{\mu}  ) = \frac{{( 1 - T_2 )}^{2}}{{( 1 + \mu^2 )}^{2}} $
implying $$ S^*_4 S_4 = S_4 S^*_4 = \frac{\mu^2 {( 1 - T_2 )}^{2}}{{( 1 + \mu^2 )}^{2}} .$$
$ M^* M = N - N^2  $ gives 
$$ \frac{S^*_2 S_2}{{( 1 + \mu^2 )}^{2}} = \frac{1 - T_2}{1 + \mu^2} - \frac{{( 1 - T_2  )}^{2}}{{( 1 + \mu^2 )}^{2}}, $$
which means $ \frac{S^*_2 S_2}{{( 1 + \mu^2 )}^{2}} = \frac{( 1 + \mu^2 ) ( 1 - T_2 )  - {( 1 - T_2 )}^{2}}{{( 1 + \mu^2 )}^{2}} $ implying
$$ S^*_2 S_2 = ( 1 - T_2 ) ( \mu^2 + T_2 )  .$$
$ M M^* = \mu^2 N - \mu^4 N^2  $ gives 
$ \frac{S_2 S^*_2}{{( 1 + \mu^2 )}^{2}} = \mu^2 ( \frac{1 - T_2}{1 + \mu^2} ) - \mu^4 {( \frac{1 - T_2 }{1 + \mu^2}  )}^{2}  ,$ which implies $$ S_2 S^*_2 = \mu^2 ( 1 - T_2 ) ( 1 + \mu^2 T_2 ).$$
$ C^* C = N - N^2  $ gives 
$ \mu T^*_3 \mu T_3 = \frac{1 - T_2}{1 + \mu^2} - {( \frac{1 - T_2}{1 + \mu^2} )}^{2}, $
which implies $$( 1 - T_2 ) ( \mu^2 + T_2 ) = \mu^2 {( 1 + \mu^2 )}^{2} T^*_3 T_3 . $$
$ C C^* = \mu^2 N - \mu^4 N^2  $ gives 
$ ( \mu T_3 )  {(  \mu T_3 )}^* = \mu^2 ( \frac{1 - T_2}{1 + \mu^2} ) - \mu^4 {( \frac{1 - T_2}{1 + \mu^2} )}^2,  $
which implies $ \mu^2 T_3 T^*_3 = \frac{\mu^2 ( 1 + \mu^2 ) ( 1 - T_2 ) - \mu^4 {( 1 - T_2 )}^{2}}{{( 1 + \mu^2  )}^{2}}  .$ 
Therefore, $ ( 1 + \mu^2 ) ( 1 - T_2 ) - \mu^2 {( 1 - T_2 )}^{2} = {( 1 + \mu^2 )}^{2} T_3 T^*_3 .$ Thus,
$$ ( 1 - T_2 ) ( 1 + \mu^2 T_2 ) = {( 1 + \mu^2 )}^{2} T_3 T^*_3 .$$
$ L N = \mu^4 N L  $ gives 
$ S_3 ( \frac{1 - T_2}{1 + \mu^2} ) = \mu^4 ( \frac{1 - T_2}{1 + \mu^2}  ) S_3, $
implying $$ S_3 ( 1 - T_2 ) = \mu^4 ( 1 - T_2 ) S_3 .$$
$ G N = N G  $ gives 
$ - \frac{S_4}{\mu} ( \frac{1 - T_2}{1 + \mu^2}  ) = ( \frac{1 - T_2}{1 + \mu^2}  ) ( - \frac{S_4}{\mu} ),  $ which implies $ S_4 ( 1 - T_2 ) = ( 1 - T_2 ) S_4 .$ Thus, $$ S_4 T_2 = T_2 S_4 .$$
$ M N = \mu^2 N M  $ gives 
$ ( - \frac{S_2}{1 + \mu^2} ) ( \frac{1 - T_2}{1 + \mu^2} ) = \mu^2 ( \frac{1 - T_2}{1 + \mu^2} ) ( - \frac{S_2}{1 + \mu^2} ) ,$ implying $$ S_2 ( 1 - T_2 ) = \mu^2 ( 1 - T_2 ) S_2 .$$
$ C N = \mu^2 N C  $ gives
$ \mu T_3 ( \frac{1 - T_2}{1 + \mu^2}  ) = \mu^2 ( \frac{1 - T_2}{1 + \mu^2}  ) \mu T_3 , $ which implies $$ T_3 ( 1 - T_2 ) = \mu^2 ( 1 - T_2 ) T_3 .$$
$ L G = \mu^4 G L  $ gives
$ S_3 ( - \mu^{ -1} S_4  ) = \mu^4 ( - \mu^{- 1} S_4  ) S_3 ,$ that is, $$ S_3 S_4 = \mu^4 S_4 S_3 .$$
$ L M = \mu^2 M L  $ gives 
$ S_3 ( - \frac{S_2}{1 + \mu^2}  ) = \mu^2 ( - \frac{S_2}{1 + \mu^2}  ) S_3 ,$ that is, $$ S_3 S_2 = \mu^2 S_2 S_3. $$
$ M G = \mu^2 G M  $ gives 
$ ( - \frac{S_2}{1 + \mu^2}  ) ( - \frac{S_4}{\mu}  ) = \mu^2 ( - \frac{S_4}{\mu} ) ( - \frac{S_2}{1 + \mu^2}  ), $
that is, $$ S_2 S_4 = \mu^2 S_4 S_2 .$$ 
$ C M = M C  $ gives 
$ ( \mu T_3 ) (  - \frac{S_2}{1 + \mu^2} ) = (  - \frac{S_2}{1 + \mu^2} ) ( \mu T_3 )  ,$ that is, $$ T_3 S_2 = S_2 T_3. $$
$ L G^* = \mu^4 G^* L  $ gives
$ S_3 ( - \frac{S^{*}_{4}}{\mu}  ) = \mu^4 ( - \frac{S^{*}_{4}}{\mu} ) S_3, $ that is, $$ S_3 S^*_4 = \mu^4 S^*_4 S_3. $$
$ M^2 = \mu^{- 1} L G $ gives
$ {( - \frac{S_2}{1 + \mu^2} )}^{2} = \mu^{- 1} S_3 ( - \frac{S_4}{\mu}  ) ,$ that is, $$ S_3 S_4 = - \frac{\mu^2}{{(  1 + \mu^2 )}^{2} } S^{2}_2. $$
$ M^* L = \mu^{- 1} ( I - N ) C  $ gives
$ - \frac{S^*_2}{1 + \mu^2} S_3 = \mu^{- 1} ( 1 - \frac{1 - T_2}{1 + \mu^2} ) \mu T_3 ,$ that is, $$ - S^*_2 S_3 = ( \mu^2 + T_2 ) T_3 .$$
$ N^* = N  $ gives 
$ \frac{{( 1 - T_2 )}^*}{1 + \mu^2} = \frac{1 - T_2}{1 + \mu^2} ,$
that is, $$ T^*_2 = T_2 .$$

Thus, we are led to prove the following lemmas, in some of which we will need $ \mu^{2} \neq 1.$ The case $ \mu = 1 $ will be dealt separately.

\blmma

\label{spherecheckinghomomorphism zero}  

$ S^{*}_{2} S_2 = ( 1 - T_2 ) ( \mu^2 + T_2 ). $

\elmma

{\it Proof :} Subtracting the equation obtained by multiplying $ c ( 1 + \mu^2 ) $ with ( \ref{spherehomomorphism B*B = A - A2 2} ) from ( \ref{spherehomomorphism B*B = A - A2 1} ), we have
 \bean  ( 1 + 2c ( 1 + \mu^2 ) ) S^{*}_{2} S_2 + c {( 1 + \mu^2 )}^{2} ( 1 - \mu^2 ) S^{*}_{4} S_4  \eean 
\be \label{spherecheckinghomomorphism zero_1} = ( 1 - T_2 ) ( \mu^2 + T_2 ) - c ( 1 + \mu^2 ) ( \mu^2 + 2 T_2 - 1 ) T_2 + c {( 1 + \mu^2 )}^{2} ( \mu^2 - 1 ) T_3 T^{*}_3 + c {( 1 + \mu^2 )}^{2}.1. \ee
Again, by adding ( \ref{spherehomomorphism B*B = A - A2 1} ) with $ c {( 1 + \mu^2 )}^{2} $ times ( \ref{spherehomomorphism B*B = A - A2 3} ) gives
\bean  ( 1 + c {( 1 + \mu^2 )}^{2} ) S^{*}_{2} S_2 +  c ( 1 - \mu^4 ) {( 1 + \mu^2 )}^{2} S^{*}_{4} S_4.  \eean
\be \label{spherecheckinghomomorphism zero_2} =  ( 1 - T_2 ) ( \mu^2 + T_2 ) - c {( 1 + \mu^2 )}^{2} T^{2}_{2} + c {( 1 + \mu^2 )}^{2} ( {\mu}^4 - 1 ) T_3 T^{*}_3 + c {( 1 + \mu^2 )}^{2}.1. \ee
Subtracting the equation obtained by multiplying $ ( \mu^2 + 1 ) $ with ( \ref{spherecheckinghomomorphism zero_1} ) from ( \ref{spherecheckinghomomorphism zero_2} ) we obtain
\bean - ( \mu^2 + c {( 1 + \mu^2 )}^{2} ) S^{*}_{2} S_2   
 = ( 1 - T_2 ) ( \mu^2 + T_2 ) - c {( 1 + \mu^2 )}^{2} T^{2}_{2} \eean 
 \bean  - ( 1 + \mu^2 ) ( 1 - T_2 ) ( \mu^2 + T_2 ) - c \mu^2 {( 1 + \mu^2 )}^{2}.1 + c {( 1 + \mu^2 )}^{2} ( \mu^2 + 2 T_2 - 1 ) T_2. \eean
The right hand side can be seen to equal $ - ( \mu^2 + c {( 1 + \mu^2 )}^{2} ) ( 1 - T_2 ) ( \mu^2 + T_2 ) .$
Thus, $ S^{*}_{2} S_2 = ( 1 - T_2 ) ( \mu^2 + T_2 ) .$ \qed

\blmma

\label{spherecheckinghomomorphism one}

\be \label{spherecheckinghomomorphism one_1} \mu^2 {( 1 + \mu^2 )}^{2} T^{*}_{3} T_3 = ( 1 - T_2 ) ( \mu^2 + T_2 ), \ee
\be \label{spherecheckinghomomorphism one_2} {( 1 + \mu^2 )}^{2} T_3 T^{*}_{3} = ( 1 - T_2 ) ( 1 + \mu^2 T_2 ), \ee
\be \label{spherecheckinghomomorphism one_3} S_2 S^{*}_{2} = \mu^2 ( 1 - T_2 ) ( 1 + \mu^2 T_2 ). \ee

\elmma

{\it Proof :} Applying $ \kappa $ on  Lemma \ref{spherecheckinghomomorphism zero}, we obtain ( \ref{spherecheckinghomomorphism one_1}  ).

Unitarity of the matrix $ Z $ ( ( 2, 2 ) position of the matrix $ Z^* Z $ ) gives $ \mu^2 ( 1 + \mu^2 ) T^{*}_{3} T_{3} + T^{2}_{2} + ( 1 + \mu^2 ) T_3 T^{*}_{3} = 1 .$

Using ( \ref{spherecheckinghomomorphism one_1} ) we deduce $ - ( 1 + \mu^2 )^{2} T_3 T^{*}_{3} = ( T_2 - 1 ) ( 1 + \mu^2 T_2 ).$ Thus we obtain ( \ref{spherecheckinghomomorphism one_2} ).

Applying $ \kappa $ on ( \ref{spherecheckinghomomorphism one_2} ), we deduce ( \ref{spherecheckinghomomorphism one_3}  ).  \qed

\blmma

\label{spherecheckinghomomorphism two}

$ S^{*}_{4} S_4 = S_4 S^{*}_{4} = {( 1 + \mu^2 )}^{-2} \mu^{2} {( 1 - T_2 )}^{2}. $

\elmma

{\it Proof :} Adding ( \ref{k on B*B = A - A2 2} ) and ( \ref{k on B*B = A - A2 3} ), we have :
$ - \mu^4 {( 1 + \mu^2 )}^3 ( 1 - \mu^2 ) T^*_3 T_3 + \mu^4 {( 1 + \mu^2 )}^2 ( 1 - \mu^2 ) S_4 S^*_4 = - \mu^2 ( 1 + \mu^2 ) ( 1 - \mu^2 ) T_2 ( 1 - T_2 ) - \mu^2 ( 1 - \mu^2 ) S_2 S^*_2 .$

Using $ \mu^2 \neq 1  ,$ we obtain,
$$ - \mu^4 {( 1 + \mu^2 )}^3 T^*_3 T_3 + \mu^4 {( 1 + \mu^2 )}^2 S_4 S^*_4 = - \mu^2 ( 1 + \mu^2 ) T_2 ( 1 - T_2 ) - \mu^2 S_2 S^*_2 .  $$ 
Now using ( \ref{spherecheckinghomomorphism one_1} ) and ( \ref{spherecheckinghomomorphism one_3} ), we reduce the above equation to 
 \bean \lefteqn{ \mu^4 {( 1 + \mu^2 )}^2 S_4 S^*_4 }\\
 & = & - \mu^2 ( 1 - T_2 ) ( T_2 + \mu^2 T_2 + \mu^2 + \mu^4 T_2 ) + \mu^2 ( 1 + \mu^2 ) ( 1 - T_2 ) ( \mu^2 + T_2 ) \\ & = & \mu^6 {( 1 - T_2 )}^2. \eean 
Thus,\bean  S_4 S^*_4 
& = & \frac{\mu^6}{\mu^4 {( 1 + \mu^2 )}^2} {( 1 - T_2 )}^2 \\
& = & \frac{\mu^2}{{( 1 + \mu^2 )}^2} ( 1 - T_2 )^2. \eean
Applying $ \kappa ,$ we have $ S^*_4 S_4 = \frac{\mu^2}{{( 1 + \mu^2 )}^2} {( 1 - T_2 )}^2 .$

Thus, $ S^*_4 S_4 = S_4 S^*_4 = \frac{\mu^2}{{( 1 + \mu^2 )}^2} ( 1 - T_2 )^2 .$ \qed

 




\blmma

\label{spherecheckinghomomorphism three}

$ \mu^{2} {( 1 + \mu^2 )}^{2} S^{*}_{3} S_3 = ( \mu^2 + T_2 ) [ \mu^2 ( 1 + \mu^2 ) - ( 1 - T_2 ) ]. $ 

\elmma

{\it Proof :} Using Lemma \ref{spherecheckinghomomorphism zero} in ( \ref{spherehomomorphism B*B = A - A2 1}  ), we have 
 \be \label{spherecheckinghomomorphism_lemma_in_expression} S^*_3 S_3 + T^*_3 T_3 + T_3 T^*_3 + S^*_4 S_4 = 1.  \ee
 The lemma is derived by substituting the expressions of $ T^{*}_{3} T_3, ~ T_3 T^{*}_{3} $ and $ S^{*}_{4} S_4 $ from   (\ref{spherecheckinghomomorphism one_1}), ( \ref{spherecheckinghomomorphism one_2} ) and Lemma \ref{spherecheckinghomomorphism two} in the equation ( \ref{spherecheckinghomomorphism_lemma_in_expression} ). \qed

\blmma

\label{spherecheckinghomomorphism four}

$ {( 1 + \mu^2 )}^{2} S_3 S^{*}_{3} = ( 1 + \mu^2 T_2 ) ( 1 + \mu^2 - \mu^4 ( 1 - T_2 ) ). $

\elmma

{\it Proof :} By unitarity of the matrix $ Z $, in particular equating the ( 1, 1 ) th  entry of $ Z Z^{*} $ to 1 we get $ S_3 S^{*}_{3} + \mu^2 ( 1 + \mu^2 ) T_3 T^{*}_{3} + \mu^2 S^{*}_{4} S_4 = 1 .$ Then the Lemma follows by using  ( \ref{spherecheckinghomomorphism one_2} ) and Lemma \ref{spherecheckinghomomorphism two}  in the above equation. \qed

\blmma

\label{spherecheckinghomomorphism five}

$ - S^{*}_{2} S_3 = ( \mu^2 + T_2 ) T_3. $

\elmma

{\it Proof :}  By applying $ \ast $ and then multiplying by $ \mu^2 $ on ( \ref{spherehomomorphism B*B = A - A2 5} ) we have $ \mu^2 S^*_2 S_3 + \mu^2 S^*_4 S_2 = - \mu^2 T_3 ( \mu^2 + T_2 ) + \mu^2 ( 1 - T_2 ) T_3 .$ Subtracting this from ( \ref{spherehomomorphism B*B = A - A2 6} ) we have $ ( 1 - \mu^2 ) S^*_2 S_3 = - T_2 T_3 - \mu^2 T_3 T_2 + \mu^2 T_3 ( \mu^2 + T_2 ) - \mu^2 ( 1 - T_2 ) T_3 $ which implies $ - S^*_2 S_3 = ( \mu^2 + T_2 ) T_3 $ as $\mu^2 \neq 1.$  \qed

\blmma

\label{spherecheckinghomomorphism six}

$ S_2 ( 1 - T_2 ) = \mu^2 ( 1 - T_2 ) S_2. $ 

\elmma 

{\it Proof :} Applying $ \kappa $  to  Lemma \ref{spherecheckinghomomorphism five} and then taking adjoint, we have


\be \label{spherecheckinghomomorphism six_1}  \mu^2 {( 1 + \mu^2 )}^{2} T^{*}_3 S_3 = - ( \mu^2 + T_2 ) S_{2}. \ee
Adding ( \ref{k on BB* = mu2 A - mu4 A2 B} ) and ( \ref{k on BB* = mu2 A - mu4 A2 AB} ) and then taking adjoint, we get ( by using $ \mu^2 \neq 1$ ) 
\be \label{spherecheckinghomomorphism six_2}  \mu^2 {( 1 + \mu^2 )}^{2} T_3 S_4 = \mu^4 ( 1 - T_2 ) S_2. \ee
Moreover, ( \ref{k on B*B = A - A2 5} ) gives
$$ \mu^2 {( 1 + \mu^2 )}^{2} S_4 T_3 = - S_2 ( \mu^2 + T_2 ) + ( 1 - T_2 )S_2 - \mu^2 {( 1 + \mu^2 )}^{2} T^{*}_{3} S_3. $$
Using ( \ref{spherecheckinghomomorphism six_1} ), the right hand side of this equation turns out to be $ S_2 ( 1 - T_2 ).$

Thus,

\be \label{spherecheckinghomomorphism six_3}  {( 1 + \mu^2 )}^{2} S_4 T_3 = \mu^{-2} S_2 ( 1 - T_2 ). \ee
Again, application of adjoint to the equation ( \ref{k on BB* = mu2 A - mu4 A2 B} ) gives :
$$ \mu^2 {( 1 + \mu^2 )}^{2} S_3 T^{*}_3 = - \mu^2 {( 1 + \mu^2 )}^{2} T_3 S_4 - \mu^2 ( 1 + \mu^2 ) S_2 + \mu^4 ( 1 - T_2 )S_2 + \mu^4 S_2 ( 1 - T_2 ). $$
Using ( \ref{spherecheckinghomomorphism six_2} ), we get
\be \label{spherecheckinghomomorphism six_4} {( 1 + \mu^2 )}^{2} S_3 T^{*}_3 = - S_2 ( 1 + \mu^2 T_2 ). \ee
Using ( \ref{spherecheckinghomomorphism six_1} ) - ( \ref{spherecheckinghomomorphism six_4} ) to the equation ( \ref{spherehomomorphism BA = mu2 AB A} ), we obtain :

$ S_2 - 2 S_2 T_2 - ( 1 + \mu^2 )^{-1} \mu^2 S_2 ( 1 + \mu^2 T_2 ) + \mu^{-2} ( 1 + \mu^2 )^{-1} S_2 ( 1 - T_2 ) = \mu^2 S_2 - 2 \mu^2 T_2 S_2 + ( 1 + \mu^2 )^{-1} \mu^6 ( 1 - T_2 )S_2 - ( 1 + \mu^2 )^{-1} ( \mu^2 + T_2 ) S_2. $

This gives

$ \mu^2 ( 1 + \mu^2 ) [ ( S_2 - S_2 T_2 ) - ( \mu^2 S_2 - \mu^2 T_2 S_2 ) ] - \mu^2 ( 1 + \mu^2 ) ( S_2 T_2 - \mu^2 T_2 S_2 ) - \mu^4 S_2 - \mu^6 S_2 T_2 + S_2 ( 1 - T_2 ) - \mu^8 ( S_2 - T_2 S_2 ) + \mu^4 S_2 + \mu^2 T_2 S_2 = 0 .$

Thus,
$ \mu^2 ( 1 + \mu^2 )[ S_2 ( 1 - T_2 ) - \mu^2 ( 1 - T_2 ) S_2 ] + S_2 ( 1 - T_2 ) - \mu^2 ( S_2 - T_2 S_2 ) + \mu^6 [ S_2 ( 1 - T_2 ) - \mu^2 ( 1 - T_2 ) S_2 ] - \mu^6 ( 1 - T_2 )S_2 + \mu^4 S_2 ( 1 - T_2 ) + \mu^2 ( S_2 ( 1 - T_2 ) - \mu^2 ( 1 - T_2 ) S_2 ) = 0 .$

On simplifying,$ ( \mu^6 + 2 \mu^4 + 2 \mu^2 + 1 )( S_2 ( 1 - T_2 ) - \mu^2 ( 1 - T_2 ) S_2 ) = 0 ,$ which proves the lemma as $ 0 < \mu < 1 .$ \qed

\blmma


\label{spherecheckinghomomorphism seven}


\be \label{spherecheckinghomomorphism seven_1} T_3 ( 1 - T_2 ) = \mu^2 ( 1- T_2 )T_3,  \ee
\be \label{spherecheckinghomomorphism seven_2} S_3 S^*_4 = \mu^4 S^*_4 S_3.  \ee

\elmma

{\it Proof :} The equation ( \ref{spherecheckinghomomorphism seven_1} ) follows by applying $ \kappa $ on Lemma \ref{spherecheckinghomomorphism six} and then taking $\ast.$

We have $ S^*_4 S_3 = - T^2_3 $ from ( \ref{spherehomomorphism B*B = A - A2 7} ). On the other hand we have $ S_3 S^*_4 = - \mu^4 T^2_3 $ from ( \ref{spherehomomorphism BB* = mu2 A - mu4 A2 B2} ). Combining these two, we get ( \ref{spherecheckinghomomorphism seven_2} ). \qed

\blmma


\label{spherecheckinghomomorphism eight}

$ S_4 T_2 = T_2 S_4. $

\elmma

{\it Proof :}  
Subtracting ( \ref{k on BA = mu2 AB B*} ) from ( \ref{spherehomomorphism BA = mu2 AB B*} )  we get the required result. \qed

\blmma


\label{spherecheckinghomomorphism nine}

$ T_3 S_2 = S_2 T_3. $

\elmma

{\it Proof :}  By applying $ \ast $ on  ( \ref{k on BA = mu2 AB B} ) and then subtracting it from ( \ref{spherehomomorphism BA = mu2 AB B} ) we obtain $ S_2 T_3 - T_3 S_2  = 0 .$
 \qed






\blmma


\label{spherecheckinghomomorphism ten}

$ S_3 ( 1 - T_2 ) = \mu^4 ( 1 - T_2 ) S_3. $ 

\elmma

{\it Proof :}  By adding ( \ref{spherehomomorphism BA = mu2 AB B} ) with ( \ref{spherehomomorphism BA = mu2 AB AB} ) we obtain 
$$ S_3 ( 1 - T_2 ) + \mu^2 S_3 ( T_2 - 1 ) = \mu^2 ( \mu^2 - 1 ) T_3 S_2. $$ 
Thus, using $ \mu^2 \neq 1,$
\be \label{spherecheckinghomomorphism ten_1}  S_3 ( 1 - T_2 ) = - \mu^2 T_3 S_2.   \ee
 Moreover, by applying $ \ast $ on ( \ref{k on BA = mu2 AB B} ), we obtain $$ \mu^2 ( 1 - T_2 ) S_3 = \mu^2 S_2 T_3 - T_3 S_2 + S_3 ( 1 - T_2 ) .$$
 
 Thus, $$ \mu^4 ( 1 - T_2 ) S_3 = \mu^4 S_2 T_3 - \mu^2 T_3 S_2 + \mu^2 S_3 ( 1 - T_2 ) .$$ 
 Hence, to prove the Lemma it suffices to prove: 
$$ S_3 ( 1 - T_2 ) = \mu^4 S_2 T_3 - \mu^2 T_3 S_2 + \mu^2 S_3 ( 1 - T_2 ). $$
After using $ T_3 S_2 = S_2 T_3 $ obtained from Lemma \ref{spherecheckinghomomorphism nine} we get this to be  the same as $ ( 1 - \mu^2 ) S_3 ( 1 - T_2 ) = \mu^2 ( \mu^2 - 1 ) T_3 S_2  .$ This is equivalent to $ S_3 ( 1 - T_2 ) = - \mu^2 T_3 S_2 $ ( as $ \mu^2 \neq 1 $ ) which follows from ( \ref{spherecheckinghomomorphism ten_1} ). \qed

\bppsn

\label{sphere_identification}

Assume $ \mu \neq 1 .$
The map $ SO_{\mu}( 3 ) \rightarrow \clq $ sending $ M, L, G, N, C $ to $ - ( 1 + \mu^2 )^{-1} S_2,~ S_3,~ - \mu^{-1} S_4,~ ( 1 +\mu^2 )^{-1} ( 1 - T_2 ),~ \mu T_3 $ respectively is a $ \ast $ homomorphism.

\eppsn

{\it Proof :} Now, we note that the proof of this Lemma reduces to verification of the  relations on $ \clq $  as derived in Lemmas \ref{spherecheckinghomomorphism zero} - \ref{spherecheckinghomomorphism ten} along with the following equations :


\be \label{spherecheckinghomomorphism eleven} S_3 S_4 = \mu^4 S_4 S_3,   \ee
\be \label{spherecheckinghomomorphism twelve} S_3 S_2 = \mu^2 S_2 S_3, \ee
\be \label{spherecheckinghomomorphism thirteen} S_2 S_4 = \mu^2 S_4 S_2,  \ee
\be \label{spherecheckinghomomorphism fourteen}  S_3 S_4 = - \frac{\mu^2}{{( 1 + \mu^2 )}^2} S^2_2,  \ee
which follow from   Remark \ref{sphere_remark_k}, ( \ref{k on BA = mu2 AB B2} ), ( \ref{k on BA = mu2 AB B*2} ), ( \ref{k on BB* = mu2 A - mu4 A2 B2} ) respectively.   \qed

\bthm

\label{sphere_final_theorem}

For $ \mu \neq 1, ~ QISO^{+}_{R} ( \clo(S^{2}_{\mu,c}), ~ \clh, ~ D ) \cong SO_{\mu}( 3 ). $

\ethm

{\it Proof :} We have seen in Theorem \ref{sphere_SU_mu_2_object} that $ SU_{\mu} ( 2 ) $ is an object in $ {\bf Q^{\prime}_{R}} ( D ) $ and  $ SO_{\mu}( 3 ) $ is the 
 corresponding maximal Woronowicz subalgebra for which the action is faithful. Thus, $SO_\mu(3)$  is a quantum subgroup of $ QISO^{+}_{R}(D).$ 
Now, Proposition \ref{sphere_identification} implies that $QISO^+_R(D)$ is a quantum subgroup of $SO_\mu(3)$, thereby completing the proof. \qed 

\brmrk

\label{sphere_final_theorem_remark3}

We observe that in the proof of Theorem \ref{sphere_final_theorem}, the only place where the structure of $ D $ was used was in Proposition 
\ref{sphere_Dabrowski_linearity_2} and there we used the fact that the unitary commutes with $ \left| D \right|.$ Thus, if we replace this 
spectral triple by $ ( \clo(S^2_{\mu,c}), \clh, \left| D \right| ) ,$ everything remains same and we 
deduce that $$ QISO^{+}_{R}( \clo(S^2_{\mu,c}), \clh, \left| D \right| ) \cong  QISO^{+}_{R} ( \clo(S^2_{\mu,c}),  \clh, D ) \cong  SO_{\mu} ( 3 ).$$    

\ermrk

\subsection{The quantum isometry group in the case $ \mu = 1$}

As we had mentioned earlier, some of the Lemmas which were required for the proof of Theorem \ref{sphere_final_theorem} needed the condition $ \mu \neq 1.$ In this subsection, we will give the proof for $ \mu = 1 $ case.

To begin with, we prove some of the Lemmas in the case $ \mu = 1 $ which needed $ \mu \neq 1 $ previously.  We note that the proofs of Lemmas \ref{spherecheckinghomomorphism zero}, \ref{spherecheckinghomomorphism one}, \ref{spherecheckinghomomorphism six}, \ref{spherecheckinghomomorphism seven}, \ref{spherecheckinghomomorphism eight}, \ref{spherecheckinghomomorphism nine} go through even in the case $ \mu = 1. $ Therefore, we can use these Lemmas here.

 \blmma

 \label{sphere_mu=1_S3T2=T2S3}

$ S_3 ( 1 - T_2 ) = (1 - T_2) S_3. $

\elmma

 {\it Proof :}  From ( \ref{k on BA = mu2 AB AB} ) , we obtain $ T_3 S_2 + S_3 T_2 = T_2 S_3 + S_2 T_3. $

Using $T_3 S_2 = S_2 T_3,$ from Lemma \ref{spherecheckinghomomorphism nine}, we have $ S_3 T_2 = T_2 S_3 $ which proves the Lemma. 

\qed

\blmma

\label{sphere_mu=1_S2T2=T2S2}

$ S_2 (1 - T_2) = (1 - T_2) S_2. $

\elmma

{\it Proof :} From (  \ref{spherehomomorphism BA = mu2 AB B2}   ) and ( \ref{spherehomomorphism BA = mu2 AB B*2}  ), we have respectively $ S_3 T_3 = T_3 S_3 $ and $ S_4 T^*_3 = T^*_3 S_4 .$ From (  \ref{spherecheckinghomomorphism one_1}   ) and (   \ref{spherecheckinghomomorphism one_2}      ), we see that in the case $\mu = 1,$ $T_3$ is normal. Hence, $S_3 T^*_3 = T^*_3 S_3 $ and $S_4 T_3 = T_3 S_4.$ Using these in ( \ref{spherehomomorphism BA = mu2 AB 1}  ), we have $S_2 (1 - T_2) = (1 - T_2) S_2.$ \qed

\blmma

\label{sphere_mu=1_S3_normal}

$S_3 S^*_3 = S^*_3 S_3.$

\elmma

{\it Proof :} 
Multiplying by $ 4 $ the equation (  \ref{spherehomomorphism BB* = mu2 A - mu4 A2 A2}  ),  we have
$$ 4 S_2 S^*_2 - 4 S_3 S^*_3 - 4 S_4 S^*_4 = - 4 T^2_2 + 4 T_3 T^*_3 + 4 T^*_3 T_3.  $$
Again, from (  \ref{k on B*B = A - A2 2}  ), we have 
$$ 16 T^*_3 T_3 - 4 S^*_3 S_3 - 4 S_4 S^*_4 = - 4 T^2_2 + S_2 S^*_2 + S^*_2 S_2. $$
Subtracting this from the previous equation, we have
$$ 4 S_2 S^*_2 - 4 ( S_3 S^*_3 - S^*_3 S_3  ) - 16 T^*_3 T_3 = 4 T_3 T^*_3 + 4 T^*_3 T_3 - S_2 S^*_2 - S^*_2 S_2.  $$
Again using Lemma \ref{spherecheckinghomomorphism zero}, (  \ref{spherecheckinghomomorphism one_1}   ), (  \ref{spherecheckinghomomorphism one_2}     ) and (  \ref{spherecheckinghomomorphism one_3}   ) in this equation, we obtain 
 $$ - 4 ( S_3 S^*_3 - S^*_3 S_3  ) = 4 ( 1 - T^2_2 ) - 6 ( 1 - T^2_2 ) + 2 ( 1 - T^2_2 ) = 0, $$
which implies $ S^*_3 S_3 = S_3 S^*_3 .$ \qed

\blmma

\label{sphere_mu=1_S4_normal}

$ S^*_4 S_4 = S_4 S^*_4. $

\elmma

{\it Proof :} 
From ( \ref{spherehomomorphism B*B = A - A2 2} ), we have
$$ - 2 S^*_2 S_2 + 2 S^*_3 S_3 + 2 S^*_4 S_4 = 2 T^2_2 - 2 T_3 T^*_3 - 2 T^*_3 T_3.  $$
Again, from  (  \ref{spherehomomorphism BB* = mu2 A - mu4 A2 A}  ), we obtain
$$ - 2 S_2 S^*_2 + 2 S_3 S^*_3 + 2 S_4 S^*_4 = 2 T_2 - 2 ( 1 - T_2 ) T_2 - 2 T_3 T^*_3 - 2 T^*_3 T_3 . $$
Subtracting this from the previous equation and using Lemma \ref{spherecheckinghomomorphism zero}, ( \ref{spherecheckinghomomorphism one_3}   ) and Lemma \ref{sphere_mu=1_S3_normal} we get $ 2 ( S^*_4 S_4 - S_4 S^*_4 ) = 0 $ implying  $ S^*_4 S_4 = S_4 S^*_4.$ \qed

\vspace{6mm}

Now we prove that $ QISO^{+}_{R} ( \clo(S^{2}_{1,c}), ~ \clh, ~ D ) $ is commutative as a $ C^* $ algebra. As $ QISO^{+}_{R} ( \clo(S^{2}_{1,c}), ~ \clh, ~ D ) $ is generated by $ T_2,~ T_3,~ S_2,~ S_3,~ S_4 ,$ it is enough to show that these elements belong to the centre of $ QISO^{+}_{R} ( \clo(S^{2}_{\mu,c}), ~ \clh, ~ D ). $

\blmma

\label{sphere_mu=1_everything_in_centre}

$ T_2,~ T_3,~ S_2,~ S_3,~ S_4 ,$ belong to the centre of $ QISO^{+}_{R} ( \clo(S^{2}_{\mu,c}), ~ \clh, ~ D ). $

\elmma

{\it Proof :}  $ T_2 $ is self adjoint. From ( \ref{spherecheckinghomomorphism one_1} ) and ( \ref{spherecheckinghomomorphism one_2} ) we note that $ T_3 $ is normal in the case $ \mu = 1.$ From Lemma \ref{spherecheckinghomomorphism zero} and (  \ref{spherecheckinghomomorphism one_3}   ), we deduce that $ S_2 $ is normal in the case $ \mu = 1. $ Similarly, from Lemma \ref{sphere_mu=1_S3_normal} and Lemma \ref{sphere_mu=1_S4_normal} we obtain that $ S_3 $ and $ S_4 $ are normal. Hence, it is enough to show that the elements $ T_2, T_3, S_2, S_3, S_4  $ commute among themselves.

Now, from ( \ref{spherecheckinghomomorphism seven_1} ), Lemma \ref{sphere_mu=1_S2T2=T2S2},  Lemma \ref{sphere_mu=1_S3T2=T2S3}, Lemma \ref{spherecheckinghomomorphism eight} we get that $ T_2 $ commutes with $ T_3, S_2, S_3, S_4 $ respectively. 
From Lemma \ref{spherecheckinghomomorphism nine} and ( \ref{spherehomomorphism BA = mu2 AB B2} ) we have that $ T_3 $ commutes with $ S_2, S_3 $ respectively. Now taking adjoint on the equation ( \ref{spherehomomorphism BA = mu2 AB B*2} )  we obtain that $ T_3 $ commutes with $ S^*_4 $ implying that $T_3 $ commutes with $ S_4. $
  From ( \ref{k on BA = mu2 AB B2} ) and ( \ref{k on BA = mu2 AB B*2} ), we have that $ S_2 $ commutes with $ S_3 $ and $ S_4 $ respectively.
           Finally, $ S_3 S_4 = S_4 S_3 $ follows from Remark \ref{sphere_remark_k}.
\qed

\vspace{4mm}

From Lemma \ref{sphere_mu=1_everything_in_centre}, we deduce the following. 

\bthm

\label{sphere_mu=1_final_theorem}

$ QISO^{+}_{R} ( \clo(S^{2}_{1,c}),  \clh,  D ) $ is commutative as a $ C^* $ algebra and hence coincides with the classical compact quantum group of orientation preserving isometries of the sphere, that is, $ C ( SO ( 3 ) ).$

\ethm

\brmrk

Our characterization of $SO_\mu(3)$ as the quantum isometry group of a noncommutative Riemannian manifold generalizes the classical description of the group
 $SO(3)$ as the group of orientation preserving isometries of the usual Riemannian structure on the $2$-sphere. It may be mentioned here that in a very recent
 article (\cite{soltan_new}), P. M. Soltan has characterized $SO_\mu(3)$ as the universal compact quantum group acting on the finite dimensional $C^*$-algebra
 $M_2(\IC)$ such that the action preserves a functional $\omega_\mu$ defined in \cite{soltan_new}. 
 In the classical case, we have three equivalent descriptions of $SO(3)$: (a) as a quotient of $SU(2)$, (b) as the group of (orientation preserving) isometries 
 of $S^2$, and (c) as the automorphism group of $M_2$. In the quantum case the definition of $SO_\mu(3)$ is an analogue of (a), so the
 characterization of $SO_\mu(3)$  obtained in this paper  as the quantum isometry 
 group, together with Soltan's characterization, completes the generalization of all three classical descriptions of $SO(3)$.  

\ermrk

\subsection{ Existence of $ \widetilde{QISO^{+}}( D ) $ }

 For the above spectral triple we have been unable to settle the issue of the existence of $ \widetilde{QISO^{+}}( D ) $ which is 
the universal object ( if it exists ) in the category $ \bf{Q^{\prime}} ( D ) $ mentioned in subsection \ref{qorient_subsection_proof_of_main_theorem}. 
Nevertheless, we now show that if it exists, the Woronowicz subalgebra $ QISO^{+}( D ) $ must be $ SO_{\mu} ( 3 ).$ In particular, the universal object in the subcategory of CQG s acting 
by orientation preserving isometries and containing $ SO_{\mu} ( 3 ) $ as a quantum subgroup exists.

\blmma

\label{sphere_unresricted}

If $ \widetilde{QISO^{+}}( D ) $ exists, its induced action on $ S^{2}_{\mu,c} $, say $ \alpha_0 $, must preserve
 the state $h$ on the subspace  spanned by $ \{1,A,B,B^*,AB,AB^*,A^2,B^2,B^{*^2} \}.$

\elmma

{\it Proof :} Let $ \clw_0 =  \IC.1,$ $\clw_{\frac{1}{2}} = {\rm Span}  \{ 1, A, B, B^* \}$,

 $ \clw_{\frac{3}{2}} ={\rm Span}
 \{ 1, A, B, B^*, AB, AB^*, A^2, B^2, B^{*2}  \}$. 

We note that the proof of Proposition \ref{sphere_Dabrowski_linearity_2} and the Lemmas preceding it do not use the assumption that the action is volume 
 preserving, 
so the proof of Proposition \ref{sphere_Dabrowski_linearity_2} goes through verbatim implying that 
$ \alpha_0 $ keeps invariant the subspace spanned by $ \{ 1, A, B, B^* \} $ and hence 
it preserves  $ \clw_{\frac{3}{2}}$ as well.   Let $\clw_{\frac{3}{2}}=\clw_{\frac{1}{2}} \oplus \clw^\prime$ be the orthogonal decomposition w.r.t. the Haar 
 state (say $h_0$) of $QISO^+(D)$. Since
 $ SO_{\mu} ( 3 ) $ is a sub-object of $QISO^+(D)$,  
 there is a CQG morphism $ \pi $ from $ QISO^{+}( D ) $ onto $ SO_{\mu} ( 3 ) $ satisfying $ ( {\rm id} \otimes \pi ) \alpha_0 = \Delta$, where 
 $\Delta$ is the $SO_\mu(3)$ action on $S^2_{\mu,c}$. It follows from this that any $QISO^+(D)$-invariant subspace 
(in particular $\clw^\prime)$ is also $SO_\mu(3)$-invariant. On the other hand, it is easily seen that on $\clw_{\frac{3}{2}}$, the $SO_\mu(3)$-action decomposes 
  as $\clw_{\frac{1}{2}} \oplus \clw^{\prime \prime}$, (orthogonality w.r.t. $h$, the Haar state of $SO_\mu(3))$), where $\clw^{\prime \prime}$ is a five 
 dimensional irreducible subspace.

We claim that $\clw^\prime=\clw^{\prime \prime}$, which will prove that the $QISO^+(D)$-action $\alpha_0$ has the same $h$-orthogonal decomposition 
 as the $SO_\mu(3)$-action on $\clw_{\frac{3}{2}}$, so preserves $\IC.1$ and its $h$-orthogonal complements. This will prove that $\alpha_0$ preserves the Haar
 state $h$ on $\clw^{\frac{3}{2}}$. 

We now prove the claim. We observe that $\clv:=\clw^\prime \bigcap \clw^{\prime \prime}$ is invariant under the $SO_\mu(3)$-action but due to the irreducibility of $ \Delta $ on the vector space $ \clw^{\prime} $ or $ \clw^{\prime \prime} ,$ it has to be zero or  $ \clw^{\prime} = ~ \clw^{\prime \prime} .$ Now,  $ {\rm dim} ( \clv ) = 0 $ implies $ {\rm dim} ( \clw^{\prime} ) + {\rm dim} ( \clw^{\prime \prime} ) = 5 + 5 > 9 = {\rm dim} ( \clw_{\frac{3}{2}} ) $ which is a contradiction unless $ \clw^{\prime} = ~ \clw^{\prime \prime} .$


\qed 

 \bthm

If $ \widetilde{QISO^{+}}( D ) $ exists, then we have $ QISO^{+} ( D ) \cong SO_{\mu} ( 3 ). $ In particular, the universal object in the subcategory of $ \bf{Q^{\prime}} ( D ) $ with objects  containing $ SO_{\mu} ( 3 ) $ as a sub- object, exists.

\ethm

{\it Proof :} In Lemma \ref{sphere_unresricted}, it was noted that Proposition \ref{sphere_Dabrowski_linearity_2} follows as before.
 We observe that the other Lemmas used  to prove Theorem \ref{sphere_final_theorem} require the conclusion of Lemma \ref{sphere_unresricted} as the only extra ingredient. \qed

\section{The spectral triple by Chakraborty and Pal on $ S^{2}_{\mu,c}, c > 0 $}

Now, we shall consider another class of spectral triples on the Podles spheres and show that they give rise to completely different quantum groups of (orientation preserving) isometries. Indeed, for these spectral triples, we have been able to prove the existence of $\widetilde{{QISO}^+}$ and identify it with the CQG $ C^* ( \IZ_2 \ast \IZ^{\infty} )$ where $ \IZ^{\infty} $ denotes countably infinite copies of the group of integers.  
  
  In this section, we will work with $ c > 0.$
  
  \subsection{The spectral triple}
  
  Let us describe the spectral triple on $S^2_{\mu,c}$  introduced and studied  in \cite{chak_pal}.
  
Let $ \clh_{+} = \clh_{-} = l^2 ( \IN \bigcup\{0\} ),\clh = \clh_{+} \oplus \clh_{-} .$

Let $\{  e_n, n \geq 0 \} $ be the canonical  orthonormal basis of $ \clh_{+} = \clh_{-} $ and $ N $ be the operator defined on it by $ N ( e_n ) = n e_n .$

We recall the irreducible  representations $ \pi_{+} $ and $ \pi_{-}  : \clh_{\pm} \rightarrow \clh_{\pm} $   as in \cite{chak_pal}.
\be \label{sphere_chakraborty_pal_pi ( A )} \pi_{\pm} ( A ) e_n = \lambda_{\pm} {\mu}^{2n} e_n,  \ee
\be \label{sphere_chakraborty_pal_pi ( B )} \pi_{\pm} ( B ) e_n = {c_{\pm}( n )}^{\frac{1}{2}} e_{n - 1},  \ee
where \be \label{sphere_chak_pal_c_pm(n)} e_{- 1} = 0,~ \lambda_{\pm} =  \frac{1}{2} \pm {( c + \frac{1}{4} )}^{\frac{1}{2}},~ c_{\pm} ( n ) = \lambda_{\pm} \mu^{2n} - {( \lambda_{\pm} \mu^{2n} )}^{2} + c .\ee
Let $ \pi = \pi_{+} \oplus \pi_{-} $ and $ D = \left ( \begin {array} {cccc}
   0 & N  \\ N & 0 \end {array} \right ) .$
   
Then $ ( S^{2}_{\mu,c}, \pi, \clh, D ) $ is a spectral triple.

We note that the eigenvalues of $ D $ are $ \{ n : n \in \IZ \} $ and eigenspace is spanned by $ \left ( \begin {array} {cccc}
   e_n \\ e_n \end {array} \right ) $ corresponding to the positive eigenvalue $ n $ and $ \left ( \begin {array} {cccc}
   e_n \\ - e_n \end {array} \right ) $ for the negative eigenvalue $ - n .$
   
 \blmma
 
 \label{sphere_chak_pal_pi_B}

  $$ \pi_{+} ( B^* ) ( e_n ) = {c_{+} ( n + 1 )}^{\frac{1}{2}} e_{ n + 1 }, $$  
  $$ \pi_{-} ( B^* ) ( e_n ) = {c_{-} ( n + 1 )}^{\frac{1}{2}} e_{ n + 1 }. $$
  
  \elmma
  
  {\it Proof :} $ \left\langle  \pi ( B ) ( \sum_n c_n \left ( \begin {array} {cccc}
   e_n \\ 0 \end {array} \right ) ) , \left ( \begin {array} {cccc}
   e_{n^{\prime}} \\ 0 \end {array} \right ) \right\rangle = \sum_n c_n {c_{+} ( n )}^{\frac{1}{2}} \left\langle  e_{n - 1}, ~ e_{n^{\prime}} \right\rangle = c_{n^{\prime} + 1} $
   
   $ {c_{+} ( n^{\prime} + 1 )}^{\frac{1}{2}} = \sum_n c_n {c_+ ( n^{\prime} + 1 )}^{\frac{1}{2}} \left\langle  e_n, ~ e_{n^{\prime} + 1} \right\rangle = \sum_n c_n $
   
   $ \left\langle e_n, ~ \overline{{c_+ ( n^{\prime} + 1 )}^{\frac{1}{2}}} e_{n^{\prime} + 1} \right\rangle = \left\langle \sum_n c_n \left ( \begin {array} {cccc}
   e_n \\ 0 \end {array} \right ) , ~ \left ( \begin {array} {cccc}
 {c_+ ( n^{\prime} + 1 )}^{\frac{1}{2}}  e_{n^{\prime} + 1} \\ 0 \end {array} \right ) \right\rangle.   $
 
 Hence, $ \pi_+ ( B^* ) ( e_n ) = {c_+ ( n + 1 )}^{\frac{1}{2}} e_{n + 1}. $
 
 Similarly, $ \pi_{-} ( B^* ) ( e_n ) = {c_{-} ( n + 1 )}^{\frac{1}{2}} e_{n + 1}. $ \qed
  
  \blmma
   
   \label{sphere_chak_pal_projections_in_algebra}
   
   If $ P_n, Q_n $ denote the projections onto the subspace generated by $ \left ( \begin {array} {cccc}
   e_n \\ 0 \end {array} \right ) $ and $ \left ( \begin {array} {cccc}
   0 \\ e_n \end {array} \right ) $ respectively, then $P_n, Q_n $ belong to $ \pi(S^2_{\mu,c}).$ 
     
  \elmma
  
  {\it Proof :} We claim that  for all $ n \neq 0,~ c_+ ( n ) $ and $ c_{-} ( n ) $ are distinct.
  
               Let $ c_+ ( n ) = c_{-} ( n ).$ Therefore, 
   $ \lambda_{+} \mu^{2n} - {( \lambda_{+} \mu^{2n} )}^{2} + c = \lambda_{-} \mu^{2n} - {( \lambda_{-} \mu^{2n} )}^{2} + c.$ This implies $ ( \lambda_+ + \lambda_{-} ) \mu^{2n} = 1.$ Thus, $ \mu^{2n} = 1 $ and so $ n $ has to be $0.$
   
   Now, for all $ n \geq 1, ~ \pi ( B^* B )\left ( \begin {array} {cccc}
   e_n \\ 0 \end {array} \right ) = {c_{+} ( n )}^{\frac{1}{2}} \pi ( B^* ) \left ( \begin {array} {cccc}
   e_{n - 1} \\ 0 \end {array} \right ) = c_+ ( n ) \left ( \begin {array} {cccc}
   e_n \\ 0 \end {array} \right ).$
   
   Similarly,  for all $ n \geq 1, ~ \pi ( B^* B ) \left ( \begin {array} {cccc}
   0 \\ e_n \end {array} \right ) = c_{-} ( n ) \left ( \begin {array} {cccc}
   0 \\ e_n \end {array} \right ).  $ 
   
   Hence,  for all $ n \geq 1, ~ c_+ ( n ) $ and $ c_{-} ( n ) $ is a discrete distinct set of eigenvalues of $ B^* B $ with eigenspace spanned by $ \left ( \begin {array} {cccc}
   e_n \\ 0 \end {array} \right ) $ and $ \left ( \begin {array} {cccc}
   0 \\ e_n \end {array} \right ) $ respectively. Hence, the eigenprojections corresponding to these eigenvalues belong to the $ C^* ( B^* B ) \subseteq \pi ( S^2_{\mu,c} ).$ Hence, $ P_n, ~ Q_n $ belong to $ \pi ( S^2_{\mu,c} ) $  for all $ n \geq 1.$
   
   Moreover, $ \pi ( A ) \left ( \begin {array} {cccc}
   e_0 \\ 0 \end {array} \right ) = \lambda_+ \left ( \begin {array} {cccc}
   e_0 \\ 0 \end {array} \right )  $ and $ \pi ( A ) \left ( \begin {array} {cccc}
   0 \\ e_0 \end {array} \right ) = \lambda_{-} \left ( \begin {array} {cccc}
   0 \\ e_0 \end {array} \right ) .$
   
   Thus, by the same reasons as above,  $ P_0, ~ Q_0 $ belong to $ \pi ( S^2_{\mu,c} ) .$ \qed

\blmma

\label{sphere_chak_pal_algebra_preservation}
 
$ \pi {( S^{2}_{\mu,c} )}^{''} = \{ \left ( \begin {array} {cccc}
   X_{11} & X_{12}  \\ X_{21} & X_{22} \end {array} \right ) $ belong to $ \clb ( \clh \oplus \clh ) : X_{12} = X_{21} = 0 \} .$
   
 \elmma
 
 {\it Proof :} It suffices to prove that the commutant $ {\pi ( S^{2}_{\mu,c}} )^{\prime} $ is the Von Neumann algebra of operators of the form $ \{ \left ( \begin {array} {cccc}
   c_1 I & 0  \\ 0 & c_2 I \end {array} \right ) $ for some $ c_1, c_2 $ in $ \IC.$ We use the fact that $ \pi_{+} $ and $ \pi_{-} $ are irreducible representations.
 
      Let $ X = \left ( \begin {array} {cccc}
   X_{11} & X_{12}  \\ X_{21} & X_{22} \end {array} \right ) \in { \pi ( S^{2}_{\mu,c} )}^{\prime} .$ Using the fact that $ X $ commutes with $ \pi ( A ),~ \pi ( B ),~ \pi ( B^* ) ,$ we have: $ X_{11} $ belongs to $ {\pi_{+} ( S^{2}_{\mu,c} )}^{\prime} \cong \IC $ and $ X_{22} $ belongs to $ { \pi_{-} ( S^{2}_{\mu,c} )}^{\prime} \cong \IC ,$ so let $ X_{11} = c_1 I, ~ X_{22} = c_2 I $ for some $ c_1,c_2 .$
   
   Moreover, \be \label{sphere_chak_pal_i} X_{12} \pi_{-} ( A ) = \pi_{+} ( A ) X_{12}, \ee   
             \be \label{sphere_chak_pal_ii} X_{12} \pi_{-} ( B ) = \pi_{+} ( B ) X_{12}, \ee             
             \be \label{sphere_chak_pal_iii} X_{12} \pi_{-} ( B^* ) = \pi_{+} ( B^* ) X_{12}. \ee             
   Now, ( \ref{sphere_chak_pal_ii} ) implies $ X_{12} e_0 $ belongs to $ {\rm Ker} ( \pi_{+} ( B ) ) = \IC e_0.$
   Let $ X_{12} e_0= p_0 e_0.$ 
   
   We have, $ \pi_{+} ( B ) ( X_{12} e_1 ) = c^{\frac{1}{2}}_{-} ( 1) X_{12} e_0 = p_0 c^{\frac{1}{2}}_{-} ( 1 ) e_0 ,$ that is, $ \pi_{+} ( B ) ( X_{12} e_1 ) $ belongs to $ \IC e_0.$
   
 Since it follows from the definition of $ \pi_{+} ( B ) $ that $ \pi_{+} ( B ) $ maps $ \overline{\text{span} ~ \{ e_i : i \geq 2 \}} $ to $ {( \IC e_0 )}^{\bot} = \overline{\text{span} \{ e_i : i \geq 1 \}} , ~ X_{12} e_1 $ must belong to $ \text{span} \{ e_0, ~ e_1 \}. $
   
   Inductively, we conclude that  for all $ n,~ X_{12} ( e_n ) $ belongs to $ {\rm span} \{ e_0, e_1,......e_n \}.$
   
   Using the definition of $ \pi_{\pm} ( B^* ) e_n $ along with ( \ref{sphere_chak_pal_iii} ), we have  $ c^{\frac{1}{2}}_{-} ( 1 ) X_{12} e_1 = p_0 c^{\frac{1}{2}}_{+} ( 1 ) e_1 ,$ that is, $ X_{12} e_1 = p_0 \frac{c^{\frac{1}{2}}_{+} ( 1)}{c^{\frac{1}{2}}_{-} ( 1 )} e_1 .$ 
   
   We argue in a similar way  by induction that $ X_{12} e_n = c^{\prime}_{n} e_n $ for some constants $ c^{\prime}_{n}.$
   
   Now we apply  ( \ref{sphere_chak_pal_iii} ) and ( \ref{sphere_chak_pal_ii} ) on the vectors $ e_n$ and $e_{n + 1}$ to get 
   
  $ c^{\prime}_{n + 1} = \frac{c^{\prime}_{n} c^{\frac{1}{2}}_{+} ( n + 1 )}{c^{\frac{1}{2}}_{-} ( n + 1 )} $ and $ c^{\prime}_{n + 1} = \frac{c^{\prime}_{n} c^{\frac{1}{2}}_{-} ( n + 1 )}{c^{\frac{1}{2}}_{+} ( n + 1 )}. $
  
 Since $ c_{+} ( n + 1 ) \neq c_{-} ( n + 1 ) $ for $ n \geq 0,$ we have $ c^{\prime}_{n} = 0 .$
  
  Hence, $ c^{\prime}_{n}  = 0 $ for all $ n $ implying $ X_{12} = 0.$
  
  It follows similarly that $ X_{21} = 0 .$ \qed
   
   \vspace{4mm}  
   
   \subsection{Computation of the quantum isometry group} 
   
   \label{quantumsphere_subsection_chakpal_computation}     
 
 Let $ ( \widetilde{\clq}, \Delta, U ) $ be an object in the category $ {\bf Q^\prime}(D) ,$ with  $ \alpha=\alpha_U$ and  the corresponding Woronowicz $ C^{*} $ subalgebra of $ \widetilde{\clq} $ generated by $\{ <(\xi \ot 1), \alpha(x)(\eta \ot 1)>_{\tilde{\clq}},~\xi, \eta \in \clh,~x \in S^2_{\mu,c} \}$ is denoted by $\clq$. Assume, without loss of generality, that the representation $U$ is faithful. 
 
 Since $ U $ commutes with $ D $, it preserves the eigenvectors  $ \left ( \begin {array} {cccc}
   e_n \\ e_n \end {array} \right ) $ and  $ \left ( \begin {array} {cccc}
   e_n \\ - e_n \end {array} \right ) .$
   
   Let $ U ~  \left ( \begin {array} {cccc}
   e_n \\ e_n \end {array} \right ) ~  = ~  \left ( \begin {array} {cccc}
   e_n \\ e_n \end {array} \right )  \otimes q^{+}_{n}, $
   
   $ U ~  \left ( \begin {array} {cccc}
   e_n \\ - e_n \end {array} \right ) ~  = ~ \left ( \begin {array} {cccc}
   e_n \\ - e_n \end {array} \right )  \otimes q^{-}_{n}, $

   for some $ q^{+}_{n}, q^{-}_{n} $ in $ \widetilde{\clq} .$
   
   \blmma
   
   \label{sphere_chak_pal_alg_0}
   
    For all $ n \geq 0,~ \alpha ( A )  \left ( \begin {array} {cccc}
   e_n \\ 0 \end {array} \right )  =  \left ( \begin {array} {cccc}
   e_n \\ 0 \end {array} \right ) \otimes \frac{1}{4}  \{ \lambda_{+} \mu^{2n} ( 1 + q^{+}_{n} q^{-*}_{n} ) + \lambda_{-} \mu^{2n} ( 1 - q^{+}_{n} q^{-*}_{n} ) + \lambda_{+} \mu^{2n}( 1 + q^{-}_{n} q^{+*}_{n} ) - \lambda_{-}\mu^{2n} (  q^{-}_{n} q^{+*}_{n} - 1 ) \} + \left ( \begin {array} {cccc}
   0 \\ e_n \end {array} \right ) \otimes \frac{1}{4} \{ \lambda_{+} \mu^{2n} ( 1 + q^{+}_{n} q^{-*}_{n} ) + \lambda_{-} \mu^{2n} ( 1 - q^{+}_{n} q^{-*}_{n} ) - \lambda_{+} \mu^{2n} ( 1 + q^{-}_{n} q^{+*}_{n} ) + \lambda_{-}\mu^{2n} (  q^{-}_{n} q^{+*}_{n} - 1 ) \} . $
   
   \vspace{4mm}
   
    For all $ n \geq 0,~ \alpha ( A )  \left ( \begin {array} {cccc}
   0 \\ e_n \end {array} \right )  =  \left ( \begin {array} {cccc}
   e_n \\ 0 \end {array} \right ) \otimes \frac{1}{4} [ \lambda_{+} \mu^{2n} ( 1 - q^{+}_n q^{- *}_{n}   ) + \lambda_{-} \mu^{2n} ( 1 + q^{+}_n q^{-*}_n ) + \lambda_+ \mu^{2n} ( q^{-}_n q^{+*}_n - 1  ) - \lambda_{-} \mu^{2n} ( q^{-}_n q^{+*}_n + 1 ) ] +  \left ( \begin {array} {cccc}
   0 \\ e_n \end {array} \right ) \otimes \frac{1}{4} [ \lambda_+ \mu^{2n} ( 1 - q^{+}_n q^{- *}_n  ) + \lambda_{-} \mu^{2n} ( 1 + q^{+}_n q^{- *}_n  ) - \lambda_{+} \mu^{2n} ( q^{-}_n q^{+*}_n - 1 ) + \lambda_{-} \mu^{2n} ( q^{-}_n q^{+*}_n + 1 )  ]. $
   
  \vspace{4mm}
  
   For all $ n  \geq 1, ~ \alpha ( B )  \left ( \begin {array} {cccc}
   e_n \\ 0 \end {array} \right )  =   \left ( \begin {array} {cccc}
   e_{n - 1} \\ 0 \end {array} \right )  \otimes \frac{1}{4} [ ( {c_{+} ( n )}^{\frac{1}{2}} + {c_{-} ( n )}^{\frac{1}{2}}  ) q^{+}_{n - 1} q^{+*}_n + ( {c_{+} ( n )}^{\frac{1}{2}} - {c_{-} ( n )}^{\frac{1}{2}}  ) q^{-}_{n - 1} q^{+*}_n + ( {c_{+} ( n )}^{\frac{1}{2}} - {c_{-} ( n )}^{\frac{1}{2}}  ) q^{+}_{n - 1} q^{- *}_n + ( {c_{+} ( n )}^{\frac{1}{2}} + {c_{-} ( n )}^{\frac{1}{2}}  ) q^{-}_{n - 1} q^{-*}_n ] +  \left ( \begin {array} {cccc}
   0 \\ e_{n - 1} \end {array} \right )  \otimes \frac{1}{4} [ ( {c_{+} ( n )}^{\frac{1}{2}} + {c_{-} ( n )}^{\frac{1}{2}}  ) q^{+}_{n - 1} q^{+*}_n - ( {c_{+} ( n )}^{\frac{1}{2}} - {c_{-} ( n )}^{\frac{1}{2}}  ) q^{-}_{n - 1} q^{+*}_n +  ( {c_{+} ( n )}^{\frac{1}{2}} - {c_{-} ( n )}^{\frac{1}{2}}  ) q^{+}_{n - 1} q^{-*}_n - ( {c_{+} ( n )}^{\frac{1}{2}} + {c_{-} ( n )}^{\frac{1}{2}}  ) q^{-}_{n - 1} q^{- *}_n  ].$
   
   \vspace{4mm}
   
    For all $ n \geq 1, ~ \alpha ( B ) \left ( \begin {array} {cccc}
   0 \\ e_n \end {array} \right ) = \left ( \begin {array} {cccc}
   e_{n - 1} \\ 0 \end {array} \right ) \otimes \frac{1}{4} [  ( {c_{+} ( n )}^{\frac{1}{2}} + {c_{-} ( n )}^{\frac{1}{2}}  ) q^{+}_{n - 1} q^{+*}_n + ( {c_{-} ( n )}^{\frac{1}{2}} - {c_{+} ( n )}^{\frac{1}{2}}  ) q^{+}_{n - 1} q^{-*}_n + ( {c_{+} ( n )}^{\frac{1}{2}} - {c_{-} ( n )}^{\frac{1}{2}}  ) q^{-}_{n - 1} q^{+*}_n - ( {c_{+} ( n )}^{\frac{1}{2}} + {c_{-} ( n )}^{\frac{1}{2}}  ) q^{-}_{n - 1} q^{-*}_n ] + \left ( \begin {array} {cccc}
   0 \\ e_{n - 1} \end {array} \right ) \otimes \frac{1}{4} [ ( {c_{+} ( n )}^{\frac{1}{2}} + {c_{-} ( n )}^{\frac{1}{2}}  ) q^{+}_{n - 1} q^{+*}_n + ( {c_{-} ( n )}^{\frac{1}{2}} - {c_{+} ( n )}^{\frac{1}{2}}  ) q^{+}_{n - 1} q^{-*}_n - ( {c_{+} ( n )}^{\frac{1}{2}} - {c_{-} ( n )}^{\frac{1}{2}}  ) q^{-}_{n - 1} q^{+*}_n +  ( {c_{+} ( n )}^{\frac{1}{2}} + {c_{-} ( n )}^{\frac{1}{2}}  ) q^{-}_{n - 1} q^{-*}_n ]. $
   
   \vspace{4mm}
   
   For all $ n \geq 0, ~ \alpha ( B^* ) \left ( \begin {array} {cccc}
   e_n \\ 0 \end {array} \right ) = \left ( \begin {array} {cccc}
   e_{n + 1} \\ 0 \end {array} \right ) \otimes \frac{1}{4} [ ( {c_{+} ( n + 1 )}^{\frac{1}{2}} + {c_{-} ( n + 1 )}^{\frac{1}{2}} ) ( q^{+}_{n + 1} q^{+*}_n + q^{-}_{n + 1} q^{- *}_n ) + ( {c_{+} ( n + 1 )}^{\frac{1}{2}} - {c_{-} ( n + 1 )}^{\frac{1}{2}} ) ( q^{+}_{n + 1} q^{-*}_n + q^{-}_{n + 1} q^{+*}_n ) ] + \left ( \begin {array} {cccc}
   0 \\ e_{n + 1} \end {array} \right ) \otimes \frac{1}{4} [ ( {c_{+} ( n + 1 )}^{\frac{1}{2}} + {c_{-} ( n + 1 )}^{\frac{1}{2}} ) ( q^{+}_{n + 1} q^{+*}_n - q^{-}_{n + 1} q^{-*}_n  ) + ( {c_{+} ( n + 1 )}^{\frac{1}{2}} - {c_{-} ( n + 1 )}^{\frac{1}{2}} ) ( q^{+}_{n + 1} q^{-*}_n - q^{-}_{n + 1} q^{+*}_n  )  ]. $

   \elmma
   
 {\it Proof :}  One has,
 $ \alpha ( A ) ( \left ( \begin {array} {cccc}
   e_n \\ 0 \end {array} \right )  ) $
   
   $ = \widetilde{U} ( A \otimes 1 ) {\widetilde{U}}^{*} ( \left ( \begin {array} {cccc}
   e_n \\ 0 \end {array} \right ) \otimes 1 ) $
   
   $ = \frac{1}{2} \widetilde{U} ( \pi ( A ) \otimes 1 ) {\widetilde{U}}^{*} ( \left ( \begin {array} {cccc}
   e_n \\ e_n \end {array} \right ) \otimes 1 + \left ( \begin {array} {cccc}
   e_n \\ - e_n \end {array} \right ) \otimes 1  ) $
   
   $ = \frac{1}{2} \widetilde{U} ( \pi ( A ) \otimes 1 ) [ \left ( \begin {array} {cccc}
   e_n \\ e_n \end {array} \right ) \otimes q^{+ *}_{n} +  \left ( \begin {array} {cccc}
   e_n \\ - e_n \end {array} \right ) \otimes q^{- *}_{n} ] $
   
   $ = \frac{1}{2} \widetilde{U} [ \left ( \begin {array} {cccc}
 \pi_+ ( A ) e_n \\ \pi_{-} ( A ) e_n \end {array} \right ) \otimes q^{+ *}_{n} +  \left ( \begin {array} {cccc}
 \pi_+ (  A ) e_n \\ - \pi_{-} ( A ) e_n \end {array} \right ) \otimes q^{- *}_{n} ]. $
 
 By using ( \ref{sphere_chakraborty_pal_pi ( A )} ), we get this to be equal to 
   
   $ = \frac{1}{2} \widetilde{U} [ \left ( \begin {array} {cccc}
  \lambda_{+} \mu^{2n} e_n \\ \lambda_{-} \mu^{2n} e_n \end {array} \right ) \otimes q^{+*}_{n} + \left ( \begin {array} {cccc}
   \lambda_{+} \mu^{2n}e_n \\ - \lambda_{-} \mu^{2n} e_n \end {array} \right ) \otimes q^{-*}_{n} ] $
   
   $ = \frac{1}{2} \widetilde{U} [ \left ( \begin {array} {cccc}
   e_n \\ 0 \end {array} \right ) \otimes \lambda_{+} \mu^{2n} ( q^{+*}_{n} + q^{- *}_{n} ) + \left ( \begin {array} {cccc}
   0 \\ e_n \end {array} \right ) \otimes \lambda_{-} \mu^{2n} ( q^{+*}_{n} - q^{- *}_{n}  )  ] $
   
   $ = \frac{1}{4} \widetilde{U} [ \left ( \begin {array} {cccc}
   e_n \\  e_n \end {array} \right ) \otimes \{ \lambda_{+} \mu^{2n} ( q^{+*}_{n} + q^{-*}_{n} ) + \lambda_{-} \mu^{2n} ( q^{+*}_{n} - q^{-*}_{n} ) \} + \left ( \begin {array} {cccc}
   e_n \\ - e_n \end {array} \right ) \otimes \{ \lambda_{+} \mu^{2n} ( q^{+*}_{n} + q^{-*}_{n} ) - \lambda_{-} \mu^{2n} ( q^{+*}_{n} - q^{-*}_{n} ) \}  ] $
   
   $ = \left ( \begin {array} {cccc}
   e_n \\ e_n \end {array} \right ) \otimes \frac{1}{4} q^{+}_n \{ \lambda_+ \mu^{2n} ( q^{+*}_n + q^{-*}_{n} ) + \lambda_{-} \mu^{2n} ( q^{+*}_{n} - q^{- *}_{n}  ) \} +  \left ( \begin {array} {cccc}
   e_n \\ - e_n \end {array} \right ) \otimes \frac{1}{4} q^{-}_{n} \{ \lambda_+ \mu^{2n} ( q^{+*}_n + q^{- *}_n  ) - \lambda_{-} \mu^{2n} ( q^{+*}_n - q^{-*}_n ) \} $
   
   $ =  \left ( \begin {array} {cccc}
   e_n \\ 0 \end {array} \right ) \otimes \frac{1}{4}  \{ \lambda_{+} \mu^{2n} ( 1 + q^{+}_{n} q^{-*}_{n} ) + \lambda_{-} \mu^{2n} ( 1 - q^{+}_{n} q^{-*}_{n} ) + \lambda_{+} \mu^{2n}( 1 + q^{-}_{n} q^{+*}_{n} ) - \lambda_{-}\mu^{2n} (  q^{-}_{n} q^{+*}_{n} - 1 ) \} + \left ( \begin {array} {cccc}
   0 \\ e_n \end {array} \right ) \otimes \frac{1}{4} \{ \lambda_{+} \mu^{2n} ( 1 + q^{+}_{n} q^{-*}_{n} ) + \lambda_{-} \mu^{2n} ( 1 - q^{+}_{n} q^{-*}_{n} ) - \lambda_{+} \mu^{2n} ( 1 + q^{-}_{n} q^{+*}_{n} ) + \lambda_{-}\mu^{2n} (  q^{-}_{n} q^{+*}_{n} - 1 ) \} . $
   
   
   Similarly, $ \alpha ( A ) \left ( \begin {array} {cccc}
   0 \\ e_n \end {array} \right )   $
   
     $ = \widetilde{U} ( \pi ( A ) \otimes 1  ) {\widetilde{U}}^{*} ( \left ( \begin {array} {cccc}
   0 \\ e_n \end {array} \right ) \otimes 1 )  $
   
    $ =  \frac{1}{2} \widetilde{U} ( \pi ( A ) \otimes 1 ) {\widetilde{U}}^{*}  ( [ \left ( \begin {array} {cccc}
   e_n \\ e_n \end {array} \right ) - \left ( \begin {array} {cccc}
   e_n \\ - e_n \end {array} \right ) ] \otimes 1  )  $
   
   $ = \frac{1}{2} \widetilde{U} [  \left ( \begin {array} {cccc}
  \pi_+ ( A ) e_n \\ \pi_{-} ( A ) e_n \end {array} \right ) \otimes q^{+*}_n - \left ( \begin {array} {cccc}
  \pi_+ ( A )  e_n \\ - \pi_{-} ( A ) e_n \end {array} \right ) \otimes q^{- *}_n ]  $
   
   $ = \frac{1}{2} \widetilde{U} [ \left ( \begin {array} {cccc}
   \lambda_+ \mu^{2n} e_n \\ \lambda_{-} \mu^{2n} e_n \end {array} \right ) \otimes q^{+*}_n  -  \left ( \begin {array} {cccc}
  \lambda_+ \mu^{2n} e_n \\ - \lambda_{-} \mu^{2n} e_n \end {array} \right ) \otimes q^{-*}_n  ]  $
  
  $ = \frac{1}{4} \widetilde{U} [ \left ( \begin {array} {cccc}
   e_n \\ e_n \end {array} \right ) \otimes \{ \lambda_+ \mu^{2n} ( q^{+*}_n - q^{- *}_n  ) + \lambda_{-} \mu^{2n} ( q^{+*}_n + q^{- *}_n  ) \} + \left ( \begin {array} {cccc}
   e_n \\ - e_n \end {array} \right ) \otimes \{ \lambda_+ \mu^{2n} ( q^{+*}_n - q^{-*}_n  ) - \lambda_{-} \mu^{2n} (  q^{+*}_n + q^{-*}_n )  \}  ] $
   
  $ = \left ( \begin {array} {cccc}
   e_n \\ 0 \end {array} \right ) \otimes \frac{1}{4} [ \lambda_+ \mu^{2n} ( 1 - q^{+}_n q^{-*}_n ) + \lambda_- \mu^{2n} (  1 + q^{+}_n q^{-*}_n ) + \lambda_+ \mu^{2n} ( q^{-}_n q^{+*}_n - 1 ) - \lambda_{-} \mu^{2n} (  q^{-}_n q^{+*}_n + 1 ) ] + \left ( \begin {array} {cccc}
   0 \\ e_n \end {array} \right ) \otimes \frac{1}{4} [ \lambda_+ \mu^{2n} ( 1 - q^{+}_n q^{-*}_n ) + \lambda_{-} \mu^{2n} ( 1 + q^{+}_n q^{-*}_n  ) - \lambda_+ \mu^{2n} ( q^{-}_n q^{+*}_n - 1  ) + \lambda_{-} \mu^{2n} ( q^{-}_n q^{+*}_n + 1  )  ]. $

  As the proof of the others are exactly similar, we omit the proofs.

 \blmma
 
 \label{sphere_chak_pal_alg}
 We have:
\be \label{sphere_chak_pal_alg1} q^{+}_{n}q^{-*}_{n} = q^{-}_{n}q^{+*}_{n} ~ {\rm ~ for ~ all}~ n ,\ee
\bean  ( {c_{+}( n )}^{\frac{1}{2}} + {c_{-}( n )}^{\frac{1}{2}} ) ( q^{+}_{n - 1} q^{+*}_{n} - q^{-}_{n - 1} q^{-*}_{n} ) +  ( {c_{+}( n )}^{\frac{1}{2}} - {c_{-}( n )}^{\frac{1}{2}} ) ( q^{+}_{n - 1} q^{-*}_{n} - q^{-}_{n - 1} q^{+*}_{n} ) = 0 \eean 
\be \label{sphere_chak_pal_alg2}  ~ {\rm ~ for ~ all}~ n \geq 1, \ee
\bean  ( {c_{+}( n )}^{\frac{1}{2}} + {c_{-}( n )}^{\frac{1}{2}} ) ( q^{+}_{n - 1} q^{+*}_{n} - q^{-}_{n - 1} q^{-*}_{n} ) +  ( {c_{+}( n )}^{\frac{1}{2}} - {c_{-}( n )}^{\frac{1}{2}} ) ( q^{-}_{n - 1} q^{+*}_{n} - q^{+}_{n - 1} q^{-*}_{n} ) = 0 \eean
\be \label{sphere_chak_pal_alg3} ~ {\rm ~ for ~ all}~ n \geq 1, \ee
\bean ( {c_{+}( n + 1 )}^{\frac{1}{2}} + {c_{-}( n + 1 )}^{\frac{1}{2}} ) ( q^{+}_{n + 1} q^{+*}_{n} - q^{-}_{n + 1} q^{-*}_{n} ) =  ( {c_{+}( n + 1 )}^{\frac{1}{2}} - {c_{-}( n + 1 )}^{\frac{1}{2}} ) ( q^{-}_{n + 1} q^{+*}_{n} -  \eean
\be  \label{sphere_chak_pal_alg4} q^{+}_{n + 1} q^{-*}_{n} )  ~ {\rm ~ for ~ all} ~ n.  \ee
 
 \elmma
 
{\it Proof :} Since $ \alpha ( A ) $ maps $ \pi( S^{2}_{\mu,c} ) $ into its double commutant, we conclude by using the description of $ {\pi ( S^{2}_{\mu,c} )}^{\prime \prime} $ given in Lemma \ref{sphere_chak_pal_algebra_preservation} that the coefficient of $ \left ( \begin {array} {cccc}
   0 \\ e_n \end {array} \right ) $ in $ \alpha ( A ) \left ( \begin {array} {cccc}
   e_n \\ 0 \end {array} \right ) $ must be $ 0 ,$ which implies ( by Lemma \ref{sphere_chak_pal_alg_0} )
   
   $ \lambda_+ [ 1 + q^{+}_n q^{-*}_n - ( 1 + q^{-}_n q^{+*}_n   ) ] + \lambda_{-} [ 1 - q^{+}_n q^{-*}_{n} + q^{-}_n q^{+*}_n - 1 ] = 0. $
   
   Hence,
   
    $ ( \lambda_+ - \lambda_{-} ) ( q^{+}_{n}q^{-*}_{n} - q^{-}_{n}q^{+*}_{n} ) = 0 .$ Hence, $ ( q^{+}_{n}q^{-*}_{n} - q^{-}_{n}q^{+*}_{n} ) = 0.$
   
 Proceeding in a similar way, ( \ref{sphere_chak_pal_alg2} ),( \ref{sphere_chak_pal_alg3} ), ( \ref{sphere_chak_pal_alg4} ) follow  from the facts that coefficients of $  \left ( \begin {array} {cccc}
   0 \\ e_{n - 1} \end {array} \right ) , \left ( \begin {array} {cccc}
   e_{n - 1} \\ 0 \end {array} \right ) $ and $ \left ( \begin {array} {cccc}
   0 \\ e_{n + 1} \end {array} \right ) $ in $ \alpha ( B )  \left ( \begin {array} {cccc}
   e_n \\ 0 \end {array} \right ) ,$
   
  $  \alpha ( B ) \left ( \begin {array} {cccc}
   0 \\ e_{n} \end {array} \right ) $ and   $ \alpha ( B^* ) \left ( \begin {array} {cccc}
   e_n \\ 0 \end {array} \right ) $ ( respectively ) are zero. \qed
   
   \bcrlre
   
   \label{sphere_chak_pal_corollary}
   
   We have 
   
   \bean \alpha ( A ) = \sum^{\infty}_{n = 0} A P_n \otimes \frac{1}{2 \lambda_{+}} \{ \lambda_{+} ( 1 + q^{+}_{n} q^{-*}_{n} ) + \lambda_{-} ( 1 - q^{+}_{n} q^{-*}_{n} ) \} \eean   
   \bean + \sum^{\infty}_{n = 0} A Q_n \otimes \frac{1}{2 \lambda_{-}} \{ \lambda_{+} ( 1 - q^{+}_{n} q^{-*}_{n} ) + \lambda_{-} ( 1 + q^{+}_{n} q^{-*}_{n} ) \}. \eean     
    \bean \alpha ( B ) = \sum^{\infty}_{n = 1} B P_n \otimes \frac{1}{4c_{+}( n )} [ ( {c_{+} ( n )}^{\frac{1}{2}} + {c_{-} ( n )}^{\frac{1}{2}} ) ( q^{+}_{n - 1} q^{+*}_{n} + q^{-}_{n - 1} q^{-*}_{n} ) \eean    
    \bean + ( {c_{+} ( n )}^{\frac{1}{2}} - {c_{-} ( n )}^{\frac{1}{2}} ) ( q^{-}_{n - 1} q^{+*}_{n} + q^{+}_{n - 1} q^{-*}_{n} ) ] + \sum^{\infty}_{n = 1} B Q_n \otimes \frac{1}{4 c_{-}( n )} [ ( {c_{+} ( n )}^{\frac{1}{2}} + {c_{-} ( n )}^{\frac{1}{2}}  ). \eean    
    \bean ( q^{+}_{n - 1} q^{+*}_{n} + q^{-}_{n - 1} q^{-*}_{n} ) - ( {c_{+} ( n )}^{\frac{1}{2}} - {c_{-} ( n )}^{\frac{1}{2}} ) ( q^{+}_{n - 1} q^{-*}_{n} + q^{-}_{n - 1} q^{+*}_{n} ) ]. \eean

   \ecrlre
   
   {\it Proof :} 
We note that $ \pi ( A ) \left ( \begin {array} {cccc}
   e_n \\ 0 \end {array} \right ) = \left ( \begin {array} {cccc}
   \pi_+ ( A ) & 0 \\ 0 & \pi_{-} ( A )  \end {array} \right ) \left ( \begin {array} {cccc}
   e_n \\ 0 \end {array} \right ) = \left ( \begin {array} {cccc}
   \pi_{+} ( A ) e_n \\ 0 \end {array} \right ) = \lambda_+ \mu^{2n} \left ( \begin {array} {cccc}
   e_n \\ 0 \end {array} \right ).$
   
   Thus, $ \left ( \begin {array} {cccc}
   e_n \\ 0 \end {array} \right ) = \frac{ \pi ( A ) \left ( \begin {array} {cccc}
   e_n \\ 0 \end {array} \right ) } {\lambda_+ \mu^{2n}}. $
   
   Similarly, $ \left ( \begin {array} {cccc}
   0 \\ e_{n} \end {array} \right ) = \frac{\pi ( A ) \left ( \begin {array} {cccc}
   0 \\ e_n \end {array} \right ) } {\lambda_{-} \mu^{2n}} .$
   
   \vspace{4mm}
   
   Now, using (  \ref{sphere_chak_pal_alg1}  ), 
   $$  \alpha ( A ) \left ( \begin {array} {cccc}
   e_n \\ 0 \end {array} \right ) = \left ( \begin {array} {cccc}
   e_n \\ 0 \end {array} \right ) \otimes \frac{1}{2} \{ \lambda_+ \mu^{2n} ( 1 + q^+_n q^{- *}_n ) + \lambda_{-} \mu^{2n} ( 1 - q^{+}_n q^{-*}_n ) \} $$
   
   $$ = \pi( A ) \left ( \begin {array} {cccc}
   e_n \\ 0 \end {array} \right ) \otimes \frac{1}{2 \lambda_{+}} \{ \lambda_+ (  1 + q^+_n q^{-*}_n ) + \lambda_{-} ( 1 - q^+_n q^{-*}_n ) \} .  $$
   
   Similarly, $$ \alpha ( A ) \left ( \begin {array} {cccc}
   0 \\ e_n \end {array} \right ) = \pi ( A ) \left ( \begin {array} {cccc}
   0 \\ e_n \end {array} \right ) \otimes \frac{1}{2 \lambda_{-}} \{ \lambda_+ ( 1 - q^+_n q^{-*}_n ) + \lambda_{-} ( 1 + q^{+}_n q^{-*}_n ) \} .  $$
   
   Thus, $ \alpha ( A ) = \sum^{\infty}_{n = 0} A P_n \otimes \frac{1}{2 \lambda_{+}} \{ \lambda_{+} ( 1 + q^{+}_{n} q^{-*}_{n} ) + \lambda_{-} ( 1 - q^{+}_{n} q^{-*}_{n} ) \} + + \sum^{\infty}_{n = 0} A Q_n \otimes \frac{1}{2 \lambda_{-}} \{ \lambda_{+} ( 1 - q^{+}_{n} q^{-*}_{n} ) + \lambda_{-} ( 1 + q^{+}_{n} q^{-*}_{n} ) \}.$

By similar considerations from   $~  \alpha ( B ) \left ( \begin {array} {cccc}
   e_n \\ 0 \end {array} \right ) , \alpha ( B ) \left ( \begin {array} {cccc}
   0 \\ e_{n} \end {array} \right ) $, the result follows.  \qed

 \blmma
 
 \label{sphere_chak_pal_alpha_Pn_tilda}
 
 Let $ \widetilde{P_n} = P_n + Q_n.$ Then for each vector $ v $ in  $ \clh, ~ \alpha ( \widetilde{P_n} ) v = \widetilde{P_n} v \otimes 1.$
 
 \elmma
 
 {\it Proof :} To start with, we recall that $ P_n $ and $ Q_n $ belong to $ \pi ( S^{2}_{\mu,c} ) $ ( Lemma \ref{sphere_chak_pal_projections_in_algebra} ). Hence, $ \widetilde{P_n} $ belongs to $ \pi ( S^{2}_{\mu,c} ).$                 \bean \lefteqn{ \alpha ( \widetilde{P_n} )\left ( \begin {array} {cccc}
   e_n \\ e_n \end {array} \right )  }\\
                       &=& \widetilde{U} (  \widetilde{P_n}  \otimes 1 ) {\widetilde{U}}^{*}   \left ( \begin {array} {cccc}
   e_n \\ e_n \end {array} \right ) \eean
     \bean &=& \widetilde{U} (  \left ( \begin {array} {cccc}
   e_n \\ e_n \end {array} \right ) \otimes q^{+*}_n )\\
   &=& \left ( \begin {array} {cccc}
   e_n \\ e_n \end {array} \right ) \otimes 1\\
   &=& ( \widetilde{P_n} \otimes 1 ) ( \left ( \begin {array} {cccc}
   e_n \\ e_n \end {array} \right ) \otimes 1 ). \eean
   
   For $ k \neq n ,$   
   \bean \lefteqn{ \widetilde{U} ( \widetilde{P_k} \otimes 1 ) {\widetilde{U}}^{*} \left ( \begin {array} {cccc}
   e_n \\ e_n \end {array} \right )}\\
      &=& \widetilde{U} ( \widetilde{P_k} \otimes 1  ) ( \left ( \begin {array} {cccc}
   e_n \\ e_n \end {array} \right ) \otimes q^{+*}_n   )\\
    &=& 0\\
    &=& ( \widetilde{P_k} \otimes 1  ) ( \left ( \begin {array} {cccc}
   e_n \\ e_n \end {array} \right ) \otimes 1 ). \eean

   Similarly, \bean \lefteqn{ \alpha ( \widetilde{P_n} ) \left ( \begin {array} {cccc}
   e_n \\ - e_n \end {array} \right )}\\
    &=& \widetilde{U} ( \left ( \begin {array} {cccc}
   e_n \\ - e_n \end {array} \right ) \otimes 	q^{-*}_n  )\\
   &=&  \left ( \begin {array} {cccc}
   e_n \\ - e_n \end {array} \right ) \otimes 1\\
   &=& ( \widetilde{P_n} \otimes 1  ) ( \left ( \begin {array} {cccc}
   e_n \\ - e_n \end {array} \right ) \otimes 1 ), \eean
   
   and for $$ k \neq n, ~ \alpha ( \widetilde{P_k} ) \left ( \begin {array} {cccc}
   e_n \\ - e_n \end {array} \right ) = 0 = ( \widetilde{P_k} \otimes 1 ) ( \left ( \begin {array} {cccc}
   e_n \\  - e_n \end {array} \right ) \otimes 1 ). $$   
   Combining all these, we get the required result. \qed

 \bppsn
   
   \label{sphere_chak_pal_description of Q'+}

   As a $C^*$-algebra, $ \widetilde{\clq}$ is  generated by the unitaries $  \{ q^{+}_{n} \}_{n \geq 0}$, and the self-adjoint unitary $y_0  = q^{-*}_{0} q^{+}_{0} .$ Moreover, $ \clq $ is generated by unitaries $ z_n = q^{+}_{n - 1} q^{+*}_{n}, ~ n \geq 1 $ along with a self adjoint unitary $ w^{\prime} .$ 
   
   \eppsn

 {\it Proof :} Replacing $ n + 1 $ by $ n $ in ( \ref{sphere_chak_pal_alg4} ) we have, 
\be  \label{sphere_chak_pal_id5} ( {c_{+}( n  )}^{\frac{1}{2}} + {c_{-}( n )}^{\frac{1}{2}} ) ( q^{+}_{n} q^{+*}_{n - 1} - q^{-}_{n} q^{-*}_{n- 1} ) - ( {c_{+}( n )}^{\frac{1}{2}} - {c_{-}( n )}^{\frac{1}{2}} ) ( q^{-}_{n} q^{+*}_{n - 1} - q^{+}_{n} q^{-*}_{n - 1} ) = 0 ~ {\rm for ~ all}~ n \geq 1. \ee
Subtracting ( \ref{sphere_chak_pal_id5} ) from the equation obtained by applying $ \ast  $ on ( \ref{sphere_chak_pal_alg2} ), we have : $ 2 ( {c_{+}( n  )}^{\frac{1}{2}} - {c_{-}( n )}^{\frac{1}{2}} ) ( q^{-}_{n} q^{+*}_{n - 1} - q^{+}_{n} q^{-*}_{n - 1} ) = 0 $ for all $ n \geq 1  .$ Now, from the proof of Lemma \ref{sphere_chak_pal_projections_in_algebra},  $ ( {c_{+}( n  )}^{\frac{1}{2}} - {c_{-}( n )}^{\frac{1}{2}} ) \neq 0$ for all $ n \geq 1.$ This implies :
\be \label{sphere_chak_pal_id6} q^{-}_{n} q^{+*}_{n - 1} = q^{+}_{n} q^{-*}_{n - 1} ~ {\rm for ~ all}~ n \geq 1. \ee
Using ( \ref{sphere_chak_pal_id6} ) in ( \ref{sphere_chak_pal_id5} ), we have 
 \be \label{sphere_chak_pal_id7} q^{+}_{n} q^{+*}_{n - 1} = q^{-}_{n} q^{-*}_{n - 1} ~ {\rm ~ for ~ all}~ n \geq 1. \ee 
 Let $ y_n = q^{-*}_{n} q^{+}_{n}.$
 
 Then, using  ( \ref{sphere_chak_pal_alg1} ),  we have $ q^{-*}_{n} q^{+}_n = q^{+*}_n q^{-}_n. $ Hence, $ y_n = y^{*}_n. $ Moreover, as $ y_n $ is a product of unitaries, it is a self adjoint unitary.
 
 Now, from ( \ref{sphere_chak_pal_id6} ), we have $q^{-}_{n} = q^{+}_n y^{*}_{n - 1} $ for all $ n \geq 1.$ Therefore, 
  \be \label{sphere_chak_pal_id8} q^{-}_{n} = q^{+}_{n} y_{n - 1} ~ {\rm ~ for ~ all}~ n \geq 1. \ee

  
  
  
 Next, from ( \ref{sphere_chak_pal_id7}  ), we obtain  $ q^{-*}_{n} q^{+}_{n} = q^{-*}_{n - 1} q^{+}_{n - 1} $ for all $n \geq 1$ implying  
  \be \label{sphere_chak_pal_id9} y_n = y_{n - 1} ~ {\rm ~ for ~ all}~ n \geq 1. \ee  
  From ( \ref{sphere_chak_pal_id8} ) and ( \ref{sphere_chak_pal_id9} ) and the faithfulness of the representation $U$, we conclude that $ \widetilde{\clq} $ is generated by $ \{ q^{+}_{n} \}_{n \geq 0} $ and $ y_{0}.$

 Now we prove the second part of the proposition.
 

 Using Lemma \ref{sphere_chak_pal_alpha_Pn_tilda},  we note that  for all $ v $ in $ \clh, ~ \alpha ( A \widetilde{P_k} )v = \alpha ( A ) ( \widetilde{P_k}v \otimes 1 ) = A P_k v \otimes \frac{1}{2 \lambda_{+}} \{ \lambda_{+} ( 1 + q^{+}_{k}q^{-*}_{k} ) + \lambda_{-} ( 1 - q^{+}_{k} q^{-*}_{k} ) \} +  A Q_k v \otimes \frac{1}{2 \lambda_{-}} \{ \lambda_{+} ( 1 - q^{+}_{k} q^{-*}_{k} ) + \lambda_{-} ( 1 + q^{+}_{k} q^{-*}_{k} ) \} .$ Therefore, $ \alpha ( A \widetilde{P_k} ) = A P_k  \otimes \frac{1}{2 \lambda_{+}} \{ \lambda_{+} ( 1 + q^{+}_{k}q^{-*}_{k} ) + \lambda_{-} ( 1 - q^{+}_{k} q^{-*}_{k} ) \} +  A Q_k  \otimes \frac{1}{2 \lambda_{-}} \{ \lambda_{+} ( 1 - q^{+}_{k} q^{-*}_{k} ) + \lambda_{-} ( 1 + q^{+}_{k} q^{-*}_{k} ) \} .$  
  
 Now, $ A P_k $ and $ A Q_k $ being distinct elements, there exist linear functional  $ \phi $
 such that $ \phi ( A P_k ) = 1, ~ \phi ( A Q_k ) = 0 $ and vice versa. Hence, $ ( \phi \otimes {\rm id} ) \alpha ( A \widetilde{P_k}  ) =  \lambda_{+} ( 1 + q^{+}_{m}q^{-*}_{m} ) + \lambda_{-} ( 1 - q^{+}_{m} q^{-*}_{m} ) $ belongs to $ \clq.$   
  Similarly, $ \lambda_{+} ( 1 - q^{+}_{m} q^{-*}_{m} ) + \lambda_{-} ( 1 + q^{+}_{m} q^{-*}_{m} ) $ belongs to $ \clq $  for all $ m.$
   
   Subtracting, we get $ q^{+}_{m} q^{-*}_{m} $ belongs to $ \clq.$ 
   
 Using the expression of $ \alpha ( B ) $ in a similar way, we have 
 
  $ ( {c_{+} ( n )}^{\frac{1}{2}} + {c_{-} ( n )}^{\frac{1}{2}} ) ( q^{+}_{n - 1} q^{+*}_{n} + q^{-}_{n - 1} q^{-*}_{n} ) + ( {c_{+} ( n )}^{\frac{1}{2}} - {c_{-} ( n )}^{\frac{1}{2}} ) ( q^{-}_{n - 1} q^{+*}_{n} + q^{+}_{n - 1} q^{-*}_{n} ) $ belongs to $ \clq $ for all $ n \geq 1.$ 
  
  and  $ ( {c_{+} ( n )}^{\frac{1}{2}} + {c_{-} ( n )}^{\frac{1}{2}} ) ( q^{+}_{n - 1} q^{+*}_{n} + q^{-}_{n - 1} q^{-*}_{n} ) - ( {c_{+} ( n )}^{\frac{1}{2}} - {c_{-} ( n )}^{\frac{1}{2}} ) ( q^{+}_{n - 1} q^{-*}_{n} + q^{-}_{n - 1} q^{+*}_{n} ) $ belongs to $ \clq $ for all $ n \geq 1.$
   
   Adding and subtracting, we have   
   \be \label{sphere_chak_pal_id10}   q^{+}_{n - 1} q^{+*}_{n} + q^{-}_{n - 1} q^{-*}_{n} \in \clq  ~ {\rm ~ for ~ all} ~ n \geq 1, \ee    
   \be \label{sphere_chak_pal_id11}  q^{-}_{n - 1} q^{+*}_{n} + q^{+}_{n - 1} q^{-*}_{n} \in \clq ~ {\rm ~ for ~ all} ~ n \geq 1. \ee   
  Recalling ( \ref{sphere_chak_pal_id9} ), we have $ q^{-}_{n} = q^{+}_{n} y_{n - 1} = q^{+}_{n} y_0 .$ Using this  in ( \ref{sphere_chak_pal_id10} ), we obtain  
  \bean \lefteqn{ q^{+}_{n - 1} q^{+*}_n + q^{-}_{n - 1} q^{-*}_n}\\
           &=& q^{+}_{n - 1} q^{+*}_n + q^{+}_{n - 1} y_0 y^{*}_0 q^{+*}_n\\
           &=& q^{+}_{n - 1} q^{+*}_n + q^{+}_{n - 1} q^{+*}_n\\
           &=& 2 q^{+}_{n - 1} q^{+*}_n. \eean  
   Similarly, using $ q^{-}_{n} = q^{+}_{n} y_0 $ in ( \ref{sphere_chak_pal_id11} ), one has    
   \bean \lefteqn{q^{-}_{n - 1} q^{+*}_n + q^{+}_{n - 1} q^{-*}_n }\\
          &=& q^{+}_{n - 1} y_0 q^{+*}_n + q^{+}_{n - 1} y_0 q^{+*}_n ~ = ~ 2 q^{+}_{n - 1} y_0 q^{+*}_n. \eean   
  Hence, we conclude that $ q^{+}_{n - 1} q^{+*}_{n} $ and $ q^{+}_{n - 1} y_0 q^{+*}_{n} $ are in $ \clq $ for all $ n \geq 1.$
  
  Let $$ z_n = q^{+}_{n - 1} q^{+*}_{n}, $$  
      $$ w_n =  q^{+}_{n - 1} y_0 q^{+*}_{n}, $$      
       for all $ n \geq 1 .$
      
      Then, we observe that \bean z^{*}_{n} w_n = q^{+}_{n} y_0 q^{+*}_{n} = q^{+}_{n} q^{-*}_{n}. \eean      
      Moreover, \bean q^{+}_0 q^{-*}_0 = q^{+}_0 {( q^{+}_0 y^{*}_0 )}^{*} = q^{+}_0 y_0 q^{+*}_0 = q^{+}_0 y_0 q^{+*}_1  q^{+}_1 q^{+*}_0 = w_1 z^{*}_1.\eean      
    Thus  for all $ n \geq 0, ~ q^{+}_{n} q^{-*}_{n} $ belong to $ C^{*} ( \{ z_n, w_n \}_{n \geq 1} ) .$     
   \bean {\rm~ For ~ all} ~ n \geq 2, ~ q^{-}_{n - 1} q^{-*}_n = q^{+}_{n - 1} y^{*}_{n - 2} {( q^{+}_n y^{*}_{n -1} )}^* = q^{+}_{n - 1} y^*_0 y_0 q^{+*}_n = q^{+}_{n - 1} q^{+*}_n = q^{+}_{n - 1} q^{+*}_n = z_n, \eean   
   \bean q^{-}_0 q^{-*}_1 = q^{+}_0 y_0 {( q^{+}_1 y_0)}^{*} = q^{+}_0 y_0 y^{*}_0 q^{+*}_1 = q^{+}_0 q^{+*}_1 = z_1.\eean
   
   Finally,\bean q^{-}_{n - 1} q^{+*}_n = q^{+}_{n - 1} y^*_0 q^{+*}_n = w_n \eean  and \bean q^{-}_0 q^{+*}_1 = q^{+}_0 y_0 q^{+*}_1 = w_1. \eean     
    Now, from the expressions of $ \alpha ( A ) $ and $\alpha( B ) ,$ it is clear that $ \clq $ is generated by $ q^{+}_{n} {q^{-}_{n}}^{*}, ~ q^{+}_{n - 1} {q^{+}_{n}}^{*} + q^{-}_{n - 1} {q^{-}_{n}}^{*}, ~ q^{-}_{n - 1} {q^{+}_{n}}^{*} + q^{+}_{n - 1} {q^{-}_{n}}^{*} .$ By the above made observations, these belong to $ C^{*} ( {\{ z_n, w_n \}}_{n \geq 1} ) $ which implies that $ \clq $ is a $ C^* $ subalgebra of $ C^{*} ( {\{ z_n, w_n \}}_{n \geq 1} ) .$ Moreover, from the definitions of $ z_n,~ w_n $ it is clear that $ C^{*} ( {\{ z_n, w_n \}}_{n \geq 1} ) $ is a $ C^* $ subalgebra of $ \clq .$ 
    
   Therefore, $ \clq \cong C^{*} ( \{ z_n, w_n \}_{n \geq 1} ) .$
    
    In fact, a simpler description is possible by noting that $ z_n w_{n + 1} = q^{+}_{n - 1} q^{+*}_n q^{+}_n y_0 q^{+*}_{n + 1} = q^{+}_{n - 1} y_0 q^{+*}_{n + 1} =  q^{+}_{n - 1} y_0 q^{+*}_n q^{+}_n q^{+*}_{n + 1} = w_n z_{n + 1} $ and so,  $ w_{n + 1} = {z_n}^{*} w_n z_{n + 1} $ which implies  $ \{ w_n \}_{n \geq 1} $ is a subset of $ C^{*} ( \{ z_n \}_{n \geq 1}, ~ w_1 ) .$
    
    Defining $ w^{\prime} = w^{*}_{1} z_1 ,$ we note that  $ z_1 = q^{+}_0 y_0 q^{+*}_1  q^{+}_1 y^*_0 q^{+*}_1 = w_1 ( z^*_1 w_1 ) $ which implies $ w^*_1 z_1 = z^*_1 w_1.$ Thus, $ w^{\prime} $ is self adjoint. It is a unitary as it is a product of unitaries.
    
    Thus $ \clq \cong C^{*} \{ \{ z_n \}_{n \geq 1}, ~ w^{\prime} \}.$ \qed







    \blmma
    
    \label{sphere_chak_pal_coproduct}
    
   $ \Delta ( q^{\pm}_{n} ) = q^{\pm}_{n} \otimes q^{\pm}_{n}, $
    
    $ \Delta ( y_1 ) = y_1 \otimes y_1. $
    
    \elmma
   
  {\it Proof :}    We use the fact that $ U $ is a unitary representation.
    
   $ ( {\rm id} \otimes \Delta ) U \left ( \begin {array} {cccc}
   e_n \\ e_n \end {array} \right ) = ( {\rm id} \otimes \Delta ) ( \left ( \begin {array} {cccc}
   e_n \\ e_n \end {array} \right ) \otimes q^{+}_{n} ) =  ( \left ( \begin {array} {cccc}
    e_n \\ e_n \end {array} \right ) \otimes \Delta ( q^{+}_{n} ). $
   
   $ U_{(12)} U_{(13)} ( \left ( \begin {array} {cccc}
   e_n \\ e_n \end {array} \right ) ) =  \left ( \begin {array} {cccc}
   e_n \\ e_n \end {array} \right ) \otimes q^{+}_{n} \otimes  q^{+}_{n}.$
   
   Hence, $ \Delta ( q^{+}_{n} ) = q^{+}_{n} \otimes q^{+}_{n}.$
   
   Similarly, $ \Delta ( q^{-}_{n} ) = q^{-}_{n} \otimes q^{-}_{n}.$
   
   Moreover, $ \Delta ( y_1 ) = \Delta ( q^{-}_{n} {q^{+}_{n}}^{*} ) = ( q^{-}_{n} \otimes q^{-}_{n} ) ( {q^{+}_{n}}^{*} \otimes {q^{+}_{n}}^{*} ) =  q^{-}_{n}  {q^{+}_{n}}^{*} \otimes q^{-}_{n} {q^{+}_{n}}^{*} = y_1 \otimes y_1.$  \qed 
   
   \vspace{4mm}

 Let us now consider the quantum group $\widetilde{\cls}  \cong C^* ( \IZ_2 \ast \IZ^\infty ),$ where $\IZ^\infty=\IZ \ast \IZ \ast \cdots$ denotes the free product of countably infinitely many copies of $\IZ$. By the Remarks \ref{preliminaries_free_product_C*_distributive}, \ref{preliminaries_amenable_reduced=full} and \ref{preliminaries_reduced_dual}, $ \widetilde{\cls}  \cong C(\IZ_2) \ast C(\IT) \ast C(\IT) \ast  \cdots,$ and let us denote by  $ r^{+}_{n} $ the generator of n th copy of $ C ( \IT ) $ and by $ y $  the generator of $ C ( \IZ_{2} ).$
 
 The coproduct  $\Delta_0$ on $\widetilde{\cls}$ is  given by $\Delta_0(r^+_n) = r^+_n \ot r^+_n,~\Delta_0(y)=y \ot y.$

Define 
 \bean V \left ( \begin {array} {cccc}
   e_n \\ e_n \end {array} \right ) = \left ( \begin {array} {cccc}
   e_n \\ e_n \end {array} \right ) \otimes r^{+}_{n}. \eean   
 \bean  V \left ( \begin {array} {cccc}
   e_n \\ - e_n \end {array} \right ) = \left ( \begin {array} {cccc}
   e_n \\ - e_n \end {array} \right ) \otimes r^{+}_{n} y. \eean

   \blmma
   
$ V $ commutes with $ D $ and $ V $ is a unitary representation of $\widetilde{\cls}$ , that is $(\widetilde{\cls}, \Delta_0, V)$  is an object in ${\bf Q}^\prime(D).$

   \elmma
   
   {\it Proof :}  As the eigenspaces corresponding to distinct eigenvalues of $ D $ are spanned by $ \left ( \begin {array} {cccc}
   e_n \\ e_n \end {array} \right ) $ and $ \left ( \begin {array} {cccc}
   e_n \\ - e_n \end {array} \right ) ,~ V $ commutes with $ D .$
   
   The fact that  $ V $ is a representation follows from the proof of Lemma \ref{sphere_chak_pal_coproduct}.
   
   To prove that $ V $ is a unitary, it is enough to check the following:
   
   \vspace{4mm}
   
   $ \left\langle V \left ( \begin {array} {cccc}
   e_n \\ e_n \end {array} \right ), V  \left ( \begin {array} {cccc}
   e_m \\ e_m \end {array} \right ) \right\rangle = \left\langle \left ( \begin {array} {cccc}
   e_n \\ e_n \end {array} \right ) , \left ( \begin {array} {cccc}
   e_m \\ e_m \end {array} \right ) \right\rangle.1, ~~ \left\langle V \left ( \begin {array} {cccc}
   e_n \\ e_n \end {array} \right ), V  \left ( \begin {array} {cccc}
   e_m \\ - e_m \end {array} \right ) \right\rangle = \left\langle \left ( \begin {array} {cccc}
   e_n \\ e_n \end {array} \right ) , \left ( \begin {array} {cccc}
   e_m \\ - e_m \end {array} \right ) \right\rangle.1, ~~ \left\langle V \left ( \begin {array} {cccc}
   e_n \\ - e_n \end {array} \right ), V  \left ( \begin {array} {cccc}
   e_m \\ -  e_m \end {array} \right ) \right\rangle = \left\langle \left ( \begin {array} {cccc}
   e_n \\ - e_n \end {array} \right ) , \left ( \begin {array} {cccc}
   e_m \\ - e_m \end {array} \right ) \right\rangle.1, ~~ $ 
    
    \vspace{4mm}
   
   The proofs being similar, we prove only the first equation.
   
   \vspace{4mm}
   
    $ \left\langle V \left ( \begin {array} {cccc}
   e_n \\ e_n \end {array} \right ), V  \left ( \begin {array} {cccc}
   e_m \\ e_m \end {array} \right ) \right\rangle = \left\langle  \left ( \begin {array} {cccc}
   e_n \\ e_n \end {array} \right ) \otimes r^{+}_n  , ~  \left ( \begin {array} {cccc}
   e_m \\ e_m \end {array} \right ) \otimes r^{+}_m \right\rangle $

    $ = \left\langle  \left ( \begin {array} {cccc}
   e_n \\ e_n \end {array} \right ), ~ \left ( \begin {array} {cccc}
   e_m \\ e_m \end {array} \right ) \right\rangle \otimes r^{+*}_n r^{+}_m, $
   
   which, if $ n \neq m $ equals $ 0 = \left\langle \left ( \begin {array} {cccc}
   e_n \\ e_n \end {array} \right )  ,  \left ( \begin {array} {cccc}
   e_m \\ e_m \end {array} \right ) \right\rangle. 1, $
   
   and if $ n = m,$ equals $ 1 = \left\langle \left ( \begin {array} {cccc}
   e_n \\ e_n \end {array} \right ) , \left ( \begin {array} {cccc}
   e_m \\ e_m \end {array} \right )  \right\rangle.1.  $ \qed

   \vspace{4mm}

   Setting $ r^{-}_{n} = r^{+}_{n} y ,$ we observe that 
    $ r^{-}_{n} $ is a unitary and satisfies :   
   \be \label{sphere_chak_pal_id1} r^{-}_{n} r^{+*}_{n - 1} = r^{+}_{n} r^{-*}_{n - 1} ~~ {\rm ~ for ~ all} ~ n \geq 1. \ee   
   \be \label{sphere_chak_pal_id2} r^{+}_{n} r^{+*}_{n - 1} = r^{-}_{n} r^{-*}_{n - 1} ~~ {\rm ~ for ~ all} ~ n \geq 1 . \ee   
   Using $ r^{-}_{n} = r^{+}_{n} r^{-*}_{n - 1} r^{+}_{n - 1} $ ( from  \ref{sphere_chak_pal_id1} ) in ( \ref{sphere_chak_pal_id2} ) we have $ r^{+}_{n - 1} = r^{-}_{n - 1} r^{+*}_{n - 1} r^{-}_{n - 1} .$
   
   This implies   
   \be \label{sphere_chak_pal_r1} r^{+}_{n}{r^{-}_{n}}^{*} = r^{-}_{n} r^{+*}_{n} ~~ {\rm ~ for ~ all} ~ n. \ee    
     Moreover, taking $ \ast $ on ( \ref{sphere_chak_pal_id1} ) and ( \ref{sphere_chak_pal_id2} )  respectively, we get the following:     
    \be \label{sphere_chak_pal_r2} r^{+}_{n - 1} r^{-*}_{n} - r^{-}_{n - 1}r^{+*}_{n} = 0 ~~ {\rm ~ for ~ all} ~ n \geq 1. \ee     
    \be   \label{sphere_chak_pal_r3}  r^{+}_{n - 1} r^{+*}_{n} - r^{-}_{n - 1} r^{-*}_{n} = 0 ~~ {\rm ~ for ~ all} ~ n \geq 1. \ee     
    Thus, the equations ( \ref{sphere_chak_pal_alg1} ) -  ( \ref{sphere_chak_pal_alg4} ) in Lemma \ref{sphere_chak_pal_alg} are satisfied with $ q_n^{\pm} $'s replaced by $ r_n^{\pm} $'s and hence it is easy to see that there is a $C^*$-homomorphism from $\widetilde{\cls}$ to $\tilde{\clq}$ sending $ y, r^{+}_{n} $ to $ y_0 $ and $ q^{+}_{n} $ respectively,  which is surjective by Proposition \ref{sphere_chak_pal_description of Q'+} and is a CQG morphism by Lemma \ref{sphere_chak_pal_coproduct}. In other words, $(\widetilde{\cls}, \Delta_0, V)$ is indeed a universal object in ${\bf Q}^\prime(D)$. It is clear that the maximal Woronowicz subalgebra of $\widetilde{\cls}$ for which the action is faithful, that is $QISO^+(D)$, is generated by $r^+_{n-1}r^{+^*}_n$, $n \geq 1$ and $r_0^+y r_1^{+^*}$, so again isomorphic with $C^{*} ( \IZ_{2}  \ast \IZ^{\infty} ).$

 This is summarized in the following:

   \bthm
   
   The universal object in the category ${\bf Q}^\prime(D)$, that is  $ \widetilde{QISO^{+}}( D ) $ exists and is isomorphic with $ C^{*} ( \IZ_{2}  \ast \IZ^{\infty} ).$  Moreover, the quantum group $QISO^+(D)$ is again isomorphic with $C^{*} ( \IZ_{2}  \ast \IZ^{\infty} ).$
   
   \ethm

\brmrk
This example shows that $QISO^+$ in general may not be matrix quantum group, that is may not have a finite dimensional fundamental unitary, even if the underlying spectral triple is of compact type. This is somewhat against the intuition derived from the classical situation, since for a classical compact Riemannian manifold the group of isometries is always a compact Lie group, hence has an embedding into the group of orthogonal matrices of some finite dimension. 
\ermrk

\vspace{4mm}

We end this chapter by noting that $ \alpha $ gives an example where the quantum group of orientation preserving isometries does not have a $ C^* $ action. Before that, we recall some useful properties of the so called Toeplitz algebra from \cite{davidson}.

\bppsn

\label{sphere_toeplitz_properties}

Let $\tau_1$ be the unilateral shift operator on $l^2( \IN )$ defined by $\tau_1 ( e_n ) = e_{n - 1}, ~ n \geq 1, ~ \tau ( e_0 ) = 0.$ Then the Toeplitz algebra $ C^* ( \tau_1 ).$ is the $C^*$ algebra generated by $\tau_1,$  on $l^2(\IN ) .$ It contains all compact operators and moreover, the commutator of any two elements of $ C^* ( \tau_1 ) $ is compact. 

\eppsn

Let $\tau$ be the operator on $\clh$ defined by $\tau = \tau_1 \otimes {\rm id}.$ 

\blmma

\label{sphere_chak_pal_not_action_1}

$ B = \tau \left| B \right| .$

\elmma

{\it Proof :} We note that \bean \lefteqn{ \left| B \right| (  \left ( \begin {array} {cccc}
   e_n \\ 0 \end {array} \right ) )
= {( A - A^2 + c I )}^{\frac{1}{2}} (  \left ( \begin {array} {cccc}
   e_n \\ 0 \end {array} \right ) )}\\
 &=& \sqrt{\lambda_+ \mu^{2n} - \lambda^2_{+} \mu^{4n} + c }  (  \left ( \begin {array} {cccc}
   e_n \\ 0 \end {array} \right ) )
 = {c_+ ( n )}^{\frac{1}{2}} \left ( \begin {array} {cccc}
   e_n \\ 0 \end {array} \right ) \eean   
   and  hence $ \tau \left| B \right| (  \left ( \begin {array} {cccc}
   e_n \\ 0 \end {array} \right ) ) = {c_+ (n)}^{\frac{1}{2}} (  \left ( \begin {array} {cccc}
   e_{n - 1} \\ 0 \end {array} \right ) ) = B (  \left ( \begin {array} {cccc}
   e_n \\ 0 \end {array} \right ) ).$

Similarly, $ \tau \left| B \right| (  \left ( \begin {array} {cccc}
   0 \\ e_n \end {array} \right ) ) = B (  \left ( \begin {array} {cccc}
   0 \\ e_n \end {array} \right ) ). $ This completes the proof of the Lemma. \qed

\blmma

\label{sphere_chak_pal_not_action_2}

 $$ \alpha(\tau)=\sum_{n \geq 1} \tau (P_n +Q_n) \otimes r^+_{n-1}{r^+_n}^*,$$ where $r_n^{\pm}$ are the elements of $\widetilde{QISO^+}(D)$
 introduced before. 

\elmma 

{\it Proof :} For all $ n \geq 1, $ we have  \bean \lefteqn{ \alpha ( \tau ) ( \left ( \begin {array} {cc}
   e_n \\ 0 \end {array} \right ) )  =  \widetilde{U} ( \tau \otimes {\rm id} ) {\widetilde{U}}^* \left ( \begin {array} {cccc}
   e_n \\ 0 \end {array} \right )}\\
&=& \frac{1}{2} \widetilde{U} ( \tau \otimes {\rm id} ) {\widetilde{U}}^* [ \left ( \begin {array} {cccc}
   e_n \\ e_n \end {array} \right ) + \left ( \begin {array} {cccc}
   e_n \\ - e_n \end {array} \right )  ]\\
&=& \frac{1}{2} \widetilde{U} ( \tau \otimes {\rm id} ) [ \left ( \begin {array} {cccc}
   e_n \\ e_n \end {array} \right ) \otimes  r^{+*}_n + \left ( \begin {array} {cccc}
   e_n \\ -  e_n \end {array} \right ) \otimes r^{-*}_n ]\\
&=&  \frac{1}{2} \widetilde{U} ( \tau \otimes {\rm id} ) [ \left ( \begin {array} {cccc}
   e_n \\ 0 \end {array} \right ) \otimes ( r^{+*}_n + r^{-*}_n ) + \left ( \begin {array} {cccc}
   0 \\ e_n \end {array} \right ) \otimes ( r^{+*}_n - r^{-*}_n ) ]\\
&=&  \frac{1}{2} \widetilde{U} [ \left ( \begin {array} {cccc}
   e_{n - 1} \\ 0 \end {array} \right )  \otimes ( r^{+*}_n + r^{-*}_n ) + \left ( \begin {array} {cccc}
   0 \\ e_{n - 1} \end {array} \right ) \otimes ( r^{+*}_n - r^{-*}_n )   ]\\
&=& \frac{1}{4} \widetilde{U} [ \left ( \begin {array} {cccc}
   e_{n - 1} \\ e_{n - 1} \end {array} \right ) \otimes ( 2 r^{+*}_n ) + \left ( \begin {array} {cccc}
   e_{n - 1} \\ - e_{n - 1} \end {array} \right ) \otimes 2 r^{-*}_n ]\\
 &=& \frac{1}{2} [ \left ( \begin {array} {cccc}
   e_{n - 1} \\ e_{n - 1} \end {array} \right ) \otimes r^+_{n - 1} r^{+*}_n + \left ( \begin {array} {cccc}
   e_{n - 1} \\ - e_{n - 1} \end {array} \right ) \otimes r^{-}_{n - 1} r^{-*}_n  ] = \left ( \begin {array} {cccc}
   e_{n - 1} \\ 0 \end {array} \right ) \otimes  r^+_{n - 1} r^{+*}_n. \eean
Similarly, $ \alpha ( \tau ) \left ( \begin {array} {cccc}
   0 \\  e_{n} \end {array} \right ) = \left ( \begin {array} {cccc}
   0 \\  e_{n - 1} \end {array} \right ) \otimes r^{+}_{n - 1} r^{+*}_n $ for all $ n \geq 1. $
   
   Moreover, $ \alpha ( \tau ) \left ( \begin {array} {cccc}
   e_0 \\  0 \end {array} \right ) = \alpha ( \tau ) \left ( \begin {array} {cccc}
   0 \\  e_{0} \end {array} \right ) = 0.   $  

Thus, $ \alpha ( \tau ) = \sum_{n \geq 1} \tau P_n \otimes r^+_{n - 1} r^{+*}_n + \sum_{n \geq 1} \tau Q_n \otimes r^{+}_{n - 1} r^{+*}_n = \sum_{n \geq 1} \tau ( P_n + Q_n ) \otimes r^+_{n - 1} r^{+*}_n$
 \qed

\bthm
 \label{no_action}
 The $\ast$-homomorphism $ \alpha $ is not a $C^*$ action.

\ethm
{\it Proof :}  We begin with the observation that each of the $C^*$ algebras $\pi_{\pm}(S^2_{\mu, c})$ is nothing but the Toeplitz algebra. For example, consider $\clc:=\pi_+(S^2_{\mu, c})$.
 Clearly, $T=\pi_+(B)$ in an invertible operator with the polar decomposition given by, $ T=\tau_1 |T|,$ hence $\tau_1$ belongs to $ \clc.$ 
Thus, $\clc$ contains the Toeplitz algebra $C^*(\tau_1)$, which by Proposition \ref{sphere_toeplitz_properties} contains all compact operators.  
In particular, $C^*(\tau_1)$ must contain $\pi_+(A)$ as well as all the eigenprojections  $P_n$ of $|\pi_+(B)|$ so it must contain the whole of $\clc$. 
Similar arguments will work for $\pi_-(S^2_{\mu,c}).$

Thus,  $\tau= \tau_1 \oplus \tau_1 = \pi(B)| \pi(B)|^{-1}  $ belongs to $ \pi(S^2_{\mu,c}).$ If $\alpha$ is a $C^*$ action, then for an arbitrary state $\phi$ on $QISO^+(D)$ we must have $\alpha_\phi(\tau) \equiv ({\rm id} \otimes \phi) \circ \alpha (\tau) $ is in $ \pi(S^2_{\mu,c}),$  hence $\alpha_\phi(\tau)P_+$ must belong to $\clc=\pi_+(S^2_{\mu,c})$, where $P_+$ denotes the projection 
onto $\clh_+$. By Proposition \ref{sphere_toeplitz_properties}, this implies that $[\alpha_\phi(\tau)P_+, \tau_1]$ must be a compact operator.
 We claim that for suitably chosen $\phi$, this compactness condition is violated, which will complete the proof of the theorem.

To this end, fix an irrational number $\theta$ and consider the sequence $\lambda_n=e^{2 \pi i n \theta}$ of complex number of unit modulus. 
 We note that the linear functionals  which send the generator of $C(Z_2) $( which is $ y $ ) to $1$ and the generator of the $n$-th copy of $C(\IT)$ ( which is $ r^+_{n - 1} r^{+*}_{n} $ by Proposition \ref{sphere_chak_pal_description of Q'+}    ) to $\lambda_n$ are evaluation maps and hence homomorphisms. Using Remark \ref{preliminaries_free_product_of_homomorphisms}, we have a unital $\ast$-homomorphism $\phi : QISO^+(D) =C(Z_2) \ast C(\IT)^{*^\infty} \rightarrow \IC$ which extends the above mentioned homomorphisms. Hence, $ \alpha_\phi(\tau)=\sum_n \lambda_n \tau (P_n +Q_n).$ Moreover,  we see that 
\bean \lefteqn{[\alpha_\phi(\tau)P_+, \tau_1]\left ( \begin {array} {cccc}
   e_{n} \\ 0 \end {array} \right )}\\
 &=& ({\rm id} \otimes \phi ) \alpha ( \tau )\left ( \begin {array} {cccc}
   e_{n - 1} \\ 0 \end {array} \right ) - \tau ( {\rm id} \otimes \phi ) [ \left ( \begin {array} {cccc}
   e_{n - 1} \\ 0 \end {array} \right ) \otimes r^+_{n - 1} r^{+*}_n ]\\
  &=& (\lambda_{n-1}-\lambda_n) \left ( \begin {array} {cccc}
   e_{n - 2} \\ 0 \end {array} \right ) ,~~n \geq 2.\eean
Similarly, $ [\alpha_\phi(\tau)P_+, \tau_1] \left ( \begin {array} {cccc}
   0 \\ e_n \end {array} \right ) = (\lambda_{n - 1} - \lambda_n) \left ( \begin {array} {cccc}
   0 \\ e_{n - 2} \end {array} \right ) .$
Hence, the above commutator cannot be compact since $\lambda_n-\lambda_{n-1}$ does not go to $0$ as $n \rightarrow \infty.$ 
\qed
   
\bcrlre

\label{sphere_chak_pal_non_existence_wang_univ}

The subcategory  of ${\bf Q}^\prime(D)$ consisting of objects $(\tilde{\clq}, U)$ where $\alpha_U$ is a $C^*$ action does not have a universal object.

\ecrlre
 {\it Proof:} By Theorem \ref{no_action}, the proof will be complete if we can show that if a universal object exists for the subcategory (say ${\bf Q}_1^\prime$) mentioned 
above, then it must be isomorphic with $\widetilde{QISO^+}(D)$. For this, consider the quantum subgroups $\widetilde{\clq}_N$, $N=1,2,...,$ 
of $\widetilde{QISO^+}(D)$ generated by $r_n^+, n=1,...,N$ and $y$. Let $\pi_N : \widetilde{QISO^+}(D) \raro \widetilde{\clq}_N$ be 
the CQG morphism given by $\pi_N(y)=y,~ \pi_N(r_n^+)=r_n^+$ for $n \leq N$ and $\pi_N(r^+_n)=1$ for $n > N.$ 

We claim that $(\widetilde{\clq}_N, U_N:=({\rm id} \ot \pi_N)\circ V)$ is an object in ${\bf Q}_1^\prime$ 
( where $V$ denotes the unitary representation of $\widetilde{QISO^+}(D)$ on $\clh$ ). 
To see this, we first note that for all $ N,$
$ ( {\rm id} \otimes \pi_N ) \alpha ( A ) = \sum^N_{n = 0} A P_n \otimes \frac{1}{2 \lambda_+} \{ \lambda_+ ( 1 + r^+_n y r^{+*}_n ) + \lambda_{-} ( 1 - r^+_n y r^{+*}_n ) \} + \sum^N_{n = 0} A Q_n \otimes \frac{1}{2 \lambda_{-}} \{ \lambda_+ ( 1 - r^+_n y r^{+*}_n ) + \lambda_{-} ( 1 + r^{+}_n y r^{+*}_n ) \} + \sum^{\infty}_{n = N + 1} A P_n \otimes \frac{1}{2 \lambda_+} \{ \lambda_+ (  1 + y ) + \lambda_{-} ( 1 - y ) \} + \sum^{\infty}_{n = N + 1} A Q_n \otimes  \frac{1}{2 \lambda_{-}} \{  \lambda_+ ( 1 - y ) + \lambda_{-} ( 1 + y ) \}.$

Among the four summands, the first two clearly belong to $ \cla \otimes \widetilde{\clq}_N.$ Moreover, the sum of the third and the fourth summand equals $ A ( 1 - \sum^{N}_{n = 1} P_n ) 
\otimes \frac{1}{2 \lambda_+} \{ \lambda_+ (  1 + y ) + \lambda_{-} ( 1 - y ) \} + A ( 1 - \sum^N_{n = 1} Q_n ) \otimes \frac{1}{2 \lambda_{-}} \{ \lambda_+ ( 1 - y ) + \lambda_{-} ( 1 + y ) \} $ which is an element of $ \cla \otimes \widetilde{\clq}_N.$

 We proceed similarly in the case of  $ B, $ to note that it is enough to show that for all $ N,$
$$  \sum^{\infty}_{n = N + 2} B P_n \otimes \frac{{c_{+}( n )}^{\frac{1}{2}} + {c_{-}( n )}^{\frac{1}{2}} }{2 c_+ ( n )}  + \sum^{\infty}_{n = N + 2} B P_n \otimes \frac{( {c_{+}( n )}^{\frac{1}{2}} - {c_{-}( n )}^{\frac{1}{2}} ) y }{2 c_+ ( n )} $$
$$ + \sum^{\infty}_{n = N + 2} B Q_n \otimes \frac{{c_{+}( n )}^{\frac{1}{2}} + {c_{-}( n )}^{\frac{1}{2}} }{2 c_- ( n )} - 
\sum^{\infty}_{n = N + 2} B Q_n \otimes \frac{( {c_{+}( n )}^{\frac{1}{2}} - {c_{-}( n )}^{\frac{1}{2}} ) y }{2 c_- ( n )} $$
 belongs to $ \cla \otimes \widetilde{\clq_{N}}.$ The norm of the second and the fourth term can be easily seen to be finite. 
The first term equals $ \frac{1}{2} B ( 1  - \sum^{N + 1}_{n = 1} P_n ) [ {( A - A^2 + c I )}^{- \frac{1}{2}} + {( A - A^2 + c I  )}^{ - 1} 
{\{ \frac{\lambda_{-}}{\lambda_+} A - {( \frac{\lambda_{-}}{\lambda_{+}} A  )}^2 + c I \}}^{\frac{1}{2}}  ] \otimes 1 $ and
 therefore belongs to $ \cla \otimes \widetilde{\clq_{N}}.$ The third term can be treated similarly.

 Thus, there is surjective  CQG morphism 
$\psi_N $ from the universal object, say $\widetilde{\clg}$, of ${\bf Q}_1^\prime$ to $\widetilde{\clq}_N$. Clearly, $ (  \widetilde{\clq}_N )_{N \geq 1}$ 
form an inductive system of objects in ${\bf Q}^\prime(D)$, with the inductive limit being $\widetilde{QISO^+}(D)$, and the morphisms 
$\psi_N$ induce a surjective CQG morphism (say $\psi$) from $\widetilde{\clg}$ to $\widetilde{QISO^+}(D)$. But $\widetilde{\clg}$ is an object
 in ${\bf Q}^\prime(D)$, so must be a quantum subgroup of the universal object in this category, that is, $\widetilde{QISO^+}(D)$. 
This gives the CQG morphism from $\widetilde{QISO^+}(D)$ onto $\widetilde{\clg}$, which is obviously the inverse of $\psi$, and hence we get the 
desired isomorphism between $\widetilde{\clg}$ and $\widetilde{QISO^+}(D).$
 \qed

\cleardoublepage
\addcontentsline{toc}{chapter} 
{\protect\numberline{Bibliography\hspace{-96pt}}} 
\cleardoublepage

\end{document}